\documentclass{report}

\usepackage{amsfonts,amsthm}
\usepackage{amsmath}   
\usepackage{latexsym}
\usepackage{amssymb}
\usepackage{mathrsfs}
\usepackage{pifont}
\usepackage{tikz}
\usepackage[all]{xy}
\usepackage[hidelinks
]{hyperref}
\usepackage{xcolor}
\usepackage{amsthm}
\usepackage{framed}

\newtheoremstyle{rem}{3pt}{3pt}{}{}
{\bfseries}{.}{.5em}{}

\newtheorem{theo}{Theorem}[chapter]

\newtheorem{prop}[theo]{Proposition}

\newtheorem{coro}[theo]{Corollary}
\newtheorem{conj}[theo]{Conjecture}
\newtheorem{conc}[theo]{Conclusion}
\theoremstyle{rem}
\newtheorem{remax}[theo]{Remark}
\newenvironment{rema}
  {\pushQED{\qed}\remax}
  {\popQED\endremax}
\theoremstyle{definition}
\newtheorem{defi}[theo]{Definition}
\newtheorem{exam}[theo]{Example}

\newtheorem*{term*}{Notation/Terminology}

\colorlet{shadecolor}{gray!15}
\newenvironment{framedconc}
  {\begin{shaded}\begin{conc}}
  {\end{conc}\end{shaded}}
\colorlet{shadecolor}{gray!15}
\newenvironment{framedtheo}
  {\begin{shaded}\begin{theo}}
  {\end{theo}\end{shaded}}
\colorlet{shadecolor}{gray!15}
\newenvironment{frameddefi}
  {\begin{shaded}\begin{defi}}
  {\end{defi}\end{shaded}}
\colorlet{shadecolor}{gray!15}
\newenvironment{framedconj}
  {\begin{shaded}\begin{conj}}
  {\end{conj}\end{shaded}}

\newcommand{\C}{\mathbb{C}}
 
\newcommand{\Mat}{\operatorname{Mat}}

\newcommand{\mS}{\mathfrak{S}}

\newcommand{\tg}{\widetilde{g}}

\newcommand{\tm}{\textbf{m}}

\def\sl{{\mathfrak{sl}}}

\def\sl{{\mathfrak{sl}}}

\def\p{{\textit{\textsf{P}}}}
\def\f{{\textit{\textsf{F}}}}

\newcommand{\bT}{\boldsymbol{t}}
\newcommand{\bbT}{\boldsymbol{T}}
\newcommand{\STab}{\mathrm{STab}}
\newcommand{\SSTab}{\mathrm{SSTab}}

\newcommand{\bk}{\mathbf{k}}
\newcommand{\bC}{\mathbb{C}}
\newcommand{\bZ}{\mathbb{Z}}
\newcommand{\CC}{\mathbb{C}}
\newcommand{\si}{\sigma}

\newcommand{\cc}{\mathrm{cc}}
\newcommand{\qc}{\mathrm{c}}

\newcommand{\AS}{\mathrm{AS}}

\newcommand{\cD}{\mathcal{D}}

\newcommand{\HJ}{H^J}

\newcommand{\diag}[3]{ \foreach \t in {1,...,#3} {\draw[thick] (#1+\t,#2-1) rectangle (#1+\t-1,#2);} }
\newcommand{\diagg}[4]{ \foreach \t in {1,...,#3} {\draw[thick] (#1+\t,#2-1) rectangle (#1+\t-1,#2);} \foreach \t in {1,...,#4} {\draw[thick] (#1+\t,#2-1) rectangle (#1+\t-1,#2-2);} }
\newcommand{\diaggg}[5]{ \foreach \t in {1,...,#3} {\draw[thick] (#1+\t,#2-1) rectangle (#1+\t-1,#2);} \foreach \t in {1,...,#4} {\draw[thick] (#1+\t,#2-1) rectangle (#1+\t-1,#2-2);}
                         \foreach \t in {1,...,#5} {\draw[thick] (#1+\t,#2-2) rectangle (#1+\t-1,#2-3);} }
\newcommand{\diagggg}[6]{ \foreach \t in {1,...,#3} {\draw[thick] (#1+\t,#2-1) rectangle (#1+\t-1,#2);} \foreach \t in {1,...,#4} {\draw[thick] (#1+\t,#2-1) rectangle (#1+\t-1,#2-2);}
                         \foreach \t in {1,...,#5} {\draw[thick] (#1+\t,#2-2) rectangle (#1+\t-1,#2-3);} \foreach \t in {1,...,#6} {\draw[thick] (#1+\t,#2-3) rectangle (#1+\t-1,#2-4);} }    
\newcommand{\diaggggg}[6]{ \foreach \t in {1,...,#3} {\draw[thick] (#1+\t,#2-1) rectangle (#1+\t-1,#2);} \foreach \t in {1,...,#4} {\draw[thick] (#1+\t,#2-1) rectangle (#1+\t-1,#2-2);}
                         \foreach \t in {1,...,#5} {\draw[thick] (#1+\t,#2-2) rectangle (#1+\t-1,#2-3);} \foreach \t in {1,...,#6} {\draw[thick] (#1+\t,#2-3) rectangle (#1+\t-1,#2-4);}
                         \foreach \t in {1,...,#6} {\draw[thick] (#1+\t,#2-4) rectangle (#1+\t-1,#2-5);} }

\usepackage{geometry}
\usepackage{minitoc}
\makeatletter
\let\original@addcontentsline\addcontentsline
\newcommand{\dummy@addcontentsline}[3]{}
\newcommand{\DeactivateToc}{\let\addcontentsline\dummy@addcontentsline}
\newcommand{\ActivateToc}{\let\addcontentsline\original@addcontentsline}
\makeatother

\renewcommand\thechapter{\Roman{chapter}}

\begin{document}

\begin{center}

\thispagestyle{empty}

\vspace{1cm}
\textbf{\hrulefill}

\huge{\textsf{Centralisers and Hecke algebras in Representation Theory}}\\
\LARGE{\textsf{with applications to Knots and Physics}}

\vspace{0.3cm}
\textbf{\hrulefill}

\vskip 1cm
{\Large L. Poulain d'Andecy\footnote{Laboratoire de Math\'ematiques de Reims UMR 9008, Universit\'e de Reims Champagne-Ardenne,
Moulin de la Housse BP 1039, 51100 Reims, France.
loic.poulain-dandecy@univ-reims.fr}}

\vskip 2cm
{\large M\'emoire d'Habilitation \`a diriger des recherches}

\end{center}

\vskip 1cm
{\normalsize \textbf{Abstract.} This document is a thesis presented for the ``Habilitation \`a diriger des recherches''. The first chapter provides some background and sketch the story of the classical Schur--Weyl duality and its quantum analogue involving the Hecke algebra. Quantum groups and their centralisers are discussed along with their connections with the braid group, the Yang--Baxter equation and the construction of link invariants. A quick review of several related notions of Hecke algebras concludes the first chapter. The second chapter presents the fused Hecke algebras and an application to the construction of solutions for the Yang--Baxter equation. The third chapter summarises recent works on diagonal centralisers of Lie algebras, and gives an application to the study of the missing label for sl(3). The last chapter discusses algebraic results on the Yokonuma--Hecke algebras and their applications to the study of link invariants.}

\setcounter{secnumdepth}{3}
\dominitoc
\setcounter{tocdepth}{0}
\tableofcontents
\addtocontents{toc}{\protect\thispagestyle{empty} 
                    \protect\pagestyle{empty}}

\thispagestyle{empty}

\renewcommand{\theequation}{0.{\arabic{equation}}}
\DeactivateToc
\chapter*{\huge{Introduction}}
\ActivateToc
\addcontentsline{toc}{chapter}{\large{Introduction\vspace{0.3cm}}}
\mtcaddchapter
\setcounter{page}{1}

\renewcommand{\theequation}{\thechapter.{\arabic{equation}}}

Representation theory was officially born at the very end of the 19th century in the work of Georg Frobenius connected to group theory and number theory. It developed rapidly and is now a very active mathematical subject on its own, which moreover has many interactions with other subjects in mathematics and in other sciences. The usual way to approach representation theory is through group theory. First one learns about abstract groups and the general theory of group actions on sets, and then one defines representations as a particular kind of group actions: they are group actions on vector spaces through linear operators. The first historical motivations for representations were to find applications in group theory and in number theory. Striking applications have indeed been found, for example some representation theory is used in the classification of finite simple groups, or in the proof of Fermat's last theorem. However, it is clear now that representation theory is not reduced to its possible applications to other areas: it is considered as a fundamental subject in mathematics on its own. For historical accounts of representation theory, one can consult \cite{Cur2,Haw1,Haw2}.

A significant increase of interest in representation theory in the first half of the 20th century is tied to the development of quantum physics (see for example \cite{Mac1,Wey3}). Since quantum physics fundamentally involves vector spaces and linear operators, the groups or algebras, describing symmetries or commutation relations between observables, appear in fact through their linear representations. In physics, representations are usually what appears concretely while group theory, or more involved algebraic theories, are the mathematical tools to understand them (this is quite visible for instance in the standard model of particle physics). The theory of ``spin'' in quantum physics is a good example. The spin is modelled using a Hilbert space of dimension 1 for spin 0 particles, of dimension 2 for spin 1/2 particles, and so on. On each of these vector spaces, there are observables, giving for example the total spin or the values of the components along chosen axes. They are matrices of various sizes depending on the value of the spin, and they are the objects of primary importance for the physical theory. It turns out that all these vector spaces are best understood as irreducible representations of the group $SU(2)$ (of which there is indeed one for each dimension) and therefore all matrices involved in spin theory are unified under the common banner of the representation theory of $SU(2)$.

The above example is meant to illustrate how fundamental representations can be in physics, and how useful it can be to realise a set of matrices as being in fact the representation of a group. This agrees with the general philosophy that \emph{a representation is a group (or an algebra) disguised as a set of matrices}. Matrices are simple and fundamental objects in mathematical sciences and this explains the usefulness and universality of representation theory.

\vskip .2cm
Issai Schur and Hermann Weyl are two of the pioneers of representation theory. Their names are forever associated in what is called the Schur--Weyl duality. The remarkable and beautiful aspect of the Schur--Weyl duality is that it relates intimately the representation theories of, probably, the two most fundamental groups that we can think about. On one hand, we have a finite group, the symmetric group $S_n$, and on the other hand, we have a matrix group, or a Lie group, the general linear group $GL_N(\mathbb{C})$ (if one prefers unitary matrices, we can formulate the Schur--Weyl duality with the special unitary group $SU(N)$ instead of $GL_N(\mathbb{C})$).

There are several possible formulations of the Schur--Weyl duality, with various degrees of ``abstractness'' and generality. One of the simplest and most concise is the following: on the vector space $\bigl(\mathbb{C}^N\bigr)^{\otimes n}$, there is a representation of $GL_N(\mathbb{C})$ and a representation of $S_n$, and the corresponding two sets of matrices mutually centralise each other. This formulation highlights the notion of centraliser of a representation.

When we have at hand a set of matrices (for example, a representation), looking at its centraliser is a very natural idea. In physics, this often corresponds to studying symmetries. In general, from a mathematical point of view, the centraliser contains a large amount of information on the set of matrices that we started with. So following this general line of thought, in representation theory, as soon as we have a representation of some groups or some algebras that we would like to study, we should desire to understand its centraliser. An especially enjoyable situation is when we recover, looking at such centralisers, algebras in which we are independently interested, like in the Schur--Weyl duality. In this case, this whole idea of centralisers becomes a bridge between two independently interesting algebraic structures.

\paragraph{Main theme of this thesis.} This thesis discusses several works which are unified under the common theme of centraliser algebras and their applications, focusing especially in their applications to knots and links theory and to mathematical physics. A picture is:
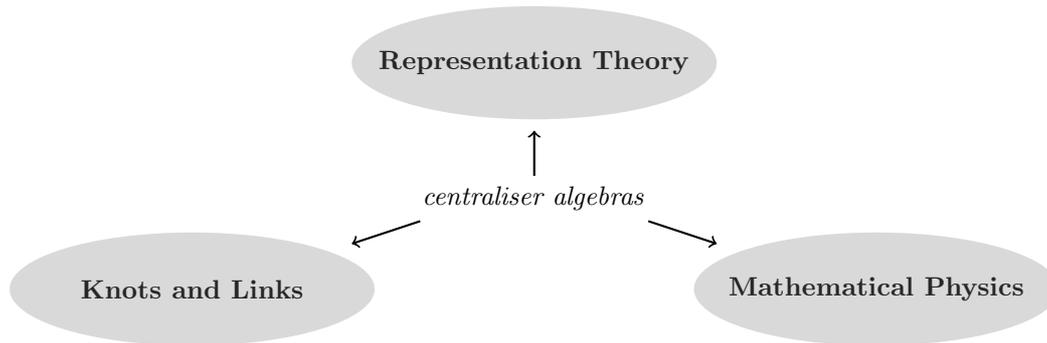
\begin{center}
\begin{tikzpicture}[scale=0.3]
\node at (0,0) {\emph{centraliser algebras}};
\node at (0,6) {\textbf{Representation Theory}};
\fill[thin, fill=gray,opacity=0.3] (0,6) ellipse (8cm and 2.5cm);
\draw[thick,->] (0,1)--(0,3);
\node at (-15,-4) {\textbf{Knots and Links}};
\fill[thin, fill=gray,opacity=0.3] (-15,-4) ellipse (8cm and 2.5cm);
\draw[thick,->] (-5,-1)--(-8,-2);
\node at (15,-4) {\textbf{Mathematical Physics}};
\fill[thin, fill=gray,opacity=0.3] (15,-4) ellipse (8cm and 2.5cm);
\draw[thick,->] (5,-1)--(8,-2);
\end{tikzpicture}
\end{center}
The relevance of representation theory to knots and links can be briefly explained as follows. First the study of knots and links, and especially of their topological invariants, is related to an algebraic theory: the braid group makes its appearance and one needs to find so-called Markov traces. The second step comes from the fact that the study of Markov traces on the braid group in general is too difficult, so we decide to study them through some representations of the braid group. This amounts to replace the braid group by simpler algebras, and this allows to actually produce concrete invariants. The best known examples are the Temperley--Lieb algebra, the Hecke algebra and the Birman--Murakami--Wenzl algebra producing, respectively, the Jones polynomial, the HOMFLYPT polynomial and the Kauffman polynomial. 

But where are centraliser algebras in this story? It turns out that the three algebras mentioned above appear as centraliser algebras. For example, the Hecke algebra is a deformation of the symmetric group and there is a suitable deformation of $GL_N(\mathbb{C})$, called a quantum group, such that the Schur--Weyl duality extends to this quantum setting. In fact, the three algebras above appear as centraliser algebras for some quantum groups and this is no coincidence. 

Quantum groups are deformations of universal enveloping algebras of simple Lie algebras which are carefully designed in order to have the following properties. On one hand they retain most of the algebraic properties of the Lie algebras, in particular their semisimple representation theory is well-understood. On the other hand they diverge from the Lie algebra with respect to tensor products of representations. From their axioms, the centraliser algebras of their tensor product representations support a representation of the braid group, and it turns out that they are sufficiently deformed from the Lie algebra situation so that they provide non-trivial representations (here non-trivial means that they do not factorize through the symmetric group in contrast with the classical Lie algebra situation). Furthermore, all these representations are shown to support a Markov trace leading to links invariants (Reshetikhin--Turaev invariants).

The brief summary in the last paragraph already provides a strong motivation, from knots and links theory, for the study of centraliser algebras of quantum groups representations. There is more motivation coming from mathematical physics. In fact, quantum groups originally appeared from mathematical physics considerations. In a nutshell, they were constructed in order to provide solutions of the Yang--Baxter equation, a master equation in the theory of integrable systems. Again, the centralisers are relevant since these solutions actually live, in a sense, in the centraliser of representations. This is best understood with the example of the Hecke algebra in mind, in which a simple formula provides an abstract solution of the Yang--Baxter equation which is then represented to recover some of the matrix solutions coming from quantum groups.

\vskip .1cm
Historically, Hecke algebras appeared first in a context seemingly very different from quantum groups: the representation theory of finite simple groups of Lie type. A remarkable fact is that they also appeared as centraliser algebras of some representations, so that both points of view are quite well adapted to our discussion. This point of view from finite group theory supplies an infinite number of possible generalisations, and some of them fit quite well in the picture involving knots and links and integrable models. One chapter of this thesis will be devoted to Yokonuma--Hecke algebras which are natural generalisations of Hecke from the point of view of finite groups of Lie type, and which recently turned out to have applications in knots and links theory. Another example are the affine Hecke algebras (coming from $p$-adic groups) which are involved with affine quantum groups, with knots and links inside the solid torus and with the reflection equation (a second master equation in integrable models). More generally, the various Hecke algebras coming from finite or $p$-adic groups of Lie type have been found to have a large number of applications besides pure representation theory, such as in categorification or in the theory of special functions and orthogonal polynomials.

\vskip .1cm
There is another class of algebras besides the Hecke-type algebras which fits as well in the above diagram. The famous examples are the Racah and the Askey--Wilson algebras. They appear originally in the theory of orthogonal polynomials, but it became clear recently that they are also closely related to centraliser algebras. One has to shift a little bit the point of view on centralisers, and consider the centraliser of the diagonal embedding of Lie algebras in $n$-fold tensor product of their enveloping algebras (and the same thing for quantum groups). We call this algebra ``diagonal centralisers''. A chapter of this thesis is devoted to the diagonal centralisers. The Racah algebra has been used, among other things, as a symmetry algebra of superintegrable systems, so the connection with mathematical physics is still found here\footnote{the Racah algebra is named after the physicist Giulo Racah.}. Regarding knot theory, the connection is less obvious (for one thing, they are not immediately quotients of the braid group algebra) but still exists. In fact, the Askey--Wilson algebra has recently been shown  to be isomorphic to a Kauffman bracket skein algebra. Besides, by construction they are closely related to the usual centralisers of representations discussed above, so they are not so far from quantum groups and link invariants. Actually, the combinations of these algebras with more braid-like algebras is an exciting perspective. To conclude, the diagonal centralisers are less well-known than Hecke algebras, but they certainly belong to the framework of this thesis as expressed by the preceding diagram.

\paragraph{Contents.} Some parts of the first chapter are sketchy to save space while still being able to tell the story we wanted to tell. Many references are given. The proofs and more details on the materials of the subsequent chapters can be found in the publications listed at the end before the bibliography.

\vskip .1cm
The first chapter discusses in more details some of the story sketched above. It is intended at providing some background, some notations and some motivation for the results presented afterwards. Centraliser algebras and the Schur--Weyl duality are discussed first. Then the braid group and the Yang--Baxter equation are introduced and the Hecke algebra, as a first example, appears here. In order to properly fit the Hecke algebra in the general theory of centralisers, some discussion about quantum groups is given leading to the remarkable properties of their centralisers mentioned above. Then the construction of Markov traces and link invariants is recalled. Finally, a last section on Hecke algebras in general is provided to connect with the historical definition of Hecke algebras and to set the stage for the Yokonuma--Hecke algebras of the last chapter. Also this last section mentions some recent results on affine and quiver Hecke algebras which are not included in the main part of the thesis since they depart a bit too much from the main theme.

\vskip .1cm
The second chapter presents a new class of algebras fitting precisely in the picture above. These algebras are analogues of the Hecke algebras for more general tensor products of $U_q(gl_N)$ representations. The Schur--Weyl duality is extended to incorporate these algebras, termed fused Hecke algebras. Partial results and conjectures about the algebraic description of the relevant centraliser algebras are presented in terms of the fused Hecke algebras. A braid-like description of the fused Hecke algebra is provided along with a complete study of its representation theory (in the semisimple case). Finally, the bridge towards mathematical physics is crossed with an explicit formula for a solution of the Yang--Baxter equation inside the fused Hecke algebra.

\vskip .1cm
The third chapter puts the centraliser story directly inside the universal enveloping algebras of Lie algebras. The notion of ``diagonal centraliser'' is introduced and the story of the Racah algebra is sketched (the $q$-deformation named Askey--Wilson algebra is also briefly evoked). Then two new examples of diagonal centralisers are discussed in details. For one of them (corresponding to $sl_3$) a surprising symmetry involving the Weyl group of type $E_6$ is discovered. Finally, the bridge towards mathematical physics is crossed again with a concrete application to the study of the missing label of $sl_3$, for which our algebraic treatment allowed us to discover new symmetries.

\vskip .1cm
The last shorter chapter is entirely devoted to Yokonuma--Hecke algebras and to the construction of link invariants. As always, we start with ``pure'' algebra: we give an isomorphism theorem relating the Yokonuma--Hecke algebras with the usual Hecke algebras. Then this is applied to the algebraic classification of Markov traces on the family of Yokonuma--Hecke algebra. In this chapter, the bridge we cross to conclude is the one towards knots and links theory, and we apply our algebraic theorems to provide an explicit formula for the topological invariants coming from the Yokonuma--Hecke algebra in terms of the well-known HOMFLYPT invariant.

\paragraph{Acknowledgements.} I warmly thank the referees for their time and their report on this work, and the member of the jury for the defense. Many thanks to my colleagues from the University of Reims for all the support and my coauthors for all the interesting scientific interactions. It is also a pleasure to warmly thank the mathematical physics team of the CRM in Montreal for the welcoming atmosphere while a part of this document was written.

\setcounter{equation}{0}
\DeactivateToc
\chapter{\huge{Preliminaries: Uses of centralisers}}\label{chap-prel}
\ActivateToc
\addcontentsline{toc}{chapter}{\large{Chapter \thechapter. \hspace{0.2cm}Preliminaries: Uses of centralisers\vspace{0.3cm}}}

\setcounter{minitocdepth}{1}
\minitoc

\section{The general philosophy of centralisers}\label{subsec-cent}

A historical reference is \cite{We1}, see also the influential book \cite{We2} or a more modern account \cite{GW}. We note that the theory of centralisers goes beyond the semisimple situation that we discuss here for simplicity (see for example \cite[\S 3]{BDK}, \cite[\S 4]{GJ} and \cite{Gr,JS,MMcA}).

\paragraph{Centraliser of one matrix.} The first example of a centraliser algebra occurs very early in linear algebra. We take a square matrix $M$. Its centraliser is the following subalgebra of the algebra of matrices:
\[Z(M)=\{\text{Matrices $N$ s.t.}\ [M,N]=0\}\ .\]
The basic example already showing some features of centraliser theory is when $M$ is a diagonal matrix. Then we have:
\[
M=\left(\begin{array}{ccccccc}
\!\!\lambda_1\!\!\!\! & & & & & & \\[-0.4em]
& \ddots\!\!\!\! & \\[-0.4em]
& & \lambda_1\!\!\!\! & \\[-0.2em]
& & & \lambda_2\!\!\!\! & \\[-0.4em]
& & & & \ddots \!\!\!\!& & \\[-0.4em]
& & & & & \lambda_2\!\!\!\! \\[-0.4em]
& & & & & & \ddots\!\!\!\!
\end{array}\right)\ \ \ \Rightarrow\ \ \ Z(M)=\{\left(\begin{array}{ccccccc}
 & & & & & & \\[-0.4em]
& \!\!\!\!A_1\!\!\!\! & \\[-0.4em]
& &  & \\[-0.2em]
& & &  & \\[-0.4em]
& & & & \!\!\!\!A_2\!\!\!\! & & \\[-0.4em]
& & & & & \\[-0.4em]
& & & & & & \ddots
\end{array}\right)\}\ ,\]
where, if the multiplicity of the eigenvalue $\lambda_i$ is $m_i$, then the matrices in $Z(M)$ are bock-diagonal matrices and the size of $A_i$ is $m_i$. So we can say somehow  that the centraliser $Z(M)$ ``knows'' about the multiplicities $m_i$.

\paragraph{centraliser of a set of matrices.}
Going to a more involved example, take a set of square matrices $S=\{M_1,\dots,M_k\}$ of the same size. Naturally, its centraliser is:
\[
Z(S)=\{\text{Matrices $N$ s.t.}\ [M_i,N]=0\,,\ \forall i\}\ .
\]
Say for example $S$ is a commutative set of diagonalizable matrices. Using the language of quantum mechanics, we say that $S$ is complete if $S$ generates a maximal commutative subalgebra. Physically, this means that the common eigenspaces of elements in $S$ are all one-dimensional and thus vectors forming a basis of the underlying vector space can be uniquely indexed (up to scalars) by eigenvalues of the elements of $S$. Now we have that $S$ is complete if and only if $Z(S)$ is commutative (and in this case, $Z(S)$ is in fact the commutative algebra generated by $S$). So we can say that the centraliser $Z(S)$ knows if the set is complete, and moreover, by construction, if $S$ is not complete, it is from $Z(S)$ that we must take elements to complete it.

\paragraph{centraliser of a representation.} Finally, we get to the more general situation in representation theory. Take an algebra $\mathcal{A}$ and let $\rho\ :\ \mathcal{A}\to\text{End}(E)$ be a representation on a finite-dimensional complex vector space $E$. The centraliser of the representation is simply the centraliser in $\text{End}(E)$ of the family of endomorphisms $\{\rho(a)\}_{a\in \mathcal{A}}$\,:
\[\text{End}_\mathcal{A}(E)=\ \text{centraliser of $\{\rho(a)\}_{a\in \mathcal{A}}$}=\{x\in\text{End}(E)\ |\ \rho(a)x=x\rho(a)\,,\ \forall a\in \mathcal{A}\}\ .\]
Here again, the centraliser ``knows'' many things about the representation. A first example is Schur's lemma, which says that: 
\[\text{$\rho$ is irreducible\ \ \ $\Leftrightarrow$\ \ \ $\text{End}_\mathcal{A}(E)=\{\lambda \text{Id}_E\}_{\lambda\in\mathbb{C}}$\ .}\]
More generally, assume that the representation decomposes as a direct sum of irreducible representations:
\[E\cong m_1 A_1\oplus m_2 A_2\oplus ......\oplus m_k A_k\,,\]
where $m_i\in\mathbb{Z}_{>0}$ are the multiplicities of the irreducible representations $A_i$ in $E$. Using Schur Lemma, the centraliser $\text{End}_\mathcal{A}(E)$ is described as follows. Each isotypic component $A_i^{\oplus m_i}$ is stable, and one can choose a basis (concatenating a basis of $A_i$) such that, on $A_i^{\oplus m_i}$, the operators of $\mathcal{A}$ and of the centralisers look as follows:
\[\mathcal{A}\ni a\mapsto\left(\begin{array}{cccc}\rho_i(a) & 0 & \dots & 0\\
0 & \rho_i(a) & \ddots & \vdots\\
\vdots & \ddots & \ddots & 0\\
0 & \ldots & 0 & \rho_i(a)\\\end{array}\right)\,,\ \ \ \ \ \ \text{End}_\mathcal{A}(E)\ni x\mapsto\left(\begin{array}{cccc}x_{11}\text{Id} & x_{12}\text{Id} & \dots & x_{1m_i}\text{Id}\\
x_{21}\text{Id} & x_{22}\text{Id} & \ddots & \vdots\\
\vdots & \ddots & \ddots & \vdots\\
x_{m_i1}\text{Id} & \ldots & \ldots & x_{m_im_i}\text{Id}\\\end{array}\right)\,,\]
where we are showing block-matrices with blocks of size $m_i$. Thus we see that the centraliser $\text{End}_\mathcal{A}(E)$ is isomorphic to a direct sum of matrix algebras:
\[\text{End}_\mathcal{A}(E)\cong \text{Mat}_{m_1}(\mathbb{C})\oplus\text{Mat}_{m_2}(\mathbb{C})\oplus ......\oplus \text{Mat}_{m_k}(\mathbb{C})\ ,\]
the sizes of which correspond to the multiplicities. Thus the centraliser ``knows'' about the multiplicities.

\paragraph{Two points of view.} Setting $\mathcal{B}:=\text{End}_\mathcal{A}(E)$, one can go further and note that $E$ can also be seen as a representation of $\mathcal{B}$. Picking the coefficients $x_{ij}$ above, we build a representation of $\mathcal{B}$ which is of dimension $m_i$ (the multiplicity of $A_i$) and which appears in $E$ a number of times equal to the dimension of $A_i$. We have the two following points of view now.

\vspace{0.2cm}
\begin{minipage}{.45\textwidth}
\underline{From the point of view of $\mathcal{A}$:}
\[\rho\cong \dim(B_1)A_1\oplus ......\oplus \dim(B_k)A_k\]
\end{minipage}
\vrule{} \ \ 
\begin{minipage}{.5\textwidth}
\underline{From the point of view of $\mathcal{B}$:}
\[\rho\cong \dim(A_1)B_1\oplus ......\oplus \dim(A_k)B_k\]
\end{minipage}
\vspace{0.2cm}

This leads to the slogan:
\[
\text{``Multiplicities for $\mathcal{A}$\ $=$\ Dimensions for $\mathcal{B}$\ \ \ (and vice-versa)''.}
\]

Finally, as the images commute, $E$ can also be seen as a representation of the algebra $\mathcal{A}\otimes \mathcal{B}$.
To summarise, we have the following decompositions of $E$ as a representation, respectively, of $\mathcal{A}$, of $\mathcal{B}$ and finally of $\mathcal{A}\otimes \mathcal{B}$:
\begin{equation}\label{centraliser-general}
E=\bigoplus_i A_i^{\oplus \dim(B_i)}\,,\quad\ \ \ E=\bigoplus_i B_i^{\oplus \dim(A_i)}\,,\quad\ \ \ E=\bigoplus_i A_i\otimes B_i\,.
\end{equation}

\paragraph{Labelling vectors.} Returning briefly to the quantum mechanics terminology in the same spirit as above, we ask how to naturally label vectors of the vector space $E$. The last decomposition in (\ref{centraliser-general}) gives a clear answer. First select a maximal commutative set in the algebra $\mathcal{A}$. This set acting on $E$ will allow to distinguish between different representations $A_i$ and $A_j$ (with $i\neq j$) and to distinguish vectors in a given $A_i$. 

However, as soon as one multiplicity $m_i$ is strictly greater than 1, this is not enough to specify uniquely a vector in $E$. This is easy to understand since there is no way that an operator from $\mathcal{A}$ will make a difference between different isomorphic copies of the same irreducible representation $A_i$. In other words, from the point of view of $\mathcal{A}$, the spaces $B_i$ are transparent, they only record the multiplicity through their dimensions. We say sometimes that we have a ``missing label'' for the vectors in $E$ (from the point of view of $\mathcal{A}$).

What is missing is actually operators acting on the spaces $B_{i}$ in order to distinguish between different copies of a given $A_i$. The solution is obvious in our formalism coming from representation theory: the centraliser algebra $\mathcal{B}$ provides such operators, and in turn, the additional labels for the vectors. At the end of the day, and quite naturally, to get a complete set of observables, one needs to concatenate a complete set from $\mathcal{A}$ with a complete set from $\mathcal{B}$. A concrete example is studied in details in Chapter \ref{chap-rac}.

\section{The Schur--Weyl duality}\label{sec-SW}

We illustrate the general theory of the previous section with the celebrated classical Schur--Weyl duality. References are \cite{Ful,FH,GW,We2}. The definition of the standard combinatorial objects (partitions, tableaux, etc.) used below are recalled in the appendix.

\subsection{Representations of $GL_N$ and $S_n$} \label{subsec-GL-Sn}

Take a complex vector space $V$ of finite dimension $N$ and take a non-negative integer $n$. The group $GL(V)$ acts naturally on $V$ and it also acts on $V^{\otimes n}$ via the diagonal action:
\[g\cdot (v_1\otimes v_2\otimes\dots\otimes v_n)=(g\cdot v_1)\otimes (g\cdot v_2)\otimes\dots\otimes (g\cdot v_n)\ \ \ \ \quad\forall g\in GL(V)\ .\]
We also have the symmetric group $S_n$ acting on $V^{\otimes n}$ by permuting the factors:
\[\pi\cdot (v_1\otimes v_2\otimes\dots\otimes v_n)=v_{\pi^{-1}(1)}\otimes v_{\pi^{-1}(2)}\otimes\dots\otimes v_{\pi^{-1}(n)}\ \ \ \ \quad\forall \pi\in S_n\ .\]
The representation $V^{\otimes n}$ decomposes into  a direct sum of irreducible representations, both for $GL(V)$ and for $S_n$. Let us start with $GL(V)$. We have:
\begin{equation}\label{SW-GL}\text{\underline{For $GL(V)$:}}\qquad\ \ \ \ \ V^{\otimes n}\cong \bigoplus _{\substack{\lambda\vdash n\\[0.2em] \ell(\lambda)\leq N}} V_{\lambda}^{\oplus m_{\lambda}}\,,
\end{equation}
The representation $V_{\lambda}$ is the irreducible representation of $GL(V)$ of highest-weight $\lambda$, and the positive integers $m_{\lambda}$ are the multiplicities.

On the other hand, we have:
\begin{equation}\label{SW-Sn}\text{\underline{For $S_n$:}}\qquad\ \ \ \ \ \ V^{\otimes n}\cong \bigoplus _{\substack{\lambda\vdash n\\[0.2em] \ell(\lambda)\leq N}} S_{\lambda}^{\oplus M_{\lambda}}\,,
\end{equation}
The representation $S_{\lambda}$ is the irreducible representation of $S_n$ associated with $\lambda$ (see Appendix), and the positive integers $M_{\lambda}$ are the multiplicities.

\paragraph{Putting $GL(V)$ and $S_n$ together.} We have two decompositions in (\ref{SW-GL}) and (\ref{SW-Sn}) which are very similar, so let us combine them in one formula. We consider the direct product group $GL(V)\times S_n$. Both groups are acting on $V^{\otimes n}$ and this induces an action of the direct product, since the two actions commute.

The last decomposition of $V^{\otimes n}$ we will give is as a representation of the direct product $GL(V)\times S_n$. Recall that an irreducible representation of this direct product is of the form $M\otimes N$ where $M$ is an irreducible representation of $GL(V)$ and $N$ is an irreducible representation of $S_n$. We have that $V^{\otimes n}$ decomposes as a representation of $GL(V)\times S_n$ as follows:
\begin{equation}\label{SW-both}
V^{\otimes n}=\bigoplus_{\substack{\lambda\vdash n\\[0.2em] \ell(\lambda)\leq N}} V_{\lambda}\otimes S_{\lambda}\,.
\end{equation}
Restricting this decomposition to, first, the group $GL(V)$, and second, the group $S_n$, we get back our two decompositions:
\begin{equation*}
V^{\otimes n}=\bigoplus_{\substack{\lambda\vdash n\\[0.2em] \ell(\lambda)\leq N}} (V_{\lambda})^{\oplus \dim(S_{\lambda})}\,,\quad\ \ \ V^{\otimes n}=\bigoplus_{\substack{\lambda\vdash n\\[0.2em] \ell(\lambda)\leq N}} (S_{\lambda})^{\oplus \dim(V_{\lambda})}\,.
\end{equation*}
The duality here is apparent. In particular, the dimensions of the irreducible representations of $S_n$ are the multiplicities of the irreducible representations of $GL(V)$, and vice versa. They are expressed in terms of a certain number of Young tableaux:
\begin{equation}\label{dim-mult}
m_{\lambda}=\dim(S_{\lambda})=|\STab(\lambda)|\ \ \ \ \ \ \text{and}\ \ \ \ \ \ M_{\lambda}=\dim(V_{\lambda})=|\SSTab_N(\lambda)|\ .
\end{equation}
The fruitful idea behind the above construction is that we can use some knowledge about the representations of $S_n$ to study the representations of $GL(V)$, or the other way around. These sorts of ideas have been present already in the work of Schur in the early days of representation theory \cite{Sch}, and were greatly popularized by Weyl \cite{We2}. This situation and its many generalisations are often referred to as the Schur--Weyl duality.

\subsection{The classical Schur--Weyl duality}\label{subsec-SW}

We summarise the previous discussion about $GL(V)$ and $S_n$ in the light of Section \ref{subsec-cent}, and we formulate the main facts in a way which will lay the ground for the generalisations we will consider.

If we want to discuss about centralisers, it is better to be talking about algebras (elements commuting with a certain set of matrices form an algebra). So first we will replace the groups $GL(V)$ and $S_n$ by algebras. We replace the group $S_n$ by its group algebra $\mathbb{C}S_n$ and the group $GL(V)$ by the universal enveloping algebra of its Lie algebra, that we denote $U(gl_N)$ to emphasize the dependence on $N$.

So the situation is that we have two algebras represented on the same vector space, and we can picture it like this:
\[U(gl_N)\ \ \stackrel{\phi}{\longrightarrow}\ \ \ \ \text{End}(V^{\otimes n})\ \ \ \ \stackrel{\rho}{\longleftarrow}\ \ \mathbb{C}S_n\ \]
To be precise, $\rho$ is simply the action of $S_n$ extended linearly to $\mathbb{C}S_n$. The map $\phi$ is the extension to $U(gl_N)$ of the action of $gl_N$ on $V^{\otimes n}$. An element of $gl_N=\text{End}(V)$ acts via its diagonal action on $V^{\otimes n}$, which is the sum of the actions on each factor:
\[g\cdot(v_1\otimes v_2\otimes \dots\otimes v_n)=(g.v_1)\otimes v_2\otimes\dots\otimes v_n+\ldots\ldots+v_1\otimes v_2\otimes \dots\otimes (g.v_n)\ .\]
Now we can state the main statements that form what is called the Schur--Weyl duality. Note that we see now the representation $V_{\lambda}$ of $GL(V)$ as a representation of $U(gl_N)$, and we will indicate from now on the dependence on $N$ in the notation: $V^N_{\lambda}$ (since the same partition indexes a representation of $U(gl_N)$ for various $N$).
\begin{framedtheo}[Schur--Weyl duality]\label{theo-SW1}$\ $
\begin{enumerate}
\item In $\text{End}(V^{\otimes n})$, the images of $U(gl_N)$ and $\mathbb{C}S_n$ are their mutual centralisers:
\begin{equation}\label{SW1}
\text{End}_{U(gl_N)}(V^{\otimes n})=\rho\bigl(\mathbb{C}S_n\bigr)\ \ \ \text{and}\ \ \ \text{End}_{\mathbb{C}S_n}(V^{\otimes n})=\phi\bigl(U(gl_N)\bigr)\ . \tag{\textbf{SW1}}
\end{equation}
\item  The decomposition of $V^{\otimes n}$, as a representation of $U(gl_N)\otimes \mathbb{C}S_n$ is:\begin{equation}\label{SW2}
V^{\otimes n}=\bigoplus_{\substack{\lambda\vdash n\\[0.2em] \ell(\lambda)\leq N}} V_{\lambda}^N\otimes S_{\lambda}\,.\tag{\textbf{SW2}}
\end{equation}
\end{enumerate}
\end{framedtheo}
As we explained before, the first part (\ref{SW1}) of the Schur--Weyl duality implies that there is a decomposition as in (\ref{SW2}). But the main point in (\ref{SW2}), which is an important second part of the Schur--Weyl duality, is to actually identify the subset of the irreducible representations of $GL(V)$ and of $S_n$ which appear in $V^{\otimes n}$. 

\subsubsection{Kernel in the Schur--Weyl duality} Clearly an important role is played above by the centralisers $\text{End}_{U(gl_N)}(V^{\otimes n})$ and $\text{End}_{\mathbb{C}S_n}(V^{\otimes n})$. So to complete our study of the Schur--Weyl duality, it is natural to ask about an explicit algebraic description of these algebras. For example, at this point, to the question: \emph{What is the centraliser of the action of $GL(V)$ on $V^{\otimes n}$?}, we can only answer that it is an homomorphic image of the group algebra $\mathbb{C}S_n$. Or in other words, that it is a quotient of $\mathbb{C}S_n$. So we would like to study the kernel of the map:
\[\mathbb{C}S_n\ \ \stackrel{\rho}{\longrightarrow}\ \ \ \ \text{End}(V^{\otimes n})\ .\]
In fact, the kernel of the above map was actually implicitly described in (\ref{SW2}) since the irreducible representations of $S_n$ appearing in $V^{\otimes n}$ were identified, and thus also those who are lost through the representation $\rho$. We will discuss this in more detail in the next paragraph, our aim being to deduce an algebraic description of the quotient (say, with explicit defining relations).

\begin{rema}
Of course, one can ask as well for the kernel of the map from $U(gl_N)$, and thus for the quotient of $U(gl_N)$ isomorphic to the centraliser of the $S_n$-action on $V^{\otimes n}$. We will not follow this direction, but we indicate that we get well-studied and interesting algebras in this way, called Schur algebras $S(n,N)$. We refer to \cite{Gr} and also to \cite{DG} for an explicit description as a quotient of $U(gl_N)$.
\end{rema}

\paragraph{Generalities on quotients of semisimple algebras.} Let $A$ be a finite-dimensional semisimple algebra over $\CC$. Let $S$ be an indexing set for a complete set of pairwise non-isomorphic irreducible representations of $A$. The Artin--Wedderburn theorem asserts that we have the following isomorphism of algebras:
\begin{equation}\label{AW}
A\cong \bigoplus_{\lambda\in S} \text{End}(S_{\lambda})\ ,
\end{equation}
where $S_{\lambda}$ is a realisation of the irreducible representation corresponding to $\lambda$. The isomorphism is given naturally by sending $a\in A$ to the endomorphism in $\text{End}(S_{\lambda})$ corresponding to the action of $a$ on $S_{\lambda}$.

From the Artin--Wedderburn decomposition (\ref{AW}), one sees immediately that ideals (and equivalently quotients) of $A$ are in correspondence with subsets $S'\subset S$ as follows:
\[I_{S'}=\bigoplus_{\lambda\in S'} \text{End}(V_{\lambda})\ \ \ \ \ \ \text{and}\ \ \ \ \ \ A/I_{S'}\cong\bigoplus_{\lambda\in S\backslash S'} \text{End}(V_{\lambda})\ .\]
Informally speaking, to obtain a quotient of a semisimple algebra $A$, one needs to choose a certain subset of irreducible representations (here $S'$) which have to disappear when going from $A$ to the quotient $A/I_{S'}$. That is a representation-theoretic description of an ideal.

Then we define the minimal central idempotent $E_{\lambda}$ as the element of $A$ corresponding under the Artin--Wedderburn decomposition of $A$ to $\text{Id}_{V_{\lambda}}$, in the component corresponding to $\lambda$, and $0$ in all other components. The set $\{E_{\lambda}\}_{\lambda\in S}$ is a complete set of minimal central orthogonal idempotents of $A$, meaning that they are central, they add up to 1, they satisfy $E_{\lambda}E_{\lambda'}=\delta_{\lambda,\lambda'}E_{\lambda}$ and they cannot be written as the sum of two non-zero central idempotents.

In any representation $W$ of $A$, the action of a minimal central idempotent $E_{\lambda}$ by definition projects onto the isotopic component of $W$ corresponding to $\lambda$. More precisely, if the decomposition of $W$ into irreducibles is
\[W\cong \bigoplus_{\lambda'\in S}V_{\lambda'}^{\oplus m_{\lambda'}}\ ,\]
then the action of $E_{\lambda}$ is the projection onto the summand $V_{\lambda}^{\oplus m_{\lambda}}$, that is:
\[{E_{\lambda}}_{|_{V_{\lambda}}}=\text{Id}_{V_{\lambda}}\ \ \ \ \text{and}\ \ \ \ E_{\lambda}(V_{\lambda'})=0\ \ \text{if $\lambda'\neq \lambda$\ .}\]
Note that one set of generators of the ideal $I_{S'}$ consists of the elements $E_{\lambda}$ with $\lambda\in S'$.

\paragraph{Back to the Schur--Weyl duality.} Now we are quite ready to describe the quotient of $\mathbb{C}S_n$ appearing in the Schur--Weyl duality. It is immediate from (\ref{SW2}) that the kernel of the map form $\mathbb{C}S_n$ to $\text{End}(V^{\otimes n})$ is the ideal $I_{S'}$, where $S'$ is the subset of partitions of $n$ with strictly more than $N$ rows. 

It is useful to start with the situation $n=N+1$. In this case, there is only one such partition, which is a single column of $N+1$ boxes. This is the sign representation of $S_{N+1}$. Let us describe its minimal central idempotent. This idempotent coincides with the so-called \textbf{antisymmetriser} of $\mathbb{C}S_{N+1}$.

For any $n$, the antisymmetriser of $\mathbb{C}S_n$  is the following element:
\begin{equation}\label{def-P'Sn}
P'_{n}=\frac{1}{n!}\sum_{w\in S_{n}}sgn(w)w .
\end{equation}
It is easy to check that $P'_n$ satisfies $wP'_n=P'_nw=sgn(w)P'_n$ and that $P'_n$ is an idempotent such that the image of the multiplication by $P'_n$ in $\mathbb{C}S_n$ is of dimension 1.

From what we have recalled above, the antisymmetriser $P'_{N+1}$ is a generator of the kernel in the Schur--Weyl duality if $n=N+1$. Now if $n\geq N+1$, we can take $P'_{N+1}$, which is an element of $\mathbb{C}S_{N+1}$ and see it as an element of $\mathbb{C}S_{n}$ thanks to the natural inclusion of $S_{N+1}$ in $S_n$. It turns out that in $\mathbb{C}S_{n}$, the antisymmetriser also generates the ideal corresponding to the subset of partitions of $n$ with strictly more than $N$ rows. This follows from the knowledge of the restriction of representations of $S_n$ to $S_{N+1}$ (see the Bratteli diagram in Appendix). Indeed, one sees that the sign representation of $S_{N+1}$ appears precisely in the restriction of those representations of $S_n$ with strictly more than $N$ rows. So the image of  $P'_{N+1}$ is non-zero exactly in the correct set of representations, and this is enough to imply the statement.

We remark that up to now the description of the centraliser did not explicitly depend on the dimension $N$ of $V$. Of course, this dependence is hidden in the map from $\mathbb{C}S_n$ to $\text{End}(V^{\otimes n})$ and in its kernel. So it is only natural that it appears explicitly now.
\begin{framedtheo}[SW2']\label{SW2b}$\ $
\begin{itemize}
\item If $n\leq N$ then the kernel of the map is $\{0\}$ and the centraliser $\text{End}_{U(gl_N)}(V^{\otimes n})$ is isomorphic to $\mathbb{C}S_n$.
\item If $n>N$ then the kernel of the map is generated in $\mathbb{C}S_n$ by the element $P'_{N+1}$. So the centraliser $\text{End}_{U(gl_N)}(V^{\otimes n})$ is the quotient of $\mathbb{C}S_n$ by the relation $P'_{N+1}=0$.
\end{itemize}  
\end{framedtheo}
In words, to obtain the complete description of the centraliser, we start from the algebra of the symmetric group $\mathbb{C}S_n$ and, if $n>N$, we cancel the antisymmetriser on $N+1$ letters. Explicitly, the centraliser $\text{End}_{U_q(sl_N)}(V^{\otimes n})$ is described as the algebra generated by $s_1,\dots,s_{n-1}$ with defining relations:
\[\begin{array}{ll}
s_is_{i+1}s_i=s_{i+1}s_is_{i+1}\,,\ \ \  & \text{for $i\in\{1,\dots,n-2\}$}\,,\\[0.2em]
s_is_j=s_js_i\,,\ \ \  & \text{for $i,j\in\{1,\dots,n-1\}$ such that $|i-j|>1$}\,,\\[0.2em]
s_i^2=1\,,\ \ \  & \text{for $i\in\{1,\dots,n-1\}$}\,,\\[0.2em]
P'_{N+1}=0 & \text{if $n>N$.} 
\end{array}
\]

\begin{exam}[Temperley--Lieb algebra]\label{exa-TL}
Let $N=2$. In this case, the centraliser $\text{End}_{U(sl_2)}(V^{\otimes n})$ is called the Temperley--Lieb algebra. The additional relation $P'_3=0$ reads:
\[1-s_1-s_2+s_1s_2+s_2s_1-s_1s_2s_1=0\ .\]
One can show that it implies the same relation with indices $i,i+1$ for all $i=1,\dots,n-2$. Then setting $t_i:=s_i-1$, one recovers the other standard presentation of the Temperley--Lieb algebra:
\[t_i^2=-2t_i\,,\ \ \ \ t_it_{i+1}t_i=t_i\,,\ \ \ \ t_{i+1}t_{i}t_{i+1}=t_{i+1}\ \ \ \text{and}\ \ \ t_it_j=t_jt_i\ \text{if $|i-j|>1$.}\]
\end{exam}

\subsection{Extensions of the Schur--Weyl duality}

In view of its fundamental importance in representation theory, the Schur--Weyl duality has numerous generalisations. For example  one could ask for more general representations of $GL(V)$ than $V$ itself. This will be the subject of one chapter of this thesis. 

In this section we shall only evoke the first historical extension of the Schur--Weyl duality from $GL(V)$ to the classical groups $O(V)$ and $SP(V)$, and a more recent one involving the symmetric group, but in a different fashion than before. In fact we will give a brief description of the following table:
\begin{center}
\begin{tabular}{c|c}
Groups acting diagonally on $V^{\otimes n}$ & Their centralisers\\
\hline & \\[-0.2em]
$GL(V)$ & Symmetric group\\
$\cup$ & $\cap$\\
$O(V)$ and $SP(V)$ & Brauer algebra\\
$\cup$ & $\cap$\\
$S_N$ where $N=\dim(V)$ & Partition algebra
\end{tabular}
\end{center}
On the left, the inclusion of $S_N$ is in the orthogonal group $O(V)$, when we see the elements of $S_N$ as permutation matrices.

\subsubsection{The classical groups and the Brauer algebra}

\paragraph{Classical groups.} For the sake of simplicity and explicitness, we will work with matrix groups. We take a vector space $V$ of finite dimension $N$ and we fix a basis $v_1,\dots,v_N$, so that elements of $\text{End}(V)$ are identified with $N\times N$ matrices. We take an invertible $N\times N$ matrix $G=(g_{ij})$ which is either symmetric or antisymmetric, and we see it as the Gram matrix of a non-degenerate bilinear form: $<v,w>_G={}^tvGw$, that is, $<v_i,v_j>_G=g_{ij}$. The bilinear form is thus either symmetric or alternating (and if it is alternating, then necessary $N$ is even).

We will consider the subgroup $\Gamma_G$ of $GL(V)$ consisting of invertible matrices $M$ such that ${}^tMGM=G$. Of course, this is equivalent to ask that $M$ preserves the bilinear form: $<Mv,Mw>_G=<v,w>_G$. The subgroup $\Gamma_G$ acts on the tensor space $V^{\otimes n}$ with the diagonal action, which is the restriction to $\Gamma_G$ of the diagonal action of $GL(V)$ studied earlier.

As can be expected in general, the subgroup $\Gamma_G$ has more invariant subspaces in $V^{\otimes n}$ than $GL(V)$. In particular, in $V\otimes V$, the main novelty is that the vector
\[v_0:=\sum_{i,j}\overline{g}_{ij}v_i\otimes v_j\ ,\]
is now invariant by $\Gamma_G$, where the inverse of $G$ is $\overline{G}=(\overline{g}_{ij})$. This allows to build an operator of $\text{End}(V\otimes V)$ commuting with $\Gamma_G$:
\[E\ :\ v\otimes w\mapsto <v,w>_G v_0\ .\]
This is a new element in the centraliser of $\Gamma_G$ compared with the situation of $GL(V)$, where we had only the permutation operator $P(v\otimes w)=w\otimes v$. It is easy to check that the elements $P$ and $E$ satisfy the relations:
\[P^2=\text{Id}_{V\otimes V}\,,\quad\ \ E^2=\pm N E\ \ \quad\text{and}\quad\ \ EP=PE=\pm E\ ,\]
where the signs $\pm$ are the same as in $^{t}G=\pm G$.

\paragraph{The Brauer algebra and the centraliser of $O(N)$, $SP(N)$.}
The first main result on the centralisers of the classical groups is that the various natural embeddings of $P$ and $E$ in $\text{End}(V^{\otimes n})$ generate the whole centraliser. We first need an algebra encoding the relations of the elements $P$ and $E$ above. So we fix a non-zero parameter $\delta\in\mathbb{C}$ and an integer $n\geq 0$. We define the Brauer algebra $Br_n(\delta)$ as the algebra generated by elements:
\[s_1,\dots,s_{n-1},\ e_1,\dots,e_{n-1}\,,\]
with the defining relations:
\[\begin{array}{c}
s_i^2=1\,,\ \ \ \ e_i^2=\delta e_i\,,\ \ \ \ e_is_i=s_ie_i=e_i\,,\ \ \ \quad i=1,\dots,n-1,\\[0.5em]
s_is_j=s_js_i\,,\ \ \ \ e_ie_j=e_je_i\,,\ \ \ \ \ s_ie_j=e_js_i\,,\ \ \ \ \quad \text{if $|i-j|>1$,}\\[0.5em]
s_is_{i+1}s_i=s_{i+1}s_is_{i+1}\,,\ \ \ \ e_is_{i+1}s_i=s_{i+1}s_ie_{i+1}\,,\ \ \ \ e_{i+1}s_is_{i+1}=s_is_{i+1}e_i\,,\\[0.5em]
e_is_{i+1}e_i=e_i\,,\ \ \ \ \ \ \quad i=1,\dots,n-2.
\end{array}\]
Consider first the case where the matrix $G$ is symmetric. We take $G$ to be the identity matrix of size $N$, and we denote the resulting group $\Gamma_G$ by $O(N)$. We consider the Brauer algebra $Br_n(N)$, for $\delta=N$, and we build a representation of it on $V^{\otimes n}$:
\[\rho_{O(N)}\ :\ \begin{array}{c}
Br_n(N)\ \to\ \text{End}(V^{\otimes n})\\[0.5em]
s_i\mapsto P_i\ \ \ \text{and}\ \ \ e_i\mapsto E_i
\end{array}\]
where $P_i$ (respectively, $E_i$) is the operator $P$ (respectively, $E$) acting on copies $i$ and $i+1$ of $V$ in $V^{\otimes n}$.

Then consider the case where the matrix $G$ is antisymmetric, which implies that $N$ is even. We take $G$ to be the matrix $\left(\begin{array}{cc} 0 & \text{-Id}\\ \text{Id} & 0\end{array}\right)$, where the blocks are of size $N/2$. We denote the resulting group $\Gamma_G$ by $SP(N)$. We consider the Brauer algebra $Br_n(-N)$, for $\delta=-N$, and we build a representation of it on $V^{\otimes n}$:
\[\rho_{SP(N)}\ :\ \begin{array}{c}
Br_n(-N)\ \to\ \text{End}(V^{\otimes n})\\[0.5em]
s_i\mapsto -P_i\ \ \ \text{and}\ \ \ e_i\mapsto E_i
\end{array}\]
where $P_i$ and $E_i$ are as above. 

The analogue of the first part of the Schur--Weyl duality, for the classical groups, is that in both cases the image of the Brauer algebra is the whole centraliser of the diagonal action of the group on $V^{\otimes n}$:
\[\rho_{O(N)}(Br_n(N))=\text{End}_{O(N)}(V^{\otimes n})\ \ \quad\text{and}\quad\ \ \rho_{SP(N)}(Br_n(-N))=\text{End}_{SP(N)}(V^{\otimes n})\ .\]
In both cases ($O(N)$ and $SP(N)$) there is also an analogue of the second part of the Schur--Weyl duality, namely a decomposition of $V^{\otimes n}$ as in (\ref{SW2}), involving the irreducible representations of the Brauer algebra. We do not give the details and refer for example to \cite{Wen}.

\paragraph{Diagrammatic description.} Consider the following set of diagrams (Brauer diagrams). We take two rows of $n$ dots, one below the other, and we draw edges connecting the dots two by two (there are $n$ edges). Note that this is a combinatorial object, where the only information that matters is which dot is connected to which other dot. In other words, such a diagram leads to a partition of the set $\{1,\dots,2n\}$ (we number the dots by $1,\dots,2n$)  into $n$ $2$-element subsets, and two diagrams are seen as equivalent if the resulting partitions are the same.

Consider the vector space of linear combinations of Brauer diagrams, so that basis elements are identified with Brauer diagrams. A multiplication is defined on the Brauer diagrams and extended linearly. Take $d$ and $d'$ two diagrams and perform the following procedure:
\begin{itemize}
\item We place the diagram of $d$ on top of the diagram of $d'$ by identifying the bottom row of dots of $d$ with the top row of dots of $d'$. We denote $k(d,d')$ the number of closed loops in the resulting diagram.
\item Then we forget the dots in the middle row and the product $dd'$ is defined to be $\delta^{k(d,d')}$ times the resulting diagram without the closed loops.
\end{itemize}
This defines an associative unital algebra on the vector space of linear combinations of Brauer diagrams which turns out to be isomorphic to $Br_n(\delta)$.

For example, considering only the diagrams where the edges always connect the top row to the bottom row, we recover the symmetric group $S_n$, and restricting to these permutation diagrams, the multiplication is the multiplication of $S_n$ (composition of permutations). The algebra $Br_n(\delta)$ is bigger, there are actually $(2n-1)!!=1\times 3\times\dots\times (2n-1)$ different Brauer diagrams ($2n-1$ choices for connecting the first dot, then $2n-3$ choices for connecting the next dots, and so on).

The generators in the algebraic presentation of $Br_n(\delta)$ above correspond to the following diagrams:
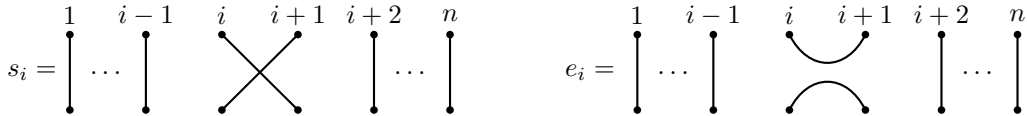
\[ \begin{tikzpicture}[scale=0.25]
\node at (0,0) {$s_i=$};
\node at (2,3) {$1$};\fill (2,2) circle (0.2);\fill (2,-2) circle (0.2);
\draw[thick] (2,2) -- (2,-2);
\node at (4,0) {$\dots$};
\draw[thick] (6,2) -- (6,-2);
\node at (6,3) {$i-1$};\fill (6,2) circle (0.2);\fill (6,-2) circle (0.2);
\node at (10,3) {$i$};\fill (10,2) circle (0.2);\fill (10,-2) circle (0.2);
\node at (14,3) {$i+1$};\fill (14,2) circle (0.2);\fill (14,-2) circle (0.2);
\draw[thick] (14,2) -- (10,-2);
\draw[thick] (10,2) -- (14,-2);
\draw[thick] (18,2) -- (18,-2);
\node at (18,3) {$i+2$};\fill (18,2) circle (0.2);\fill (18,-2) circle (0.2);
\node at (20,0) {$\dots$};
\draw[thick] (22,2) -- (22,-2);\fill (22,2) circle (0.2);\fill (22,-2) circle (0.2);
\node at (22,3) {$n$};
\end{tikzpicture}\ \ \ \ \ \ \ \ \ \ 
 \begin{tikzpicture}[scale=0.25]
\node at (-0.5,0) {$e_i=$};
\node at (2,3) {$1$};\fill (2,2) circle (0.2);\fill (2,-2) circle (0.2);
\draw[thick] (2,2) -- (2,-2);
\node at (4,0) {$\dots$};
\draw[thick] (6,2) -- (6,-2);
\node at (6,3) {$i-1$};\fill (6,2) circle (0.2);\fill (6,-2) circle (0.2);
\node at (10,3) {$i$};\fill (10,2) circle (0.2);\fill (10,-2) circle (0.2);
\node at (14,3) {$i+1$};\fill (14,2) circle (0.2);\fill (14,-2) circle (0.2);
\draw[thick] (10,2)..controls (11,0) and (13,0) .. (14,2);
\fill[white] (12,0) circle (0.4);
\draw[thick] (10,-2)..controls (11,0) and (13,0) .. (14,-2);
\draw[thick] (18,2) -- (18,-2);
\node at (18,3) {$i+2$};\fill (18,2) circle (0.2);\fill (18,-2) circle (0.2);
\node at (20,0) {$\dots$};
\draw[thick] (22,2) -- (22,-2);\fill (22,2) circle (0.2);\fill (22,-2) circle (0.2);
\node at (22,3) {$n$};
\end{tikzpicture}\]

\paragraph{Remarks and references.} The chain of Brauer algebras $\{Br_n(\delta)\}_{n\geq 0}$ has a very interesting representation theory. When semisimple, it is summarised in its Bratteli diagram, which is obtained by the ``up-down'' process from the Young diagram (the Bratteli diagram of the symmetric group found in the appendix): to go from $n-1$ to $n$, one can add a box to a partition or remove a box from a partition. The simple modules of the Brauer algebras $Br_n(\delta)$ are indexed by the partitions of $k$, for all $k\leq n$ such that $n-k$ is even. Some references for the definitions, the structure and the semisimple representation theory of the Brauer algebras are \cite{Br,Bro,HW,IMO,Na1,Ra,Wen}.

The Brauer algebras $Br_n(\delta)$ are semisimple for generic $\delta$, but not for particular specialisations of $\delta$. Nonsemisimplicity occurs when $\delta$ is an integer with $|\delta|<n$. Recall that the Brauer--Schur--Weyl duality for a classical group $\Gamma$ is expressed as:
\[\text{End}_{\Gamma}(V^{\otimes n})=\rho\bigl(Br_n(\pm N)\bigr)\,,\]
where $N=\dim(V)$ and $\rho\bigl(Br_n(\pm N)\bigr)$ is a certain homomorphic image of $Br_n(\pm N)$. Thus, a somewhat intriguing feature is that the centralisers of classical groups (intrinsically semisimple) are in general expressed as quotients of a non semisimple algebra. The representation theory of the Brauer algebras outside of the semisimple situation turns out to show some interesting properties, see for example \cite{CDVM1,CDVM2,Mar5}. 

\subsubsection{The symmetric group and the partition algebra}

To conclude this section, we briefly indicate another generalisation of the Schur--Weyl duality. This time, we consider the symmetric group $S_N$ acting in its natural permutation representation. Namely, we take a vector space $V$ with given basis $v_1,\dots,v_N$, and $S_N$ acts on $V$ by permuting the vectors $v_1,\dots,v_N$ (not to be confused with the action of $S_n$ on $V^{\otimes n}$ by permuting the copies of $V$). Of course permutation matrices are orthogonal so $S_N$ acts on $V$ through a subgroup of the orthogonal group $O(N)$. The centraliser of $O(N)$ was described in the preceding section in terms of the Brauer algebra, so we expect to find here an algebra containing and generalising the Brauer algebra.

The algebra that appears in this setting is the partition algebra $P_n(\delta)$. Instead of the Brauer diagrams (that we recall corresponded to partitions of $2n$ points into $2$-element subsets), we consider any possible partition of $2n$ points. We draw a corresponding diagram by connecting the points in the same subset of the given partition. The multiplication on this set of diagrams is defined as for the Brauer diagrams: first we concatenate the two diagrams $d$ and $d'$ and we count how many connected components we have involving only middle dots, say $k(d,d')$. Then we forget the dots in the middle row and the product $dd'$ is defined to be $\delta^{k(d,d')}$ times the resulting diagram without the connected components in the middle.

Clearly the partition algebra $P_n(\delta)$ contains the Brauer algebra $Br(\delta)$ which contains the symmetric group algebra $\mathbb{C}S_n$. This is the chain of inclusions mentioned in the table at the beginning of this section.

The action of the partition algebra $P_n(\delta)$ on a tensor space $V^{\otimes n}$ is very simple. Denote the dots in the top row from left to right by $1,\dots,n$ and similarly for the bottom row with $1',\dots,n'$. A diagram $d$ corresponds to a partition of $\{1,\dots,n,1',\dots,n'\}$. We set:
\[(d)^{i_1,\dots,i_n}_{i_{1'},\dots,i_{n'}}=\left\{\begin{array}{ll} 1 & \text{if $i_a=i_b$ whenever $a$ and $b$ are in the same block,}\\
0 & \text{otherwise}\end{array}\right.\]
This defines an operator on $V^{\otimes n}$ by:
\[\rho(d)\cdot v_{i_1}\otimes \dots\otimes v_{i_n}=\sum_{i_{1'},\dots,i_{n'}}(d)^{i_1,\dots,i_n}_{i_{1'},\dots,i_{n'}}v_{i_{1'}}\otimes \dots\otimes v_{i_{n'}}\,,\]
where all indices run from $1$ to $N=\dim(V)$. One can check that this generalises the action of $Br(N)$ on $V^{\otimes n}$ discussed for the centraliser of $O(N)$. The analogue of the Schur--Weyl duality here is:
\[\text{End}_{S_N}(V^{\otimes n})=\rho\bigl(P_n(\pm N)\bigr)\,,\]
or in words, the centraliser of the diagonal action of $S_N$ on $V^{\otimes n}$ is the image of the partition algebra.

\paragraph{Remarks and references.} The partition algebra was introduced in the context of the Potts model in statistical physics. The reader is referred to \cite{HR,Mar1,Mar2,Mar3,Mar4,MR,Jo4} and references therein. The partition algebra is also used for the study of the Kronecker coefficients for the symmetric group \cite{BDVO}.

\section{The Braid Group and the Yang--Baxter Equation}\label{sec-BrYB}

We move on to a (seemingly) different subject: the braid group and the Yang--Baxter equation. The connections with what precedes will appear through the examples in the last subsection.

\subsection{The Braid group}

For precise definitions and more details on the braid group, we refer to \cite{Bir,BZ,Kam,KT}. Here is an example of a braid (with $5$ strands):
\vskip .4cm
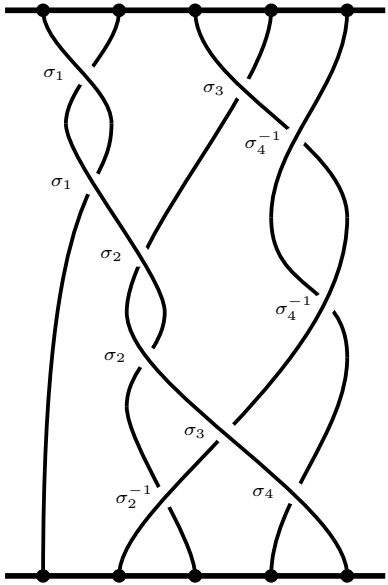
\begin{minipage}{0.3\linewidth}
\begin{center}
\begin{tikzpicture}[scale=0.5]
\draw[line width=0.75mm] (1,0)--(11,0);
\fill[black] (2,0) circle (0.18);\fill[black] (4,0) circle (0.18);
\fill[black] (6,0) circle (0.18);\fill[black] (8,0) circle (0.18);
\fill[black] (10,0) circle (0.18);
\draw[line width=0.5mm] (4,0)..controls +(0,-1) and +(0,+1) .. (2.6,-3);
\node at (2.3,-1.7) {$\scriptstyle{\si_1}$};\fill[white] (3.1,-1.7) circle (0.3);
\draw[line width=0.5mm] (2,0)..controls +(0,-1) and +(0,+1) .. (3.8,-3);
\draw[line width=0.5mm] (3.8,-3)..controls +(0,-2) and +(0,11) .. (2,-15);
\node at (2.5,-4.6) {$\scriptstyle{\si_1}$};\fill[white] (3.3,-4.6) circle (0.3);
\draw[line width=0.5mm] (8,0)..controls +(0,-2) and +(0,+2) .. (4.2,-8);
\node at (3.8,-6.5) {$\scriptstyle{\si_2}$};\fill[white] (4.5,-6.5) circle (0.3);
\draw[line width=0.5mm] (2.6,-3)..controls +(0,-1) and +(0,+1) .. (5.2,-8);
\node at (6.5,-2.1) {$\scriptstyle{\si_3}$};\fill[white] (7.3,-2.1) circle (0.3);
\draw[line width=0.5mm] (6,0)..controls +(0,-2) and +(0,+2) .. (10,-5.5);
\node at (7.8,-3.5) {$\scriptstyle{\si_4^{-1}}$};\fill[white] (8.6,-3.4) circle (0.3);
\draw[line width=0.5mm] (10,0)..controls +(0,-2) and +(0,+2) .. (8,-5.5);
\draw[line width=0.5mm] (5.2,-8)..controls +(0,-1) and +(0,+1) .. (4.2,-10.5);
\node at (3.9,-9.2) {$\scriptstyle{\si_2}$};\fill[white] (4.7,-9.2) circle (0.3);
\draw[line width=0.5mm] (4.2,-10.5)..controls +(0,-1) and +(0,+1) .. (6,-15);
\node at (4.4,-12.9) {$\scriptstyle{\si_2^{-1}}$};\fill[white] (5.2,-12.9) circle (0.3);
\draw[line width=0.5mm] (8,-5.5)..controls +(0,-2) and +(0,+2) .. (10,-9.2);
\node at (8.6,-7.9) {$\scriptstyle{\si_4^{-1}}$};\fill[white] (9.4,-7.8) circle (0.3);
\draw[line width=0.5mm] (10,-9.2)..controls +(0,-2) and +(0,+2) .. (8,-15);
\node at (7.8,-12.8) {$\scriptstyle{\si_4}$};\fill[white] (8.6,-12.8) circle (0.3);
\draw[line width=0.5mm] (10,-5.5)..controls +(0,-4) and +(0,+2) .. (4,-15);
\node at (6,-11.2) {$\scriptstyle{\si_3}$};\fill[white] (6.8,-11.2) circle (0.3);
\draw[line width=0.5mm] (4.2,-8)..controls +(0,-2) and +(0,+2) .. (10,-15);
\draw[line width=0.75mm] (1,-15)--(11,-15);
\fill[black] (2,-15) circle (0.18);\fill[black] (4,-15) circle (0.18);
\fill[black] (6,-15) circle (0.18);\fill[black] (8,-15) circle (0.18);
\fill[black] (10,-15) circle (0.18);
\end{tikzpicture}
\end{center}
\end{minipage}
\hspace{2cm}
\begin{minipage}{0.5\linewidth}
We fix two horizontal parallel lines each having $5$ fixed points, and we connect bijectively each point on the top line to a point on the bottom line by a ``strand''. \\[0.4em]
This is an object in a three-dimensional space, so the strands can pass ``over'' or ``under'' other strands. The labels $\sigma_i^{\pm1}$ next to the crossings indicate the standard generators of the braid group.\\[0.4em]
Braids are considered up to isotopy, meaning that we can move continuously the strands while leaving their end points fixed, and this is still the same braid.
\end{minipage}

\vskip .4cm
As shown by the example, a braid with $n$ strands is drawn in a rectangular strip with a top line of $n$ fixed dots and a bottom line of $n$ fixed dots. We connect bijectively each top dot to a bottom dot by a strand inside the strip. A strand is a continuous line going from a top point to a bottom point. The vertical coordinate is always decreasing along the line (no strand is allowed to go back towards the top). 

At each point of the strip at most two strands are crossing each other, and at each crossing, we indicate which strand pass over the other one. We call a crossing  positive (resp. negative) when the strand coming from the left passes over (resp. under) the strand coming from the right. We assume that there are only a finite number of crossings. Such diagram is called a braid with $n$ strands and braids are considered up to isotopy (continuous moves of the strands with fixed end points).

An associative multiplication is defined on the set $B_n$ of all braids with $n$ strands, the multiplication being simply the vertical concatenation of diagrams. If $\alpha,\beta\in B_n$, to perform the product $\alpha\beta$, we place the diagram of $\alpha$ on top of the diagram of $\beta$, we identify the bottom dots of $\alpha$ with the top dots of $\beta$, thereby connecting the strands, and we delete these middle dots.

With this multiplication, the set $B_n$ of braids with $n$ strands forms a group. The identity element $1$ of $B_n$ is the braid where the $n$ strands are vertical and parallel. The following elements are called elementary braidings:
\begin{center}
 \begin{tikzpicture}[scale=0.3]
\node at (0,0) {$\si_i=$};
\node at (2,3) {$1$};\fill (2,2) circle (0.2);\fill (2,-2) circle (0.2);
\draw[thick] (2,2) -- (2,-2);
\node at (4,0) {$\dots$};
\draw[thick] (6,2) -- (6,-2);
\node at (6,3) {$i-1$};\fill (6,2) circle (0.2);\fill (6,-2) circle (0.2);
\node at (10,3) {$i$};\fill (10,2) circle (0.2);\fill (10,-2) circle (0.2);
\node at (14,3) {$i+1$};\fill (14,2) circle (0.2);\fill (14,-2) circle (0.2);
\draw[thick] (14,2)..controls +(0,-2) and +(0,+2) .. (10,-2);
\fill[white] (12,0) circle (0.4);
\draw[thick] (10,2)..controls +(0,-2) and +(0,+2) .. (14,-2);
\draw[thick] (18,2) -- (18,-2);
\node at (18,3) {$i+2$};\fill (18,2) circle (0.2);\fill (18,-2) circle (0.2);
\node at (20,0) {$\dots$};
\draw[thick] (22,2) -- (22,-2);\fill (22,2) circle (0.2);\fill (22,-2) circle (0.2);
\node at (22,3) {$n$};
\node at (38,0) {$i\in\{1,\dots,n-1\}$.};
\end{tikzpicture}
\end{center}
The inverse $\si_i^{-1}$ is the same diagram as $\si_i$ except that a negative crossing is used instead of a positive one. The elementary braidings $\si_1,\dots,\si_{n-1}$ are generators of the braid group $B_n$. The following relations are satisfied:
\[\si_i\si_j=\si_j\si_i\ \ \ \ \ \text{if $|i-j|>1$,}\] 
since if $|i-j|>1$, then the pairs of strands $i,i+1$ and $j,j+1$ can be manipulated independently. It is also immediate that the following relations are satisfied:
\vskip .1cm
\begin{minipage}{0.3\linewidth}
\begin{center}
\begin{tikzpicture}[scale=0.3]
\node at (0,3) {\scriptsize{$i$}};
\node at (1.8,3) {\scriptsize{$i\!+\!1$}};
\node at (4,3) {\scriptsize{$i\!+\!2$}};
\draw[thick] (4,2)..controls +(0,-4) and +(0,+4) .. (0,-6);
\draw[thick] (2,2)..controls +(0,-3) and +(0,+3) .. (0,-2);
\fill[white] (0.6,-0.2) circle (0.3);
\fill[white] (2,-2) circle (0.3);
\draw[thick] (0,2)..controls +(0,-4) and +(0,+4) .. (4,-6);
\fill[white] (0.6,-3.8) circle (0.3);
\draw[thick] (0,-2)..controls +(0,-3) and +(0,+3) .. (2,-6);
\node at (6,-2) {$=$};
\node at (8,3) {\scriptsize{$i$}};
\node at (9.8,3) {\scriptsize{$i\!+\!1$}};
\node at (12,3) {\scriptsize{$i\!+\!2$}};
\draw[thick] (12,2)..controls +(0,-4) and +(0,+4) .. (8,-6);
\fill[white] (11.4,-0.2) circle (0.3);
\fill[white] (10,-2) circle (0.3);
\draw[thick] (10,2)..controls +(0,-3) and +(0,+3) .. (12,-2);
\draw[thick] (12,-2)..controls +(0,-3) and +(0,+3) .. (10,-6);
\fill[white] (11.4,-3.8) circle (0.3);
\draw[thick] (8,2)..controls +(0,-4) and +(0,+4) .. (12,-6);
\end{tikzpicture}
\end{center}
\end{minipage}
\hspace{0.5cm}
\begin{minipage}{0.6\linewidth}
Indeed, one simply has to move the middle strand. In algebraic terms, we have:
\[\si_i\si_{i+1}\si_i=\si_{i+1}\si_i\si_{i+1}\,,\ \ \ \ \text{for $i=1,\dots,n-2$.}\]
\end{minipage}
\vskip .1cm
\vskip .2cm
It turns out that only the two sorts of relations above are enough to completely characterise the braid group. This is the well known following theorem.
\begin{theo}[Artin \cite{Ar1,Ar2,Bo}]\label{theo-artin}
The braid group $B_n$ is generated by $\si_1,\dots,\si_{n-1}$ with defining relations:
\begin{equation}\label{rel-Artin}
\begin{array}{ll}
\si_i\si_{i+1}\si_i=\si_{i+1}\si_i\si_{i+1}\,,\ \ \  & \text{for $i\in\{1,\dots,n-2\}$}\,,\\[0.2em]
\si_i\si_j=\si_j\si_i\,,\ \ \  & \text{for $i,j\in\{1,\dots,n-1\}$ such that $|i-j|>1$}\,.
\end{array}
\end{equation} 
\end{theo}

\paragraph{Local representations of $B_n$.} 
Let $V$ be a vector space and form the tensor product $V^{\otimes n}$ of $n$ copies of $V$. We take an invertible element $R\in\text{End}(V\otimes V)$ and we define elements $R_1,\dots,R_{n-1}$ like this:
\[\text{Operators $R_1,\dots,R_{n-1}$ on $V^{\otimes n}$:}\qquad\ \ \ \rlap{$\overbrace{\phantom{V\otimes V}}^{R_1}$}V\otimes\overbrace{V\otimes V}^{R_2}\otimes\dots\overbrace{\ldots\otimes V}^{R_{n-1}}\]
More precisely, for $i=1,\dots,n-1$, the operator $R_{i}$ acts on $V^{\otimes n}$ non-trivially only in the copies $i,i+1$ of $V$, where it acts by $R$. That is, it is defined by $R_{i}=\text{Id}_{V^{\otimes i-1}}\otimes R\otimes\text{Id}_{V^{\otimes n-i-1}}$.

We say that $R$ defines a local representation of the braid group $B_n$, or that the following map
\[\begin{array}{lcrcl}\rho\ &:&\ B_n & \to & \text{End}(V^{\otimes n})\\[0.5em]
 & & \si_i & \mapsto & R_i \end{array}\]
is a local representation of $B_n$ on the tensor space $V^{\otimes n}$, if the braid relation is satisfied
\begin{equation}\label{braid}
R_1R_{2}R_1=R_{2}R_1R_{2}\ \quad \ \ \ \ \text{on $V\otimes V\otimes V$.} 
\end{equation}
This immediately implies all the other defining relations of $B_n$, thereby making the map $\rho$ a representation of the group $B_n$. A local representation is thus a representation of $B_n$ on a tensor space $V^{\otimes n}$ with the special (local) form explained above for the action of the generators $\si_i$. It depends only on the invertible element $R\in\text{End}(V\otimes V)$ satisfying (\ref{braid}). 

\begin{exam}
The permutation operator $P\in\text{End}(V\otimes V)$, sending $x\otimes y$ to $y\otimes x$, provides the first example of a local representation of the braid group. In fact, it factors through the representation $S_n$ of the symmetric group on $V^{\otimes n}$.
\end{exam}

\subsection{The Yang--Baxter equation}\label{subsec-YB}

Let $V$ be a finite-dimensional vector space. We introduce some more notations for operators on $V^{\otimes n}$. If $S\in \text{End}(V\otimes V)$, we use the notation $S_{ij}$ for the operator on $V^{\otimes n}$ acting as $S$ on copies $i$ and $j$ and trivially otherwise. We consider a function 
\[\begin{array}{lcrcl}R\ &:&\ \bC^2 & \to & \text{End}(V\otimes V)\\[0.5em]
 & & (u,v) & \mapsto & R(u,v)\ . \end{array}\]
The Yang--Baxter equation (YB equation for short) is a functional equation for such a function $R$:
\begin{equation}\label{YBprime}
R_{12}(u_1,u_2)R_{13}(u_1,u_3)R_{23}(u_2,u_3)=R_{23}(u_2,u_3)R_{13}(u_1,u_3)R_{12}(u_1,u_2)\ \quad \ \ \ \ \text{on $V\otimes V\otimes V$.} 
\end{equation}
In this context, the variables $u_1,u_2,u_3$ are often called \emph{spectral parameters}. In our perspective, the so-called braided version of this equation is more relevant. From $R$ define another function:
\[\begin{array}{lcrcl}\check R\ &:&\ \bC^2 & \to & \text{End}(V\otimes V)\\[0.5em]
 & & (u,v) & \mapsto & PR(u,v) \end{array}\ ,\] 
 where $P$ is the permutation operator sending $x\otimes y$ to $y\otimes x$. To see better the connections with the braid relation, set $\check R_{1}(u,v):=\check R_{12}(u,v)$ and $\check R_{2}(u,v):=\check R_{23}(u,v)$. Then an easy manipulation shows that the YB equation is equivalent to the following equation only involving the function $\check R$:
\begin{equation}\label{YB}
\check R_{1}(u_1,u_2)\check R_{2}(u_1,u_3)\check R_{1}(u_2,u_3)=\check R_{2}(u_2,u_3)\check R_{1}(u_1,u_3)\check R_{2}(u_1,u_2)\ \quad \ \ \ \ \text{on $V\otimes V\otimes V$.} 
\end{equation}
When precision is needed, we will refer to this version as the braided YB equation.

It happens in many situations that the function $\check R$ (or equivalently $R$) depends only on the ratio $u/v$. In this particular case (set $u=u_1/u_2$ and $v=u_2/u_3$), the braided YB equation becomes:
\begin{equation}\label{YBmult}
\check R_1(u)\check R_{2}(uv)\check R_1(v)=\check R_{2}(v)\check R_1(uv)\check R_{2}(u)\ \quad \ \ \ \ \text{on $V\otimes V\otimes V$.} 
\end{equation}
Instead of a multiplicative version, we sometimes also have an additive version, where the function $\check R$ depends only on the difference $u-v$:
\begin{equation}\label{YBadd}
\check R_1(u)\check R_{2}(u+v)\check R_1(v)=\check R_{2}(v)\check R_1(u+v)\check R_{2}(u)\ \quad \ \ \ \ \text{on $V\otimes V\otimes V$.} 
\end{equation}

\begin{exam}[\textbf{constant solution, Yang solution}]\label{ex-solYB}$\ $\\
$\bullet$ A constant solution (that is, not depending on the spectral parameters) of the braided YB equation is equivalent to a local representation of the braid group.

\noindent $\bullet$ One can check by hand that the following function satisfies (\ref{YBprime}):
\[R(u,v)=\text{Id}_{V\otimes V}+\frac{P}{u-v}\,,\]
or equivalently that the braided YB equation (\ref{YB}) is satisfied by:
\[\check R(u,v)=PR(u,v)=P+\frac{\text{Id}_{V\otimes V}}{u-v}\,.\] 
This solution is called Yang solution.
\end{exam}

\paragraph{References and a diagram.} The Yang--Baxter equation is certainly one of the fundamental equations in theoretical and mathematical physics. Some references among the vast literature on this subject are \cite{AYP,Fad,FST,GRS,Ji-int,Ji-ed,JM,KBI,PAY}. Historical references are \cite{McG,Ons,Yang1,Yang2,ZZ} and the standard book \cite{Bax}. Though the YB equation belongs to the theoretical and mathematical side of Physics, it is interesting to note that connections with experiments can be investigated today, see for example \cite{Bat,BF}.

In one of its interpretation in physics, the YB equation is represented by the following picture: 
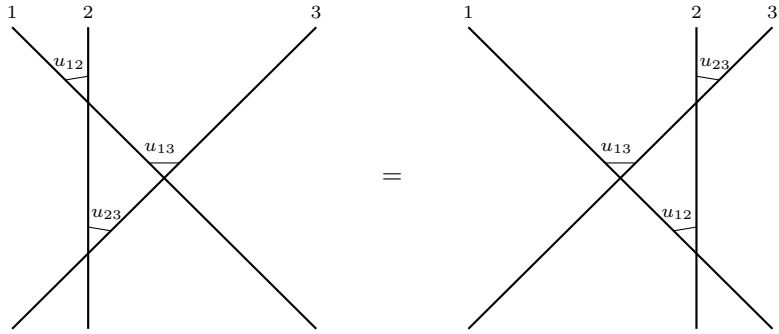
\begin{center}
\begin{tikzpicture}[scale=0.5]
\node[above] at (0,2) {\scriptsize{$1$}};
\node[above] at (2,2) {\scriptsize{$2$}};
\node[above] at (8,2) {\scriptsize{$3$}};
\draw[thick] (0,2)--(8,-6);
\draw[thick] (2,2)--(2,-6);
\draw[thick] (8,2)--(0,-6);
\draw (1.4,0.6)--(2,0.7);
\node[above] at (1.5,0.65) {\scriptsize{$u_{12}$}};
\draw (3.6,-1.6)--(4.4,-1.6);
\node[above] at (3.9,-1.6) {\scriptsize{$u_{13}$}};
\draw (2,-3.3)--(2.6,-3.4);
\node[above] at (2.5,-3.35) {\scriptsize{$u_{23}$}};
\node at (10,-2) {$=$};
\node[above] at (12,2) {\scriptsize{$1$}};
\node[above] at (18,2) {\scriptsize{$2$}};
\node[above] at (20,2) {\scriptsize{$3$}};
\draw[thick] (12,2)--(20,-6);
\draw[thick] (18,2)--(18,-6);
\draw[thick] (20,2)--(12,-6);
\draw (18.6,0.6)--(18,0.7);
\node[above] at (18.5,0.65) {\scriptsize{$u_{23}$}};
\draw (15.6,-1.6)--(16.4,-1.6);
\node[above] at (15.9,-1.6) {\scriptsize{$u_{13}$}};
\draw (18,-3.3)--(17.4,-3.4);
\node[above] at (17.5,-3.35) {\scriptsize{$u_{12}$}};
\end{tikzpicture}
\end{center}
where $u_{ij}$ are some physical parameters. Here, the YB equation is a compatibility condition for the factorization of the $n$-body interaction in terms of the two-body interaction. This originates in the papers \cite{McG,Yang1,Yang2,ZZ}. Note that the YB equation only involves three particles. This interpretation of the YB equation is very similar to the following classical fact about the braid group and the symmetric group. The longest element of the symmetric group (reversing the order of $1,\dots,n$) can be factorised in different ways as a minimal-length product of simple transpositions, and all these different factorisations are seen to be equal just by assuming the braid relation (which involves only three strands). The YB equation above plays the same role as the braid equation.

\subsection{Motivating examples}\label{subsec-ex}

\subsubsection{The symmetric group}

To start to unfold the connections between this section and the preceding one, the first example should be the symmetric group. A presentation by generators and relations of the symmetric group $S_n$ is with generators $s_1,\dots,s_{n-1}$ satisfying the defining relations:
\[\begin{array}{ll}
s_is_{i+1}s_i=s_{i+1}s_is_{i+1}\,,\ \ \  & \text{for $i\in\{1,\dots,n-2\}$}\,,\\[0.2em]
s_is_j=s_js_i\,,\ \ \  & \text{for $i,j\in\{1,\dots,n-1\}$ such that $|i-j|>1$}\,,\\[0.2em]
s_i^2=1\,,\ \ \  & \text{for $i\in\{1,\dots,n-1\}$}\,.
\end{array}
\]
In this presentation, the generators $s_i$ are identified with the elementary transposition $(i,i+1)$. This presentation makes it apparent that the symmetric group $S_n$ is a quotient of the braid group $B_n$, since we add the defining relations $s_i^2=1$ to the Artin presentation of the braid group (Theorem \ref{theo-artin}). 

It turns out that the symmetric group has also something to do with the YB equation, as can be seen in Example \ref{ex-solYB}. To formalise it, we consider the following functions:
\begin{equation}\label{Baxt-Sn}
s_i(u):=s_i+\frac{1}{u}\,,\ \ \ \ \ \text{where $i=1,\dots,n-1$}.
\end{equation}
These are functions taking values in the group algebra $\mathbb{C}S_n$. An easy direct verification shows that the YB equation is satisfied:
\[
s_i(u)s_{i+1}(u+v)s_i(v)=s_{i+1}(v)s_i(u+v)s_{i+1}(u)\ . 
\]
So we have a solution of the YB equation inside the group algebra $\mathbb{C}S_n$. However it is not a \emph{genuine} solution of the YB equation since it is not yet an operator on a tensor product of vector spaces. It will give a matrix solution each time we choose a local representation of $S_n$. To be more precise, say we have a representation of $S_n$ of the form:
\[s_i\mapsto \check R_{i}=\text{Id}_{V^{\otimes i-1}}\otimes \check R\otimes\text{Id}_{V^{\otimes n-i-1}}\,,\] 
for some operator $\check R\in\text{End}(V\otimes V)$. Then we only have to apply the formula (\ref{Baxt-Sn}) inside the representation, which means that we set:
\[\check R(u)=\check R+\frac{\text{Id}_{V\otimes V}}{u}\ ,\]
in order to obtain solutions of the YB equation on $V\otimes V$. Of course, for any vector space $V$, we have the natural permutation representation of $S_n$ on $V^{\otimes n}$. This simply amounts to take $\check R=P_{V\otimes V}$. Thus we recover the Yang solution as coming from the ``abstract solution'' in $\mathbb{C}S_n$ applied in the tensor representation of $S_n$.

A formula such as (\ref{Baxt-Sn}) is often called a Baxterisation formula. It builds  a solution of the YB equation with spectral parameters (here $s_i(u)$) out of a constant solution (here $s_i$) satisfying a certain algebra (here $\mathbb{C}S_n$). We will give more examples of Baxterisation formulas in this section.

\subsubsection{The Hecke algebra}

\paragraph{Definition.} We give the definition directly related to the braid group. The idea is as follows. The braid group is a fairly complicated group algebraically since we can make an arbitrary number of crossings between strands. So the naive idea is to add an algebraic relation allowing to deal with all these crossings. The simplest way to do this is to add the relation $\si_i^2=1$, yielding the symmetric group. This is too simple, and the next to simplest way is the one resulting in the Hecke algebra.

More precisely, we add to the definition of the braid group the following relation:
\begin{equation}\label{rel-skeinH}
 \begin{tikzpicture}[scale=0.25]
\draw[thick] (0,2)..controls +(0,-2) and +(0,+2) .. (4,-2);
\fill[white] (2,0) circle (0.4);
\draw[thick] (4,2)..controls +(0,-2) and +(0,+2) .. (0,-2);
\node at (6,0) {$=$};
\draw[thick] (12,2)..controls +(0,-2) and +(0,+2) .. (8,-2);
\fill[white] (10,0) circle (0.4);
\draw[thick] (8,2)..controls +(0,-2) and +(0,+2) .. (12,-2);
\node at (17,0) {$-\,(q-q^{-1})$};
\draw[thick] (21,2) -- (21,-2);\draw[thick] (25,2) -- (25,-2);
\end{tikzpicture}
\end{equation}
where $q$ is a parameter, which can be seen as an indeterminate or as a non-zero complex number. This relation is a local relation, meaning that it is applied only to a small part of a braid (around a crossing) leaving the rest of the braid unmodified. 

Algebraically, the Hecke algebra $H_n(q)$ is defined as the algebra generated by elements $\si_1,\dots,\si_{n-1}$ with defining relations:
\begin{equation}\label{rel-H}
\begin{array}{ll}
\si_i\si_{i+1}\si_i=\si_{i+1}\si_i\si_{i+1}\,,\ \ \  & \text{for $i\in\{1,\dots,n-2\}$}\,,\\[0.2em]
\si_i\si_j=\si_j\si_i\,,\ \ \  & \text{for $i,j\in\{1,\dots,n-1\}$ such that $|i-j|>1$}\,,\\[0.2em]
\si_i^2=1+(q-q^{-1})\si_i\,,\ \ \  & \text{for $i\in\{1,\dots,n-1\}$}\,.
\end{array}
\end{equation} 
It is an algebra over $\bC$ if $q$ is a complex number, and can also be defined as an algebra over $\mathbb{Z}[q,q^{-1}]$ if $q$ is an indeterminate. The conventional choice of normalisation of the generators was made such that the characteristic equation of degree 2 is:
\[(\si_i-q)(\si_i+q^{-1})=0\ .\]
The most general characteristic equation of order 2 would be with two independent eigenvalues, but if we impose $\si_i$ to be invertible (and we do), then we can renormalise $\si_i$ such that the product of the eigenvalues is $-1$, as above (note that renormalising $\si_i$ does not change the braid relations).

If $q^2=1$ the relations (\ref{rel-H}) are defining relations for the symmetric group $S_n$, and the Hecke algebra $H_n(1)$ is therefore the group algebra $\bC S_n$. The Hecke algebra is thus a deformation of the group algebra of the symmetric group $S_n$.

\paragraph{The YB equation.} For any $i\in\{1,\dots,n-1\}$, define the following function taking values in the algebra $H_n(q)$:
\begin{equation}\label{Baxt-H}
\si_i(u)=\si_i+(q-q^{-1})\frac{1}{u-1}\ .
\end{equation}
Then a straightforward calculation using the relations in $H_n(q)$ shows that:
\[\si_i(u)\si_{i+1}(uv)\si_i(v)=\si_{i+1}(v)\si_i(uv)\si_{i+1}(u)\ ,\]
namely, the braided YB equation is satisfied in $H_n(q)$ by these functions. So we have a solution of the braided YB equation inside the Hecke algebra $H_n(q)$, given by a rather simple formula. Reproducing  the same discussion as for the symmetric group, after the Baxterisation formula (\ref{Baxt-Sn}), we see that we need local representations of the Hecke algebra to obtain genuine matrix solutions of the YB equation.

It turns out that, for any vector space $V$, we can construct a local representation of $H_n(q)$ on $V^{\otimes n}$, as a generalisation of the tensor representation of $S_n$. More precisely, fix a basis $(e_1,\dots,e_N)$ of $V$ and define a linear operator on $V\otimes V$ by:
\begin{equation}\label{rep-Hn}\check R(e_a\otimes e_b):=\left\{\begin{array}{ll}
q\,e_a\otimes e_b\ \  & \text{if $a=b$,}\\[0.8em]
e_b\otimes e_a+(q-q^{-1})\,e_a\otimes e_b & \text{if $a<b$,}\\[0.4em]
e_b\otimes e_a & \text{if $a>b$.}
\end{array}\right.\ \ \ \text{where $a,b=1,\dots,N$.}
\end{equation}
One can check by hand that 
\[\si_i\mapsto \check R_{i}=\text{Id}_{V^{\otimes i-1}}\otimes \check R\otimes\text{Id}_{V^{\otimes n-i-1}}\,,\]
is a representation of the Hecke algebra. So finally, combining the Baxterisation formula (\ref{Baxt-H}) of $H_n(q)$ with this representation, we get another solution of the YB equation for any vector space $V$.
\begin{exam}\label{ex-solYBq}
If $\dim(V)=2$, in the lexicographic ordering of the basis of $V\otimes V$, we find:
\[\check R=\left(\begin{array}{cccc}
q & \cdot & \cdot & \cdot \\
\cdot  & q-q^{-1} & 1& \cdot \\
\cdot  & 1 & 0& \cdot \\
\cdot & \cdot & \cdot & q
\end{array}\right)\ ,\ \ \ \ \ \ 
\check R(u)=\left(\begin{array}{cccc}
\displaystyle\frac{qu-q^{-1}}{u-1} & \cdot & \cdot & \cdot \\
\cdot   & \displaystyle\frac{(q-q^{-1})u}{u-1}& 1& \cdot \\
\cdot  & 1 & \displaystyle\frac{q-q^{-1}}{u-1} & \cdot \\
\cdot & \cdot & \cdot & \displaystyle\frac{qu-q^{-1}}{u-1}
\end{array}\right)\ .\]
The solution $\check R(u)$ coming from the Hecke algebra is a deformation of the Yang solution $\check R(\alpha)=P+\frac{\text{Id}}{\alpha}$ from Example \ref{ex-solYB} in the following sense: set $u=q^{2\alpha}$ and take the limit $q\to 1$.
\end{exam}

\subsubsection{The Brauer algebra and the Birman--Murakami--Wenzl algebra}

The Birman--Murakami--Wenzl (BMW) algebra $BMW_n(a,q)$, originating from \cite{BW,Mu1}, is the quotient of the braid group algebra $\bC B_n$ by the relations:
\[
\begin{array}{cl}
\si_i\varepsilon_i =  a \varepsilon_i & \text{for $i=1,\dots, n-1$,}\\[0.4em]
\varepsilon_i \si_{i+1} \varepsilon_i = a^{-1} \varepsilon_{i}\quad\text{and}\quad \varepsilon_i \si_{i+1}^{-1} \varepsilon_i = a \varepsilon_{i} & \text{for $i=1,\dots, n-2$,}
\end{array}\]
where we have set 
$$\varepsilon_i = 1-\frac{\si_i -\si_i^{-1}}{q - q^{-1}}\ .$$ 
Here $q$ and $a$ are two non-zero complex numbers and $q^2\neq 1$. Or alternatively, we work with indeterminates $q$ and $a$ over the ring $\mathbb{Z}[q^{\pm1},a^{\pm1},(q-q^{-1})^{-1}]$.

The first relation, when written only in terms of $\si_i$, is:
\[(\si_i-q)(\si_i+q^{-1})(\si_i-a)=0\ ,\]
so that we have a cubic characteristic equation for the generators. It is clear from the algebraic presentation that if we impose morever $\varepsilon_i=0$ then we recover the Hecke algebra $H_n(q)$.

\begin{rema}
The BMW algebra has also a topological definition, using the notion of tangles instead of braids, and a certain notion of regular isotopy. In this presentation, the generators are the usual elementary braidings, while the elements $\varepsilon_i$ correspond to the elementary tangle \begin{tikzpicture}[scale=0.1]
\fill (10,2) circle (0.2);\fill (10,-2) circle (0.2);
\fill (14,2) circle (0.2);\fill (14,-2) circle (0.2);
\draw[thick] (10,2)..controls (11,0) and (13,0) .. (14,2);
\fill[white] (12,0) circle (0.4);
\draw[thick] (10,-2)..controls (11,0) and (13,0) .. (14,-2);
\end{tikzpicture}
connecting the $i$-th and $i+1$-st dots. The relation between $\si_i$ and $\varepsilon_i$ amounts to the Kauffman bracket relation:
\begin{equation}\label{rel-skeinBMW}
 \begin{tikzpicture}[scale=0.25]
\draw[thick] (0,2)..controls +(0,-2) and +(0,+2) .. (4,-2);
\fill[white] (2,0) circle (0.4);
\draw[thick] (4,2)..controls +(0,-2) and +(0,+2) .. (0,-2);
\node at (6,0) {$-$};
\draw[thick] (12,2)..controls +(0,-2) and +(0,+2) .. (8,-2);
\fill[white] (10,0) circle (0.4);
\draw[thick] (8,2)..controls +(0,-2) and +(0,+2) .. (12,-2);
\node at (17,0) {$=-(q-q^{-1})$};
\node at (22,0) {$\Bigl($};
\draw[thick] (23.5,2) -- (23.5,-2);\draw[thick] (27.5,2) -- (27.5,-2);
\node at (29.5,0) {$-$};
\draw[thick] (31.5,2)..controls (32.5,0) and (34.5,0).. (35.5,2);
\draw[thick] (31.5,-2)..controls (32.5,0) and (34.5,0).. (35.5,-2);
\node at (37,0) {$\Bigr)$};
\end{tikzpicture}
\end{equation}
which is a local relation, allowing in some sense to ``resolve'' the crossings, in a way similar to the Hecke algebra situation, though more involved. The Hecke algebra leads to the Jones/HOMFLYPT polynomial and the BMW algebra leads to the Kauffman polynomial \cite{Kau}. We refer to \cite{BW,Mor,MT} for more details.
\end{rema}

\paragraph{The YB equation.} For any $i\in\{1,\dots,n-1\}$, define the following function taking values in the algebra $BMW_n(q,a)$:
\begin{equation}\label{Baxt-BMW}
\si_i(u)=\si_i+(q-q^{-1})\frac{1}{u-1}+(q-q^{-1})\frac{1}{a^{-1}qu+1}\varepsilon_i\ .
\end{equation}
It turns out that the braided YB equation is satisfied in the BMW algebra by these functions \cite{Jo2}:
\[\si_i(u)\si_{i+1}(uv)\si_i(v)=\si_{i+1}(v)\si_i(uv)\si_{i+1}(u)\ .\]
Following the same steps as for the symmetric group and the Hecke algebra above, it remains to discuss whether we can find some local representation of the BMW algebra, in order to obtain genuine matrix solutions of the YB equation. Again, it turns out to be possible. We do not give the details here and refer to \cite{Mu2,Tur} for explicit formulas. In both cases (Hecke and BMW), the existence of a local representation has a far-reaching significance in the context of Schur--Weyl dualities for quantum groups. We will discuss this later in Section \ref{sec-QG}.

\paragraph{Connection with the Brauer algebra.} The BMW algebra is not defined above for a value of $q$ such that $q^2=1$. The issue was the division by $(q-q^{-1})$ in the definition of $\epsilon_i$. So we will give an alternative presentation of the BMW algebra, and this will show explicitly the connections with the Brauer algebra. First we have to reinforce $\varepsilon_i$ as generators, so the generators are now $\si_i$ and $\varepsilon_i$, for $i=1,\dots,n-1$. We take as defning relations:
\[\begin{array}{c}
\displaystyle\si_i-\si_i^{-1}=(q-q^{-1})(1-\varepsilon_i)\,,\ \ \ \ \varepsilon_i^2=(1-\frac{a-a^{-1}}{q-q^{-1}})\varepsilon_i\,,\ \ \ \ \varepsilon_i\si_i=\si_i\varepsilon_i=a\varepsilon_i\,,\ \ \ \quad i=1,\dots,n-1,\\[0.5em]
\si_i\si_j=\si_j\si_i\,,\ \ \ \ \varepsilon_i\varepsilon_j=\varepsilon_j\varepsilon_i\,,\ \ \ \ \ \si_i\varepsilon_j=\varepsilon_j\si_i\,,\ \ \ \ \quad \text{if $|i-j|>1$,}\\[0.5em]
\si_i\si_{i+1}\si_i=\si_{i+1}\si_i\si_{i+1}\,,\ \ \ \ \varepsilon_i\si_{i+1}\si_i=\si_{i+1}\si_i\varepsilon_{i+1}\,,\ \ \ \ \varepsilon_{i+1}\si_i\si_{i+1}=\si_i\si_{i+1}\varepsilon_i\,,\\[0.5em]
\varepsilon_i \si_{i+1} \varepsilon_i = a^{-1} \varepsilon_{i}\quad\text{and}\quad \varepsilon_i \si_{i+1}^{-1} \varepsilon_i = a \varepsilon_{i}\,,\ \ \ \ \ \ \quad i=1,\dots,n-2.
\end{array}\]
The slight nuance compared to before is that the algebra is defined as long as $\frac{a-a^{-1}}{q-q^{-1}}$ is defined. If $(q-q^{-1})$ is invertible, then the first relation allows to express $\varepsilon_i$ in terms of $\si_i$ and $\si_i^{-1}$, and all the relations we have added compared to before are easily implied. So if $q^2\neq 1$ then  this is the same algebra as before.

Now consider the specialisation where $a=\pm q^{\pm (1-\delta)}$ (any choice of sign will do). We see immediately that a value of $q$ such that $q^2=1$ is allowed, and moreover, that for such values of $a$ and $q$, one recovers the Brauer algebra. In this sense, the BMW algebra is a deformation of the Brauer algebra.

The specialisation from the BMW algebra to the Brauer algebra allows to obtain from the formula (\ref{Baxt-BMW}) a corresponding formula for the Brauer algebra. Replacing the spectral parameter $u$ by $q^{2u}$, one finds that:
\[s_i(u)=s_i+\frac{1}{u}-\frac{1}{u+1-\delta/2}e_i\ ,\]
satisfies the YB equation with additive spectral parameters (to be precise, one must use the specialisation with $a=-q^{\delta-1}$; the other one recovers the Yang solution $s_i+1/u$).

\subsubsection{More quotients of the algebra of the braid group?}

To pursue the list of examples, it should seem natural to try to consider other quotients of the algebra of the braid group. It turns out to be very difficult to find some by explicit generators and relations (and reasonable to study). The world of cyclotomic Hecke algebras associated to complex reflection groups provides only a small number of possibilities, see \cite{Ma1,Ch}. Excluding the usual Hecke algebras and the commutative algebras generated by a single element ($n=2$), we only get 5 other algebras and, for example, the number of strands $n$ has to be smaller than 5. Even the question of the most generic finite-dimensional quotient of the braid group algebra including a cubic characteristic equation for the generators is difficult, see for example \cite{Ma2}.

\section{Quantum groups and their centralisers}\label{sec-QG}

In the preceding section, in our discussion of the braid group and the YB equation, the algebras that were our main examples appeared to be some deformations of the algebras playing a role in the classical Schur--Weyl dualities of Section \ref{sec-SW}. This suggests that maybe they are also centralisers of certain objects which should be deformations of the classical groups. This is where quantum groups enter the picture.

The quantum groups $U_q(\mathfrak{g})$ for simple Lie algebras $\mathfrak{g}$ were introduced in the mid 80s by Drinfeld and Jimbo \cite{Dr-QG,Ji-QG}. Some constructions in special cases were already found in \cite{KR,Skl}. General references on quantum groups are \cite{CP,Dri,Kas,KS}. The quantum groups $U_q(\mathfrak{g})$ are examples of quasi-triangular Hopf algebras. We will not give the definition of quasi-triangular Hopf algebras, nor of the quantum groups and their affine versions. So this part will be quite sketchy and is only meant to indicate that there exists a general and interesting algebraic setting which provides us with representations of the braid group and even more, with solutions of the Yang--Baxter equation. Our primary goal is to motivate the study of the centralisers from the point of view of the braid group and of the YB equation. The examples of $U_q(sl_N)$ will be treated with more explicit details.

\subsection{centralisers of tensor products}\label{subsec-QG-braid}

Let $\mathfrak{g}$ be a simple Lie algebra. To get an algebraic feeling of what we are talking about, the quantum group that we denote $U_q(\mathfrak{g})$ is an associative algebra, which can be defined explicitly by generators and relations. As a vector space, it looks like the universal enveloping algebra $U(\mathfrak{g})$ of the Lie algebra $\mathfrak{g}$. However, the way of multiplying elements is different. It is a deformation of the multiplication in $U(\mathfrak{g})$, in the sense that there is a certain way to send the parameter $q$ to $1$ which recovers the algebra $U(\mathfrak{g})$. 

The properties we are going to discuss are\footnote{The existence of a trivial representation and of contragredient (or dual) representations are some of the important properties of quantum groups that we omit from our discussion.}:
\begin{itemize}
\item the fact that we can make tensor products of representations of $U_q(\mathfrak{g})$;
\item the fact that $V\otimes W$ and $W\otimes V$ are isomorphic as representations of $U_q(\mathfrak{g})$, and that the isomorphisms provide local representations of the braid group;
\item and finally, the fact that we can upgrade this picture, using affine quantum groups, to get solutions of the Yang--Baxter equation.
\end{itemize}

Let $(\rho_1,V_1)$ and $(\rho_2,V_2)$ be two representations of $U_q(\mathfrak{g})$. Then $V_1\otimes V_2$ can be made a representation of $U_q(\mathfrak{g})$. Note that the vector space $V_1\otimes V_2$ naturally carries a representation of $U_q(\mathfrak{g})\otimes U_q(\mathfrak{g})$. Then the way to construct a representation of $U_q(\mathfrak{g})$ on $V_1\otimes V_2$ is simply by precomposing with the so-called coproduct $\Delta$:
\[\Delta\ :\ \ U_q(\mathfrak{g})\ \longrightarrow\ U_q(\mathfrak{g})\otimes U_q(\mathfrak{g})\quad\ \ \ \leadsto\quad\ \ \ \rho_{V_1\otimes V_2}=(\rho_1\otimes \rho_2)\circ \Delta\ .\]
Some properties are required (such as the coassociativity) for the coproduct and we refer to any book on Hopf algebras. It might be worth emphasizing that, in general, for an arbitrary associative algebra, there is no natural coproduct (take the Hecke algebra for example) and no natural way of performing tensor product of representations.

\begin{rema}
For a group algebra $\bC G$, the map $\Delta$ from $\bC G$ to $\bC G\otimes \bC G$ defined by $\Delta(g)=g\otimes g$ for any $g\in G$ extends to a morphism of algebras and leads to the standard way of performing tensor products of representations of a group. 

For a Lie algebra $\mathfrak{g}$, the map $\Delta$ from $U(\mathfrak{g})$ to $U(\mathfrak{g})\otimes U(\mathfrak{g})$ defined by $\Delta(x)=x\otimes 1+1\otimes x$ for any $x\in \mathfrak{g}$ extends to a morphism of algebras and leads to the standard way of performing tensor products of representations of a Lie algebra.
\end{rema}

\vskip .2cm
For us, the crucial property of tensor products of representations of $U_q(\mathfrak{g})$ is the following: for any two representations $V,W$ of $U_q(\mathfrak{g})$, the representations $V\otimes W$ and $W\otimes V$ are isomorphic, that is, we have an invertible linear operator $\check R_{V,W}$ between $V\otimes W$ and $W\otimes V$ which commutes with the action of $U_q(\mathfrak{g})$:
\[ \check R_{V,W}\ :\ \ V\otimes W\ \stackrel{\sim}{\to}\ W\otimes V\ ,\ \ \ \ \ \ \ \text{isomorphism of $U_q(\mathfrak{g})$-representations.}\]
One of the main interests of quantum groups is that their coproduct is not cocommutative, which means that the images $\Delta(a)$, with $a\in U_q(\mathfrak{g})$, are not in general invariant under the transposition of the two components. In other words, 
if we denote $\Delta^{op}$ the composition of $\Delta$ with the permutation of $U_q(\mathfrak{g})\otimes U_q(\mathfrak{g})$ ($a\otimes b\mapsto b\otimes a$), then we have $\Delta\neq\Delta^{op}$.

However, in $U_q(\mathfrak{g})$, there is a property weakening the cocommutativity which is satisfied and which ensures the isomorphism of representations $V\otimes W$ and $W\otimes V$. This is the existence of an invertible element $\mathcal{R}$ in $U_q(\mathfrak{g})\otimes U_q(\mathfrak{g})$ such that:
\begin{equation}\label{univ-R}
\mathcal{R}\cdot\Delta(a)=\Delta^{op}(a)\cdot\mathcal{R}\ \ \ \ \forall a\in U_q(\mathfrak{g})\ .
\end{equation}
It follows immediately from this that the following element provides an isomorphism $\check R_{V,W}$ of representations between $V\otimes W$ and $W\otimes V$:
\[\check R_{V,W}=P_{V\otimes W}\cdot \rho_{V\otimes W}(\mathcal{R})\,,\]
where $\rho_{V\otimes W}(\mathcal{R})$ is the image of the element $\mathcal{R}$ in the representation on $V\otimes W$, and $P_{V\otimes W}$ is the permutation operator from $V\otimes W$ to $W\otimes V$.

There is more. The mere existence of the isomorphisms $\check R_{V,W}$ is not enough for our purpose, we need moreover a compatibility condition for these isomorphims. Consider a tensor product of three representations:
\begin{center}
 \begin{tikzpicture}[scale=0.3]
\node at (0,0) {$V_1\otimes V_2\otimes V_3$};
\draw[thick][->] (4,0.1) -- (7,2.1);
\draw[thick][->] (4,-0.1) -- (7,-2.1);
\node at (11,2.1) {$V_2\otimes V_1\otimes V_3$};
\node at (11,-2.1) {$V_1\otimes V_3\otimes V_2$};
\draw[thick][->] (15,2.1) -- (18,2.1);
\draw[thick][->] (15,-2.1) -- (18,-2.1);
\node at (22,2.1) {$V_2\otimes V_3\otimes V_1$};
\node at (22,-2.1) {$V_3\otimes V_1\otimes V_2$};
\draw[thick][->] (26,2.1) -- (29,0.1);
\draw[thick][->] (26,-2.1) -- (29,-0.1);
\node at (33,0) {$V_3\otimes V_2\otimes V_1$};
\end{tikzpicture}
\end{center}
This diagram shows the two possible paths from $V_1\otimes V_2\otimes V_3$ to $V_3\otimes V_2\otimes V_1$ by applying the isomorphisms $\check R_{V_i,V_j}$. The compatibility condition is that these two paths coincide. This condition writes as an equality of linear operator from $V_1\otimes V_2\otimes V_3$ to $V_3\otimes V_2\otimes V_1$:
\begin{equation}\label{comp}
\bigl(\check R_{V_1,V_2}\otimes \text{Id}_{V_3}\bigr)\circ\bigl(\text{Id}_{V_2}\otimes \check R_{V_1,V_3}\bigr)\circ\bigl(\check R_{V_2,V_3}\otimes \text{Id}_{V_1}\bigr)=\bigl(\text{Id}_{V_1}\otimes \check R_{V_2,V_3}\bigr)\circ\bigl(\check R_{V_1,V_3}\otimes \text{Id}_{V_2}\bigr)\circ\bigl(\text{Id}_{V_3}\otimes \check R_{V_1,V_2}\bigr)\ .
\end{equation}
At this point, there is no reason for the compatibility condition (\ref{comp}) to be satisfied for representations of $U_q(\mathfrak{g})$. The fact that it is indeed satisfied follows from the so-called quasitriangularity conditions, which read as:
\begin{equation}\label{quasiT}
(\Delta\otimes\text{Id})(\mathcal{R})=\mathcal{R}_{13}\mathcal{R}_{23}\ \ \ \text{and}\ \ \ (\text{Id}\otimes\Delta)(\mathcal{R})=\mathcal{R}_{13}\mathcal{R}_{12}\ .
\end{equation}
These two conditions imply in particular the constant YB equation $\mathcal{R}_{12}\mathcal{R}_{13}\mathcal{R}_{23}=\mathcal{R}_{23}\mathcal{R}_{13}\mathcal{R}_{12}$ directly in $U_q(\mathfrak{g})\otimes U_q(\mathfrak{g})\otimes U_q(\mathfrak{g})$. From this, the compatibility condition for all representations as in (\ref{comp}) is ensured.

\vskip .2cm
Such an element $\mathcal{R}$ satisfying (\ref{univ-R}) and (\ref{quasiT}) is called a universal $R$-matrix. It turns out that the quantum group $U_q(\mathfrak{g})$ admits a universal $R$-matrix satisfying all the required properties. We note that an appropriate definition of $U_q(\mathfrak{g})$ and $U_q(\mathfrak{g})\otimes U_q(\mathfrak{g})$ is needed, involving formal power series and suitable completions.

\paragraph{The braid relation.} Now if we consider the situation where $V_1=V_2=V_3=V$ and we denote simply by $\check R$ the isomorphism $\check R_{V,V}$ which is an invertible element in $\text{End}(V\otimes V)$, then the compatibility condition above becomes simply the braid relation:
\[
\check R_1\check R_{2}\check R_1=\check R_{2}\check R_1\check R_{2}\ \quad \ \ \ \ \text{on $V\otimes V\otimes V$.} 
\]
The operator $\check R$ on $V\otimes V$ is called the $R$-matrix associated to the representation $V$ of the quantum group $U_q(\mathfrak{g})$. Thus we have elements $\check R_1,\dots,\check R_{n-1}$ in $\text{End}(V^{\otimes n})$, which are constructed from the element $\check R\in\text{End}(V\otimes V)$ like this:
\[\rlap{$\overbrace{\phantom{V\otimes V}}^{\check R_1}$}V\otimes\overbrace{V\otimes V}^{\check R_2}\otimes\dots\overbrace{\ldots\otimes V}^{\check R_{n-1}}\]
and which satisfy the braid relations. By construction, they commute with the action of $U_q(\mathfrak{g})$ so they actually live in the centraliser, and we have our first conclusion of this discussion.
\begin{framedconc}\label{conc-QGbr}
For any representation $V$ of $U_q(\mathfrak{g})$, we have a local representation of the braid group on $V^{\otimes n}$. Moreover, it takes values in the centraliser:
\[\begin{array}{rcl} B_n\ & \to\ & \text{End}_{U_q(\mathfrak{g})}(V^{\otimes n})\\[0.2em]
\sigma_i & \mapsto & \check R_i
\end{array}\]
\end{framedconc}
It is important to note that the image of the braid group does not necessarily generate the whole centraliser $\text{End}_{U_q(\mathfrak{g})}(V^{\otimes n})$. Indeed in general the matrices $\check R_i$ do not generate the centraliser. It is satisfied for any representation of $U_q(sl_2)$, and for the standard (vector) representations of $U_q(\mathfrak{g})$ with $\mathfrak{g}$ a classical Lie algebra or the one of type $G_2$ \cite{LZ1}. In general, the structure of the centralisers is not well known, outside of the usual cases treated by the quantum analogue of the Schur--Weyl duality, to be discusses later.

\paragraph{The YB equation.} The quantum groups were originally designed for being useful in the study of the YB equation, but so far we have only discussed how they are related to local representations of the braid group. In short we must explain how to add the spectral parameters. This requires increasing by one (rather large) step the technical difficulties by going to the so-called affine quantum groups. We will only try to convey the idea, in an informal way, that the existence of solutions of the YB equation comes from the existence of affine quantum groups with good properties. We note that the use of affine quantum groups for the YB equation was present since the origin of the theory of quantum groups \cite{Ji86}. References for this part are \cite{CP1,CP,EFK,FR,GRS,He}, see also references therein.

\vskip .2cm
Fix a representation $V$ of $U_q(\mathfrak{g})$. Let us imagine that we have a ``larger'' algebra $U_q(\hat{\mathfrak{g}})$ containing $U_q(\mathfrak{g})$. Moreover, imagine that we have a family of representations $V(a)$ of $U_q(\hat{\mathfrak{g}})$, depending on a complex parameter $a$, such that the restrictions to $U_q(\mathfrak{g})$ always give the representation $V$ we started with. We will call the parameter $a$ an ``evaluation'' parameter.

Finally, last assumption, imagine that we have similar properties for tensor products of representations of $U_q(\hat{\mathfrak{g}})$ as the ones we have been discussing for $U_q(\mathfrak{g})$. That is, at least for generic values of the evaluation parameters, we can find an operator intertwining the representations $V(a)\otimes V(b)$ and $V(b)\otimes V(a)$. This operator, denoted $\check R(a,b)$, thus lives in $\text{End}(V\otimes V)$ since the vector space underlying both representations is $V$. Then we can see it as a function of the parameters $a,b$:
\[\begin{array}{lcrcl}\check R\ &:&\ \bC^2 & \to & \text{End}(V\otimes V)\\[0.5em]
 & & (a,b) & \mapsto & \check R(a,b) \end{array}\ ,\]
and the compatibility condition applied to $V(a)\otimes V(b)\otimes V(c)$ looks like this:
\begin{center}
 \begin{tikzpicture}[scale=0.3]
\node at (0,0) {$V(a)\otimes V(b)\otimes V(c)$};
\draw[thick][->] (6,0.1) -- (8,2.1);
\draw[thick][->] (6,-0.1) -- (8,-2.1);
\node at (14,2.1) {$V(b)\otimes V(a)\otimes V(c)$};
\node at (14,-2.1) {$V(a)\otimes V(c)\otimes V(b)$};
\draw[thick][->] (20,2.1) -- (22,2.1);
\draw[thick][->] (20,-2.1) -- (22,-2.1);
\node at (28,2.1) {$V(b)\otimes V(c)\otimes V(a)$};
\node at (28,-2.1) {$V(c)\otimes V(a)\otimes V(b)$};
\draw[thick][->] (34,2.1) -- (36,0.1);
\draw[thick][->] (34,-2.1) -- (36,-0.1);
\node at (42,0) {$V(c)\otimes V(b)\otimes V(a)$};
\end{tikzpicture}
\end{center}
For the function $\check R(a,b)$, this reads:
\[
\check R_{12}(b,c)\check R_{23}(a,c)\check R_{12}(a,b)=\check R_{23}(a,b)\check R_{12}(a,c)\check R_{23}(b,c)\ \quad \ \ \ \ \text{on $V\otimes V\otimes V$.} 
\]
This is the braided YB equation. From the point of view of $U_q(\mathfrak{g})$, the dependance on the ``evaluation'' parameter disappears and this equality lives in $V^{\otimes 3}$. Moreover, the intertwining property for $U_q(\hat{\mathfrak{g}})$, restricted to $U_q(\mathfrak{g})$, says that all of this is actually happening in the centraliser. Let us draw our second conclusion at this point.
\begin{framedconc}\label{conc-QGYB}
In such a situation as above, we have a solution $\check R(a,b)$ of the YB equation on $V\otimes V$. Moreover, 
the functions $\check R_1(a,b),\dots,\check R_{n-1}(a,b)$ take values in the centraliser:
\[\begin{array}{rcl} \mathbb{C}^2\ & \to\ & \text{End}_{U_q(\mathfrak{g})}(V^{\otimes n})\\[0.2em]
(a,b) & \mapsto & \check R_i(a,b)
\end{array}\]
\end{framedconc}
Note that in the end the spectral parameters come from a family of representations of $U_q(\hat{\mathfrak{g}})$ extending the given representation of $U_q(\mathfrak{g})$. So, very roughly speaking, the YB equation may be seen as the braid relation applied at the level of $U_q(\hat{\mathfrak{g}})$.

\vskip .2cm
Of course we still need to discuss if such a situation as above exists at least for some representations $V$ of some quantum group $U_q(\mathfrak{g})$. It turns out that the answer is yes. First there is an algebra $U_q(\hat{\mathfrak{g}})$ with good properties. It can be seen as the quantum group associated to the affine Lie algebra $\hat{\mathfrak{g}}$, or alternatively, as a certain affinization of $U_q(\mathfrak{g})$. It has all the nice properties we need about tensor products of representations. Moreover, generalising the natural embedding of $\mathfrak{g}$ in $\hat{\mathfrak{g}}$, one may regard $U_q(\mathfrak{g})$ as a subalgebra of $U_q(\hat{\mathfrak{g}})$.

More subtle is the situation concerning the assumption that we had a family of ``evaluation'' representations $V(a)$ of $U_q(\hat{\mathfrak{g}})$ extending a representation $V$ of $U_q(\mathfrak{g})$. This assumption is not always satisfied and we will discuss briefly its validity.

First, for the simplest situation $\mathfrak{g}=sl_N$, this assumption is valid for any representation $V$ of $U_q(\mathfrak{g})$. It follows from the existence of a morphism from $U_q(\hat{\mathfrak{g}})$ to $U_q(\mathfrak{g})$ depending on a parameter $a$, which is the analogue of the evaluation morphism from $\hat{\mathfrak{g}}$ to $\mathfrak{g}$. 

For other simple Lie algebras $\mathfrak{g}$, such an analogue of evaluation morphism does not exist. However, what is true in general is that there is a family of automorphisms of $U_q(\hat{\mathfrak{g}})$ depending on a parameter $a$. So assume that the representation $V$ of $U_q(\mathfrak{g})$ can be extended to a representation of $U_q(\hat{\mathfrak{g}})$. Then by twisting this representation by the automorphisms, one gets the family of representations $V(a)$ of $U_q(\hat{\mathfrak{g}})$ extending $V$. 

This leaves us with the following conclusion: we have a solution of the YB equation on representations $V$ of $U_q(\mathfrak{g})$ which can be extended to $U_q(\hat{\mathfrak{g}})$. To insist once more, for $\mathfrak{g}=sl_N$, all representations can be extended but for arbitrary $\mathfrak{g}$ this is not true. This leads to the natural concept of ``minimal affinization'' of a representation \cite{CP2}. Note that reversing the point of view, one could also say: we have a solution of the YB equation on representations $V$ of $U_q(\mathfrak{g})$ which can be obtained from restrictions of representations of $U_q(\hat{\mathfrak{g}})$. At the end, we do have a rather large set of solutions of the YB equation obtained through this construction. As a example, we refer to \cite{ZJ2} where we can see an instance of these subtleties associated to affinization of representations.

\begin{rema}
Example \ref{ex-solYBq} is associated to $U_q(sl_2)$ and the fundamental representation. Before that, we had Example \ref{ex-solYB}, the Yang solution, which can be obtained as a certain limit $q\to 1$ of the one from $U_q(sl_2)$. In general a certain limit $q\to 1$ of the solutions associated to the centraliser of $U_q(\mathfrak{g})$ will produce interesting solutions, which live in the centraliser of $U(\mathfrak{g})$. 
These solutions also come from an algebra, sharing many properties with $U_q(\hat{\mathfrak{g}})$, which is called the Yangian $Y(\mathfrak{g})$ (see \cite{CP}). This algebra, without parameter $q$, can be seen a certain limit of $U_q(\hat{\mathfrak{g}})$. 
\end{rema}

\subsection{The quantum Schur--Weyl duality}\label{subsec-qSW}

Here we come back to a more concrete discussion to explain that our examples from earlier, the Hecke and the BMW algebra, come indeed from this centraliser story.

\paragraph{Quantum group $U_q(sl_2)$.} We start with $\mathfrak{g}=sl_2$. As a vector space, $sl_2$ is generated by the following matrices, and it has the following natural Lie algebra structure:
\begin{equation}\label{rel-sl2}h=\left(\begin{array}{cc} 1 & 0 \\ 0 & -1 \end{array}\right)\,,\ \ \ x=\left(\begin{array}{cc} 0 & 1 \\ 0 & 0 \end{array}\right)\,,\ \ \ y=\left(\begin{array}{cc} 0 & 0 \\ 1 & 0 \end{array}\right)\ \ \ \ \text{and}\ \ \ \ [h,x]=2x\,,\ \ [h,y]=-2y\,,\ \ [x,y]=h\ .
\end{equation}
The definition of $U_q(sl_2)$ that we are going to use is the following one. As a vector space, we define $U_q(sl_2)$ to be
\[U(sl_2)[[\alpha]]=\{c_0+c_1\alpha+\dots\,,\ c_i\in U(sl_2)\}\,,\]
the vector space of formal power series in $\alpha$ with coefficients in $U(sl_2)$. We abuse notation and keep the names $x,y,h$ for the generators of $U(sl_2)$ but we insist that here we mean only the vector space $U(sl_2)$; the multiplication will be different. The multiplication of elements of $U_q(sl_2)$ is the usual multiplication of formal series, together with the following defining relations between the elements $x,y,h$:
\begin{equation}\label{def-sl2-q2}
hx-xh=2x\,,\ \ hy-yh=-2y\,,\ \ \ xy-yx=\frac{e^{\alpha h}-e^{-\alpha
h}}{e^{\alpha}-e^{-\alpha}}\ .
\end{equation}
One has to notice that the right hand side of the last relation is indeed a power series in $\alpha$. We should insist that as a vector space, $U(sl_2)[[\alpha]]$ is not the same as $U(sl_2)$ with
coefficients in $\mathbb{C}[[\alpha]]$. This is because $U(sl_2)$ is not finite dimensional. For
example, consider $e^{\alpha h}$.

It is usual to denote $q=e^{\alpha}$ and $K=e^{\alpha h}$. And using $K$, one can avoid the formal series in the defining relations, but then the $R$-matrix below is harder to introduce, and the connections with $U(sl_2)$ are not as transparent. Here, the limit $\alpha=0$ is well-defined and one recovers the algebra $U(sl_2)$ as can be seen immediately in the defining relations.

\paragraph{Coproduct and $R$-matrix.} 
We consider the completed tensor product $U_q(sl_2)\hat{\otimes} U_q(sl_2)$, which means that as a vector space, this is the space of power series in $\alpha$ with coefficients in the usual tensor product $U(sl_2)\otimes U(sl_2)$. The coproduct is defined on the generators by:
\begin{equation}\label{coprod-sl2}\Delta(x)=x\otimes e^{-\alpha h/2}+e^{\alpha h/2}\otimes x\,,\ \ \ \Delta(y)=y\otimes
e^{-\alpha h/2}+e^{\alpha h/2}\otimes y\,,\ \ \ \Delta(h)=h\otimes 1+1\otimes h\,,
\end{equation}
One can check that $\Delta$ extends to an algebra homomorphism from $U_{q}(sl_2)$ to $U_q(sl_2)\hat{\otimes} U_q(sl_2)$, and it is immediate to see that the limit $\alpha=0$ gives back the usual coproduct on $U(sl_2)$.

The completion of the tensor product becomes important when talking about the universal $R$-matrix $\mathcal{R}$. Indeed an explicit formula is:
\begin{equation}\label{R-sl2}\mathcal{R}=e^{\alpha(h\otimes h)/2}\sum_{n\geq 0}\frac{(q-q^{-1})^n}{[n]_q!}q^{n(n-1)/2}\Bigl(e^{-\alpha h/2}y\otimes x e^{\alpha h/2}\Bigr)^n\,,
\end{equation}
where we have set $q=e^{\alpha h}$, $[n]_q!=[2]_q\cdot [3]_q\dots [n]_q$ where $[k]_q=\frac{q^k-q^{-k}}{q-q^{-1}}$. The term in the sum being a multiple of $\alpha^n$, it follows that the element $\mathcal{R}$ is a well-defined element of the completed tensor product.

\paragraph{The quantum group $U_q(sl_N)$.} Now we take $\mathfrak{g}=sl_N$ for some $N>0$. Natural generators of $sl_N$ are the elements $h_i,x_i,y_i$, with $i=1,\dots,N-1$. The triplet of elements $h_i,x_i,y_i$ corresponds to the following matrices of $sl_N$: take the matrices $\left(\begin{array}{cc} 1 & 0 \\ 0 & -1 \end{array}\right)$, $\left(\begin{array}{cc} 0 & 1 \\ 0 & 0 \end{array}\right)$, $\left(\begin{array}{cc} 0 & 0 \\ 1 & 0 \end{array}\right)$ as for $sl_2$ and plug them in lines and columns $i$ and $i+1$.

The multiplication is deformed as follows. First each triplet $h_i,x_i,y_i$ satisfies the relations of $U_q(sl_2)$ in (\ref{def-sl2-q2}). Then we keep the following usual relations of matrices:
\[[h_i,x_{i\pm1}]=-x_{i\pm1}\,,\ \ \ \ \ [h_i,y_{i\pm1}]=y_{i\pm1}\,,\ \ \ \ \ [h_i,x_{j}]=[h_i,y_{j}]=0\ \ (|i-j|>1)\,,\]
\[[h_a,h_b]=0\ \ (\forall a,b)\,,\ \ \ \ \ \ [x_i,y_j]=0\ \ (i\neq j)\ .\]
And finally we add the deformation of the so-called Serre relations:
\[x_jx_i^2-(q+q^{-1})x_ix_jx_i+x_i^2x_{j}=y_jy_i^2-(q+q^{-1})y_iy_jy_i+y_i^2y_{j}=0\,,\ \ \ \ \text{ if $|i-j|=1$\,.}\]
For each triplet $h_i,x_i,y_i$, the coproduct is given by the same formulas as for $U_q(sl_2)$ in (\ref{coprod-sl2}).

\paragraph{Vector representation of $U_q(sl_N)$.} The fundamental (or vector) representation $V$ of $U_q(sl_N)$ is simply given by assigning to each $h_i,x_i,y_i$ its natural $N\times N$ matrix. All defining relations of $U_q(sl_N)$ are indeed satisfied. It looks as if we were dealing with $U(sl_N)$, but of course the ``quantum'' novelty of $U_q(sl_N)$ is that we have a different way to construct a representation on $V^{\otimes n}$. In particular, the permutations of the factors do not belong to the centraliser $\text{End}_{U_q(sl_N)}(V^{\otimes n})$. This is where we meet again with the Hecke algebra.

There is an explicit formula for the universal $R$-matrix of $U_q(sl_N)$ but we will not give it. It turns out that when we consider its image in the representation $V\otimes V$ and we multiply by the permutation operator, we find:
\[\check R(e_a\otimes e_b):=\left\{\begin{array}{ll}
q\,e_a\otimes e_b\ \  & \text{if $a=b$,}\\[0.4em]
e_b\otimes e_a+(q-q^{-1})\,e_a\otimes e_b & \text{if $a<b$,}\\[0.4em]
e_b\otimes e_a & \text{if $a>b$.}
\end{array}\right.\ \ \ \text{where $a,b=1,\dots,N$,}\]
where $(e_1,\dots,e_N)$ is the basis used to construct the representation $V$ above. The crucial point for what follows is that we recover the representation of the Hecke algebra $H_n(q)$ on the tensor product $V^{\otimes n}$ that was given explicitly in (\ref{rep-Hn}).

\paragraph{Hecke algebra $H_n(q)$ and quantum Schur--Weyl duality.} Information on the quantum Schur--Weyl duality can be found in \cite{CP,Ji86,KS}. Here we consider $U_q(gl_N)$ instead of $U_q(sl_N)$ for a better comparison with the first section. Recall that we have a representation of the Hecke algebra $H_n(q)$ on the tensor product $V^{\otimes n}$. So our situation is that we have two algebras represented on the same vector space, and we can picture it like this:
\[U_q(gl_N)\ \ \stackrel{\rho}{\longrightarrow}\ \ \ \ \text{End}(V^{\otimes n})\ \ \ \ \stackrel{\pi}{\longleftarrow}\ \ H_n(q)\ \]
Now we can state the (first part of) the Schur--Weyl duality.
\begin{framedtheo}[quantum Schur--Weyl I]\label{theo-qSW1}
The centraliser $\text{End}_{U_q(gl_N)}(V^{\otimes n})$ is the image of the Hecke algebra $H_n(q)$:
\begin{equation}\label{qSW1}
\text{End}_{U_q(gl_N)}(V^{\otimes n})=\pi\bigl(H_n(q)\bigr)\ .
\tag{\textbf{qSW1}}
\end{equation}
\end{framedtheo}
Here we are in the semisimple situation ($q$ is generic) and from the general consideration on centralisers sketched above, we have actually that $\pi\bigl(H_n(q)\bigr)$ and $\rho\bigl(U_q(gl_N)\bigr)$ are the mutual centralisers of each other.

\vskip .2cm
The second part of the quantum Schur--Weyl duality is going to be also very similar to its classical counterpart, but first, we briefly discuss the irreducible representations of $U_q(gl_N)$ and of $H_n(q)$. We have:
\begin{itemize}
\item The quantum group $U_q(gl_N)$ has finite-dimensional irreducible representations in bijection with those of $U(gl_N)$. Thus, they are also indexed by highest weights, which are here again identified with partitions $\lambda$ such that $\ell(\lambda)\leq N$. Let us denote them by $V_{\lambda}^N$, as in the classical situation. The vector representation $V$ corresponds to $\lambda=(1)$. 
\item The Hecke algebra $H_n(q)$ has irreducible representations in bijection with those of the symmetric group $S_n$. Thus they are indexed by partitions $\lambda$ of size $n$. Let us denote them by $S_{\lambda}$ (see Appendix).
\end{itemize}
\begin{framedtheo}[quantum Schur--Weyl II]\label{theo-qSW2}$\ $
The decomposition of $V^{\otimes n}$, as a representation of $U_q(gl_N)\otimes H_n(q)$ is:
\begin{equation}\label{qSW2}
V^{\otimes n}=\bigoplus_{\substack{\lambda\vdash n\\[0.2em] \ell(\lambda)\leq N}} V_{\lambda}^N\otimes S_{\lambda}\,.\tag{\textbf{qSW2}}
\end{equation}
\end{framedtheo}

The complete description of the centraliser of $U_q(gl_N)$, that is, the description of the kernel of the map from $H_n(q)$ to $\text{End}(V^{\otimes n})$ follows the same lines as in the classical case. We only need the analogue of the antisymmetriser, and for this we need to go a little bit into the algebraic structure of the Hecke algebra $H_n(q)$, whose definition was given in (\ref{rel-H}).

\vskip .2cm
For any element $w$ of the symmetric group $S_n$, let $w=s_{a_1}\dots s_{a_k}$ be a reduced (\emph{i.e.} minimal length) expression for $w$ in terms of the generators $s_i=(i,i+1)$, and denote $\ell(w)=k$. Then define $\si_w:=\si_{a_1}\dots \si_{a_k}\in H_n(q)$. This definition does not depend on the reduced expression for $w$ and the set $\{\si_w\}_{w\in S_n}$ forms a basis of $H_n(q)$. Then the \textbf{$q$-antisymmetriser} in $H_n(q)$ is the following element:
\begin{equation}\label{def-P'}
P'_n=\frac{\sum_{w\in S_n}(-q^{-1})^{\ell(w)}\si_w}{\sum_{w\in S_n}q^{-2\ell(w)}}=\frac{q^{n(n-1)/2}}{[n]_q!}\sum_{w\in\mS_n}(-q^{-1})^{\ell(w)}\si_w\ ,
\end{equation}
where the $q$-numbers are defined as follows:
\begin{equation}\label{quantum-numbers}
[L]_q:=\frac{q^L-q^{-L}}{q-q^{-1}}=q^{L-1}+q^{L-3}+\dots+q^{-(L-1)}\ \ \ \ \ \text{and}\ \ \ \ \ \ \ [L]_q!:=[1]_q[2]_q\dots[L]_q\ .
\end{equation}
Clearly, if $q=1$ the $q$-antisymmetriser becomes the usual antisymmetriser. 

For the statement below, note that we can see the element $P'_{N+1}$ as an element of $H_n(q)$ by the natural inclusion of $H_{N+1}(q)$ in $H_n(q)$ (if $n>N$).
\begin{framedtheo}[quantum Schur--Weyl II']\label{theo-qSW2b}$\ $
\begin{itemize}
\item If $n\leq N$ then the kernel of the map is $\{0\}$ and the centraliser $\text{End}_{U_q(gl_N)}(V^{\otimes n})$ is isomorphic to $H_n(q)$.
\item If $n>N$ then the kernel of the map is generated in $H_n(q)$ by the element $P'_{N+1}$. So the centraliser $\text{End}_{U_q(gl_N)}(V^{\otimes n})$ is the quotient of $H_n(q)$ by the relation $P'_{N+1}=0$.
\end{itemize}  
\end{framedtheo}
Explicitly, the centraliser $\text{End}_{U_q(gl_N)}(V^{\otimes n})$ is described as the algebra generated by $s_1,\dots,s_{n-1}$ with defining relations:
\[\begin{array}{ll}
\si_i\si_{i+1}\si_i=\si_{i+1}\si_i\si_{i+1}\,,\ \ \  & \text{for $i\in\{1,\dots,n-2\}$}\,,\\[0.2em]
\si_i\si_j=\si_j\si_i\,,\ \ \  & \text{for $i,j\in\{1,\dots,n-1\}$ such that $|i-j|>1$}\,,\\[0.2em]
\si_i^2=1+(q-q^{-1})\si_i\,,\ \ \  & \text{for $i\in\{1,\dots,n-1\}$}\,,\\[0.2em]
P'_{N+1}=0 & \text{if $n>N$.} 
\end{array}
\]
\begin{exam}[Temperley--Lieb algebra]\label{exa-TL2}
Let $N=2$. In this case, the centraliser $\text{End}_{U_q(gl_2)}(V^{\otimes n})$ coincides with the Temperley--Lieb algebra. The additional relation $P'_3=0$ reads:
\[1-q^{-1}\si_1-q^{-1}\si_2+q^{-2}\si_1\si_2+q^{-2}\si_2\si_1-q^{-3}\si_1\si_2\si_1=0\ .\]
One can show that it implies the same relation with indices $i,i+1$ for all $i=1,\dots,n-2$. Then setting $e_i:=\si_i-q$, one recovers the other standard presentation of the Temperley--Lieb algebra:
\[e_i^2=-[2]_qe_i\,,\ \ \ \ e_ie_{i+1}e_i=e_i\,,\ \ \ \ e_{i+1}e_{i}e_{i+1}=e_{i+1}\ \ \ \text{and}\ \ \ e_ie_j=e_je_i\ \text{if $|i-j|>1$.}\]
\end{exam}

\paragraph{And the BMW algebra?}
We explained how the Hecke algebra fits in the story of centralisers of quantum groups representations. This nice story has its counterpart for the BMW algebra. Instead of $gl_N$ we can consider the other classical Lie algebras $so_N$ or $sp_N$. We can form the quantum groups $U_q(so_N)$ or $U_q(sp_N)$, and we can consider the tensor product $V^{\otimes n}$, where $V$ is the (analogue of the) vector representation for $U_q(so_N)$ or $U_q(sp_N)$. Then the centralisers of $V^{\otimes n}$ are described in a way similar than above, where the Hecke algebra is replaced by the BMW algebra (with the parameter $a$ specialised to a power of $q$ depending on whether we have $so_N$ or $sp_N$). For details about that along the same lines as above, we refer for example to \cite{KS}. This concludes our discussion on the interpretation of Hecke and BMW algebras as centralisers of representations of quantum groups.

\section{Markov traces on the braid group and link invariants}\label{sec-markov}

We briefly review the well-known procedure to produce invariants for links using the braid group, relying on Alexander and Markov theorems. Then we state that such a procedure can be applied to centralisers of quantum groups, leading to Reshetikhin--Turaev invariants. We give more details, as an example, for the Markov trace on the family of Hecke algebras, which will be used in Chapter \ref{chap-yok}.

\subsection{Alexander and Markov theorems}

References here are \cite{Bir,BZ,Kam,KT}. The closure of a braid is obtained by connecting the dots on top with the dots on the bottom in a braid. We connect the leftmost dot on top with the leftmost dot at the bottom, then we do the same for the second dots, and so on. The closure of a braid results in a link in three-dimensional space. We consider links up to ambient isotopy, and we have the two fundamental results.
\begin{framedtheo}$\ $
\begin{itemize}
\item (Alexander theorem, \cite{Ale}) Any link can be obtained as the closure of some braid.
\item (Markov theorem, \cite{Mark}) Two braids give the same closure up to isotopy if and only if they are equivalent under the equivalence relation generated by:
\begin{equation}\label{markov-move}\alpha\beta\sim\beta\alpha\ \ (\alpha,\beta\in B_n,\ n\geq1)\ \ \ \ \quad\text{and}\quad\ \ \ \ \alpha\sigma_n^{\pm1}\sim\alpha\ \ (\alpha\in B_n,\ n\geq1)\ .
\end{equation}
\end{itemize}
\end{framedtheo}
It is easy to see, drawing braids, that the two conditions in (\ref{markov-move}) are necessary. The main result of Markov theorem is that they are also sufficient to ensure that two equivalent braids have isotopic closures. Note that in the second equivalence, $\alpha$ is a braid on $n$ strands, while $\alpha\sigma_n^{\pm1}$ is a braid on $n+1$ strands (implicit is the use of the natural inclusion $B_n\subset B_{n+1}$ to see $\alpha$ also as an element of $B_{n+1}$). The equivalence relation is on the set of braids on any number of strands, that is, on $\bigcup_{n\geq 1} B_n$.

\paragraph{Definition of Markov traces.} From the above results, a link, up to isotopy, corresponds to an equivalence class of braids. Then invariants of links can be defined using the equivalence classes of braids, by assigning a certain value (a number, a polynomial, or anything else) to any equivalence class. Of course, one would like to calculate it directly on a braid (a representative of the class) and one is led to the notion of a Markov-compatible map, which is a map $\tau$ on the set of all braids $\bigcup_{n\geq 1} B_n$, such that:
\[\tau(\alpha\beta)\sim\tau(\beta\alpha)\ \ (\alpha,\beta\in B_n,\ n\geq1)\ \ \ \ \quad\text{and}\quad\ \ \ \ \tau(\alpha\sigma_n^{\pm1})\sim\tau(\alpha)\ \ (\alpha\in B_n,\ n\geq1)\ .\]
From now on we will see the map $\tau$ as a family of maps $\{\tau_n\}_{n\geq1}$, each $\tau_n$ acting on the braid group $B_n$. 

Of course, the perfect Markov-compatible map is the map which assigns to a braid its equivalence class. But one would like to be able to actually calculate the value of the invariants, and for this purpose, it is useful to consider algebraic structures less complicated than the braid group. These are quotients of the braid group algebra (such as the Hecke algebra).

So assume that we have, for each $n\geq1$, an algebra $\mathcal{A}_n$ which is a quotient of the braid group algebra $RB_n$, where $R$ is the ring of coefficients. For simplicity here, we will assume (as it will be the case in our examples) that the algebras $\mathcal{A}_n$ form a chain of algebras for natural inclusion maps:
\[\mathcal{A}_1\subset \mathcal{A}_2\subset \dots\dots \mathcal{A}_n\subset \mathcal{A}_{n+1}\subset\dots\dots\ .\]
We abuse notations and keep the names $\si_i$ for the images of the braid group generators in $\mathcal{A}_n$.
\begin{frameddefi} A Markov trace on $\{\mathcal{A}_n\}_{n\geq1}$ is a family of linear maps $\{\tau_n\,:\,\mathcal{A}_n\to R\}_{n\geq1}$ satisfying:
\begin{equation}\label{Markov-tau}
\begin{array}{llr}
\tau_n(xy)=\tau_n(yx)\ ,& \text{for $n\geq1$ and $x,y\in\mathcal{A}_n$;} &\quad \text{\emph{(Trace condition)}}\\[0.2em]
\tau_{n+1}(x\si_n)=\tau_{n+1}(x\si_n^{-1})=\tau_n(x)\ , & \text{for $n\geq1$ and $x\in\mathcal{A}_n$.} &\quad\text{\emph{(Markov condition)}}
\end{array}
\end{equation}
\end{frameddefi}
It may be worth emphasizing the following. An element $x\in\mathcal{A}_n$ can be seen as an element of $\mathcal{A}_{n+1}$, but we do not require that $\tau_{n+1}(x)=\tau_n(x)$. The dependence on $n$ for the maps $\{\tau_n\}_{n\geq 1}$ is important. This is a place of possible confusion since many constructions of invariants in the literature use different sort of traces (sometimes also called Markov traces). We will unravel this explicitly in the last paragraph below.

\paragraph{Invariants for links.} Whenever we have a Markov trace as above, by construction, we are very close to having a link invariant. We only need to apply, successively, Alexander and Markov theorems. Namely, start from a link $L$ and apply Alexander theorem to find a braid $\beta_L$ whose closure is $L$. Then send $\beta_L$ in the algebra $\mathcal{A}_n$ (where $n$ is the number of strands of $\beta_L$) and apply the map $\tau_n$. From Markov theorem, this is an invariant of the link $L$.

To formalize it a bit more, with a formula, we give a name to the map sending braids into algebras $\mathcal{A}_n$:
\[\pi_n\ :\ RB_n\to \mathcal{A}_n\ \ \ \qquad (n\geq 1)\ .\]
The maps $\pi_n$ are morphisms of algebras sending the braid generators $\si_i\in B_n$ to the elements $\si_i\in\mathcal{A}_n$. Note that we do not require that $\pi_n$ are surjective, namely, we do not require that the image of the braid group $B_n$ generates $\mathcal{A}_n$. The surjectivity is satisfied in many examples, such as the Hecke algebras, but we will see that in general it is not: it is not always true for centralisers of quantum groups, as already mentioned, and it will not be true for the Yokonuma--Hecke algebras discussed in Chapter \ref{chap-yok}.

The formula for the value of the invariant of the link $L$ is simply:
\[P(L)=\tau_n\bigl(\pi_n(\beta_L)\bigr)\ \ \ \ \ \quad\text{where $\beta_L\in B_n$ closes to $L$.}\]

\paragraph{A slightly different construction.} It is possible to obtain invariants via a slightly different sort of traces (see \cite{Jo1}). We can define an ``Ocneanu trace'' $\psi$ on $\bigcup_{n\geq 1} \mathcal{A}_n$ as a linear map satisfying the trace condition and the modified Markov conditions:
\[\psi(x\si_n)=z\psi(x)\ \ \ \text{and}\ \ \ \psi(x\si_n^{-1})=\overline{z}\psi(x)\ , \text{for $n\geq1$ and $x\in\mathcal{A}_n$\,,}\]
where $z,\overline{z}$ are invertible elements of $R$. Here we assume that $\psi(x)$ does not depend on whether we see $x$ as an element of $\mathcal{A}_n$ or of $\mathcal{A}_{n+1}$ (or of $\mathcal{A}_{N}$ with $N\geq n$). With these conditions, one can check that the following formula gives an invariant:
\[P(L)=(z\sqrt{\theta})^{-(n-1)}(\sqrt{\theta})^{w(\beta_L)}\psi\bigl(\pi_n(\beta_L)\bigr)\ ,\ \ \ \ \quad\text{with $\theta=\frac{\overline{z}}{z}$\ ,}\]
where $\beta_L\in B_n$ closes to $L$, and $w(\beta_L)$ is the number of crossings in $\beta_L$ counted with signs. We have extended the ring $R$ by adding a square root of $\theta$.

Actually this invariant can also be obtained through a Markov trace, in our sense. First, one has to renormalize the braid generators in $\mathcal{A}_n$. We set $\si'_i=\frac{1}{\sqrt{\theta}}\si_i$. Of course, the braid relations, being homogeneous, are still satisfied by $\si'_i$. Then the family of maps $\{\tau_n=(z\sqrt{\theta})^{-(n-1)}\psi\}_{n\geq 1}$ is a Markov trace, with respect to the generators $\si'_i$, and they produce precisely the same invariants as above.

\subsection{Reshetikhin--Turaev invariants}

We know, abstractly, a large number of chains of algebras $\{\mathcal{A}_n\}_{n\geq 1}$, such that $\mathcal{A}_n$ is a quotient of the algebra of the braid group $B_n$. As discussed in the preceding section, they are the centraliser algebras of quantum groups representations, namely,
\[\mathcal{A}_n=\text{End}_{U_q(\mathfrak{g})}(V^{\otimes n})\ ,\]
where $V$ is a representation of a quantum group $U_q(\mathfrak{g})$. The images of the generators of the braid group are given by the $R$-matrix. So one question is whether we can find a Markov trace on $\{\mathcal{A}_n\}_{n\geq 1}$. It turns out that this is possible\footnote{A more general setting is the context of ribbon Hopf algebras, see \emph{e.g.} \cite{LR}.}.
\begin{framedtheo}[\cite{Res,Tur}, see also \cite{KRT,ZGB}]$\ $\\
The family of centraliser algebras $\{\mathcal{A}_n\}_{n\geq 1}$ supports a Markov trace.
\end{framedtheo}
We will make two comments on this remarkably general result.\\[0.1cm]
$\bullet$ First, on the conceptual side, there is an interesting, abstract way to understand how to produce the Markov trace. First one has to remember that there is the notion of a dual, or contragredient, representation $V^{\star}$ of a representation $V$ of a quantum group. Moreover, there is also a ``trivial'' one-dimensional representation $V_0=\mathbb{C}$. The axioms of quantum group ensure that the trivial representation $V_0$ appears with multiplicity one in the tensor product $V\otimes V^{\star}$. This leads to a unique (up to normalization) morphism commuting with the quantum group action:
\[\text{qTr}_V\ :\ V\otimes V^{\star}\to \mathbb{C}\ .\]
The map $\text{qTr}_V$ is a quantum group analogue of the usual trace of a matrix, and is called the ``quantum trace''. Then consider the following procedure:
\[\begin{array}{ccccc} 
\text{End}(V^{\otimes n}) & \to & \text{End}(V^{\otimes (n-1)}\otimes V\otimes V^{\star}) & \to &  \text{End}(V^{\otimes (n-1)})\\[0.4em]
M & \mapsto & M\otimes \text{Id}_{V^{\star}} & \mapsto & (\text{Id}_{\text{End}(V^{\otimes (n-1)})} \otimes \text{qTr}_V)\bigl(M\otimes \text{Id}_{V^{\star}}\bigr)\ .
\end{array}\]
In simple words, we take the quantum trace over the $n$-th factor in the tensor product. It turns out that, restricted to the subalgebra $\mathcal{A}_{n}\subset \text{End}(V^{\otimes n})$, this procedure takes values in the centraliser $\mathcal{A}_{n-1}\subset\text{End}(V^{\otimes (n-1)})$. Thus we have a map from $\mathcal{A}_{n}$ to $\mathcal{A}_{n-1}$. It is sometimes called ``conditional expectation value''.

Finally, we can repeat the procedure $n$ times and obtain thus a map from $\mathcal{A}_{n}$ to $\mathbb{C}$. This gives, ignoring renormalization issues, the value of the Markov trace.\\[0.1cm]
$\bullet$ This is rather abstract, but there is also a concrete way to perform this procedure. There is a matrix $\mu_V\in\text{End}(V)$, which is in fact the image of a certain element of the quantum group in its representation, which allows to calculate the Markov trace. We will not explain how to find $\mu_V$, we only indicate that the value of the Markov trace is calculated by using $\mu_V$ and the usual trace:
\[\mathcal{A}_n\ni M\mapsto \text{Tr}\bigl(M\cdot(\mu_V\otimes\dots\otimes\mu_V)\bigr)\ .\]

\subsection{Example of the Hecke algebra: HOMFLYPT polynomial}

For our purposes, it is convenient to add another parameter in the characteristic equation for the generators of the Hecke algebra. The approach below to the Markov trace on the Hecke algebras is described for example in \cite[section 4.5]{GP}.

\vskip .2cm 
In this section, the Hecke algebra $H_n$ is defined over $\mathbb{C}[u^{\pm1},q^{\pm1},(q-q)^{-1}]$ with:
\[
\begin{array}{ll}
\si_i\si_{i+1}\si_i=\si_{i+1}\si_i\si_{i+1}\,,\ \ \  & \text{for $i\in\{1,\dots,n-2\}$}\,,\\[0.2em]
\si_i\si_j=\si_j\si_i\,,\ \ \  & \text{for $i,j\in\{1,\dots,n-1\}$ such that $|i-j|>1$}\,,\\[0.2em]
\si_i^2=u^2+u(q-q^{-1})\si_i\,,\ \ \  & \text{for $i\in\{1,\dots,n-1\}$}\,.
\end{array}
\]
From the definition, the set of Markov traces on a given family of algebras forms a vector space. The statement below says that for the Hecke algebras, this vector space is one-dimensional.
\begin{framedtheo}
Fixing the value of $\tau_1(1)$, there is a unique Markov trace $\{\tau_n\}_{n\geq 1}$ on the family of Hecke algebra $\{H_n\}_{n\geq 1}$.
\end{framedtheo}
The map $\tau_n$ can be given explicitly in the following sense. First we choose a convenient basis of $H_n$, for example the following (where the product of sets $A.B$ is $\{a.b\ |\ a\in A,\ b\in B\}$):
\[
\left\{\begin{array}{c} 1,\\ \si_1 \end{array}\right\}\cdot \left\{\begin{array}{c} 1,\\ \si_2, \\ \si_2\si_1 \end{array}\right\}\cdot \left\{\begin{array}{c} 1,\\ \si_3, \\ \si_3\si_2, \\ \si_3\si_2\si_1 \end{array}\right\}\cdot\ \ldots\ \cdot \left\{\begin{array}{c} 1,\\ \si_{n-1}, \\ \vdots \\ \si_{n-1}\dots\si_1 \end{array}\right\}\ .
\]
Note that the generator $\si_{n-1}$ appears at most once. Consider the case where $\si_{n-1}$ does not appear, that is, take $x\in H_{n-1}$. The trick is to use the formula $q-q^{-1}=u^{-1}\si_{n-1}-u\si_{n-1}^{-1}$. Using the invertibility of $(q-q^{-1})$, we have:
\[\tau_n(x)=\tau_n(\frac{u^{-1}\si_{n-1}-u\si_{n-1}^{-1}}{q-q^{-1}}x)=\frac{u^{-1}-u}{q-q^{-1}}\tau_{n-1}(x)\ .\]
The remaining elements of the basis above can be written $x\si_{n-1}\dots\si_k$ with $x\in H_{n-1}$ and $k\in\{1,\dots,n-1\}$, and for them, we simply have:
\[\tau_n(x\si_{n-1}\dots\si_k)=\tau_{n-1}(x\si_{n-2}\dots\si_k)\ .\]
The two preceding formulas uniquely fix, recursively, the family of maps $\{\tau_n\}_{n\geq 1}$ up to the choice of $\tau_1(1)$.

\paragraph{HOMFLYPT and Jones polynomials.} Now that we have a Markov trace, we can define an invariant, as explained in a preceding subsection. The resulting invariant depends on the two parameters $q$ and $u$ and is called HOMFLYPT polynomial (it is a Laurent polynomial in $u$ and $q-q^{-1}$); see \cite{HOMFLY,Jo1,PT}.

Once we have a Markov trace on the Hecke algebras, we can ask whether it passes to some quotients, for example to the Temperley--Lieb algebras. With our normalization, the additional relation giving the Temperley--Lieb algebra $TL_n$ is:
\[1-(uq)^{-1}\si_1-(uq)^{-1}\si_2+(uq)^{-2}\si_1\si_2+(uq)^{-2}\si_2\si_1-(uq)^{-3}\si_1\si_2\si_1=0\ .\]
It is not too difficult to show that $\tau_n$ passes to the quotient $TL_n$ (for any $n$) if and only if $\tau_3$ applied to the left-hand-side of the above relation gives $0$. This happens for example if $u=q^2$ and this gives a one-variable invariant, which is the Jones polynomial.

\section{A jungle of Hecke algebras}\label{sec-Hecke}

The terminology ``Hecke algebra'' refers to many different algebras in different areas of mathematics. Algebras termed Hecke algebras appear in several parts of representation theory, in number theory, in the study of Coxeter groups and reflection groups, of braid groups and Artin groups, in the study of quantum groups and in categorification, in pure deformation theory, and probably in other fields as well. We discuss some definitions of Hecke algebras, the ones we use in this thesis, and some connections between them, with no intent at being exhaustive.

\paragraph{Hecke algebras as centraliser algebras.} Let $G$ be a group with a subgroup $B$. Assume for simplicity that $G$ is finite. The Hecke algebra $H(G,B)$ is the centraliser algebra of some representation of $G$. This representation is the induced representation from the trivial representation $\textbf{1}_B$ of dimension 1 of $B$. Equivalently, it is the permutation representation of $G$ associated to the action of $G$ on its coset space $G/B$:
\begin{equation}\label{defHGB} H(G,B)=\text{End}_G(\mathbb{C}G/B)=\text{End}_G(\text{Ind}_B^G\textbf{1}_B)\ .
\end{equation}
A basis of $H(G,B)$ is indexed by double cosets of $G$ with respect to $B$ and a multiplication is given in combinatorial terms, see for example \cite[\S 8.4]{GP}. Of course following the general centraliser philosophy, the Hecke algebra $H(G,B)$ contains information on the decomposition of the permutation representation into irreducible representations. One can consider more general pairs $(G,H)$ of group-subgroup, such as a topological group $G$ with a subgroup $H$ satisfying some conditions, but we skip these details.

The name ``Hecke algebra'' originates from the work of Erich Hecke \cite{Hec} on modular forms, see \cite{Bum,Iwa,Leh,Ogg,Shi}. Roughly speaking, he used a certain commutative set of operators to find a distinguished basis for some space of modular forms. These operators are closely related, in our modern terminology, to a centraliser algebra as above.

\paragraph{Hecke algebras as deformations of Coxeter groups.} A Coxeter group $W$ is defined by a set of generators $S=\{s_1,\dots,s_n\}$ and with defining relations:
\[s_i^2=1\ \ \ \ \quad \text{and}\quad\  \ \ \ \underbrace{s_is_js_i\dots}_{m_{ij}\ \text{terms}}=\underbrace{s_js_is_j\dots}_{m_{ij}\ \text{terms}}\ \ \ \ (i\neq j)\ ,\]
for a choice of $m_{ij}\in\mathbb{Z}_{\geq 2}\cup\{\infty\}$ with $m_{ij}=m_{ji}$ (if $m_{ij}=\infty$, we remove the defining relation with $s_i$ and $s_j$). The finite Coxeter groups classify the finite reflection groups \cite{Hum}.

The associated Hecke algebra $H(W)$ is defined over $\mathbb{Z}[q,q^{-1}]$ by generators $\si_1,\dots,\si_n$ and relations:
\[\si_i^2=1+(q-q^{-1})\si_i\ \ \ \ \quad \text{and}\quad\  \ \ \ \underbrace{\si_i\si_j\si_i\dots}_{m_{ij}\ \text{terms}}=\underbrace{\si_j\si_i\si_j\dots}_{m_{ij}\ \text{terms}}\ \ \ \ (i\neq j)\ .\]
We can allow for different parameters $q_i$, as long as they are constant on conjugacy classes of the generators.
For any element $w\in W$, let $w=s_{a_1}\dots s_{a_k}$ be a reduced (\emph{i.e.} minimal length) expression for $w$ in terms of the generators. Then define $\si_w:=\si_{a_1}\dots \si_{a_k}\in H(W)$. This definition does not depend on the reduced expression for $w$ and the set 
\[\{\si_w\}_{w\in W}\ ,\]
forms a basis of $H(W)$. Thus the Hecke algebra $H(W)$ is a flat deformation of the group algebra $\mathbb{C}W$. A standard reference is \cite{GP}. 

\begin{exam}
We recover the Hecke algebra $H_n$ of the previous sections, associated to the symmetric group $S_n$, by taking $S=\{s_1,\dots,s_{n-1}\}$, the integers $m_{i,i+1}=m_{i+1,i}=3$ and $m_{ij}=2$ otherwise. In this case, $W=S_n$ and the generators are $s_i=(i,i+1)$.
\end{exam}

The name Hecke algebras (also often called Iwahori--Hecke algebras, or finite Iwahori--Hecke algebras, to distinguish them from their affine versions below) for these algebras is justified through the following fact \cite{Iwah}. Roughly, if we take $G$ a reductive group over a finite field and $B$ its Borel subgroup, then the centraliser algebra as in (\ref{defHGB}) coincides with the Hecke algebra of this paragraph when $W$ is the Weyl group of $G$. More precisely, the parameter $q$ becomes the cardinal of the finite field\footnote{thus explaining the choice of letter ``$q$'': it is the letter after ``p'' and a finite field is usually denoted $\mathbb{F}_q$, where $q$ is a power of $p$. Luckily, the letter after $p$ is also the first letter of the word ``quantum'' so there is no conflict with quantum group theory.}. Thus the Iwahori--Hecke algebras allow to understand at least the piece of the representation theory of $G$ which is included in $\text{Ind}_B^G\textbf{1}_B$ (representations of $G$ with a vector fixed by $B$).

More generally, Iwahori--Hecke algebras are fundamental for the representation theory of reductive groups over finite fields. The key idea is ``Harish--Chandra philosophy'' which asserts that the study of these representation theories can be split into two steps. Roughly, one should first classify the basic blocks of the representation theory: the cuspidal, or discrete series, representations. Then it remains to understand some induced representations constructed from these building blocks. Hecke algebras, as centraliser algebras, appear in the second step. The permutation representation associated to the Borel subgroup is the first example. Remarkably, it turns out that all centraliser algebras of the induced representations (from cuspidal ones) are closely related to Iwahori--Hecke algebras associated to reflections groups \cite{HL,Lus2}. The Harish--Chandra philosophy, coming from representations of real Lie groups and harmonic analysis, is very fruitful in several other parts of representation theory. One we just mentioned, another one is for $p$-adic groups and the Langlands program, and it also involves Iwahori--hecke algebras.

\paragraph{Affine Hecke algebras.} Among the Coxeter groups, we have all the finite Weyl groups associated to finite root systems. It turns out that we also have the affine Weyl groups associated to the affine extension of root systems. These groups can be seen as semi-direct products of the finite Weyl group with the root lattice. It turns out that they also have a presentation as Coxeter groups. Therefore they have associated Hecke algebras, that are called affine Hecke algebras, or affine Iwahori--Hecke algebras.

Actually, one can consider ``extended'' affine Weyl groups by replacing the root lattice by a larger lattice (for instance, the weight lattice, but it can also be a larger lattice). One can associate an extended affine Weyl group to a root datum, which appears in the study of reductive groups (\emph{e.g.} \cite{GM}). It turns out that we can also define an affine Iwahori--Hecke algebra associated to a root datum \cite{Lus1,Sol2}. For example, we have an affine Hecke algebra of type GL, which is probably the most well-known. It will appear in Chapter \ref{chap-yok} as a particular case ($d=1$) of the affine Yokonuma--Hecke algebra.

Affine Iwahori--Hecke algebras play a similar role in the representation theory of reductive groups over $p$-adic fields as their finite version for finite fields. First if one takes $G$ a reductive group over a $p$-adic field and $B$ its Iwahori subgroup then the centraliser algebra in (\ref{defHGB}) is the affine Iwahori--Hecke algebra associated to the Weyl group of $G$ \cite{IM}. So at least part of the representation theory of $p$-adic group (representations with fixed vectors for the Iwahori subgroup) is controlled by affine Hecke algebras. It seems that much more is true again since affine Hecke algebras appear in the whole representation theory of $p$-adic groups in the same sort of Harish--Chandra philosophy as above, see \emph{e.g.} \cite{Sol1}.

We cannot possibly discuss here all the applications of finite or affine Iwahori--Hecke algebras. We can at least mention the Kazhdan--Lusztig polynomials \cite{KL} for their remarkable ubiquity in representation theory. More closely related to the themes of this thesis, let us restrict to the affine Hecke algebra of GL. First for integrable systems, it is the algebraic framework to incorporate the reflection equation into the story of the YB equation in the Hecke algebra \cite{IO2}. Then, for knots and links, it plays the same role as the Hecke algebra but for knots and links inside the solid torus \cite{Lam}. Finally, there is a generalisation of the Schur--Weyl duality relating the affine Hecke algebra and the affine quantum group of type A \cite{CP}.

\paragraph{Quiver Hecke algebras.} The quiver Hecke algebra, or KLR Hecke algebras, have been introduced by Khovanov--Lauda and Rouquier \cite{KhLa1,KhLa2,Rou} in the context of categorification of quantum groups. Roughly, the structure of their representation theories mimics the algebraic structure of (some subalgebras of) quantum groups. Even more, the cyclotomic quotients of KLR algebras categorify highest-weight integrable modules of quantum groups \cite{KK}. Such categorification statements, relating an algebra $A$ with a category of representations of another algebra $B$, are very fruitful in two directions. Of course, this produces naturally a lot of information on the representation theory of $B$. On the other hand, this often has striking applications regarding the algebra $A$ as well. For example, the positivity of the coefficients in the  decomposition of elements of $A$ in a certain basis follows immediately if this decomposition corresponds to decomposing a representation of $B$ into a direct sum (or a Jordan--Holder composition series).

The origin of these developments lies in the particularly interesting story for type A. Here a conjecture of Lascoux--Leclerc--Thibon, then proved by Ariki, is concerned with the modular representation theory of the symmetric group and of the usual Hecke algebra\footnote{It may be worth recalling that the modular (\emph{i.e.} over fields of positive characteristic) representation theory of the symmetric group is far less well understood than the characteristic zero case.}. They are related with the affine quantum group of type A, in the same spirit as in the preceding paragraph (see for example the introduction of \cite{KKK} and the references therein). When we follow the historical approach of quantum groups from integrable systems and/or from the braid group, it is a remarkable surprise to realise that they became also of fundamental importance for modular representations of the symmetric group.

The Brundan--Kleschev isomorphism explains the connections between the two paragraphs above: they prove that cyclotomic quotients of the KLR algebra of type A are in fact isomorphic to cyclotomic quotients of the affine Hecke algebra of type A, including the usual Hecke algebra and the symmetric group \cite{BK}. The isomorphism is explicit, even though quite intricated. The KLR algebras are defined by presentations with generators and relations which are not so appealing. Thanks to the Brundan--Kleshchev isomorphism, one can see them (in type A) as generalising for the non-semisimple blocks the easy explicit description of the irreducible representations in characteristic zero as in the Appendix. See \cite[\S 5]{BK}.

So, affine quantum groups of all type are categorified by KLR algebras of all type. In type A, KLR algebras are explicitly related to the affine Hecke algebra. This raises the question whether we can do something similar for other affine Hecke algebras. Conjectures relating affine Hecke algebras of types B and D to quantum groups in a similar spirit were made by M. Kashiwara and his collaborators \cite{EK,KM} and then proved by P. Shan, M. Varagnolo and E. Vasserot \cite{VV,SVV}. In their proofs, some geometrically defined algebras appear which play a role similar to the role played by the KLR algebra of type A. Based on their work, we have been able in \cite{PRo,PRu2,PRu1} to define algebras generalising the KLR algebras and to prove an explicit isomorphism with cyclotomic quotients of affine Hecke algebras of type B and D.

It turns out that the algebras defined in \cite{PRo,PRu2,PRu1} have been used recently for a Schur--Weyl duality statement relating them to quantum symmetric pairs of quantum groups \cite{AP}. These mathematical structures are related to the reflection equation in the same way as quantum groups are related to the Yang--Baxter equation \cite{BK,AV}. So the loop (for this thesis) is closed since Schur--Weyl duality, quantum groups and integrable systems also appear in these developments.

\section*{Appendix \thechapter.A $\ $Representation theory of the symmetric groups and the Hecke algebras}
\addcontentsline{toc}{section}{Appendix \thechapter.A $\ $Representation theory of the symmetric groups and the Hecke algebras}

\subsection*{\thechapter.A.1 $\ $Compositions, partitions and tableaux}
\addcontentsline{toc}{subsection}{\thechapter.A.1 $\ $Compositions, partitions and tableaux}

A tuple of non-negative integers $\nu=(\nu_1,\dots,\nu_l)$ such that $\nu_1+\dots+\nu_l=n$ is called a composition of $n$. The integer $\nu_1,\dots,\nu_l$ are called the \emph{parts} of $\nu$. We make no difference between $\nu$ and the same tuple with some parts equal to $0$ added at the end. The length of a composition $\ell(\nu)$ is the largest integer $k$ such that $\nu_k\neq 0$. We note $\nu\models n$ for ``$\nu$ is a composition of $n$''.

A partition $\lambda$ of $n$ is a composition of $n$ with non-increasing parts. So it is a tuple of integers $\lambda=(\lambda_1,\dots,\lambda_l)$ such that $\lambda_1\geq\lambda_2\geq\dots\geq\lambda_l\geq 0$ and $\lambda_1+\dots+\lambda_l=n$. For partitions, the length $\ell(\lambda)$ is simply the number of non-zero parts. We note $\lambda\vdash n$ for ``$\lambda$ is a partition of $n$''.

A pair $(x,y)\in\mathbb{Z}^2$ is called a {\em node}. The Young diagram (or Ferrer diagram) of a partition $\lambda=(\lambda_1,\dots,\lambda_l)$ is the set of nodes $(x,y)$ such that $x\in\{1,\dots,l\}$ and $y\in\{1,\dots,\lambda_x\}$. The Young diagram of $\lambda$ will be seen as a left-justified array of $l$ rows such that the $j$-th row contains $\lambda_j$ nodes for all $j=1,\dots,l$ (note that we count rows from top to bottom). A node will often be pictured by an empty box.

For a partition $\lambda$, a Young tableau of shape $\lambda$ is a map from the set of nodes of $\lambda$ to $\mathbb{Z}_{\geq1}$. It is represented by filling the nodes of the diagram of $\lambda$ by numbers in $\mathbb{Z}_{\geq1}$. The size of a Young tableau is $n$ if its shape is a partition of $n$ (the size is the number of nodes). For brevity, we will sometimes call a Young tableau of shape $\lambda$ a $\lambda$-tableau.

A Young tableau of size $n$ is called standard if the map from the set of nodes is a bijection with $\{1,\dots,n\}$ and if moreover the numbers are strictly ascending along rows and down columns of the Young diagram. We set
\[\STab(\lambda)=\{\text{standard $\lambda$-tableaux}\}\ \quad\text{and}\ \quad \STab(n)=\bigcup_{\lambda\vdash n}\STab(\lambda)\ .\]

A Young tableau $\bT$ is called semistandard if the numbers are weakly ascending along rows and strictly ascending down columns of the Young diagram. For $a\in\mathbb{Z}_{\geq1}$, let $\nu_a$ be the number of times the integer $a$ appears in the tableau $\bT$. The sequence $\nu=(\nu_1,\nu_2,\dots)$ forms a composition of $n$. We say that $\bT$ is a semistandard tableau of weight $\nu$. we set
\[\SSTab(\lambda,\nu):=\{\text{semistandard $\lambda$-tableaux of weight $\nu$}\}\ .\]
For example, a standard Young tableau is a semistandard Young tableau of weight $(1,\dots,1)$.

An important subset of semistandard tableaux are those where the numbers in the nodes are all in $\{1,\dots,N\}$ for a certain integer $N$.  So we set:
\[\SSTab_N(\lambda)=\{\text{semistandard $\lambda$-tableaux with numbers $\leq N$}\}\ \quad\text{and}\ \quad \SSTab_N(n)=\bigcup_{\lambda\vdash n}\SSTab_N(\lambda)\ .\]
Equivalently, this means that the weight $\nu$ of the tableau is not fixed, but is required to satisfy $\ell(\nu)\leq N$. That is, we have
$\SSTab_N(\lambda)=\bigcup\SSTab(\lambda,\nu)$ where the union is over the compositions $\nu$ of $n$ with $\ell(\nu)\leq N$

Here are two examples:
\[\bT=\begin{array}{ccc}
\fbox{\scriptsize{$1$}} & \hspace{-0.35cm}\fbox{\scriptsize{$2$}} & \hspace{-0.35cm}\fbox{\scriptsize{$4$}} \\[-0.2em]
\fbox{\scriptsize{$3$}} & &
\end{array}\ \qquad \ \text{and}\ \qquad\ \bbT=\begin{array}{ccc}
\fbox{\scriptsize{$1$}} & \hspace{-0.35cm}\fbox{\scriptsize{$2$}} & \hspace{-0.35cm}\fbox{\scriptsize{$3$}} \\[-0.2em]
\fbox{\scriptsize{$2$}} & &
\end{array}\,.\]
Both are $\lambda$-tableaux with $\lambda=(3,1)$, and $\bT$ is standard while $\bbT$ is semistandard. We have $\bbT\in\SSTab_3(\lambda)$ (and in fact belongs in $\SSTab_N(\lambda)$ for any $N\geq 3$). More precisely, the weight of $\bbT$ is the composition $\nu=(1,2,1)$, which is of length $3$.

\paragraph{Dominance order and Kostka numbers.} For a partition $\lambda=(\lambda_1,\dots,\lambda_l)$, we use the convention that $\lambda_{l+1}=\lambda_{l+2}=\dots=0$. For two partitions $\lambda,\mu$ of the same size, we denote 
$$\lambda\geq\mu\ \ \ \Longleftrightarrow\ \ \ \lambda_1+\dots+\lambda_i\geq \mu_1+\dots+\mu_i\,,\ \ \ \forall i.$$
This is the dominance ordering of partitions. 

For a partition $\lambda$ of $n$ and a composition $\nu$ of $n$, the number of semistandard Young tableaux of shape $\lambda$ and of weight $\nu$ is called a \emph{Kostka number}:
\[K_{\lambda,\nu}:=|\SSTab(\lambda,\nu)|\ .\]
A fundamental property is that $K_{\lambda,\nu}$ does not depend on the ordering of the parts of $\nu$ (see \cite[Theorem 7.10.2]{St} for a direct combinatorial proof). As a consequence, one can check that this implies another fundamental property of Kostka numbers, relating them with the dominance order:
$$K_{\lambda,\nu}\neq 0\ \ \ \ \ \Longleftrightarrow\ \ \ \ \ \lambda\geq \nu^{\text{ord}}\ ,$$
where $\nu^{\text{ord}}$ is the partition obtained from $\nu$ by reordering the parts in decreasing order.

\begin{exam}
Let $\lambda=(4,4)$ and $\nu=(3,2,2,1)$ then 
\[\SSTab(\lambda,\nu):=\Bigl\{ \begin{array}{cccc}
\fbox{\scriptsize{$1$}} & \hspace{-0.35cm}\fbox{\scriptsize{$1$}} & \hspace{-0.35cm}\fbox{\scriptsize{$1$}} & \hspace{-0.35cm}\fbox{\scriptsize{$2$}} \\[-0.2em]
\fbox{\scriptsize{$2$}} & \hspace{-0.35cm}\fbox{\scriptsize{$3$}} & \hspace{-0.35cm}\fbox{\scriptsize{$3$}} & \hspace{-0.35cm}\fbox{\scriptsize{$4$}}
\end{array}\ ,\qquad \begin{array}{cccc}
\fbox{\scriptsize{$1$}} & \hspace{-0.35cm}\fbox{\scriptsize{$1$}} & \hspace{-0.35cm}\fbox{\scriptsize{$1$}} & \hspace{-0.35cm}\fbox{\scriptsize{$3$}} \\[-0.2em]
\fbox{\scriptsize{$2$}} & \hspace{-0.35cm}\fbox{\scriptsize{$2$}} & \hspace{-0.35cm}\fbox{\scriptsize{$3$}} & \hspace{-0.35cm}\fbox{\scriptsize{$4$}}
\end{array}\Bigr\}\ .\]
If $\nu'=(2,3,1,2)$ (we have permuted the parts of $\nu$) then:
\[\SSTab(\lambda,\nu'):=\Bigl\{ \begin{array}{cccc}
\fbox{\scriptsize{$1$}} & \hspace{-0.35cm}\fbox{\scriptsize{$1$}} & \hspace{-0.35cm}\fbox{\scriptsize{$2$}} & \hspace{-0.35cm}\fbox{\scriptsize{$2$}} \\[-0.2em]
\fbox{\scriptsize{$2$}} & \hspace{-0.35cm}\fbox{\scriptsize{$3$}} & \hspace{-0.35cm}\fbox{\scriptsize{$4$}} & \hspace{-0.35cm}\fbox{\scriptsize{$4$}}
\end{array}\ ,\qquad \begin{array}{cccc}
\fbox{\scriptsize{$1$}} & \hspace{-0.35cm}\fbox{\scriptsize{$1$}} & \hspace{-0.35cm}\fbox{\scriptsize{$2$}} & \hspace{-0.35cm}\fbox{\scriptsize{$3$}} \\[-0.2em]
\fbox{\scriptsize{$2$}} & \hspace{-0.35cm}\fbox{\scriptsize{$2$}} & \hspace{-0.35cm}\fbox{\scriptsize{$4$}} & \hspace{-0.35cm}\fbox{\scriptsize{$4$}}
\end{array}\Bigr\}\ .\]
Finally, we take a composition with a part equal to $5$, say for example $\nu=(2,5,1)$. Here $\nu^{\text{ord}}=(5,2,1)$ and we do not have $\lambda\geq\nu^{\text{ord}}$. Therefore the set $\SSTab(\lambda,\nu)$ is empty. Here it is easy to see that one cannot put $5$ times the same numbers in the diagram $\lambda$ without having twice the same in a column.
\end{exam}

\subsection*{\thechapter.A.2 $\ $Seminormal representations of $S_n$ and $H_n(q)$}
\addcontentsline{toc}{subsection}{\thechapter.A.2 $\ $Seminormal representations of $S_n$ and $H_n(q)$}

We refer for example to \cite{Ho,IO1,OV,Ra2,Rut,Yo}. 

\paragraph{The Hecke algebra.} First we define the sequence of contents of a standard Young tableau $\bT$. The parameter $q$ is such that $q^2$ is not a root of unity. Let $\theta_i$ be the node of $\bT$ with number $i$ and denote the coordinate of $\theta_i$ by $(x_i,y_i)$. That is, in $\bT$, the number $i$ is in the line $x_i$ and the column $y_i$. We define the classical content $\cc_i(\bT)=y_i-x_i$ and the $q$-content $\qc_i(\bT)=q^{2\cc_i(\bT)}$. Here is an example of a standard Young tableau with its sequence of $q$-contents:
\[\bT=\begin{array}{ccc}
\fbox{\scriptsize{$1$}} & \hspace{-0.35cm}\fbox{\scriptsize{$2$}} & \hspace{-0.35cm}\fbox{\scriptsize{$4$}} \\[-0.2em]
\fbox{\scriptsize{$3$}} & &
\end{array}\ : \qquad \qc_1(\bT)=1\,,\ \ \ \qc_2(\bT)=q^2\,,\ \ \ \qc_3(\bT)=q^{-2}\,,\ \ \ \qc_4(\bT)=q^4\,.\]

Now let $\lambda\vdash n$ and let $S_{\lambda}$ be a vector space with a basis $\{v_{\bT}\}_{\bT\in\STab(\lambda)}$ indexed by the standard Young tableaux of shape $\lambda$. The following formula for the generators $\si_1,\dots,\si_{n-1}$ defines an irreducible representation of the Hecke algebra $H_{n}(q)$ on the space $S_{\lambda}$:
\begin{equation}\label{rep-si}
\si_{i}(v_{\bT})=\frac{(q-q^{-1})\qc_{i+1}(\bT)}{\qc_{i+1}(\bT)-\qc_{i}(\bT)}\,v_{\bT}+\frac{q\,\qc_{i+1}(\bT)-q^{-1}\qc_{i}(\bT)}{\qc_{i+1}(\bT)-\qc_{i}(\bT)}\,v_{s_i(\bT)}\ ,\ \ \ \quad\text{for $i\in\{1,\dots,n-1\}$.}
\end{equation}
where $s_i(\bT)$ is the Young tableau obtained from $\bT$ by exchanging $i$ and $i+1$. The tableau $s_i(\bT)$ is not necessarily standard and we set $v_{\bT'}:=0$ for any non-standard Young tableau $\bT'$. Note that since $q^2$ is not a root of unity, then $\qc_{i}(\bT)\neq\qc_{i+1}(\bT)$ for any $\bT\in\STab(\lambda)$ so that Formula (\ref{rep-si}) is valid.

As $\lambda$ runs over the set of partitions of $n$, the representations $S_{\lambda}$ are pairwise non-isomorphic and exhaust the set of irreducible representations of $H_{n}(q)$. The irreducible representations $S_{\lambda}$ are the ones appearing in the quantum Schur--Weyl duality. In particular, if $\lambda$ is a single row of boxes, the representation $S_\lambda$ is the one-dimensional representation given by $\si_i\mapsto q$, and if $\lambda$ is a single column of boxes, the representation $S_\lambda$ is the one-dimensional representation given by $\si_i\mapsto -q^{-1}$.

\paragraph{The symmetric group.} If we denote by $d_{i,j}(\bT):=\cc_{j}(\bT)-\cc_i(\bT)$ the axial distance between the nodes with number $j$ and $i$ in the Young tableau $\bT$, then Formula (\ref{rep-si}) can be written in terms of $q$-numbers $[d]_q=\frac{q^d-q^{-d}}{q-q^{-1}}$ as follows:
\begin{equation}\label{rep-si2}
\si_{i}(v_{\bT})=\frac{q^{d_{i,i+1}(\bT)}}{[d_{i,i+1}(\bT)]_q}\,v_{\bT}+\frac{[d_{i,i+1}(\bT)+1]_q}{[d_{i,i+1}(\bT)]_q}\,v_{s_i(\bT)}\ ,\ \ \ \quad\text{for $i\in\{1,\dots,n-1\}$.}
\end{equation}
This shows that this formula can be extended to the case where $q^2=1$, since the limit when $q^2=1$ of a $q$-number $[d]_q$ is simply $d$. The formula (\ref{rep-si2}) when $q^2=1$ gives the irreducible representations $S_{\lambda}$ of the symmetric group. They are pairwise non-isomorphic and exhaust the set of irreducible representations of $H_{n}(q)$. Again, these representations $S_{\lambda}$ are the ones appearing in the Schur--Weyl duality. For $\lambda$ a single row of boxes, the representation $S_\lambda$ is the trivial representation, while for $\lambda$ a single column of boxes, the representation $S_\lambda$ is the sign representation.

\subsection*{\thechapter.A.3 $\ $Branching rules and Bratteli diagram}
\addcontentsline{toc}{subsection}{\thechapter.A.3 $\ $Branching rules and Bratteli diagram}

The branching rules for the Hecke algebras and for the symmetric groups are the same (again, $q^2$ is not a root of unity when we consider $H_n(q)$). They are given by:
\begin{equation}\label{BR}
\text{Res}_{H_{n-1}(q)}(V_{\lambda})\cong \bigoplus_{\begin{array}{c}
\\[-1.6em]
\scriptstyle{\mu\,\vdash n-1} \\[-0.4em]
\scriptstyle{\mu\subset \lambda}
\end{array}} V_{\mu}\ .
\end{equation}
The first levels of the Bratteli diagram are given below. The vertices in level $n$ are indexed by the (isomorphism classes of) irreducible representations of $H_n(q)$ (equivalently, of $S_n$), and the edges express the branching rules: they connect $\lambda$ from level $n$ to $\mu$ from level $n-1$ if and only if $S_{\mu}$ appears in the restriction of $S_{\lambda}$. There is no multiple edges since the multiplicity are always $0$ or $1$. The resulting diagram is simply the Hasse diagram of the poset of partitions, ordered by inclusion (the so called Young lattice).

The shaded areas are meant to express visually the full statement of the Schur--Weyl duality. Namely, by deleting the vertices included in the shaded area labelled $gl_N$, the remaining vertices correspond to the representations appearing in the decomposition of the tensor product $V^{\otimes n}$, where $\dim(V)=N$. In other words, removing the shaded area (and the edges connected to vertices in this area) we obtain the Bratteli diagram of the centralisers of $GL(V)$ or of $U_q(gl_N)$. For example, if $N=2$, the resulting diagram, involving only partitions with no more than two lines, is the Bratteli diagram of the Temperley--Lieb algebras. 

\begin{center}
 \begin{tikzpicture}[scale=0.3]
\node at (0.5,4) {$\emptyset$};
\draw ( 0.5,3) -- (0.5, 1);
\diag{0}{0}{1};\node at (-1,-0.5) {$1$};

\draw (-0.5,-1.5) -- (-3,-3.5);\draw (1.5,-1.5) -- (4,-3.5);

\diag{-4}{-4}{2};\node at (-5,-4.5) {$1$}; \diagg{4}{-4}{1}{1};\node at (3,-5) {$1$};

\draw (-3.5,-5.5) -- (-6.5,-8.5);\draw (-2.5,-5.5) -- (0.5,-8.5);\draw (3.5,-6.5) -- (1.5,-8.5);\draw (5.5,-6.5) -- (8.5,-8.5);

\diag{-8}{-9}{3};\node at (-9,-9.5) {$1$};\diagg{0}{-9}{2}{1};\node at (-1,-10) {$2$};\diaggg{8}{-9}{1}{1}{1};\node at (7,-10.5) {$1$};

\draw (-7.5,-10.5) -- (-14,-14.5);\draw (-6.5,-10.5) -- (-6.5,-14.5);\draw (-0.5,-11.5) -- (-4.5,-14.5);\draw (1,-11.5) -- (1,-14.5);\draw (2.5,-11.5) -- (7,-14.5);
\draw (8.5,-12.5) -- (8.5,-14.5);\draw (9.5,-12.5) -- (15.5,-14.5);

\diag{-16}{-15}{4};\node at (-17,-15.5) {$1$};\diagg{-8}{-15}{3}{1};\node at (-9,-16) {$3$};\diagg{0}{-15}{2}{2};\node at (-1,-16) {$2$};\diaggg{8}{-15}{2}{1}{1};\node at (7,-16.5) {$3$};\diagggg{15}{-15}{1}{1}{1}{1};\node at (14,-17) {$1$};

\draw (-15,-16.5)--(-17.5,-20.5);\draw (-13,-16.5)--(-10,-20.5);\draw (-7,-17.5)--(-8.5,-20.5);\draw (-6,-17.5)--(-3,-20.5);\draw (-5,-17.5)--(9,-20.5);
\draw (0,-17.5)--(-1,-20.5);\draw (2,-17.5)--(4,-20.5);
\draw (8,-18.5)--(6,-20.5);\draw (9,-18.5)--(10,-20.5);\draw (10,-18)--(15,-20.5);\draw (15.5,-19.5)--(16,-20.5);\draw (16.5,-19)--(20,-20.5);

\diag{-20}{-21}{5};\node at (-21,-21.5) {$1$};\diagg{-11}{-21}{4}{1};\node at (-12,-22) {$4$};\diagg{-3}{-21}{3}{2};\node at (-4,-22) {$5$};
\diaggg{4}{-21}{2}{2}{1};\node at (3,-22.5) {$5$};\diaggg{9}{-21}{3}{1}{1};\node at (8,-22.5) {$6$};\diagggg{15}{-21}{2}{1}{1}{1};\node at (14,-23) {$4$};\diaggggg{20}{-21}{1}{1}{1}{1}{1};\node at (19,-23.5) {$1$};

\draw[thin, fill=gray,opacity=0.2] (2.5,-18.5)..controls +(0,17) and +(0,17) .. (27,-18.5) .. controls +(0,-17) and +(0,-17) .. (2.5,-18.5);\node at (21.5,-6) {$gl_2$};
\draw[thin, fill=gray,opacity=0.2] (13,-21)..controls +(0,12) and +(0,12) .. (24,-21) .. controls +(0,-10) and +(0,-10) .. (13,-21);\node at (19,-11) {$gl_3$};
\draw[thin, fill=gray,opacity=0.2] (18.5,-23.5)..controls +(0,5) and +(0,5) .. (22.5,-23) .. controls +(0,-5) and +(0,-5) .. (18.5,-23.5);\node at (21,-18.5) {$gl_4$};

\node at (-27,-0.5) {$n=1$};\node at (-27,-4.5) {$n=2$};\node at (-27,-9.5) {$n=3$};\node at (-27,-15.5) {$n=4$};\node at (-27,-21.5) {$n=5$};

\end{tikzpicture}
\end{center}

\setcounter{equation}{0}
\DeactivateToc
\chapter{\huge{Fused Hecke algebras}}\label{chap-fus}
\ActivateToc
\addcontentsline{toc}{chapter}{\large{Chapter \thechapter. \hspace{0.2cm}Fused Hecke algebras\vspace{0.3cm}}}

\setcounter{minitocdepth}{1}
\minitoc

\section*{Introduction}
\addcontentsline{toc}{section}{Introduction}

In the first chapter, we have encountered algebras with the following features:
\begin{itemize}
\item[$\bullet$] they are quotients of the braid group algebra;
\item[$\bullet$] they contain abstract solutions of the Yang--Baxter equation;
\item[$\bullet$] they admit representations on vector spaces of the form $V^{\otimes n}$.
\end{itemize}
The first and famous example was the Hecke algebra\footnote{Actually, the very first example should be the group algebra of the symmetric group, whose deformation is the Hecke algebra. From the point of view of the braid group, the symmetric group is a bit too simple.}, but we also had other explicit examples, namely, the Temperley--Lieb algebra and the Birman--Murakami--Wenzl algebra. These examples turn out to be the algebras behind the following invariants of links: the Jones polynomial, the HOMFLYPT polynomial and the Kauffman polynomial. 

The third item in the list above was in fact the key to unifying these examples. Indeed, through their representations on tensor spaces $V^{\otimes n}$, these algebras can be seen as centralisers of representations of some interesting algebras: quantum groups. Then we saw that indeed the quantum groups are cleverly designed algebras which produce, through their centralisers, algebras fitting the list above.

\vskip .2cm
The goal of this chapter is to describe more general centralisers of quantum groups representations, in order to obtain other algebras fitting the wish list above. The algebras will go through the name of fused Hecke algebras. Their classical limit, which are algebras with a combinatorial flavour, are called the algebras of fused permutations and are simpler and also interesting. 

The name ``fused'' comes from the fact that more or less well-hidden behind our construction is the so-called fusion procedure for the Yang--Baxter equation. The fusion procedure was designed at the matrix level to produce new solutions of the Yang--Baxter equation, starting from a known one. Our starting point for this work was actually the desire to understand some solutions of the Yang--Baxter equation (the solutions corresponding to higher spin representations of $su_2$) as coming from an abstract formula in some algebras. The relevant algebras are the fused Hecke algebras.

Once we have at hand a centraliser algebra such as the fused Hecke algebra, we know from the general picture discussed in Chapter \ref{chap-prel} that a solution of the Yang--Baxter equation should be hiding somewhere. Of course, finding the explicit formula giving this solution is another matter. It turns out that we were able to find this formula for the fused Hecke algebra, and it has a nice form using $q$-numbers (it is given in the last section of this chapter).

Along our way to the solution of the Yang--Baxter equation, we have been largely diverted from this primary goal. The representation-theoretic features of the fused Hecke algebra $H_{\bk,n}$ were too tempting to be left aside. Naturally from our starting point, the algebra $H_{\bk,n}$ is in Schur--Weyl duality with some representations of $U_q(gl_N)$. One important point is that the algebra $H_{\bk,n}$ does not depend on $N$ and in some sense interpolates between, or lives above, the centralisers for any $N$. This raises the question of the quotient which has to be made given a specific $N$. For this we have conjectures and partial results.

Besides, it is also natural to ask if the representation theory of $H_{\bk,n}$ can be derived from scratch (that is, without reference to the quantum groups $U_q(gl_N)$). In the Schur--Weyl duality spirit, this would provide an alternative point of view on some results on $U_q(gl_N)$-representation, such as the decomposition of some tensor products. We succeeded in describing the representation theory of $H_{\bk,n}$, in the semisimple situation, using solely the known representation theory of the usual Hecke algebra.

Finally, with an eye to the connections with the braid group and the link invariants, we provide a topological definition of the fused Hecke algebra using the notion of fused braids and a certain way of multiplying them. This definition allows to relate naturally the objects of our algebra with the combinatorial structures appearing for $q=1$, and also allows to give a nice and rather simple (conjectural) description of the ideal allowing to obtain the genuine centraliser. The fused Hecke algebra is a natural generalisation of the Hecke algebra and the quotient by the afore-mentioned ideal is a natural generalisation of the Temperley--Lieb algebra.

\paragraph{References.} The publications relevant for this chapter, where more details can be found, are \cite{P5,publi-CP2,publi-CP1}.

\section{Extending the Schur--Weyl duality to symmetrised powers}

The integer $N>1$ is fixed: it is the dimension of the fundamental representation $V$ of $U_q(gl_N)$. All what is below can be read in particular for $q=1$.

\subsection{Fused Hecke algebra and $q$-symmetrised powers}

\paragraph{The $q$-symmetrized powers.} Recall that for any partition $\lambda$ with $\ell(\lambda)\leq N$, we have an irreducible representation $V^N_{\lambda}$ of $U_q(gl_N)$. Here we will consider the special case of a single-row partition, $\lambda=(k)$ for some integer $k\geq 1$. In this case, we have the following description of $V^N_{(k)}$, using the quantum Schur--Weyl duality. Recall that the $k$-fold tensor product of the fundamental representation $V$ of $U_q(gl_N)$ decomposes as follows:
\[V^{\otimes n}=\bigoplus_{\substack{\lambda\vdash k\\[0.2em] \ell(\lambda)\leq N}} V_{\lambda}^N\otimes S_{\lambda}\,.\]
For $\lambda=(k)$, the representation $S_{(k)}$ of the Hecke algebra $H_n(q)$ is one-dimensional ($\si_i\mapsto q$). So looking at the above decomposition from the point of view of $H_n(q)$, we see that we can have access to the subspace $V_{(k)}^N$ using the projector of $H_n(q)$ onto the irreducible representation $S_{(k)}$. This is what the $q$-symmetriser does.

\vskip .2cm
The \textbf{$q$-symmetriser} in $H_k(q)$ is the following element (where all notations were defined around (\ref{quantum-numbers})):
\begin{equation}\label{def-P}
P_k=\frac{\sum_{w\in S_k}q^{\ell(w)}\si_w}{\sum_{w\in S_k}q^{2\ell(w)}}=\frac{q^{-k(k-1)/2}}{[k]_q!}\sum_{w\in S_k}q^{\ell(w)}\si_w\ .
\end{equation}
We take this opportunity to emphasize the nice formula $\sum_{w\in S_k}q^{2\ell(w)}=q^{k(k-1)/2}[k]_q!$\,.

This element is the orthogonal central idempotent of $H_n(q)$ corresponding to the one-dimensional representation $S_{(k)}$. Thus, we define:
\[S^kV=P_k(V^{\otimes k})\,,\]
we call the space $S^kV$ the $k$-th ``$q$-symmetrised power'' of $V$. This is a representation of $U_q(gl_N)$ which coincides with the representation $V_{(k)}^N$ discussed above.

\vskip .2cm
If $q=1$ then $P_k$ is the usual symmetriser (the normalized sum over all elements of $S_k$), and the space $S^kV$ is the usual $k$-th symmetrised power of $V$.

\paragraph{The tensor product.} Now let $\bk=(k_1,k_2,...)\in \mathbb{Z}_{\geq0}^{\infty}$ be an infinite sequence of non-negative integers, and let $n\in\mathbb{Z}_{\geq0}$ (we consider an infinite sequence $\bk$ so that $n$ can vary as we want). We consider here the following tensor product of representations of $U_q(gl_N)$:
\[S^{k_1}V\otimes\dots\otimes S^{k_n}V\ .\]
The constant sequence $\bk=(1,1,\dots)$ corresponds to the usual case $V\otimes\dots\otimes V$. The constant sequences $\bk=(k,k,\dots)$ for a fixed integer $k$ are of special interest, but we will deal with the general case of mixed tensor products.

\paragraph{The fused Hecke algebra.} Let us consider the ``big'' Hecke algebra $H_{k_1+\dots+k_n}(q)$. There is a natural (parabolic) subalgebra which is isomorphic to $H_{k_1}(q)\otimes\dots\otimes H_{k_n}(q)$. Thus we can define the following element of $H_{k_1+\dots+k_n}(q)$ by plugging the $q$-symmetriser in each factor:
\begin{equation}\label{def-Pkn} P_{\bk,n}=P_{k_1}\otimes\dots\otimes P_{k_n}\ .
\end{equation}
Alternatively, an explicit expression for $P_{\bk,n}$ is:
\[P_{\bk,n}=\prod_{i=1}^n\frac{q^{-k_i(k_i-1)/2}}{[k_i]_q!}\sum_{w\in S_{k_1}\times\dots\times S_{k_n}}q^{\ell(w)}\si_w\ ,\]
where this time we refer to the parabolic subgroup $S_{k_1}\times\dots\times S_{k_n}$ of the symmetric group $S_{k_1+\dots+k_n}$. The element $P_{\bk,n}$ is an idempotent of $H_{k_1+\dots+k_n}(q)$, the prefactor ensuring the correct normalization in order to have $P_{\bk,n}^2=P_{\bk,n}$.

Now we can define, algebraically, the fused Hecke algebra.
\begin{frameddefi}\label{def-Hkn}
The fused Hecke algebra $H_{\bk,n}(q)$ is the algebra:
\[H_{\bk,n}(q)=P_{\bk,n}H_{k_1+\dots+k_n}(q)P_{\bk,n}\ .\]
\end{frameddefi}
Some comments may be needed. For any algebra $A$ and any idempotent $e$ in $A$, the subset $eAe$ forms an algebra. It may be seen as a subalgebra of $A$, in the sense that the multiplication coincides, but one has to be careful that the unit element of $eAe$ is $e$, and not the unit element of $A$. Such idempotent subalgebras are very frequent in representation theory, see for example \cite[\S 6]{Gr} where one can find the following standard results. First, by construction, if the vector space $M$ is a representation of $A$, then $eM$ (which is possibly $\{0\}$) is a representation of $eAe$. Moreover, if $A$ is semisimple and $M_1,\dots,M_k$ is a complete set of irreducible representations, then the set of non-zero representations among $eM_1,\dots,eM_k$ is a complete set of irreducible representations of $eAe$.

\vskip .2cm
If $q=1$ then the idempotent $P_{\bk,n}$ is simply the normalized sum over all elements of the subgroup $S_{k_1}\times\dots\times S_{k_n}$ of $S_{k_1+\dots+k_n}$, and the fused Hecke algebra $H_{\bk,n}(1)$ is the algebra:
\[H_{\bk,n}(1)=P_{\bk,n}S_{k_1+\dots+k_n}P_{\bk,n}\ .\]
In this case, we will call the algebra $H_{\bk,n}(1)$ the algebra of fused permutations.

\paragraph{Generalisation of the Schur--Weyl duality.} Almost by construction, we have the fused Hecke algebra $H_{\bk,n}(q)$ acting on the desired tensor product of $U_q(gl_N)$-representations. Namely we have the two representations:
\[U_q(gl_N)\ \ \stackrel{\rho}{\longrightarrow}\ \ \ \ \text{End}(S^{k_1}V\otimes\dots\otimes S^{k_n}V)\ \ \ \ \stackrel{\pi}{\longleftarrow}\ \ H_{\bk,n}(q)\ .\]
Now we can state the generalisation of (the first part of) the Schur--Weyl duality.
\begin{framedtheo}[\cite{publi-CP1}]\label{theo-SWk1}
The centraliser $\text{End}_{U_q(gl_N)}(S^{k_1}V\otimes\dots\otimes S^{k_n}V)$ is the image of the fused Hecke algebra $H_{\bk,n}(q)$:
\begin{equation}\label{gen-qSW1}
\text{End}_{U_q(gl_N)}(S^{k_1}V\otimes\dots\otimes S^{k_n}V)=\pi\bigl(H_{\bk,n}(q)\bigr)\ .
\end{equation}
\end{framedtheo}
We are in the semisimple situation ($q$ is generic) so we have actually that $\pi\bigl(H_n(q)\bigr)$ and $\rho\bigl(U_q(gl_N)\bigr)$ are the mutual centralisers of each other.

\subsection{A small digression on parabolic Hecke algebras}

Let $W$ be a Coxeter group with set of generators $S$. Let $J$ be a subset of $S$ and $W_J$ the subgroup, called a parabolic subgroup of $W$, generated by the subset $J$ of the generators. We assume that $W_J$ is finite. All that is below is valid for $q=1$. We refer to \cite{APV,Cur}.

\paragraph{Definition of the parabolic Hecke algebra $\HJ(W)$.} 
We consider the following element of the Hecke algebra $H(W)$ (see Chapter \ref{chap-prel}), associated to $W$:
\[P_{J}=\frac{1}{W_J(q^2)}\sum_{w\in W_J}q^{\ell(w)}\si_w\ ,\]
where the normalization factor $W_J(q^2)=\sum_{w\in W_J}q^{2\ell(w)}$ is chosen such that $P_J$ is an idempotent. The element $P_J$ is the $q$-symmetriser of the Hecke algebra $H(W_J)$, which we see here as the subalgebra of $H(W)$ generated by $T_s$ with $s\in J$.
\begin{defi}
The parabolic Hecke algebra $\HJ(W)$ is the algebra defined by:
\[\HJ(W)=P_JH(W)P_J\ .\]
\end{defi}

\paragraph{Standard basis of $\HJ(W)$.} To obtain a basis of $\HJ(W)$ from the standard basis of $\{T_w\}$ of $H(W)$, we need to choose a representative of each double coset in $W_J\backslash W/W_J$.
\begin{prop}
Let $R_J$ be any set of representatives for the double cosets in $W_J\backslash W/W_J$. Then $\{P_J\si_xP_J\}_{x\in R_J}$ is a basis of $\HJ(W)$.
\end{prop}
There is a most standard choice for coset representatives which is as follows. For each double coset in $W_J\backslash W/W_J$, there is a unique element of minimal length in the coset \cite[chap. 2]{GP}. We denote $X_{JJ}$ the set of minimal-length representatives for the double cosets of $W_J$ in $W$. These distinguished coset representatives are characterized by their descent sets:
\[x\in X_{JJ}\ \ \ \ \Leftrightarrow\ \ \ \ \ \mathcal{L}(x)\cap J=\mathcal{R}(x)\cap J=\emptyset\ ,\]
where the left descent $\mathcal{L}(x)$ of an element $x\in W$ consists of the simple transpositions $s\in J$ such that $\ell(sx)<\ell(x)$, and similarly for the right descent $\mathcal{R}(x)$.

The standard basis of $\HJ(W)$ is thus
\[\{P_J\si_xP_J\}_{x\in X_{JJ}}\ .\]

\section{Fused permutations and fused braids}

The fused Hecke algebra involves ``topological objects'' having strands which can cross by passing over or below each other. For $q=1$, all topological information disappears and we obtain an algebra defined purely combinatorially. This is another example of a classical situation as shown in the following table.
\begin{center}
\begin{tabular}{c|c}
Topological Algebras & Combinatorial Algebras\\
\hline & \\[-0.2em]
Hecke algebra & Symmetric group\\[0.5em]
BMW algebra & Brauer algebra\\[0.5em]
Fused Hecke algebra & Algebra of fused permutations
\end{tabular}
\end{center}
The connections between the two columns of the table is that the algebras on the left are flat deformations of the algebras on the right. Here ``flat'' means that the vector space remains the same while the multiplication is modified (with a parameter $q$). One interest is that we obtain a natural, or standard, basis for algebras on the left indexed by combinatorial objects (coming from the algebras on the right). Let us see how it works for the fused Hecke algebras.

\subsection{The algebra of fused permutations} 

\paragraph{Objects.} We place two horizontal rows of $n$ dots, one on top of another. And we connect the dots of the top row to the dots on the bottom row. We require the following: for each $a\in\{1,\dots,n\}$, there are $k_a$ edges which start from the $a$-th dot on top and there are $k_a$ edges which arrive at the $a$-th dot at the bottom. The total number of edges is then $k_1+\dots+k_n$.

The only information that matters is the following: which dots are connected to which other dots and by how many edges. More rigorously, take a diagram as above and let $a\in\{1,\dots,n\}$. There are $k_a$ edges starting from the $a$-th dot on top and we denote by $I_a$ the multiset of bottom dots reached by these edges. We consider two diagrams equivalent if their sequences of multisets $(I_1,\dots,I_n)$ coincide.

\begin{defi}
A fused permutation is an equivalence class of diagrams as explained above. We denote by $S^{fus}_{\bk,n}$ the set of fused permutations.
\end{defi}

\paragraph{Examples.} $\bullet$ If $\bk=(1,1,1,\dots)$ consists only of 1's then the set of fused permutations $\cD_{\bk,n}$ coincides with the set of permutations of $\{1,\dots,n\}$.

$\bullet$ Let $n=3$ and  take $k_1=2$ and $k_2=k_3=1$. There are 7 distinct fused permutations in $\cD_{\bk,3}$ and here is an example (we give a diagram and the corresponding sequence of multisets):
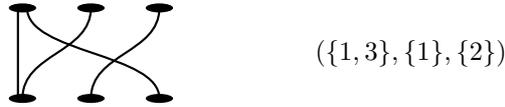
\begin{center}
\begin{tikzpicture}[scale=0.3]

\fill (21,2) ellipse (0.6cm and 0.2cm);\fill (21,-2) ellipse (0.6cm and 0.2cm);
\draw[thick] (20.8,2) -- (20.8,-2);\draw[thick] (21.2,2)..controls +(0,-2) and +(0,+2) .. (27,-2);  
\fill (24,2) ellipse (0.6cm and 0.2cm);\fill (24,-2) ellipse (0.6cm and 0.2cm);
\draw[thick] (24,2)..controls +(0,-2) and +(0,+2) .. (21,-2); 
\fill (27,2) ellipse (0.6cm and 0.2cm);\fill (27,-2) ellipse (0.6cm and 0.2cm);
\draw[thick] (27,2)..controls +(0,-2) and +(0,+2)..(24,-2);
\node at (38,0) {$(\{1,3\},\{1\},\{2\})$};

\end{tikzpicture}
\end{center}

$\bullet$ Let $n=2$ and  take $k_1=k_2=2$. We give below the three distinct fused permutations of $\cD_{\bk,3}$ (for each, we give a diagram and the corresponding sequence of multisets):

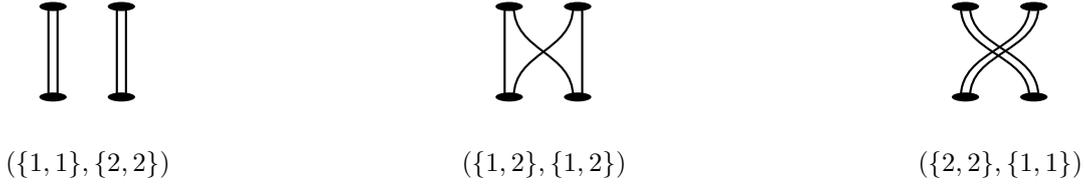
\begin{center}
\begin{tikzpicture}[scale=0.3]
\fill (1,2) ellipse (0.6cm and 0.2cm);\fill (1,-2) ellipse (0.6cm and 0.2cm);
\draw[thick] (0.8,2)..controls +(0,-2) and +(0,+2) .. (0.8,-2);\draw[thick] (1.2,2)..controls +(0,-2) and +(0,+2) .. (1.2,-2);  
\fill (4,2) ellipse (0.6cm and 0.2cm);\fill (4,-2) ellipse (0.6cm and 0.2cm);
\draw[thick] (3.8,2)..controls +(0,-2) and +(0,+2) .. (3.8,-2);\draw[thick] (4.2,2)..controls +(0,-2) and +(0,+2) .. (4.2,-2);
\node at (2.5,-5) {$(\{1,1\},\{2,2\})$};

\fill (21,2) ellipse (0.6cm and 0.2cm);\fill (21,-2) ellipse (0.6cm and 0.2cm);
\draw[thick] (20.8,2) -- (20.8,-2);\draw[thick] (21.2,2)..controls +(0,-2) and +(0,+2) .. (23.8,-2);  
\fill (24,2) ellipse (0.6cm and 0.2cm);\fill (24,-2) ellipse (0.6cm and 0.2cm);
\draw[thick] (23.8,2)..controls +(0,-2) and +(0,+2) .. (21.2,-2); \draw[thick] (24.2,2)..controls +(0,-2) and +(0,+2) .. (24.2,-2);
\node at (22.5,-5) {$(\{1,2\},\{1,2\})$};

\fill (41,2) ellipse (0.6cm and 0.2cm);\fill (41,-2) ellipse (0.6cm and 0.2cm);
\draw[thick] (40.8,2)..controls +(0,-2.1) and +(0,+1.9) .. (43.8,-2);\draw[thick] (41.2,2)..controls +(0,-1.9) and +(0,+2.1) .. (44.2,-2);
\fill (44,2) ellipse (0.6cm and 0.2cm);\fill (44,-2) ellipse (0.6cm and 0.2cm);
\draw[thick] (43.8,2)..controls +(0,-1.9) and +(0,+2.1) .. (40.8,-2);\draw[thick] (44.2,2)..controls +(0,-2) and +(0,+2) .. (41.2,-2); 
\node at (42.5,-5) {$(\{2,2\},\{1,1\})$};
\end{tikzpicture}
\end{center}

$\bullet$ Let $n=3$ and  take $k_1=k_2=k_3=2$. There are 21 distinct fused permutations in ${\cD}_{\bk,3}$  and here are three examples (for each, we give a diagram and the corresponding sequence of multisets):
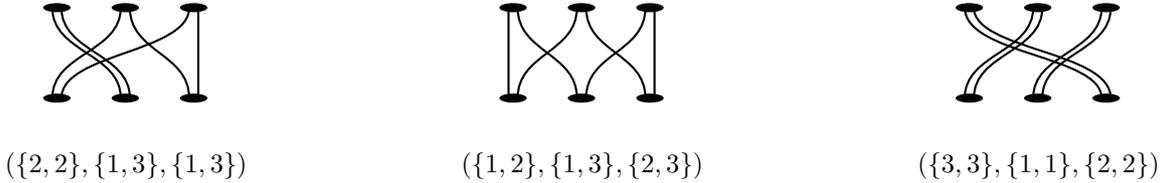
\begin{center}
\begin{tikzpicture}[scale=0.3]
\fill (1,2) ellipse (0.6cm and 0.2cm);\fill (1,-2) ellipse (0.6cm and 0.2cm);
\draw[thick] (0.8,2)..controls +(0,-2) and +(0,+2) .. (3.8,-2);\draw[thick] (1.2,2)..controls +(0,-2) and +(0,+2) .. (4.2,-2);  
\fill (4,2) ellipse (0.6cm and 0.2cm);\fill (4,-2) ellipse (0.6cm and 0.2cm);
\draw[thick] (3.8,2)..controls +(0,-2) and +(0,+2) .. (0.8,-2);\draw[thick] (4.2,2)..controls +(0,-2) and +(0,+2) .. (6.8,-2);
\fill (7,2) ellipse (0.6cm and 0.2cm);\fill (7,-2) ellipse (0.6cm and 0.2cm);
\draw[thick] (6.8,2)..controls +(0,-2) and +(0,+2) .. (1.2,-2);\draw[thick] (7.2,2) -- (7.2,-2);
\node at (4,-5) {$(\{2,2\},\{1,3\},\{1,3\})$};

\fill (21,2) ellipse (0.6cm and 0.2cm);\fill (21,-2) ellipse (0.6cm and 0.2cm);
\draw[thick] (20.8,2) -- (20.8,-2);\draw[thick] (21.2,2)..controls +(0,-2) and +(0,+2) .. (23.8,-2);  
\fill (24,2) ellipse (0.6cm and 0.2cm);\fill (24,-2) ellipse (0.6cm and 0.2cm);
\draw[thick] (23.8,2)..controls +(0,-2) and +(0,+2) .. (21.2,-2); \draw[thick] (26.8,2)..controls +(0,-2) and +(0,+2) .. (24.2,-2);
\fill (27,2) ellipse (0.6cm and 0.2cm);\fill (27,-2) ellipse (0.6cm and 0.2cm);
\draw[thick] (24.2,2)..controls +(0,-2) and +(0,+2) .. (26.8,-2);\draw[thick] (27.2,2) -- (27.2,-2);
\node at (24,-5) {$(\{1,2\},\{1,3\},\{2,3\})$};

\fill (41,2) ellipse (0.6cm and 0.2cm);\fill (41,-2) ellipse (0.6cm and 0.2cm);
\draw[thick] (40.8,2)..controls +(0,-2.1) and +(0,+1.9) .. (46.8,-2);\draw[thick] (41.2,2)..controls +(0,-1.9) and +(0,+2.1) .. (47.2,-2);
\fill (44,2) ellipse (0.6cm and 0.2cm);\fill (44,-2) ellipse (0.6cm and 0.2cm);
\draw[thick] (43.8,2)..controls +(0,-2) and +(0,+2) .. (40.8,-2);\draw[thick] (44.2,2)..controls +(0,-2) and +(0,+2) .. (41.2,-2); 
\fill (47,2) ellipse (0.6cm and 0.2cm);\fill (47,-2) ellipse (0.6cm and 0.2cm);
\draw[thick] (46.8,2)..controls +(0,-2) and +(0,+2) .. (43.8,-2);\draw[thick] (47.2,2)..controls +(0,-2) and +(0,+2) .. (44.2,-2);
\node at (44,-5) {$(\{3,3\},\{1,1\},\{2,2\})$};
\end{tikzpicture}
\end{center}

\paragraph{Multiplication.} We define the associative $\mathbb{C}$-algebra $H_{\bk,n}(1)$ as the $\bC$-vector space with basis indexed by the fused permutations in $S^{fus}_{\bk,n}$, and with the multiplication given as follows. Let $d,d'\in \cD_{\bk,n}$ and we identify $d$, respectively, $d'$, with a diagram representing it.
\begin{itemize}
\item \emph{(Concatenation)} We place the diagram of $d$ on top of the diagram of $d'$ by identifying the bottom row of dots of $d$ with the top row of dots of $d'$.
\item \emph{(Removal of middle dots)} For each $a\in\{1,\dots,n\}$, there are $k_a$ edges arriving and $k_a$ edges leaving the $a$-th dot in the middle row. We delete the $a$-th dot in the middle row and sum over all possibilities of connecting the $k_a$ edges arriving at it to the $k_a$ edges leaving from it (at the $a$-th edge, there are thus $k_a!$ possibilities).
 \item  \emph{(Normalisation)} We divide the resulting sum by $k_1!\dots k_n!$.
\end{itemize}
At the end of the procedure described above, we obtain a sum of diagrams representing a sum of fused permutations (with rational coefficients). This is what we define to be $dd'$ in $H_{\bk,n}(1)$. This diagrammatic multiplication is well-defined since the result clearly depends only on the equivalences classes of the diagrams.

The algebra $H_{\bk,n}(1)$ is an associative algebra with unit, the unit element is the fused permutation corresponding to the diagram with only vertical edges.

\paragraph{Examples.} $\bullet$ If $\bk=(1,1,1,\dots)$ consists only of 1's then the algebra $H_{\bk,n}(1)$ obviously coincides with the complex group algebra $\CC S_n$ of the symmetric group $S_n$ on $n$ letters.

$\bullet$ Here is an example of a product of two elements of $H_{\bk,2}(1)$ with $k_1=k_2=2$:
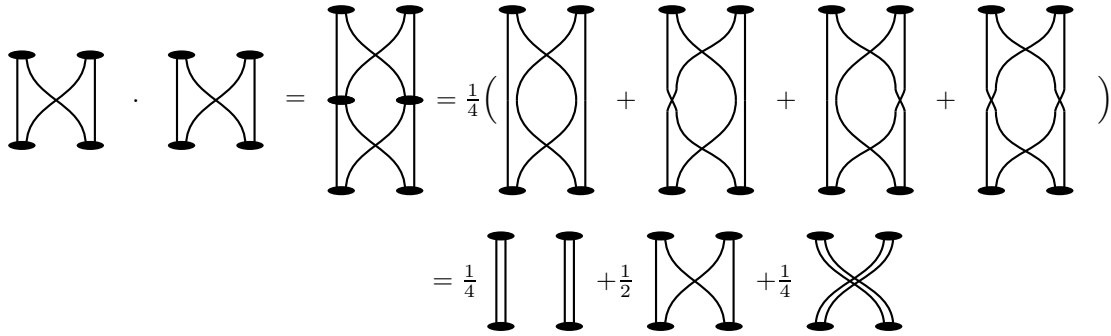
\begin{center}
\begin{tikzpicture}[scale=0.3]
\fill (1,2) ellipse (0.6cm and 0.2cm);\fill (1,-2) ellipse (0.6cm and 0.2cm);
\draw[thick] (0.8,2) -- (0.8,-2);\draw[thick] (1.2,2)..controls +(0,-2) and +(0,+2) .. (3.8,-2);  
\fill (4,2) ellipse (0.6cm and 0.2cm);\fill (4,-2) ellipse (0.6cm and 0.2cm);
\draw[thick] (3.8,2)..controls +(0,-2) and +(0,+2) .. (1.2,-2); \draw[thick] (4.2,2) -- (4.2,-2);
\node at (6,0) {$.$};
\fill (8,2) ellipse (0.6cm and 0.2cm);\fill (8,-2) ellipse (0.6cm and 0.2cm);
\draw[thick] (7.8,2) -- (7.8,-2);\draw[thick] (8.2,2)..controls +(0,-2) and +(0,+2) .. (10.8,-2);  
\fill (11,2) ellipse (0.6cm and 0.2cm);\fill (11,-2) ellipse (0.6cm and 0.2cm);
\draw[thick] (10.8,2)..controls +(0,-2) and +(0,+2) .. (8.2,-2); \draw[thick] (11.2,2) -- (11.2,-2);
\node at (13,0) {$=$};
\fill (15,4) ellipse (0.6cm and 0.2cm);\fill (15,0) ellipse (0.6cm and 0.2cm);
\draw[thick] (14.8,4) -- (14.8,0);\draw[thick] (15.2,4)..controls +(0,-2) and +(0,+2) .. (17.8,0);  
\fill (18,4) ellipse (0.6cm and 0.2cm);\fill (18,0) ellipse (0.6cm and 0.2cm);
\draw[thick] (17.8,4)..controls +(0,-2) and +(0,+2) .. (15.2,0); \draw[thick] (18.2,4) -- (18.2,0);

\fill (15,0) ellipse (0.6cm and 0.2cm);\fill (15,-4) ellipse (0.6cm and 0.2cm);
\draw[thick] (14.8,0) -- (14.8,-4);\draw[thick] (15.2,0)..controls +(0,-2) and +(0,+2) .. (17.8,-4);  
\fill (18,0) ellipse (0.6cm and 0.2cm);\fill (18,-4) ellipse (0.6cm and 0.2cm);
\draw[thick] (17.8,0)..controls +(0,-2) and +(0,+2) .. (15.2,-4); \draw[thick] (18.2,0) -- (18.2,-4);

\node at (20.5,0) {$=\frac{1}{4}\Bigl($};
\fill (22.5,4) ellipse (0.6cm and 0.2cm);
\draw[thick] (22.3,4) -- (22.3,0);\draw[thick] (22.7,4)..controls +(0,-2) and +(0,+2) .. (25.3,0);  
\fill (25.5,4) ellipse (0.6cm and 0.2cm);
\draw[thick] (25.3,4)..controls +(0,-2) and +(0,+2) .. (22.7,0); \draw[thick] (25.7,4) -- (25.7,0);

\fill (22.5,-4) ellipse (0.6cm and 0.2cm);
\draw[thick] (22.3,0) -- (22.3,-4);\draw[thick] (22.7,0)..controls +(0,-2) and +(0,+2) .. (25.3,-4);  
\fill (25.5,-4) ellipse (0.6cm and 0.2cm);
\draw[thick] (25.3,0)..controls +(0,-2) and +(0,+2) .. (22.7,-4); \draw[thick] (25.7,0) -- (25.7,-4);

\node at (27.5,0) {$+$};
\fill (29.5,4) ellipse (0.6cm and 0.2cm);
\draw[thick] (29.3,4) -- (29.3,0.5);\draw[thick] (29.7,4)..controls +(0,-2) and +(0,+2) .. (32.3,0);  
\fill (32.5,4) ellipse (0.6cm and 0.2cm);
\draw[thick] (32.3,4)..controls +(0,-2) and +(0,+2) .. (29.7,0.5); \draw[thick] (32.7,4) -- (32.7,0);

\draw[thick] (29.3,0.5)..controls +(0,-0.2) and +(0,+0.2) .. (29.7,-0.5);
\draw[thick] (29.7,0.5)..controls +(0,-0.2) and +(0,+0.2) .. (29.3,-0.5);

\fill (29.5,-4) ellipse (0.6cm and 0.2cm);
\draw[thick] (29.3,-0.5) -- (29.3,-4);\draw[thick] (29.7,-0.5)..controls +(0,-2) and +(0,+2) .. (32.3,-4);  
\fill (32.5,-4) ellipse (0.6cm and 0.2cm);
\draw[thick] (32.3,0)..controls +(0,-2) and +(0,+2) .. (29.7,-4); \draw[thick] (32.7,0) -- (32.7,-4);

\node at (34.5,0) {$+$};
\fill (36.5,4) ellipse (0.6cm and 0.2cm);
\draw[thick] (36.3,4) -- (36.3,0);\draw[thick] (36.7,4)..controls +(0,-2) and +(0,+2) .. (39.3,0.5);  
\fill (39.5,4) ellipse (0.6cm and 0.2cm);
\draw[thick] (39.3,4)..controls +(0,-2) and +(0,+2) .. (36.7,0); \draw[thick] (39.7,4) -- (39.7,0.5);

\draw[thick] (39.3,0.5)..controls +(0,-0.2) and +(0,+0.2) .. (39.7,-0.5);
\draw[thick] (39.7,0.5)..controls +(0,-0.2) and +(0,+0.2) .. (39.3,-0.5);

\fill (36.5,-4) ellipse (0.6cm and 0.2cm);
\draw[thick] (36.3,0) -- (36.3,-4);\draw[thick] (36.7,0)..controls +(0,-2) and +(0,+2) .. (39.3,-4);  
\fill (39.5,-4) ellipse (0.6cm and 0.2cm);
\draw[thick] (39.3,-0.5)..controls +(0,-2) and +(0,+2) .. (36.7,-4); \draw[thick] (39.7,-0.5) -- (39.7,-4);

\node at (41.5,0) {$+$};
\fill (43.5,4) ellipse (0.6cm and 0.2cm);
\draw[thick] (43.3,4) -- (43.3,0.5);\draw[thick] (43.7,4)..controls +(0,-2) and +(0,+2) .. (46.3,0.5);  
\fill (46.5,4) ellipse (0.6cm and 0.2cm);
\draw[thick] (46.3,4)..controls +(0,-2) and +(0,+2) .. (43.7,0.5); \draw[thick] (46.7,4) -- (46.7,0.5);

\draw[thick] (43.3,0.5)..controls +(0,-0.2) and +(0,+0.2) .. (43.7,-0.5);
\draw[thick] (43.7,0.5)..controls +(0,-0.2) and +(0,+0.2) .. (43.3,-0.5);
\draw[thick] (46.3,0.5)..controls +(0,-0.2) and +(0,+0.2) .. (46.7,-0.5);
\draw[thick] (46.7,0.5)..controls +(0,-0.2) and +(0,+0.2) .. (46.3,-0.5);

\fill (43.5,-4) ellipse (0.6cm and 0.2cm);
\draw[thick] (43.3,-0.5) -- (43.3,-4);\draw[thick] (43.7,-0.5)..controls +(0,-2) and +(0,+2) .. (46.3,-4);  
\fill (46.5,-4) ellipse (0.6cm and 0.2cm);
\draw[thick] (46.3,-0.5)..controls +(0,-2) and +(0,+2) .. (43.7,-4); \draw[thick] (46.7,-0.5) -- (46.7,-4);

\node at (48.5,0) {$\Bigr)$};

\node at (20,-8) {$=\frac{1}{4}$};
\fill (22,-6) ellipse (0.6cm and 0.2cm);\fill (22,-10) ellipse (0.6cm and 0.2cm);
\draw[thick] (21.8,-6) -- (21.8,-10);\draw[thick] (22.2,-6) -- (22.2,-10);
\fill (25,-6) ellipse (0.6cm and 0.2cm);\fill (25,-10) ellipse (0.6cm and 0.2cm);
\draw[thick] (24.8,-6) -- (24.8,-10);\draw[thick] (25.2,-6) -- (25.2,-10);

\node at (27,-8) {$+\frac{1}{2}$};

\fill (29,-6) ellipse (0.6cm and 0.2cm);\fill (29,-10) ellipse (0.6cm and 0.2cm);
\draw[thick] (28.8,-6) -- (28.8,-10);\draw[thick] (29.2,-6)..controls +(0,-2) and +(0,+2) .. (31.8,-10);  
\fill (32,-6) ellipse (0.6cm and 0.2cm);\fill (32,-10) ellipse (0.6cm and 0.2cm);
\draw[thick] (31.8,-6)..controls +(0,-2) and +(0,+2) .. (29.2,-10);  \draw[thick] (32.2,-6) -- (32.2,-10);

\node at (34,-8) {$+\frac{1}{4}$};

\fill (36,-6) ellipse (0.6cm and 0.2cm);\fill (36,-10) ellipse (0.6cm and 0.2cm);
\draw[thick] (35.8,-6)..controls +(0,-2) and +(0,+2) .. (38.8,-10);  
\draw[thick] (36.2,-6)..controls +(0,-2) and +(0,+2) .. (39.2,-10);  
\fill (39,-6) ellipse (0.6cm and 0.2cm);\fill (39,-10) ellipse (0.6cm and 0.2cm);
\draw[thick] (38.8,-6)..controls +(0,-2) and +(0,+2) .. (35.8,-10);  
\draw[thick] (39.2,-6)..controls +(0,-2) and +(0,+2) .. (36.2,-10);  

\end{tikzpicture}
\end{center}

\paragraph{Coherence of the two definitions and standard basis.} As the notation suggests, the algebra of fused permutations coincides with the algebra $H_{\bk,n}(1)$ as defined in the previous section with the help of the idempotent $P_{\bk,n}$.
\begin{prop}[\cite{publi-CP1}]
The two definitions of the algebra $H_{\bk,n}(1)$ give isomorphic algebras.
\end{prop}
Without details, the identification of the two definitions goes as follows. Recall that the idempotent subalgebra $P_{\bk,n}S_{k_1+\dots+k_n}P_{\bk,n}$ has a basis:
\begin{equation}\label{standardbasis1}\{P_{\bk,n}wP_{\bk,n}\}_{w\in X_{min}}\,,\end{equation}
where $X_{min}$ is the set of minimal-length representatives of double cosets of $S_{k_1}\times\dots \times S_{k_n}$ in $S_{k_1+\dots+k_n}$. So these $w$'s are permutations of $k_1+\dots+k_n$ and as such can be seen as permutation diagrams connecting one-to-one two rows of $k_1+\dots+k_n$. In such a diagram, in each of the two rows of dots, we glue (or ``fuse'') the $k_1$ first dots together, then the $k_2$ next dots and so on. We obtain thus a diagram as in the preceding subsection which represents a fused permutation in $S^{fus}_{\bk,n}$.

Note that in the diagrammatic point of view, the basis above, using the minimal-length representatives, is very natural. It corresponds to choosing, for a given fused representation, the unique diagram representing it which has a minimal number of intersections between edges. In particular, such a diagram satisfies that the edges leaving the same dot on top do not cross each other and similarly for the edges arriving at the same dot at the bottom. All the diagrams of fused permutations drawn above were drawn like this.

\subsection{Fused braids and the fused Hecke algebras}

\paragraph{Objects: fused braids.}  We consider the following objects, which are generalisations of usual braids. One can also see them as ``topologised'' fused permutations, in the sense that crossings will matter now. We start with the rectangular strip with a top line of $n$ fixed dots and a bottom line of $n$ fixed dots, as for usual braids. And we put strands connecting dots but now we require the following: Each strand connects a dot from the top line to a dot on the bottom line, but now, for the $a$-th dot on top and the $a$-th dot at the bottom, there are $k_a$ strands attached to each of them. So in total there are $k_1+\dots+k_n$ strands. 

To be more precise, it is better if we replace the dots by small ellipses, and we will do so from now on. Then the strands which are attached to the same ellipse are not really attached to the same point of the ellipse. Instead they are attached next to each other at the same ellipse. We use the same terminology as for usual braids, and as before we consider such diagrams up to isotopy, namely up to continuously moving the strands while leaving their end points fixed.

Such an object we call a fused braid. Needless to insist that for $\bk=(1,1,\dots)$, a fused braid is a usual braid. Some examples of fused braids can be found below.

\paragraph{The vector space.} We consider the vector space $Vect^{fus}_{\bk,n}$ of formal linear combinations of fused braids.
\begin{defi}\label{vector-fused-braids}
 The vector space $H_{\bk,n}(q)$ is the quotient of $Vect^{fus}_{\bk,n}$ by the following relations:
 \begin{itemize}
   \item[(i)] The Hecke relation for all crossings:  
   \begin{center}
 \begin{tikzpicture}[scale=0.25]
\draw[thick] (0,2)..controls +(0,-2) and +(0,+2) .. (4,-2);
\fill[white] (2,0) circle (0.4);
\draw[thick] (4,2)..controls +(0,-2) and +(0,+2) .. (0,-2);
\node at (6,0) {$=$};
\draw[thick] (12,2)..controls +(0,-2) and +(0,+2) .. (8,-2);
\fill[white] (10,0) circle (0.4);
\draw[thick] (8,2)..controls +(0,-2) and +(0,+2) .. (12,-2);
\node at (17,0) {$-\,(q-q^{-1})$};
\draw[thick] (21,2) -- (21,-2);\draw[thick] (25,2) -- (25,-2);
\end{tikzpicture}
\end{center}
  \item[(ii)]  The idempotent relations for crossings near the ellipses:
 \begin{center}
 \begin{tikzpicture}[scale=0.4]
\fill (2,2) ellipse (0.8cm and 0.2cm);
\draw[thick] (2.2,2)..controls +(0,-1.5) and +(1,1) .. (1.2,0);
\fill[white] (2,0.7) circle (0.2);
\draw[thick] (1.8,2)..controls +(0,-1.5) and +(-1,1) .. (2.8,0);
\node at (4.5,1) {$=$};
\node at (7,1) {$q$};
\fill (9,2) ellipse (0.8cm and 0.2cm);
\draw[thick] (8.8,2)..controls +(0,-1.5) and +(0.5,0.5) .. (8.2,0);
\draw[thick] (9.2,2)..controls +(0,-1.5) and +(-0.5,0.5) .. (9.8,0);

\node at (13,1) {and};

\fill (18,2) ellipse (0.8cm and 0.2cm);
\draw[thick] (17.8,2)..controls +(0,-1.5) and +(-1,1) .. (18.8,0);
\fill[white] (18,0.7) circle (0.2);
\draw[thick] (18.2,2)..controls +(0,-1.5) and +(1,1) .. (17.2,0);
\node at (20.5,1) {$= $};
\node at (23,1) {$q^{-1}$};
\fill (25,2) ellipse (0.8cm and 0.2cm);
\draw[thick] (24.8,2)..controls +(0,-1.5) and +(0.5,0.5) .. (24.2,0);
\draw[thick] (25.2,2)..controls +(0,-1.5) and +(-0.5,0.5) .. (25.8,0);
\end{tikzpicture}
\end{center}
\begin{center}
 \begin{tikzpicture}[scale=0.4]
\fill (2,0) ellipse (0.8cm and 0.2cm);
\draw[thick] (1.8,0)..controls +(0,1.5) and +(-1,-1) .. (2.8,2);
\fill[white] (2,1.25) circle (0.2);
\draw[thick] (2.2,0)..controls +(0,1.5) and +(1,-1) .. (1.2,2);
\node at (4.5,1) {$=$};
\node at (7,1) {$q$};
\fill (9,0) ellipse (0.8cm and 0.2cm);
\draw[thick] (9.2,0)..controls +(0,1.5) and +(-0.5,-0.5) .. (9.8,2);
\draw[thick] (8.8,0)..controls +(0,1.5) and +(0.5,-0.5) .. (8.2,2);

\node at (13,1) {and};

\fill (18,0) ellipse (0.8cm and 0.2cm);
\draw[thick] (18.2,0)..controls +(0,1.5) and +(1,-1) .. (17.2,2);
\fill[white] (18,1.25) circle (0.2);
\draw[thick] (17.8,0)..controls +(0,1.5) and +(-1,-1) .. (18.8,2);
\node at (20.5,1) {$=$};
\node at (23,1) {$q^{-1}$};
\fill (25,0) ellipse (0.8cm and 0.2cm);
\draw[thick] (25.2,0)..controls +(0,1.5) and +(-0.5,-0.5) .. (25.8,2);
\draw[thick] (24.8,0)..controls +(0,1.5) and +(0.5,-0.5) .. (24.2,2);
\end{tikzpicture}
\end{center}
 \end{itemize}
\end{defi}
The first relation (the Hecke relation) is valid locally for all crossings as in the situation of classical braids and Hecke algebra. The idempotents relations are also local relations. In words, they impose the following: if two strands start from the same ellipse and their first crossing is crossing each other, then the original fused braid is equal to the fused braid obtained by removing this crossing and multiplying by $q^{\pm1}$; and similarly for two strands arriving at the same ellipse.

One should note that the Hecke relation for a crossing near an ellipse is the sum of two idempotent relations, so that everything is quite compatible.

\paragraph{Multiplication.} Now we define a product on the vector space $H_{\bk,n}(q)$, which makes it an associative unital algebra. Let $b,b'$ be two fused braids. We define $bb'$ as the result of the following procedure:
\begin{itemize}
\item \emph{(Concatenation)} We place the diagram of $b$ on top of the diagram of $b'$ by identifying the ellipses at the bottom of $b$ with the ellipses on top of $b'$
\item \emph{(Removal of middle ellipses)} For each $a\in\{1,\dots,n\}$, there are $k_a$ strands incoming and $k_a$ strands leaving the $a$-th ellipse in the middle row. We remove this ellipse and replace it by the $q$-symmetriser $P_{k_a}$ of the Hecke algebra $H_{k_a}(q)$.
\end{itemize}
For the last step, it is understood that we see any element of $H_{k_a}(q)$ as a linear combination of braids on $k_a$ strands.

\begin{exam}\label{ex:rem} We illustrate the procedure to remove a middle ellipse when $k=3$. Each middle ellipse has three strands arriving and leaving:  \begin{tikzpicture}[scale=0.3]
\fill (0,0) ellipse (0.8cm and 0.2cm);
\draw[thick] (-0.3,0)--(-0.5,1);
\draw[thick] (0,0)--(0,1);
\draw[thick] (0.3,0)--(0.5,1);
\draw[thick] (-0.3,0)--(-0.5,-1);
\draw[thick] (0,0)--(0,-1);
\draw[thick] (0.3,0)--(0.5,-1);
\end{tikzpicture}; each one is replaced by the following sum normalized by $1+2q^2+2q^4+q^6$:
\begin{center}
 \begin{tikzpicture}[scale=0.3]
\fill (1,12) circle (0.2cm);\fill (1,8) circle (0.2cm);
\draw[thick] (1,12) -- (1,8);
\fill (3,12) circle (0.2cm);\fill (3,8) circle (0.2cm);
\draw[thick] (3,12) -- (3,8);
\fill (5,12) circle (0.2cm);\fill (5,8) circle (0.2cm);
\draw[thick] (5,12) -- (5,8);
\node at (7,10) {$+$};
\node at (9,10) {$q$};
\fill (10,12) circle (0.2cm);\fill (10,8) circle (0.2cm);
\draw[thick] (12,12)..controls +(0,-2) and +(0,+2) .. (10,8);\fill[white] (11,10) circle (0.4);
\draw[thick] (10,12)..controls +(0,-2) and +(0,+2) .. (12,8);
\fill (12,12) circle (0.2cm);\fill (12,8) circle (0.2cm);
\fill (14,12) circle (0.2cm);\fill (14,8) circle (0.2cm);
\draw[thick] (14,12) -- (14,8);
\node at (16,10) {$+$};
\node at (18,10) {$q$};
\fill (19,12) circle (0.2cm);\fill (19,8) circle (0.2cm);
\draw[thick] (19,12) -- (19,8);
\fill (21,12) circle (0.2cm);\fill (21,8) circle (0.2cm);
\draw[thick] (23,12)..controls +(0,-2) and +(0,+2) .. (21,8);\fill[white] (22,10) circle (0.4);
\fill (23,12) circle (0.2cm);\fill (23,8) circle (0.2cm);
\draw[thick] (21,12)..controls +(0,-2) and +(0,+2) .. (23,8);
\node at (25,10) {$+$};
\node at (27,10) {$q^{2}$};
\draw[thick] (32,12)..controls +(0,-2) and +(0,+2) .. (28,8);\fill[white] (30.6,10.4) circle (0.4);\fill[white] (29.4,9.6) circle (0.4);
\fill (28,12) circle (0.2cm);\fill (28,8) circle (0.2cm);
\draw[thick] (28,12)..controls +(0,-2) and +(0,+2) .. (30,8);
\fill (30,12) circle (0.2cm);\fill (30,8) circle (0.2cm);
\draw[thick] (30,12)..controls +(0,-2) and +(0,+2) .. (32,8);
\fill (32,12) circle (0.2cm);\fill (32,8) circle (0.2cm);
\node at (34,10) {$+$};
\node at (36,10) {$q^{2}$};
\fill (37,12) circle (0.2cm);\fill (37,8) circle (0.2cm);
\draw[thick] (41,12)..controls +(0,-2) and +(0,+2) .. (39,8);
\draw[thick] (39,12)..controls +(0,-2) and +(0,+2) .. (37,8);
\fill[white] (38.4,10.3) circle (0.4);\fill[white] (39.5,9.7) circle (0.4);
\draw[thick] (37,12)..controls +(0,-2) and +(0,+2) .. (41,8);
\fill (39,12) circle (0.2cm);\fill (39,8) circle (0.2cm);
\fill (41,12) circle (0.2cm);\fill (41,8) circle (0.2cm);
\node at (43,10) {$+$};
\node at (45,10) {$q^{3}$};
\fill (46,12) circle (0.2cm);\fill (46,8) circle (0.2cm);
\fill (48,12) circle (0.2cm);\fill (48,8) circle (0.2cm);
\fill (50,12) circle (0.2cm);\fill (50,8) circle (0.2cm);
\draw[thick] (50,12)..controls +(0,-3) and +(0,+1) .. (46,8);
\fill[white] (48,9.2) circle (0.3);
\draw[thick] (48,12) -- (48,8);
\fill[white] (48,10.8) circle (0.3);\fill[white] (49.1,10) circle (0.3);
\draw[thick] (46,12)..controls +(0,-1) and +(0,+3) .. (50,8);
\end{tikzpicture}
\end{center}

$\bullet$ Here is an example of a product of two elements of $H^{fus}_{2,2}(q)$:
\begin{center}
\begin{tikzpicture}[scale=0.3]
\fill (1,2) ellipse (0.6cm and 0.2cm);\fill (1,-2) ellipse (0.6cm and 0.2cm);
\fill (4,2) ellipse (0.6cm and 0.2cm);\fill (4,-2) ellipse (0.6cm and 0.2cm);
\draw[thick] (3.8,2)..controls +(0,-2) and +(0,+2) .. (1.2,-2); \draw[thick] (4.2,2) -- (4.2,-2);
\fill[white] (2.5,0) circle (0.4);
\draw[thick] (0.8,2) -- (0.8,-2);\draw[thick] (1.2,2)..controls +(0,-2) and +(0,+2) .. (3.8,-2);  

\node at (6,0) {$.$};

\fill (8,2) ellipse (0.6cm and 0.2cm);\fill (8,-2) ellipse (0.6cm and 0.2cm);
\fill (11,2) ellipse (0.6cm and 0.2cm);\fill (11,-2) ellipse (0.6cm and 0.2cm);
\draw[thick] (10.8,2)..controls +(0,-2) and +(0,+2) .. (8.2,-2); \draw[thick] (11.2,2) -- (11.2,-2);
\fill[white] (9.5,0) circle (0.4);
\draw[thick] (7.8,2) -- (7.8,-2);\draw[thick] (8.2,2)..controls +(0,-2) and +(0,+2) .. (10.8,-2);  

\node at (17.5,0) {$=\displaystyle\frac{1}{(1+q^2)^2}\Bigg($};
\fill (22.5,4) ellipse (0.6cm and 0.2cm);
\fill (25.5,4) ellipse (0.6cm and 0.2cm);
\draw[thick] (25.3,4)..controls +(0,-2) and +(0,+2) .. (22.7,0); \draw[thick] (25.7,4) -- (25.7,0);
\fill[white] (24,2) circle (0.4);
\draw[thick] (22.3,4) -- (22.3,0);\draw[thick] (22.7,4)..controls +(0,-2) and +(0,+2) .. (25.3,0);  

\fill (22.5,-4) ellipse (0.6cm and 0.2cm);
\fill (25.5,-4) ellipse (0.6cm and 0.2cm);
\draw[thick] (25.3,0)..controls +(0,-2) and +(0,+2) .. (22.7,-4); \draw[thick] (25.7,0) -- (25.7,-4);
\fill[white] (24,-2) circle (0.4);
\draw[thick] (22.3,0) -- (22.3,-4);\draw[thick] (22.7,0)..controls +(0,-2) and +(0,+2) .. (25.3,-4);  

\node at (27.5,0) {$+ q$};
\fill (29.5,4) ellipse (0.6cm and 0.2cm);
\fill (32.5,4) ellipse (0.6cm and 0.2cm);
\draw[thick] (32.3,4)..controls +(0,-2) and +(0,+2) .. (29.7,0.5); \draw[thick] (32.7,4) -- (32.7,0);
\fill[white] (30.9,2) circle (0.4);
\draw[thick] (29.3,4) -- (29.3,0.5);\draw[thick] (29.7,4)..controls +(0,-2) and +(0,+2) .. (32.3,0);  

\draw[thick] (29.7,0.5)..controls +(0,-0.2) and +(0,+0.2) .. (29.3,-0.5);
\fill[white] (29.5,0) circle (0.2);
\draw[thick] (29.3,0.5)..controls +(0,-0.2) and +(0,+0.2) .. (29.7,-0.5);

\fill (29.5,-4) ellipse (0.6cm and 0.2cm);
\fill (32.5,-4) ellipse (0.6cm and 0.2cm);
\draw[thick] (32.3,0)..controls +(0,-2) and +(0,+2) .. (29.7,-4); \draw[thick] (32.7,0) -- (32.7,-4);
\fill[white] (30.8,-2.1) circle (0.4);
\draw[thick] (29.3,-0.5) -- (29.3,-4);\draw[thick] (29.7,-0.5)..controls +(0,-2) and +(0,+2) .. (32.3,-4);

\node at (34.5,0) {$+q$};
\fill (36.5,4) ellipse (0.6cm and 0.2cm);
\fill (39.5,4) ellipse (0.6cm and 0.2cm);
\draw[thick] (39.3,4)..controls +(0,-2) and +(0,+2) .. (36.7,0); \draw[thick] (39.7,4) -- (39.7,0.5);
\fill[white] (38.2,2.1) circle (0.4);
\draw[thick] (36.3,4) -- (36.3,0);\draw[thick] (36.7,4)..controls +(0,-2) and +(0,+2) .. (39.3,0.5);  

\draw[thick] (39.7,0.5)..controls +(0,-0.2) and +(0,+0.2) .. (39.3,-0.5);
\fill[white] (39.5,0) circle (0.2);
\draw[thick] (39.3,0.5)..controls +(0,-0.2) and +(0,+0.2) .. (39.7,-0.5);

\fill (36.5,-4) ellipse (0.6cm and 0.2cm);
\fill (39.5,-4) ellipse (0.6cm and 0.2cm);
\draw[thick] (39.3,-0.5)..controls +(0,-2) and +(0,+2) .. (36.7,-4); \draw[thick] (39.7,-0.5) -- (39.7,-4);
\fill[white] (38,-2) circle (0.4);
\draw[thick] (36.3,0) -- (36.3,-4);\draw[thick] (36.7,0)..controls +(0,-2) and +(0,+2) .. (39.3,-4);  

\node at (41.5,0) {$+q^2$};
\fill (43.5,4) ellipse (0.6cm and 0.2cm);
\fill (46.5,4) ellipse (0.6cm and 0.2cm);
\draw[thick] (46.3,4)..controls +(0,-2) and +(0,+2) .. (43.7,0.5); \draw[thick] (46.7,4) -- (46.7,0.5);
\fill[white] (45.1,2) circle (0.4);
\draw[thick] (43.3,4) -- (43.3,0.5);\draw[thick] (43.7,4)..controls +(0,-2) and +(0,+2) .. (46.3,0.5);  

\draw[thick] (43.7,0.5)..controls +(0,-0.2) and +(0,+0.2) .. (43.3,-0.5);
\fill[white] (43.5,0) circle (0.2);
\draw[thick] (43.3,0.5)..controls +(0,-0.2) and +(0,+0.2) .. (43.7,-0.5);

\draw[thick] (46.7,0.5)..controls +(0,-0.2) and +(0,+0.2) .. (46.3,-0.5);
\fill[white] (46.5,0) circle (0.2);
\draw[thick] (46.3,0.5)..controls +(0,-0.2) and +(0,+0.2) .. (46.7,-0.5);

\fill (43.5,-4) ellipse (0.6cm and 0.2cm);
\fill (46.5,-4) ellipse (0.6cm and 0.2cm);
\draw[thick] (46.3,-0.5)..controls +(0,-2) and +(0,+2) .. (43.7,-4); \draw[thick] (46.7,-0.5) -- (46.7,-4);
\fill[white] (44.9,-2) circle (0.4);
\draw[thick] (43.3,-0.5) -- (43.3,-4);\draw[thick] (43.7,-0.5)..controls +(0,-2) and +(0,+2) .. (46.3,-4);  

\node at (48.5,0) {$\Bigg)$};

\node at (17.5,-8) {$=\displaystyle\frac{1}{(1+q^2)^2}\Bigg($};
\fill (22,-6) ellipse (0.6cm and 0.2cm);\fill (22,-10) ellipse (0.6cm and 0.2cm);
\draw[thick] (21.8,-6) -- (21.8,-10);\draw[thick] (22.2,-6) -- (22.2,-10);
\fill (25,-6) ellipse (0.6cm and 0.2cm);\fill (25,-10) ellipse (0.6cm and 0.2cm);
\draw[thick] (24.8,-6) -- (24.8,-10);\draw[thick] (25.2,-6) -- (25.2,-10);

\node at (31,-8) {$+(q-q^{-1}+2q^3)$};

\fill (37,-6) ellipse (0.6cm and 0.2cm);\fill (37,-10) ellipse (0.6cm and 0.2cm);
\fill (40,-6) ellipse (0.6cm and 0.2cm);\fill (40,-10) ellipse (0.6cm and 0.2cm);
\draw[thick] (39.8,-6)..controls +(0,-2) and +(0,+2) .. (37.2,-10);  \draw[thick] (40.2,-6) -- (40.2,-10);
\fill[white] (38.5,-8) circle (0.4);
\draw[thick] (36.8,-6) -- (36.8,-10);\draw[thick] (37.2,-6)..controls +(0,-2) and +(0,+2) .. (39.8,-10);  

\node at (42,-8) {$+q^2$};

\fill (44,-6) ellipse (0.6cm and 0.2cm);\fill (44,-10) ellipse (0.6cm and 0.2cm);
\fill (47,-6) ellipse (0.6cm and 0.2cm);\fill (47,-10) ellipse (0.6cm and 0.2cm);
\draw[thick] (46.8,-6)..controls +(0,-2) and +(0,+2) .. (43.8,-10);  
\draw[thick] (47.2,-6)..controls +(0,-2) and +(0,+2) .. (44.2,-10);  
\fill[white] (45.5,-8) circle (0.4);
\draw[thick] (43.8,-6)..controls +(0,-2) and +(0,+2) .. (46.8,-10);  
\draw[thick] (44.2,-6)..controls +(0,-2) and +(0,+2) .. (47.2,-10);  

\node at (49,-8) {$\Bigg)$};

\end{tikzpicture}
\end{center}
\end{exam}

We illustrate the use of the idempotent relations by explaining why the fused braid with only non-crossing vertical strands is the unit element of $H_{\bk,n}(q)$. For example, for $k=2$, it is:
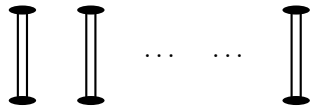
\begin{center}
\begin{tikzpicture}[scale=0.3]
\fill (1,2) ellipse (0.6cm and 0.2cm);\fill (1,-2) ellipse (0.6cm and 0.2cm);
\draw[thick] (0.8,2) -- (0.8,-2);\draw[thick] (1.2,2) -- (1.2,-2);
\fill (4,2) ellipse (0.6cm and 0.2cm);\fill (4,-2) ellipse (0.6cm and 0.2cm);
\draw[thick] (3.8,2) -- (3.8,-2);\draw[thick] (4.2,2) -- (4.2,-2);
\node at (7,0) {$\ldots$};
\node at (10,0) {$\ldots$};
\fill (13,2) ellipse (0.6cm and 0.2cm);\fill (13,-2) ellipse (0.6cm and 0.2cm);
\draw[thick] (12.8,2) -- (12.8,-2);\draw[thick] (13.2,2) -- (13.2,-2);
\end{tikzpicture}
\end{center}
Say we plug this element on top of another fused braid, and we replace the $a$-th middle ellipse by a certain $\si_w$ in $H_{k_a}(q)$. Then the whole lot of crossings added by $\si_w$ can be moved up along the parallel strands, and will thus hit an ellipse of the top row. Thus we can use the idempotent relation. Therefore, we see that replacing a middle ellipse by $\si_w$ will be the same as replacing the middle ellipse by the trivial braid, except for a factor $q^{\ell(w)}$. It remains only to perform the (normalised) sum in the $q$-symmetriser at the $a$-th ellipse and we conclude that multiplying by the above element does not modify the fused braid.

\paragraph{Coherence of the two definitions and standard basis.} As for the algebra of fused permutations, we can prove that this algebra $H_{\bk,n}(q)$ defined here coincides with the one in the previous section defined with the help of the idempotent $P_{\bk,n}$.
\begin{prop}[\cite{publi-CP1}]
The two definitions of the algebra $H_{\bk,n}(q)$ give isomorphic algebras.
\end{prop}
Recall that the standard basis of $H_{\bk,n}(q)$ (in the point of view of the idempotent subalgebra) is:
\begin{equation}\label{standardbasisq}
\{P_J\si_wP_J\}_{w\in X_{min}}\ ,
\end{equation}
where $X_{min}$ is the set of minimal-length representatives as in (\ref{standardbasis1}). Let us interpret this basis in terms of fused braids. Recall that to $w\in X_{min}$, we have fixed a standard diagram representing it (the one with a minimal number of intersections). Now in this diagram, promote each intersection between edges to a crossing, and we decide that all these crossings are positive (the edges coming from the left passes over the other one). We obtain the diagram of a fused braid, and thus an element of the fused Hecke algebra $H_{\bk,n}(q)$. Varying $w\in X_{min}$, this produces the diagrammatic version of the basis (\ref{standardbasisq}).

\begin{exam} If $n=3$ and $k_1=k_2=k_3=2$, there are 21 distinct fused permutations and here are three examples. In each case, we draw a diagram in a canonical form, and the associated standard basis element of $H_{\bk,n}(q)$ below it:
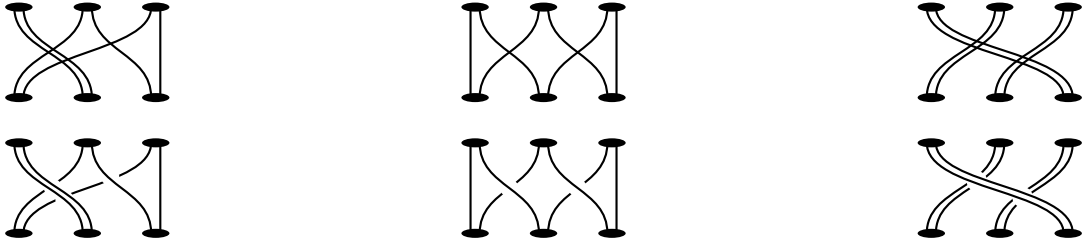
\begin{center}
\begin{tikzpicture}[scale=0.3]
\fill (1,2) ellipse (0.6cm and 0.2cm);\fill (1,-2) ellipse (0.6cm and 0.2cm);
\draw[thick] (0.8,2)..controls +(0,-2) and +(0,+2) .. (3.8,-2);\draw[thick] (1.2,2)..controls +(0,-2) and +(0,+2) .. (4.2,-2);  
\fill (4,2) ellipse (0.6cm and 0.2cm);\fill (4,-2) ellipse (0.6cm and 0.2cm);
\draw[thick] (3.8,2)..controls +(0,-2) and +(0,+2) .. (0.8,-2);\draw[thick] (4.2,2)..controls +(0,-2) and +(0,+2) .. (6.8,-2);
\fill (7,2) ellipse (0.6cm and 0.2cm);\fill (7,-2) ellipse (0.6cm and 0.2cm);
\draw[thick] (6.8,2)..controls +(0,-2) and +(0,+2) .. (1.2,-2);\draw[thick] (7.2,2) -- (7.2,-2);

\fill (1,-4) ellipse (0.6cm and 0.2cm);\fill (1,-8) ellipse (0.6cm and 0.2cm);
\fill (4,-4) ellipse (0.6cm and 0.2cm);\fill (4,-8) ellipse (0.6cm and 0.2cm);
\draw[thick] (3.8,-4)..controls +(0,-2) and +(0,+2) .. (0.8,-8);
\fill (7,-4) ellipse (0.6cm and 0.2cm);\fill (7,-8) ellipse (0.6cm and 0.2cm);
\draw[thick] (6.8,-4)..controls +(0,-2) and +(0,+2) .. (1.2,-8);\draw[thick] (7.2,-4) -- (7.2,-8);
\fill[white] (5,-5.5) circle (0.4);\fill[white] (2.5,-6) circle (0.4);\fill[white] (3,-6.5) circle (0.4);
\draw[thick] (4.2,-4)..controls +(0,-2) and +(0,+2) .. (6.8,-8);
\draw[thick] (0.8,-4)..controls +(0,-2) and +(0,+2) .. (3.8,-8);\draw[thick] (1.2,-4)..controls +(0,-2) and +(0,+2) .. (4.2,-8);  

\fill (21,2) ellipse (0.6cm and 0.2cm);\fill (21,-2) ellipse (0.6cm and 0.2cm);
\draw[thick] (20.8,2) -- (20.8,-2);\draw[thick] (21.2,2)..controls +(0,-2) and +(0,+2) .. (23.8,-2);  
\fill (24,2) ellipse (0.6cm and 0.2cm);\fill (24,-2) ellipse (0.6cm and 0.2cm);
\draw[thick] (23.8,2)..controls +(0,-2) and +(0,+2) .. (21.2,-2); \draw[thick] (26.8,2)..controls +(0,-2) and +(0,+2) .. (24.2,-2);
\fill (27,2) ellipse (0.6cm and 0.2cm);\fill (27,-2) ellipse (0.6cm and 0.2cm);
\draw[thick] (24.2,2)..controls +(0,-2) and +(0,+2) .. (26.8,-2);\draw[thick] (27.2,2) -- (27.2,-2);

\fill (21,-4) ellipse (0.6cm and 0.2cm);\fill (21,-8) ellipse (0.6cm and 0.2cm);
\draw[thick] (20.8,-4)-- (20.8,-8);
\fill (24,-4) ellipse (0.6cm and 0.2cm);\fill (24,-8) ellipse (0.6cm and 0.2cm);
\draw[thick] (23.8,-4)..controls +(0,-2) and +(0,+2) .. (21.2,-8); \draw[thick] (26.8,-4)..controls +(0,-2) and +(0,+2) .. (24.2,-8);
\fill (27,-4) ellipse (0.6cm and 0.2cm);\fill (27,-8) ellipse (0.6cm and 0.2cm);
\draw[thick] (27.2,-4) -- (27.2,-8);
\fill[white] (22.5,-6) circle (0.4);\fill[white] (25.5,-6) circle (0.4);
\draw[thick] (21.2,-4)..controls +(0,-2) and +(0,+2) .. (23.8,-8); 
\draw[thick] (24.2,-4)..controls +(0,-2) and +(0,+2) .. (26.8,-8);

\fill (41,2) ellipse (0.6cm and 0.2cm);\fill (41,-2) ellipse (0.6cm and 0.2cm);
\draw[thick] (40.8,2)..controls +(0,-2.1) and +(0,+1.9) .. (46.8,-2);\draw[thick] (41.2,2)..controls +(0,-1.9) and +(0,+2.1) .. (47.2,-2);
\fill (44,2) ellipse (0.6cm and 0.2cm);\fill (44,-2) ellipse (0.6cm and 0.2cm);
\draw[thick] (43.8,2)..controls +(0,-2) and +(0,+2) .. (40.8,-2);\draw[thick] (44.2,2)..controls +(0,-2) and +(0,+2) .. (41.2,-2); 
\fill (47,2) ellipse (0.6cm and 0.2cm);\fill (47,-2) ellipse (0.6cm and 0.2cm);
\draw[thick] (46.8,2)..controls +(0,-2) and +(0,+2) .. (43.8,-2);\draw[thick] (47.2,2)..controls +(0,-2) and +(0,+2) .. (44.2,-2);

\fill (41,-4) ellipse (0.6cm and 0.2cm);\fill (41,-8) ellipse (0.6cm and 0.2cm);
\fill (44,-4) ellipse (0.6cm and 0.2cm);\fill (44,-8) ellipse (0.6cm and 0.2cm);
\draw[thick] (43.8,-4)..controls +(0,-2) and +(0,+2) .. (40.8,-8);\draw[thick] (44.2,-4)..controls +(0,-2) and +(0,+2) .. (41.2,-8); 
\fill (47,-4) ellipse (0.6cm and 0.2cm);\fill (47,-8) ellipse (0.6cm and 0.2cm);
\draw[thick] (46.8,-4)..controls +(0,-2) and +(0,+2) .. (43.8,-8);\draw[thick] (47.2,-4)..controls +(0,-2) and +(0,+2) .. (44.2,-8);
\fill[white] (43,-5.7) circle (0.45);\fill[white] (45,-6.4) circle (0.45);
\draw[thick] (40.8,-4)..controls +(0,-2.1) and +(0,+1.9) .. (46.8,-8);\draw[thick] (41.2,-4)..controls +(0,-1.9) and +(0,+2.1) .. (47.2,-8);
\end{tikzpicture}
\end{center}
\end{exam}

\begin{framedconc}
The fused Hecke algebra $H_{\bk,n}(q)$ is a flat deformation of the algebra of fused permutations $H_{\bk,n}(1)$: they both have bases indexed by the same set, and the multiplication in $H_{\bk,n}(1)$ is recovered as a particular case ($q=1$) of the one in $H_{\bk,n}(q)$.
\end{framedconc}

\section{Representation theory of the fused Hecke algebra}

\subsection{Irreducible representations and branching rules}

Below, we denote by $\bk_{\vert n}$ the composition $(k_1,\dots,k_n)$ and by $\bk^{\text{ord}}_{\vert n}$ the partition obtained by decreasingly ordering the numbers $k_1,\dots,k_n$. All combinatorial objects below are defined in the appendix of the first chapter.

Let $\bT$ be a standard tableau of shape $\lambda\vdash k_1+\dots+k_n$. We denote by $\overline{\bT}$ the Young tableau obtained from $\bT$ by the following map from $\{1,\dots,k_1+\dots+k_n\}$ to $\{1,\dots,n\}$:
\[1,\dots,k_1 \mapsto 1\,,\ \ \ \ \ k_1+1,\dots,k_1+k_2\mapsto 2\,,\ \ \ \ \dots\ \ \ \ k_1+\dots+k_{n-1}+1,\dots,k_1+\dots+k_n\mapsto n\,,\]
that is, we replace in $\bT$ the first $k_1$ integers by 1, the next $k_2$ ones by 2, and so on. We obtain this way a Young tableau $\overline{\bT}$ of weight $\bk_{\vert n}$ (note that $\overline{\bT}$ does not have to be semistandard, as shown in the example below).

Now, let $\bbT\in\SSTab(\lambda,\bk_{\vert n})$ a semistandard Young tableau of shape $\lambda$ and of weight $\bk_{\vert n}$. We define in the irreducible representation $S_{\lambda}$ of $H_{k_1+\dots+k_n}(q)$ the following vector:
\begin{equation}\label{def-wT}
w_{\bbT}:=\sum_{\begin{array}{c}
\\[-1.6em]
\scriptstyle{\bT\in\STab(\lambda)} \\[-0.4em]
\scriptstyle{\overline{\bT}=\bbT}
\end{array}} v_{\bT}\ \in V_{\lambda}\ .
\end{equation}
where we are using the seminormal basis of $S_{\lambda}$ from Chapter \ref{chap-prel}, Appendix.
\begin{exam}
Let $n=2$, $\bk_{\vert 2}=(2,2)$ and $\lambda=(3,1)$. There is only one semistandard Young tableau of shape $\lambda$ with weight $\bk_{\vert 2}$, and that is $\begin{array}{ccc}
\fbox{\scriptsize{$1$}} & \hspace{-0.35cm}\fbox{\scriptsize{$1$}} & \hspace{-0.35cm}\fbox{\scriptsize{$2$}} \\[-0.2em]
\fbox{\scriptsize{$2$}} & &
\end{array}$. We have then : $w_{\begin{array}{ccc}
\fbox{\scriptsize{$1$}} & \hspace{-0.35cm}\fbox{\scriptsize{$1$}} & \hspace{-0.35cm}\fbox{\scriptsize{$2$}} \\[-0.2em]
\fbox{\scriptsize{$2$}} & &
\end{array}}=v_{\begin{array}{ccc}
\fbox{\scriptsize{$1$}} & \hspace{-0.35cm}\fbox{\scriptsize{$2$}} & \hspace{-0.35cm}\fbox{\scriptsize{$3$}} \\[-0.2em]
\fbox{\scriptsize{$4$}} & &
\end{array}}+v_{\begin{array}{ccc}
\fbox{\scriptsize{$1$}} & \hspace{-0.35cm}\fbox{\scriptsize{$2$}} & \hspace{-0.35cm}\fbox{\scriptsize{$4$}} \\[-0.2em]
\fbox{\scriptsize{$3$}} & &
\end{array}}$. The remaining standard Young tableau $\bT=\begin{array}{ccc}
\fbox{\scriptsize{$1$}} & \hspace{-0.35cm}\fbox{\scriptsize{$3$}} & \hspace{-0.35cm}\fbox{\scriptsize{$4$}} \\[-0.2em]
\fbox{\scriptsize{$2$}} & &
\end{array}$ of shape $\lambda$ gives a Young tableau $\overline{\bT}=\begin{array}{ccc}
\fbox{\scriptsize{$1$}} & \hspace{-0.35cm}\fbox{\scriptsize{$2$}} & \hspace{-0.35cm}\fbox{\scriptsize{$2$}} \\[-0.2em]
\fbox{\scriptsize{$1$}} & &
\end{array}$ which is not semistandard.
\end{exam}

Here is what we can prove using only the representation theory of the usual Hecke algebra.
\begin{framedtheo}[\cite{publi-CP1}]\label{thm-rep}
For any $\lambda\vdash k_1+\dots+k_n$, set $W_{\bk,\lambda}:=P_{\bk,n}(S_{\lambda})$.
\begin{enumerate}
\item The space $W_{\bk,\lambda}$ is spanned by the vectors $w_{\bbT}$, where $\bbT\in\SSTab(\lambda,\bk_{\vert n})$.
\item A complete set of pairwise non-isomorphic irreducible (non-zero) representations of $H_{\bk,n}(q)$ is 
\[\{W_{\bk,\lambda}\}_{\lambda\in S_{\bk,n}}\ \ \ \ \ \text{with $S_{\bk,n}:=\{\lambda\vdash k_1+\dots+k_n\ |\ \lambda\geq\bk^{\text{ord}}_{\vert n}\}$}.\]
The dimension of $W_{\bk,\lambda}$ is the Kostka number $K_{\lambda,\bk_{\vert_n}}=|\SSTab(\lambda,\bk_{\vert n})|$.
\item For $\lambda\in S_{\bk,n}$, the restriction of $W_{\bk,\lambda}$ to $H_{\bk,n-1}(q)$ decomposes as:
\begin{equation}\label{BR-PHP}
\text{Res}_{H_{\bk,n-1}(q)}(W_{\bk,\lambda})\,\cong\,\bigoplus_{\mu\in\text{Res}_\bk(\lambda)} W_{\bk,\mu}\,,
\end{equation}
where we have set:
$$\text{Res}_{\bk}(\lambda):=\bigl\{\mu\in S_{\bk,n-1}\ |\  \text{$\mu\subset\lambda$ and $\lambda/\mu$ contains at most one box in each column}\bigr\}.$$
\end{enumerate}
\end{framedtheo}

\begin{rema}
A numerical consequence of item 2 of the preceding theorem is that:
\[\dim H_{\bk,n}(q)=\sum_{\lambda\vdash k_1+\dots+k_n} |\SSTab(\lambda,\bk_{\vert n})|^2\ .\]
Indeed recall that $|\SSTab(\lambda,\bk_{\vert n})|=0$ if we do not have $\lambda\geq \bk_{\vert n}^{ord}$. This is one representation-theoretic interpretation of the Robinson--Schensted--Knuth correspondence \cite{Knu}. 
\end{rema}

\paragraph{The situation of a constant sequence $\bk=(k,k,...)$.}
We single out the situation of a constant sequence $\bk=(k,k,\dots)$ for an integer $k\geq 1$. In this situation, the parametrisation of irreducible representations in item 2 of Theorem \ref{thm-rep} is much simpler.
\begin{coro}[\cite{publi-CP1}]\label{coro-rep-const}
Let $\bk=(k,k,\dots)$ for an integer $k\geq 1$. A complete set of pairwise non-isomorphic irreducible (non-zero) representations of $H_{\bk,n}$ is 
\[\{W_{\bk,\lambda}\}_{\lambda\in S_{\bk,n}}\ \ \ \ \ \text{with $S_{\bk,n}:=\{\lambda\vdash kn\ |\ l(\lambda)\leq n\}$}.\]
\end{coro}

\begin{exam}
The first levels of the Bratteli diagrams of the chains $\{H_{\bk,n}\}_{n\geq 0}$ for $\bk=(2,2,2,2,...)$ and for $\bk=(3,1,1,1,...)$ are given in Appendix of this chapter.
\end{exam}

\subsection{Back to the centraliser}

Recall that we have a representation:
\[\pi\ :\ H_{\bk,n}\ \to\ \text{End}\bigl(S^{k_1}V\otimes\dots\otimes S^{k_n}V\bigr)\,,\]
whose image gives the centraliser of $U_q(sl_N)$ (this was our generalisation of the first part of the Schur--Weyl duality). It remains to understand the kernel of this map. We will give its representation-theoretic description.

The identification between irreducible representations of $H_{\bk,n}(q)$ and certain partitions refers to the classification obtained in Theorem \ref{thm-rep}.
\begin{framedtheo}[\cite{publi-CP1}]\label{theo-quot-rep}
The kernel of the representation $\pi$ is the ideal of $H_{\bk,n}(q)$ corresponding to the following subset of irreducible representations:
\[\{\ \lambda\in S_{\bk,n}\ |\ l(\lambda)>N\ \}\ .\]
In particular, the kernel is $\{0\}$ and the centraliser coincides with the algebra $H_{\bk,n}(q)$ if and only if $n\leq N$.
\end{framedtheo}

\paragraph{Decomposition of the tensor product $S^{k_1}V\otimes\dots\otimes S^{k_n}V$.} Note that we have not used the knowledge of the decomposition rules for tensor products of representations of $U_q(gl_N)$ (or simply of $GL(N)$). Here the relevant rules are known as the Pieri rule, a particular case of the Littlewood--Richardson rule. So we actually recover the Pieri rule.

More precisely, as a $U_q(gl_N)\otimes H_{\bk,n}(q)$-module, the decomposition is
\begin{equation}
S^{k_1}V\otimes\dots\otimes S^{k_n}V=\bigoplus_{\begin{array}{c}
\\[-1.6em]
\scriptstyle{\lambda\in S_{\bk,n}} \\[-0.4em]
\scriptstyle{l(\lambda)\leq N}
\end{array}} L^N_{\lambda}\otimes W_{\bk,\lambda}\ ,
\end{equation}
where we recall that the set $S_{\bk,n}=\{\lambda\vdash k_1+\dots+k_n\ |\ \lambda\geq \bk_{\vert n}^{\text{ord}}\}$ parametrises the irreducible representations of $H_{\bk,n}(q)$ and $W_{\bk,\lambda}$ is the corresponding irreducible representation constructed in Theorem \ref{thm-rep}. 

In particular, adding the information on the dimension of $W_{\bk,\lambda}$, we obtain that as a $U_q(gl_N)$-module:
\[S^{k_1}V\otimes\dots\otimes S^{k_n}V=\bigoplus_{\begin{array}{c}
\\[-1.6em]
\scriptstyle{\lambda\in S_{\bk,n}} \\[-0.4em]
\scriptstyle{l(\lambda)\leq N}
\end{array}} \bigl(L^N_{\lambda}\bigr)^{\oplus K_{\lambda,\bk_{\vert n}}}\ ,
\]
where $K_{\lambda,\bk_{\vert n}}$ is the Kostka number counting the number of semistandard Young tableaux of shape $\lambda$ and of weight $(k_1,\dots,k_n)$.

\section{Algebraic presentation of the centraliser}\label{sec-quot}

Here we assume that the sequence $\bk=(k_1,k_2,\dots)$ does not contain $0$, and is in decreasing order.

\paragraph{Generalisation of the $q$-antisymmetriser.} We define an element $H_{\bk,n}(q)$ by a simple diagrammatic procedure. For convenience, we first give the definition in the situation $q^2=1$ and then treat the general case (examples are given below). We will give an equivalent more algebraic definition just after the examples.

Let $w\in S_n$. We define a fused permutation denoted $|w|_{\bk}$ by the following procedure. Start from the permutation diagram of $w$ and add vertical edges (if necessary) between dots to form the diagram of a fused permutation corresponding to $\bk$. More precisely, for each $a\in\{1,\dots,n\}$, we add $k_a-1$ vertical edges connecting the $a$-th dot on top to the $a$-th dot at the bottom. Then we set:
\begin{equation}\label{sym1}
\AS_{\bk,n}(1)=\sum_{w\in S_n}(-1)^{\ell(w)}|w|_{\bk}\ \in H_{\bk,n}(1)\ .
\end{equation}

Similarly, we define an element $|\si_w|_{\bk}\in H_{\bk,n}(q)$ as follows. We start with the braid diagram of $\si_w$ with minimal number of crossing. We promote all dots into ellipses, and for each $a\in\{1,\dots,n\}$, we add $k_a-1$ vertical edges connecting the $a$-th top ellipse to the $a$-th bottom ellipse. The rule is a follows: at each ellipse, the new strands are attached to the right of the one strand already present; the added strands do not cross each other; the new strands are ``above'' the original ones forming $\si_w$ (above in the natural sense, as shown in the examples below). 
Then we set:
\begin{equation}\label{symq}
\AS_{\bk,n}(q)=\sum_{w\in\mS_n}(-q^{-1})^{\ell(w)}|\si_w|_{\bk}\ \in H_{\bk,n}(q)\ .
\end{equation}

\begin{exam} Let $n=3$ and $\bk=(2,2,2,\dots)$. Here is depicted the procedure to obtain $\AS_{\bk,n}(1)$:
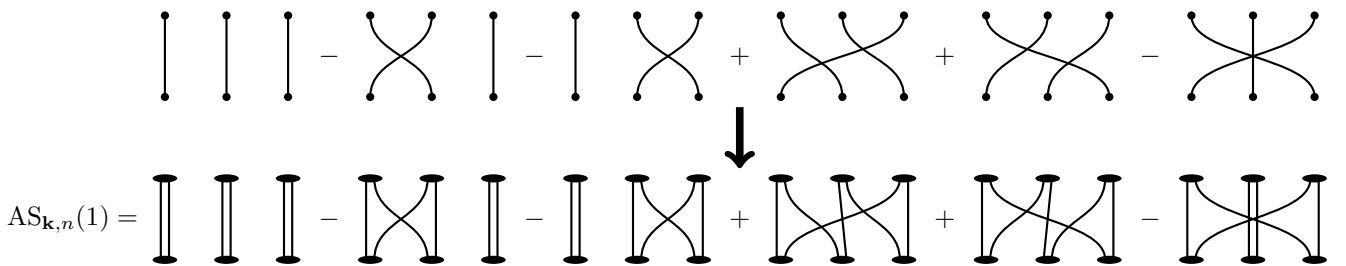
\begin{center}
 \begin{tikzpicture}[scale=0.27]
\fill (1,12) circle (0.2cm);\fill (1,8) circle (0.2cm);
\draw[thick] (1,12) -- (1,8);
\fill (4,12) circle (0.2cm);\fill (4,8) circle (0.2cm);
\draw[thick] (4,12) -- (4,8);
\fill (7,12) circle (0.2cm);\fill (7,8) circle (0.2cm);
\draw[thick] (7,12) -- (7,8);
\node at (9,10) {$-$};
\fill (11,12) circle (0.2cm);\fill (11,8) circle (0.2cm);
\draw[thick] (11,12)..controls +(0,-2) and +(0,+2) .. (14,8);
\fill (14,12) circle (0.2cm);\fill (14,8) circle (0.2cm);
\draw[thick] (14,12)..controls +(0,-2) and +(0,+2) .. (11,8);
\fill (17,12) circle (0.2cm);\fill (17,8) circle (0.2cm);
\draw[thick] (17,12) -- (17,8);
\node at (19,10) {$-$};
\fill (21,12) circle (0.2cm);\fill (21,8) circle (0.2cm);
\draw[thick] (21,12) -- (21,8);
\fill (24,12) circle (0.2cm);\fill (24,8) circle (0.2cm);
\draw[thick] (27,12)..controls +(0,-2) and +(0,+2) .. (24,8);
\fill (27,12) circle (0.2cm);\fill (27,8) circle (0.2cm);
\draw[thick] (24,12)..controls +(0,-2) and +(0,+2) .. (27,8);
\node at (29,10) {$+$};
\fill (31,12) circle (0.2cm);\fill (31,8) circle (0.2cm);
\draw[thick] (31,12)..controls +(0,-2) and +(0,+2) .. (34,8);
\fill (34,12) circle (0.2cm);\fill (34,8) circle (0.2cm);
\draw[thick] (34,12)..controls +(0,-2) and +(0,+2) .. (37,8);
\fill (37,12) circle (0.2cm);\fill (37,8) circle (0.2cm);
\draw[thick] (37,12)..controls +(0,-2) and +(0,+2) .. (31,8);
\node at (39,10) {$+$};
\fill (41,12) circle (0.2cm);\fill (41,8) circle (0.2cm);
\draw[thick] (41,12)..controls +(0,-2) and +(0,+2) .. (47,8);
\fill (44,12) circle (0.2cm);\fill (44,8) circle (0.2cm);
\draw[thick] (44,12)..controls +(0,-2) and +(0,+2) .. (41,8);
\fill (47,12) circle (0.2cm);\fill (47,8) circle (0.2cm);
\draw[thick] (47,12)..controls +(0,-2) and +(0,+2) .. (44,8);
\node at (49,10) {$-$};
\fill (51,12) circle (0.2cm);\fill (51,8) circle (0.2cm);
\draw[thick] (51,12)..controls +(0,-2) and +(0,+2) .. (57,8);
\fill (54,12) circle (0.2cm);\fill (54,8) circle (0.2cm);
\draw[thick] (54,12) -- (54,8);
\fill (57,12) circle (0.2cm);\fill (57,8) circle (0.2cm);
\draw[thick] (57,12)..controls +(0,-2) and +(0,+2) .. (51,8);
 
\draw[line width=1mm,->] (29,7.5) -- (29,4.5); 
 
\node at (-3.5,2) {$\AS_{\bk,n}(1)=$};
\fill (1,4) ellipse (0.6cm and 0.2cm);\fill (1,0) ellipse (0.6cm and 0.2cm);
\draw[thick] (0.8,4) -- (0.8,0);\draw[thick] (1.2,4) -- (1.2,0);
\fill (4,4) ellipse (0.6cm and 0.2cm);\fill (4,0) ellipse (0.6cm and 0.2cm);
\draw[thick] (3.8,4) -- (3.8,0);\draw[thick] (4.2,4) -- (4.2,0);
\fill (7,4) ellipse (0.6cm and 0.2cm);\fill (7,0) ellipse (0.6cm and 0.2cm);
\draw[thick] (6.8,4) -- (6.8,0);\draw[thick] (7.2,4) -- (7.2,0);
\node at (9,2) {$-$};
\fill (11,4) ellipse (0.6cm and 0.2cm);\fill (11,0) ellipse (0.6cm and 0.2cm);
\draw[thick] (10.8,4) -- (10.8,0);\draw[thick] (11.2,4)..controls +(0,-2) and +(0,+2) .. (13.8,0);
\fill (14,4) ellipse (0.6cm and 0.2cm);\fill (14,0) ellipse (0.6cm and 0.2cm);
\draw[thick] (13.8,4)..controls +(0,-2) and +(0,+2) .. (11.2,0);\draw[thick] (14.2,4) -- (14.2,0);
\fill (17,4) ellipse (0.6cm and 0.2cm);\fill (17,0) ellipse (0.6cm and 0.2cm);
\draw[thick] (16.8,4) -- (16.8,0);\draw[thick] (17.2,4) -- (17.2,0);
\node at (19,2) {$-$};
\fill (21,4) ellipse (0.6cm and 0.2cm);\fill (21,0) ellipse (0.6cm and 0.2cm);
\draw[thick] (20.8,4) -- (20.8,0);\draw[thick] (21.2,4) -- (21.2,0);
\fill (24,4) ellipse (0.6cm and 0.2cm);\fill (24,0) ellipse (0.6cm and 0.2cm);
\draw[thick] (23.8,4) -- (23.8,0);\draw[thick] (26.8,4)..controls +(0,-2) and +(0,+2) .. (24.2,0);
\fill (27,4) ellipse (0.6cm and 0.2cm);\fill (27,0) ellipse (0.6cm and 0.2cm);
\draw[thick] (24.2,4)..controls +(0,-2) and +(0,+2) .. (26.8,0);\draw[thick] (27.2,4) -- (27.2,0);
\node at (29,2) {$+$};
\fill (31,4) ellipse (0.6cm and 0.2cm);\fill (31,0) ellipse (0.6cm and 0.2cm);
\draw[thick] (30.8,4) -- (30.8,0);\draw[thick] (31.2,4)..controls +(0,-2) and +(0,+2) .. (33.8,0);
\fill (34,4) ellipse (0.6cm and 0.2cm);\fill (34,0) ellipse (0.6cm and 0.2cm);
\draw[thick] (33.8,4) -- (34.2,0);\draw[thick] (34.2,4)..controls +(0,-2) and +(0,+2) .. (36.8,0);
\fill (37,4) ellipse (0.6cm and 0.2cm);\fill (37,0) ellipse (0.6cm and 0.2cm);
\draw[thick] (36.8,4)..controls +(0,-2) and +(0,+2) .. (31.2,0);\draw[thick] (37.2,4) -- (37.2,0);
\node at (39,2) {$+$};
\fill (41,4) ellipse (0.6cm and 0.2cm);\fill (41,0) ellipse (0.6cm and 0.2cm);
\draw[thick] (40.8,4) -- (40.8,0);\draw[thick] (41.2,4)..controls +(0,-2) and +(0,+2) .. (46.8,0);
\fill (44,4) ellipse (0.6cm and 0.2cm);\fill (44,0) ellipse (0.6cm and 0.2cm);
\draw[thick] (43.8,4)..controls +(0,-2) and +(0,+2) .. (41.2,0);\draw[thick] (44.2,4) -- (43.8,0);
\fill (47,4) ellipse (0.6cm and 0.2cm);\fill (47,0) ellipse (0.6cm and 0.2cm);
\draw[thick] (46.8,4)..controls +(0,-2) and +(0,+2) .. (44.2,0);\draw[thick] (47.2,4) -- (47.2,0);
\node at (49,2) {$-$};
\fill (51,4) ellipse (0.6cm and 0.2cm);\fill (51,0) ellipse (0.6cm and 0.2cm);
\draw[thick] (50.8,4) -- (50.8,0);\draw[thick] (51.2,4)..controls +(0,-2) and +(0,+2) .. (56.8,0);
\fill (54,4) ellipse (0.6cm and 0.2cm);\fill (54,0) ellipse (0.6cm and 0.2cm);
\draw[thick] (53.8,4) -- (53.8,0);\draw[thick] (54.2,4) -- (54.2,0);
\fill (57,4) ellipse (0.6cm and 0.2cm);\fill (57,0) ellipse (0.6cm and 0.2cm);
\draw[thick] (56.8,4)..controls +(0,-2) and +(0,+2) .. (51.2,0);\draw[thick] (57.2,4) -- (57.2,0);
\end{tikzpicture}
\end{center}
Here is depicted the procedure to obtain to obtain $\AS_{\bk,n}(q)$:
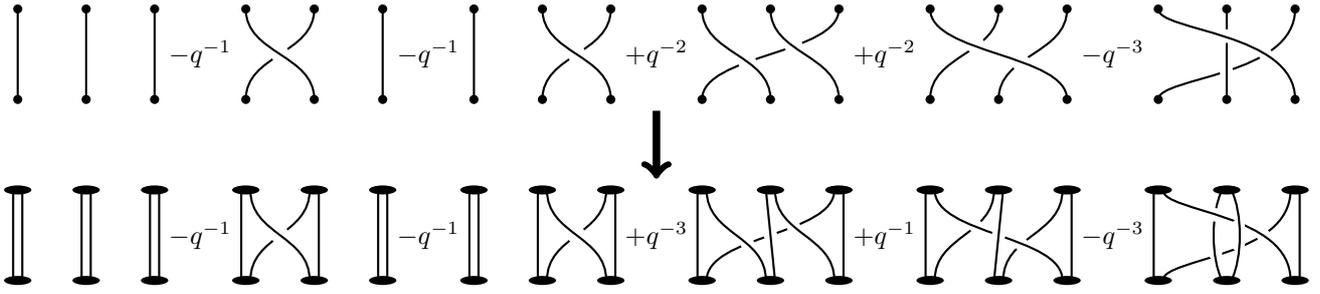
\begin{center}
 \begin{tikzpicture}[scale=0.3]
\fill (1,12) circle (0.2cm);\fill (1,8) circle (0.2cm);
\draw[thick] (1,12) -- (1,8);
\fill (4,12) circle (0.2cm);\fill (4,8) circle (0.2cm);
\draw[thick] (4,12) -- (4,8);
\fill (7,12) circle (0.2cm);\fill (7,8) circle (0.2cm);
\draw[thick] (7,12) -- (7,8);
\node at (9,10) {$-q^{-1}$};
\fill (11,12) circle (0.2cm);\fill (11,8) circle (0.2cm);
\draw[thick] (14,12)..controls +(0,-2) and +(0,+2) .. (11,8);\fill[white] (12.5,10) circle (0.4);
\draw[thick] (11,12)..controls +(0,-2) and +(0,+2) .. (14,8);
\fill (14,12) circle (0.2cm);\fill (14,8) circle (0.2cm);
\fill (17,12) circle (0.2cm);\fill (17,8) circle (0.2cm);
\draw[thick] (17,12) -- (17,8);
\node at (19,10) {$-q^{-1}$};
\fill (21,12) circle (0.2cm);\fill (21,8) circle (0.2cm);
\draw[thick] (21,12) -- (21,8);
\fill (24,12) circle (0.2cm);\fill (24,8) circle (0.2cm);
\draw[thick] (27,12)..controls +(0,-2) and +(0,+2) .. (24,8);\fill[white] (25.5,10) circle (0.4);
\fill (27,12) circle (0.2cm);\fill (27,8) circle (0.2cm);
\draw[thick] (24,12)..controls +(0,-2) and +(0,+2) .. (27,8);
\node at (29,10) {$+q^{-2}$};
\draw[thick] (37,12)..controls +(0,-2) and +(0,+2) .. (31,8);\fill[white] (33,9.6) circle (0.4);\fill[white] (35,10.4) circle (0.4);
\fill (31,12) circle (0.2cm);\fill (31,8) circle (0.2cm);
\draw[thick] (31,12)..controls +(0,-2) and +(0,+2) .. (34,8);
\fill (34,12) circle (0.2cm);\fill (34,8) circle (0.2cm);
\draw[thick] (34,12)..controls +(0,-2) and +(0,+2) .. (37,8);
\fill (37,12) circle (0.2cm);\fill (37,8) circle (0.2cm);
\node at (39,10) {$+q^{-2}$};
\fill (41,12) circle (0.2cm);\fill (41,8) circle (0.2cm);
\draw[thick] (47,12)..controls +(0,-2) and +(0,+2) .. (44,8);
\draw[thick] (44,12)..controls +(0,-2) and +(0,+2) .. (41,8);
\fill[white] (43,10.4) circle (0.4);\fill[white] (45,9.6) circle (0.4);
\draw[thick] (41,12)..controls +(0,-2) and +(0,+2) .. (47,8);
\fill (44,12) circle (0.2cm);\fill (44,8) circle (0.2cm);
\fill (47,12) circle (0.2cm);\fill (47,8) circle (0.2cm);
\node at (49,10) {$-q^{-3}$};
\fill (51,12) circle (0.2cm);\fill (51,8) circle (0.2cm);
\fill (54,12) circle (0.2cm);\fill (54,8) circle (0.2cm);
\fill (57,12) circle (0.2cm);\fill (57,8) circle (0.2cm);
\draw[thick] (57,12)..controls +(0,-3) and +(0,+1) .. (51,8);
\fill[white] (54,9.2) circle (0.3);
\draw[thick] (54,12) -- (54,8);
\fill[white] (54,10.8) circle (0.3);\fill[white] (55.7,10) circle (0.3);
\draw[thick] (51,12)..controls +(0,-1) and +(0,+3) .. (57,8);
 
\draw[line width=1mm,->] (29,7.5) -- (29,4.5); 
 
\fill (1,4) ellipse (0.6cm and 0.2cm);\fill (1,0) ellipse (0.6cm and 0.2cm);
\draw[thick] (0.8,4) -- (0.8,0);\draw[thick] (1.2,4) -- (1.2,0);
\fill (4,4) ellipse (0.6cm and 0.2cm);\fill (4,0) ellipse (0.6cm and 0.2cm);
\draw[thick] (3.8,4) -- (3.8,0);\draw[thick] (4.2,4) -- (4.2,0);
\fill (7,4) ellipse (0.6cm and 0.2cm);\fill (7,0) ellipse (0.6cm and 0.2cm);
\draw[thick] (6.8,4) -- (6.8,0);\draw[thick] (7.2,4) -- (7.2,0);
\node at (9,2) {$-q^{-1}$};
\fill (11,4) ellipse (0.6cm and 0.2cm);\fill (11,0) ellipse (0.6cm and 0.2cm);!
\draw[thick] (10.8,4) -- (10.8,0);\draw[thick] (13.8,4)..controls +(0,-2) and +(0,+2) .. (11.2,0);\fill[white] (12.5,2) circle (0.4);
\fill (14,4) ellipse (0.6cm and 0.2cm);\fill (14,0) ellipse (0.6cm and 0.2cm);
\draw[thick] (14.2,4) -- (14.2,0);\draw[thick] (11.2,4)..controls +(0,-2) and +(0,+2) .. (13.8,0);
\fill (17,4) ellipse (0.6cm and 0.2cm);\fill (17,0) ellipse (0.6cm and 0.2cm);
\draw[thick] (16.8,4) -- (16.8,0);\draw[thick] (17.2,4) -- (17.2,0);
\node at (19,2) {$-q^{-1}$};
\fill (21,4) ellipse (0.6cm and 0.2cm);\fill (21,0) ellipse (0.6cm and 0.2cm);
\draw[thick] (20.8,4) -- (20.8,0);\draw[thick] (21.2,4) -- (21.2,0);
\fill (24,4) ellipse (0.6cm and 0.2cm);\fill (24,0) ellipse (0.6cm and 0.2cm);
\draw[thick] (23.8,4) -- (23.8,0);\draw[thick] (26.8,4)..controls +(0,-2) and +(0,+2) .. (24.2,0);\fill[white] (25.5,2) circle (0.4);
\fill (27,4) ellipse (0.6cm and 0.2cm);\fill (27,0) ellipse (0.6cm and 0.2cm);
\draw[thick] (24.2,4)..controls +(0,-2) and +(0,+2) .. (26.8,0);\draw[thick] (27.2,4) -- (27.2,0);
\node at (29,2) {$+q^{-3}$};
\draw[thick] (36.8,4)..controls +(0,-2) and +(0,+2) .. (31.2,0);\draw[thick] (37.2,4) -- (37.2,0);
\fill[white] (33,1.6) circle (0.3);\fill[white] (35,2.4) circle (0.3);;\fill[white] (34,2) circle (0.3);
\fill (31,4) ellipse (0.6cm and 0.2cm);\fill (31,0) ellipse (0.6cm and 0.2cm);
\draw[thick] (30.8,4) -- (30.8,0);\draw[thick] (31.2,4)..controls +(0,-2) and +(0,+2) .. (33.8,0);
\fill (34,4) ellipse (0.6cm and 0.2cm);\fill (34,0) ellipse (0.6cm and 0.2cm);
\draw[thick] (33.8,4) -- (34.2,0);\draw[thick] (34.2,4)..controls +(0,-2) and +(0,+2) .. (36.8,0);
\fill (37,4) ellipse (0.6cm and 0.2cm);\fill (37,0) ellipse (0.6cm and 0.2cm);
\node at (39,2) {$+q^{-1}$};
\draw[thick] (46.8,4)..controls +(0,-2) and +(0,+2) .. (44.2,0);
\draw[thick] (43.8,4)..controls +(0,-2) and +(0,+2) .. (41.2,0);
\fill[white] (43,2.4) circle (0.3);\fill[white] (45,1.6) circle (0.3);
\draw[thick] (41.2,4)..controls +(0,-2) and +(0,+2) .. (46.8,0);
\fill[white] (44,2) circle (0.3);
\fill (41,4) ellipse (0.6cm and 0.2cm);\fill (41,0) ellipse (0.6cm and 0.2cm);
\draw[thick] (40.8,4) -- (40.8,0);
\fill (44,4) ellipse (0.6cm and 0.2cm);\fill (44,0) ellipse (0.6cm and 0.2cm);
\draw[thick] (44.2,4) -- (43.8,0);
\fill (47,4) ellipse (0.6cm and 0.2cm);\fill (47,0) ellipse (0.6cm and 0.2cm);
\draw[thick] (47.2,4) -- (47.2,0);
\node at (49,2) {$-q^{-3}$};
\draw[thick] (56.8,4)..controls +(0,-3) and +(0,+1) .. (51.2,0);
\fill[white] (53.5,1.2) circle (0.3);\fill[white] (54.5,1.6) circle (0.3);
\draw[thick] (53.8,4)..controls +(-0.5,-1) and +(-0.5,1) .. (53.8,0);
\fill[white] (53.5,2.9) circle (0.3);
\fill[white] (54,2.8) circle (0.3);\fill[white] (55.6,2) circle (0.3);
\draw[thick] (51.2,4)..controls +(0,-1) and +(0,+3) .. (56.8,0);
\fill[white] (54.5,2.4) circle (0.3);
\draw[thick] (54.2,4)..controls +(0.5,-1) and +(0.5,1) .. (54.2,0);
\draw[thick] (50.8,4) -- (50.8,0);
\fill (51,4) ellipse (0.6cm and 0.2cm);\fill (51,0) ellipse (0.6cm and 0.2cm);
\fill (54,4) ellipse (0.6cm and 0.2cm);\fill (54,0) ellipse (0.6cm and 0.2cm);
\fill (57,4) ellipse (0.6cm and 0.2cm);\fill (57,0) ellipse (0.6cm and 0.2cm);
\draw[thick] (57.2,4) -- (57.2,0);
\end{tikzpicture}
\end{center}
The added vertical strands are indeed above all others. Initially, we attached them at each ellipse to the right of the existing strands. In the above picture, we used some idempotent relations to suppress some crossings near the ellipses. This accounts for the modifications in the powers of $q$. Note that, at the end, the coefficient is always equal to $(-1)^{\ell(w)}$ times $q^{-1}$ to the power the signed sum of the crossings in the diagram.
\end{exam}

\paragraph{Algebraic interpretation.} Consider the following element of the usual Hecke algebra $H_{k_1+\dots+k_n}(q)$:
\begin{equation}\label{Gamma}
\Gamma=\si_{k_1}\dots\si_2\cdot \si_{k_1+k_2}\dots \si_3 \cdot \ldots\ldots\cdot \si_{k_1+\dots+k_{n-1}}\dots\si_n\ , 
\end{equation}
Here is the diagrammatic representation of $\Gamma$ in the example $n=3$:
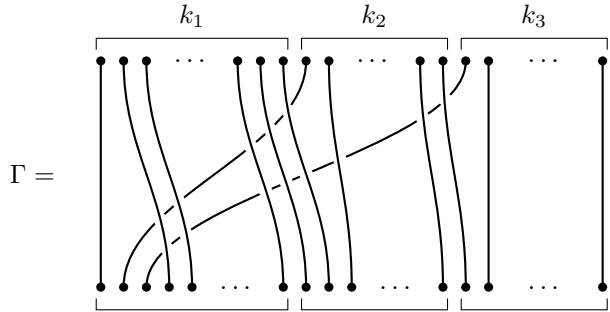
\begin{center}
 \begin{tikzpicture}[scale=0.3]
\node at (-2,0) {$\Gamma=$};
\draw (0.8,5.5) -- (0.8,6) -- (9.2,6) -- (9.2,5.5);\node at (5,7) {$k_1$};
\fill (1,5) circle (0.2);\fill (1,-5) circle (0.2);
\fill (2,5) circle (0.2);\fill (2,-5) circle (0.2);
\fill (3,5) circle (0.2);\fill (3,-5) circle (0.2);
\node at (5,5) {$\dots$};
\fill (4,-5) circle (0.2);
\fill (5,-5) circle (0.2);\node at (7,-5) {$\dots$};
\fill (7,5)  circle (0.2);
\fill (8,5)  circle (0.2);
\fill (9,5) circle (0.2);\fill (9,-5) circle (0.2);
\draw (0.8,-5.5) -- (0.8,-6) -- (9.2,-6) -- (9.2,-5.5);

\draw (9.8,5.5) -- (9.8,6) -- (16.2,6) -- (16.2,5.5);\node at (13,7) {$k_2$};
\fill (10,5) circle (0.2);\fill (10,-5) circle (0.2);
\fill (11,5) circle (0.2);\fill (11,-5) circle (0.2);
\node at (13,5) {$\dots$};\fill (12,-5) circle (0.2);
\node at (14,-5) {$\dots$};
\fill (15,5) circle (0.2);
\fill (16,5) circle (0.2);\fill (16,-5) circle (0.2);
\draw (9.8,-5.5) -- (9.8,-6) -- (16.2,-6) -- (16.2,-5.5);

\draw (16.8,5.5) -- (16.8,6) -- (23.2,6) -- (23.2,5.5);\node at (20,7) {$k_3$};
\fill (17,5) circle (0.2);\fill (17,-5) circle (0.2);
\fill (18,5) circle (0.2);\fill (18,-5) circle (0.2);
\node at (20.5,5) {$\dots$};\node at (20.5,-5) {$\dots$};
\fill (23,5) circle (0.2);\fill (23,-5) circle (0.2);
\draw (16.8,-5.5) -- (16.8,-6) -- (23.2,-6) -- (23.2,-5.5);

\draw[thick] (1,5) -- (1,-5);
\draw[thick] (10,5)..controls +(0,-4) and +(0,+4) .. (2,-5);
\draw[thick] (17,5)..controls +(0,-4) and +(0,+4) .. (3,-5);

\fill[white] (9.2,2.8) circle (0.3);\fill[white] (8.4,2) circle
(0.3);\fill[white] (7.6,1.2) circle (0.3);
\fill[white] (4.4,-1.2) circle (0.3);\fill[white] (3.6,-2) circle
(0.3);\fill[white] (4.7,-2.6) circle (0.3);\fill[white] (3.9,-3.3)
circle (0.3);
\fill[white] (16.1,3.3) circle (0.3);\fill[white] (15.1,2.5) circle (0.3);
\fill[white] (11.5,0.6) circle (0.3);\fill[white] (10,0) circle (0.3);
\fill[white] (9.1,-0.3) circle (0.3);\fill[white] (8.2,-0.7) circle (0.3);

\draw[thick] (2,5)..controls +(0,-4) and +(0,+4) .. (4,-5);
\draw[thick] (3,5)..controls +(0,-4) and +(0,+4) .. (5,-5);
\draw[thick] (7,5)..controls +(0,-4) and +(0,+4) .. (9,-5);
\draw[thick] (8,5)..controls +(0,-4) and +(0,+4) .. (10,-5);
\draw[thick] (9,5)..controls +(0,-4) and +(0,+4) .. (11,-5);
\draw[thick] (11,5)..controls +(0,-4) and +(0,+4) .. (12,-5);
\draw[thick] (15,5)..controls +(0,-4) and +(0,+4) .. (16,-5);
\draw[thick] (16,5)..controls +(0,-4) and +(0,+4) .. (17,-5);
\draw[thick] (18,5)..controls +(0,-4) and +(0,+4) .. (18,-5);
\draw[thick] (23,5)..controls +(0,-4) and +(0,+4) .. (23,-5);

\end{tikzpicture}
\end{center}

Let $w\in \mS_n$. Using the identification between $H_{\bk,n}$ and $P_{\bk,n}H_{k_1+\dots+k_n}(q)P_{\bk,n}$, we can show that an algebraic definition of $|\si_w|_{\bk}$ is
\[|\si_w|_{\bk}=P_{\bk,n}\Gamma \si_w \Gamma^{-1} P_{\bk,n}\ .\]

\paragraph{Generator of the kernel.} We make the following two conjectures. They generalise for example the description of the Temperley--Lieb algebra as a quotient of the Hecke algebra, see Example \ref{exa-TL2}. Below we see the element $\AS_{\bk,N+1}(q)$ as an element of $H_{\bk,n}(q)$ for any $n\geq N+1$, by the natural inclusion of algebras (namely, in $H_{\bk,n}(q)$, the element $\AS_{\bk,N+1}(q)$ involves only the strands attached to the $N+1$ first ellipses).
\begin{framedconj}[\cite{publi-CP1}]\label{conj1}
Let $n>N$. The kernel is generated by the element $\AS_{\bk,N+1}(q)$ and the centraliser is therefore isomorphic to the quotient of the algebra $H_{\bk,n}(q)$ by the relation:
\[\AS_{\bk,N+1}(q)=0\ .\]
\end{framedconj}
What we can show quite easily in general is that the element $\AS_{\bk,N+1}(q)$ indeed belongs to the kernel. The difficult part of the conjecture is that it is a generator. Without going into too much technicality, the difficulty is that at the first non-trivial level, $n=N+1$, the ideal already contains several irreducible representations (see the Bratteli diagram in Appendix). This is the big conceptual difference with the Hecke/Temperley--Lieb situation, where the ideal at level $N+1$ is of dimension 1. 

In particular, in our cases, there are several possible choices for a generator. For mysterious reasons, the diagrammatic construction above seems to produce a natural choice.
\begin{framedconj}[\cite{publi-CP1}]\label{conj2}
The element $\AS_{\bk,N+1}(q)$ is central in $H_{\bk,N+1}(q)$.
\end{framedconj}

\paragraph{Some known cases.} Both conjectures are supported by their verifications in some special cases below. Also, explicit (computer-aided) calculations have allowed to check their validity for all $\bk$ and $N$ such that $k_1+\dots+k_{N+1}\leq 7$. 

Recall for the statements below that $q$ is either generic (an indeterminate) or $q=1$.
\begin{prop}[\cite{publi-CP1}]\label{verif-conj}$\ $
\begin{enumerate}
\item If $\bk=(k,1,1,1,\dots)$ with $k$ arbitrary, then Conjectures \ref{conj1}-\ref{conj2} are true for any $N$.
\item If $N=2$, then for any $\bk$ Conjecture \ref{conj1} is true.
\item If $\bk$ consists only of $1$'s and $2$'s, then Conjecture \ref{conj1} is true for any $N$.
\item If $\bk$ consists only of $1$'s and $2$'s then Conjecture \ref{conj2} is true for $q^2=1$ for any $N$.
\end{enumerate}
\end{prop}

\section{Where is the braid group?}

Here we consider the constant sequence $\bk=(k,k,k,...)$ for an integer $k$ (so that all representations in the tensor product are the same). From our general discussion in the first chapter on centralisers of quantum groups, we know that that the fused Hecke algebra contains a quotient of the braid group algebra. In this section, we will identify the generators of the braid group inside $H_{\bk,n}(q)$.

We denote by $\Sigma_i$ the fused braid in $H_{\bk,n}(q)$ for which all strands starting from ellipse $i$ pass over the strands starting from ellipse $i+1$ and all other strands are vertical. A picture is better here, so for example, for $k=2$:
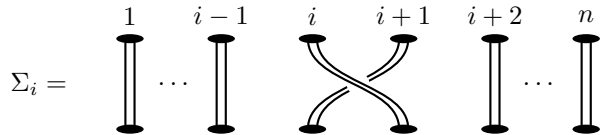
\begin{center}
 \begin{tikzpicture}[scale=0.3]
\node at (-2,0) {$\Sigma_i=$};
\node at (2,3) {$1$};\fill (2,2) ellipse (0.6cm and 0.2cm);\fill (2,-2) ellipse (0.6cm and 0.2cm);
\draw[thick] (1.8,2) -- (1.8,-2);\draw[thick] (2.2,2) -- (2.2,-2);
\node at (4,0) {$\dots$};
\draw[thick] (5.8,2) -- (5.8,-2);\draw[thick] (6.2,2) -- (6.2,-2);
\node at (6,3) {$i-1$};\fill (6,2) ellipse (0.6cm and 0.2cm);\fill (6,-2) ellipse (0.6cm and 0.2cm);
\node at (10,3) {$i$};\fill (10,2) ellipse (0.6cm and 0.2cm);\fill (10,-2) ellipse (0.6cm and 0.2cm);
\node at (14,3) {$i+1$};\fill (14,2) ellipse (0.6cm and 0.2cm);\fill (14,-2) ellipse (0.6cm and 0.2cm);

\draw[thick] (13.8,2)..controls +(0,-2) and +(0,+2) .. (9.8,-2);
\draw[thick] (14.2,2)..controls +(0,-2) and +(0,+2) .. (10.2,-2);
\fill[white] (12,0) circle (0.5);
\draw[thick] (10.2,2)..controls +(0,-2) and +(0,+2) .. (14.2,-2);
\draw[thick] (9.8,2)..controls +(0,-2) and +(0,+2) .. (13.8,-2);

\draw[thick] (17.8,2) -- (17.8,-2);\draw[thick] (18.2,2) -- (18.2,-2);
\node at (18,3) {$i+2$};\fill (18,2) ellipse (0.6cm and 0.2cm);\fill (18,-2) ellipse (0.6cm and 0.2cm);
\node at (20,0) {$\dots$};
\draw[thick] (21.8,2) -- (21.8,-2);\draw[thick] (22.2,2) -- (22.2,-2);\fill (22,2) ellipse (0.6cm and 0.2cm);\fill (22,-2) ellipse (0.6cm and 0.2cm);
\node at (22,3) {$n$};
\end{tikzpicture}
\end{center}
The multiplication in the fused Hecke algebra  $H_{\bk,n}(q)$ may seem quite intricate due to the appearances of the $q$-symmetriser replacing the middle ellipses each time we concatenate two fused braids. However, after some thoughts and with an argument in the same vein as the one after Example \ref{ex:rem}, one can show rather easily that these elements satisfy the braid relation:
\[\Sigma_i\Sigma_{i+1}\Sigma_i=\Sigma_{i+1}\Sigma_i\Sigma_{i+1}\ .\]
At the end, it follows from the following kind of obvious relations in the usual braid group\footnote{The picture looks like some kind of tubes, but we really mean a relation for usual braids (here with 6 strands)} (here $k=2$, but we can add more parallel strands for arbitrary $k$):
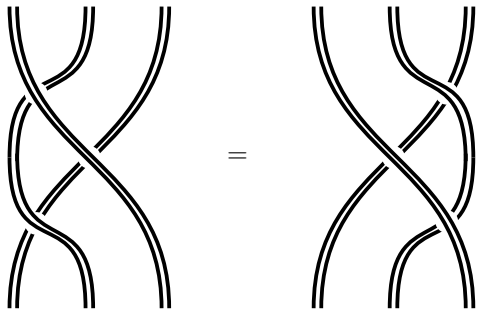
\begin{center}
\begin{tikzpicture}[scale=0.5]
\draw[line width=0.5mm] (4,2)..controls +(0,-4) and +(0,+4) .. (0,-6);
\draw[line width=0.5mm] (4.2,2)..controls +(0,-4) and +(0,+4) .. (0.2,-6);
\draw[line width=0.5mm] (2,2)..controls +(0,-3) and +(0,+3) .. (0,-2);
\draw[line width=0.5mm] (2.2,2)..controls +(0,-3) and +(0,+3) .. (0.2,-2.1);
\fill[white] (0.7,-0.3) circle (0.3);
\fill[white] (2.1,-2) circle (0.3);
\draw[line width=0.5mm] (0,2)..controls +(0,-4) and +(0,+4) .. (4,-6);
\draw[line width=0.5mm] (0.2,2)..controls +(0,-4) and +(0,+4) .. (4.2,-6);
\fill[white] (0.7,-3.75) circle (0.3);
\draw[line width=0.5mm] (0,-2)..controls +(0,-3) and +(0,+3) .. (2,-6);
\draw[line width=0.5mm] (0.2,-1.9)..controls +(0,-3) and +(0,+3) .. (2.2,-6);
\node at (6,-2) {$=$};
\draw[line width=0.5mm] (12,2)..controls +(0,-4) and +(0,+4) .. (8,-6);
\draw[line width=0.5mm] (12.2,2)..controls +(0,-4) and +(0,+4) .. (8.2,-6);
\fill[white] (11.5,-0.3) circle (0.3);
\fill[white] (10.1,-2) circle (0.3);
\draw[line width=0.5mm] (10,2)..controls +(0,-3) and +(0,+3) .. (12,-2.1);
\draw[line width=0.5mm] (10.2,2)..controls +(0,-3) and +(0,+3) .. (12.2,-2);
\draw[line width=0.5mm] (12,-1.9)..controls +(0,-3) and +(0,+3) .. (10,-6);
\draw[line width=0.5mm] (12.2,-2)..controls +(0,-3) and +(0,+3) .. (10.2,-6);
\fill[white] (11.5,-3.7) circle (0.3);
\draw[line width=0.5mm] (8,2)..controls +(0,-4) and +(0,+4) .. (12,-6);
\draw[line width=0.5mm] (8.2,2)..controls +(0,-4) and +(0,+4) .. (12.2,-6);
\end{tikzpicture}
\end{center}
In addition to the braid relations, the elements $\Sigma_1,\dots,\Sigma_{n-1}$ must also satisfy some characteristic equation since the algebra $H_{\bk,n}(q)$ is finite-dimensional. To calculate the eigenvalues directly from the defining multiplication in  $H_{\bk,2}(q)$ is not an easy exercise, but we can use the representation of $H_{\bk,n}(q)$ on the tensor product $S^kV\otimes S^kV$ of representations of $U_q(sl_N)$. In this representation, the element $\Sigma_i$ corresponds to the $R$-matrix, and from this, one can obtain the characteristic equation. 

Skipping details and summarising, the elements $\Sigma_1,\dots,\Sigma_{n-1}$ of $H_{\bk,n}(q)$ satisfy the following relations:
\[\begin{array}{ll}
\Sigma_i\Sigma_{i+1}\Sigma_i=\Sigma_{i+1}\Sigma_i\Sigma_{i+1}\,,\ \ \  & \text{for $i\in\{1,\dots,n-2\}$}\,,\\[0.2em]
\Sigma_i\Sigma_j=\Sigma_j\Sigma_i\,,\ \ \  & \text{for $i,j\in\{1,\dots,n-1\}$ such that $|i-j|>1$}\,,\\[0.2em]
\prod_{l=0}^k\Bigl(\Sigma_i-(-1)^{k+l} q^{-k+l(l+1)}\Bigr)=0\,,\ \ \  & \text{for $i\in\{1,\dots,n-1\}$}\,.
\end{array}
\]
It is important to note that these elements do not generate the whole algebra $H_{\bk,n}(q)$ as soon as $n>2$ and $k>1$. This is a striking difference compared to the usual Hecke algebra. We also note that other relations than the ones above (and not implied by the ones above) must be satisfied by the elements $\Sigma_1,\dots,\Sigma_{n-1}$ since the algebra they generate is finite-dimensional, and we know that relations above are not enough, in general, to define a finite-dimensional algebra. A complete set of defining relations for the subalgebra generated by $\Sigma_1,\dots,\Sigma_{n-1}$ in $H_{\bk,n}(q)$ is an open question.

\section{Application: Solutions of the Yang--Baxter equation}

Having identified the realisation of the braid group inside the fused Hecke algebra, now we go on in our wish list of the introduction, and ask for a solution of the YB equation. Again from our general discussion in the first chapter on centralisers of quantum groups, we know that the fused Hecke algebra must contain such a solution.

First, we introduce some natural elements in the algebra $H_{\bk,n}(q)$: the partial elementary braiding is denoted $\Sigma_i^{(p)}$ and corresponds to the fused braid for which the rightmost $p$ strands starting from ellipse $i$ pass over the $p$ leftmost strands starting from ellipse $i+1$, and all other strands are vertical. Again a picture is better (here $k=3$):
\begin{center}
 \begin{tikzpicture}[scale=0.25]
\node at (16,0) {$\Sigma_i^{(1)}=$};
\node at (20,0) {$\dots$};
\node at (22,3) {$i$};\fill (22,2) ellipse (0.8cm and 0.2cm);\fill (22,-2) ellipse (0.8cm and 0.2cm);
\node at (26,3) {$i+1$};\fill (26,2) ellipse (0.8cm and 0.2cm);\fill (26,-2) ellipse (0.8cm and 0.2cm);
\draw[thick] (21.7,2) -- (21.7,-2);
\draw[thick] (22,2) -- (22,-2);
\draw[thick] (25.7,2)..controls +(0,-2) and +(0,+2) .. (22.3,-2);
\fill[white] (24,0) circle (0.4);
\draw[thick] (22.3,2)..controls +(0,-2) and +(0,+2) .. (25.7,-2);
\draw[thick] (26,2) -- (26,-2);
\draw[thick] (26.3,2) -- (26.3,-2);
\node at (28,0) {$\dots$};
\node at (29.5,0) {$,$};

\node at (34,0) {$\Sigma_i^{(2)}=$};
\node at (38,0) {$\dots$};
\node at (40,3) {$i$};\fill (40,2) ellipse (0.8cm and 0.2cm);\fill (40,-2) ellipse (0.8cm and 0.2cm);
\node at (44,3) {$i+1$};\fill (44,2) ellipse (0.8cm and 0.2cm);\fill (44,-2) ellipse (0.8cm and 0.2cm);
\draw[thick] (39.7,2) -- (39.7,-2);
\draw[thick] (43.7,2)..controls +(0,-2) and +(0,+2) .. (40,-2);
\draw[thick] (44,2)..controls +(0,-2.5) and +(0,+2) .. (40.3,-2);
\fill[white] (42,-0.1) ellipse (0.6cm and 0.4cm);
\draw[thick] (40,2)..controls +(0,-2.5) and +(0,+2) .. (43.7,-2);
\draw[thick] (40.3,2)..controls +(0,-2) and +(0,+2) .. (44,-2);
\draw[thick] (44.3,2) -- (44.3,-2);
\node at (46,0) {$\dots$};
\node at (47.5,0) {$,$};

\node at (52,0) {$\Sigma_i^{(3)}=$};
\node at (56,0) {$\dots$};
\node at (58,3) {$i$};\fill (58,2) ellipse (0.8cm and 0.2cm);\fill (58,-2) ellipse (0.8cm and 0.2cm);
\node at (62,3) {$i+1$};\fill (62,2) ellipse (0.8cm and 0.2cm);\fill (62,-2) ellipse (0.8cm and 0.2cm);
\draw[thick] (61.7,2)..controls +(0,-2) and +(0,+2) .. (57.7,-2);
\draw[thick] (62,2)..controls +(0,-2.5) and +(0,+2) .. (58,-2);
\draw[thick] (62.3,2)..controls +(0,-3) and +(0,+2) .. (58.3,-2);
\fill[white] (60,-0.2) ellipse (1cm and 0.6cm);
\draw[thick] (57.7,2)..controls +(0,-3) and +(0,+2) .. (61.7,-2);
\draw[thick] (58,2)..controls +(0,-2.5) and +(0,+2) .. (62,-2);
\draw[thick] (58.3,2)..controls +(0,-2) and +(0,+2) .. (62.3,-2);
\node at (64,0) {$\dots$};
\node at (65.5,0) {$,$};

\end{tikzpicture}
\end{center}
where all strands starting from ellipses $1,\dots,i-1,i+2,\dots,n$ are vertical. The element $\Sigma_i^{(0)}$ is the identity element, while $\Sigma_i^{(k)}$ are the elements called $\Sigma_i$ satisfying the braid relations in the previous paragraph.

In addition to the $q$-numbers $[L]_q:=\frac{q^L-q^{-L}}{q-q^{-1}}$ and the $q$-factorial $[L]_q!:=[1]_q[2]_q\dots[L]_q$, we also define the $q$-binomials and the $q$-Pochhammer symbol:
\[
  \left[\begin{array}{c}L \\p\end{array}\right]_q:=\frac{[L]_q!}{[L-p]_q![p]_q!}\ ,\ \ \ \ \ \  (a\, ;\, q)_p=\prod_{r=0}^{p-1}(1-aq^r)\ .
\]
By convention, we have $[0]_q!=\left[\begin{array}{c}L \\0\end{array}\right]_q=(a\, ;\, q)_0=1$. 

\begin{framedtheo}[\cite{publi-CP2}]\label{th:bax}
The following function taking values in $H_{\bk,n}(q)$
 \begin{equation}
 \check R_i( u )=\sum_{p=0}^k (-q)^{k-p}\left[\begin{array}{c}k \\p\end{array}\right]^2_q\ 
 \frac{ (q^{-2}\, ;\, q^{-2})_{k-p}  }{ (uq^{-2p}\,;\, q^{-2})_{k-p}}\ \Sigma_i^{(p)}\ , \label{eq:R}
 \end{equation}
satisfies the braided Yang--Baxter equation:
\[
\check  R_i(u)\check R_{i+1}(uv)\check R_i(v)=\check R_{i+1}(v) \check R_i(uv) \check R_{i+1}(v)\ .
\]
\end{framedtheo}
The formula above can be called a ``Baxterization'' formula. The terminology is due to V. Jones \cite{Jo3}. Roughly, a Baxterization formula is a way to obtain a solution of the YB equation starting with a solution of the braid group relation. In other words, this is a way to add the spectral parameters to the braid relation. Of course, the Baxterization formulas depend on the additional relations satisfied by the generators of the braid group (unfortunately, there is no such thing as a Baxterization formula valid at the level of the braid group itself). In Chapter \ref{chap-prel}, we have given examples for Hecke and BMW algebras. The formula above is the Baxterization formula for the fused Hecke algebra.
\begin{coro}
Combining the formula above with the representations of $H_{\bk,n}(q)$ on vector spaces $S_q^k(V)^{\otimes n}$, we get matrix solutions of the Yang--Baxter equation on $q$-symmetrized powers.
\end{coro}
We emphasize that for each value of $k$, we have an infinite family of matrix solutions of the Yang--Baxter equation: one for each choice of $N$, the dimension of $V$. If $N=2$, this solution acts on the spin $k/2$ representations of $U_q(sl_2)$.

\begin{rema}$\ $\\
$\bullet$ Note that the (rational) dependence on the spectral parameter $u$ is clearly visible from the formula in $H_{\bk,n}(q)$ (\emph{i.e.} before the representations).\\[0.1cm]
$\bullet$ The formula above is already interesting for $q=1$. Note that the limit is immediate to perform (all quantum numbers become usual numbers). The solutions live in the algebras of fused permutations and their matrix versions act on usual symmetrized powers.
\end{rema}

\begin{exam}
For $k=1$, the formula above is the usual Baxterization formula of the Hecke algebra. For $k=2$, the formula is:
\[
\check R_i(u)=\Sigma_i^{(2)} - (q+q^{-1})\frac{(q^2-q^{-2}) }{1-uq^{-2}}\Sigma_i^{(1)} +q^2 \frac{(1-q^{-2})(1-q^{-4})}{(1-u)(1-uq^{-2})} \ .
\]
\end{exam}

\section{Conclusion and outlook}

We can summarise the main achievements in this chapter with the following conclusion.
\begin{framedconc}
We have constructed an algebra $H_{\bk,n}(q)$ called the fused Hecke algebra:
\begin{itemize}
\item[$\bullet$] it contains a quotient of the braid group algebra;
\item[$\bullet$] it contains an abstract solution of the Yang--Baxter equation;
\item[$\bullet$] it has representations on vector spaces $W^{\otimes n}$, where $W=S^k(V)$ is the $q$-symmetrized power of $V$, and through such representations it is in Schur--Weyl duality with $U_q(gl_N)$, where $N=\dim(V)$.
\end{itemize}
\end{framedconc}

$\bullet$ Due to their description as centralisers of quantum group representations, the Hecke and Temperley--Lieb algebras can be seen as containing Hamiltonians of spin chain models (with $U_q(sl_N)$-symmetry). In fact, the whole construction of transfer matrices and of the Bethe subalgebras of commuting elements can be performed directly inside these algebras. Of course, to perform the ``algebraic Bethe Ansatz'' machinery, one needs the solution of the Yang--Baxter equation inside the algebra. This is exactly where we are with the fused Hecke algebra and thus all the ingredients are here. From this point of view, it seems also desirable to consider an affine version of the fused Hecke algebras, in order to add into the picture non-trivial boundary conditions and the reflection equation.

\vskip .2cm
$\bullet$ The particular situation when only the first symmetric power is non-trivial of the fused Hecke algebra is the simplest one. It corresponds to the idea of ``fusing'' only the first representation in the tensor product, and thus to study particular boundary conditions. When further restricted to $sl(2)$, our definition corresponds to the so-called seam algebra introduced recently \cite{MRR} in a context related to conformal field theory. For the seam algebra the non-semisimple representation theory is quite well-understood. A similar study for the fused Hecke algebra in general should surely prove also very interesting. The point of view of parabolic Hecke algebra puts forward the idea of using the Kazhdan--Lusztig theory. This could also be useful for our questions about ideals and centralisers (see for example \cite[Chap. 22]{EMTW}).

\vskip .2cm
$\bullet$ The formulation of the fused Hecke algebra that we used makes quite obvious how to generalise it to other settings. For example, one can consider different representations of $U_q(gl_N)$ and/or other quantum (super)groups. The starting point of the approach clearly generalises as follows. One can still consider fused braids, but use a different procedure for multiplying them. Indeed one can replace the Hecke algebras by other quotients of the braid group algebras and/or one can replace the $q$-symmetrisers by other idempotents. The case of the Temperley--Lieb algebra was briefly considered in \cite{LZ2}. Another example is to replace the Hecke algebra by the BMW algebra and consider analogues of the $q$-symmetriser. The diagrammatic point of view on these algebras (namely, the classical limit for $q=1$) is also very interesting, and in this direction, it also seems natural to include the partition algebra in the algebras to be potentially considered.

\section*{Appendix \thechapter.A $\ $Bratteli diagrams for fused Hecke algebras}
\addcontentsline{toc}{section}{Appendix \thechapter.A $\ $Bratteli diagram for fused Hecke algebras}

\paragraph{The chain of algebras $H_{\bk,n}(q)$ when $\bk=(2,2,2,2,\dots)$}\label{app3}

When $\bk=(2,2,2,\dots)$ is the infinite sequence of $2$'s, the Bratteli diagram for the chain of algebras $\{H_{\bk,n}(q)\}_{n\geq0}$ begins as: 
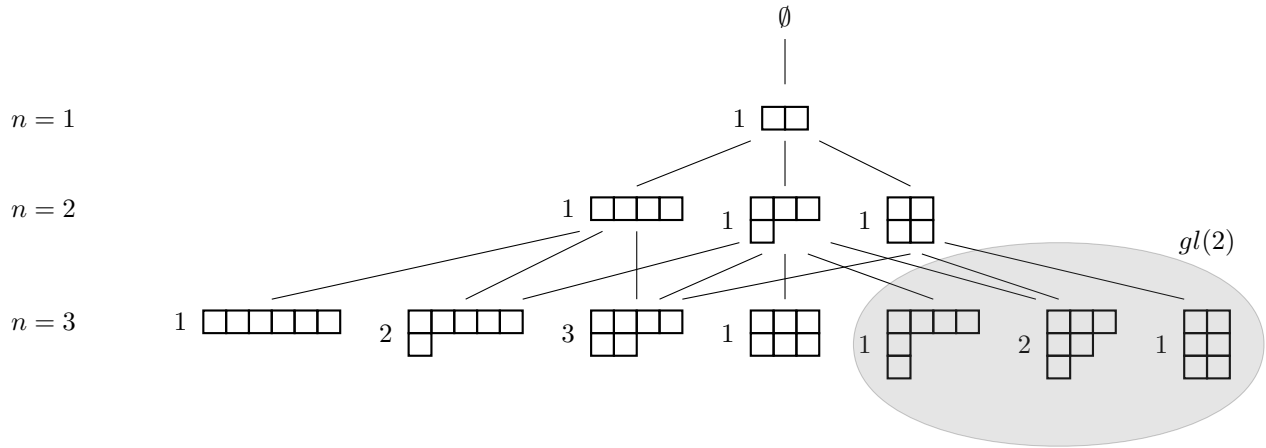
\begin{center}
 \begin{tikzpicture}[scale=0.3]
\node at (0.5,4) {$\emptyset$};
\draw ( 0.5,3) -- (0.5, 1);
\diag{-0.5}{0}{2};\node at (-1.5,-0.5) {$1$};
\draw (-1,-1.5) -- (-6,-3.5);\draw (0.5,-1.5)--(0.5,-3.5);\draw (2,-1.5) -- (6,-3.5);
\diag{-8}{-4}{4};\node at (-9,-4.5) {$1$};\diagg{-1}{-4}{3}{1};\node at (-2,-5) {$1$};\diagg{5}{-4}{2}{2};\node at (4,-5) {$1$};

\draw (-8.5,-5.5) -- (-22,-8.5);\draw (-7.5,-5.5) -- (-13.5,-8.5);\draw (-6,-5.5) -- (-6,-8.5);    \draw (-1.5,-6) -- (-11,-8.5); \draw (-0.5,-6.5) -- (-5,-8.5);\draw (0.5,-6.5) -- (0.5,-8.5);
\draw (1.5,-6.5) -- (7,-8.5);\draw (2.5,-6) -- (11.5,-8.5);\draw (6,-6.5) -- (-4,-8.5); \draw (6.5,-6.5) -- (12.5,-8.5);\draw (7.5,-6) -- (18,-8.5);

\diag{-25}{-9}{6};\node at (-26,-9.5) {$1$};\diagg{-16}{-9}{5}{1};\node at (-17,-10) {$2$};\diagg{-8}{-9}{4}{2};\node at (-9,-10) {$3$};
\diagg{-1}{-9}{3}{3};\node at (-2,-10) {$1$};\diaggg{5}{-9}{4}{1}{1};\node at (4,-10.5) {$1$}; \diaggg{12}{-9}{3}{2}{1};\node at (11,-10.5) {$2$};\diaggg{18}{-9}{2}{2}{2};\node at (17,-10.5) {$1$};

\draw[thin, fill=gray,opacity=0.2] (3.5,-10.5)..controls +(0,6) and +(0,6) .. (21.5,-10.5) .. controls +(0,-6) and +(0,-6) .. (3.5,-10.5);\node at (19,-6) {$gl(2)$};

\node at (-32,-0.5) {$n=1$};\node at (-32,-4.5) {$n=2$};\node at (-32,-9.5) {$n=3$};
\end{tikzpicture}
\end{center}

Note that there is no arrow from $\mu=\begin{array}{cc}
\fbox{\phantom{\scriptsize{$2$}}} &\hspace{-0.35cm}\fbox{\phantom{\scriptsize{$2$}}}\\[-0.2em]
\fbox{\phantom{\scriptsize{$2$}}} &\hspace{-0.35cm}\fbox{\phantom{\scriptsize{$2$}}}
\end{array}$ to $\lambda=\begin{array}{ccc}
\fbox{\phantom{\scriptsize{$2$}}} &\hspace{-0.35cm}\fbox{\phantom{\scriptsize{$2$}}} & \hspace{-0.35cm}\fbox{\phantom{\scriptsize{$2$}}}\\[-0.2em]
\fbox{\phantom{\scriptsize{$2$}}} &\hspace{-0.35cm}\fbox{\phantom{\scriptsize{$2$}}} & \hspace{-0.35cm}\fbox{\phantom{\scriptsize{$2$}}}
\end{array}$ even if $\mu\subset\lambda$ since $\lambda/\mu$ contains two boxes in the same column.
The shaded area indicates the connections between the fused Hecke algebras $H_{\bk,n}(q)$ and the centraliser of the representations of $U_q(gl_N)$, as in Chapter \ref{chap-prel}, Appendix.

\paragraph{The chain of algebras $H_{\bk,n}(q)$ when $\bk=(3,1,1,1,\dots)$}

When $\bk=(3,1,1,1\dots)$, the Bratteli diagram for the chain of
algebras $\{H_{\bk,n}(q)\}_{n\geq0}$ begins as (the shaded areas have a similar meaning as in the preceding examples):
\begin{center}
 \begin{tikzpicture}[scale=0.3]
\node at (0.5,4) {$\emptyset$};
\draw ( 0.5,3) -- (0.5, 1);
\diag{-1}{0}{3};\node at (-2,-0.5) {$1$};
\draw (-0.5,-1.5) -- (-3,-3.5);\draw (1.5,-1.5) -- (3.5,-3.5);
\diag{-5}{-4}{4};\node at (-6,-4.5) {$1$};\diagg{2}{-4}{3}{1};\node at
(1,-5) {$1$};

\draw (-5,-5.5) -- (-10.5,-8.5);\draw (-3,-5.5) -- (-3,-8.5);\draw
(1.7,-6.3) -- (-1,-8.5);\draw (3.5,-6.3) -- (3.5,-8.5);\draw (5,-5.5)
-- (9.5,-8.5);

\diag{-13}{-9}{5};\node at (-14,-9.5) {$1$};\diagg{-5}{-9}{4}{1};\node
at (-6,-10) {$2$};\diagg{2}{-9}{3}{2};\node at (1,-10)
{$1$};\diaggg{8}{-9}{3}{1}{1};\node at (7,-10.5) {$1$};

\draw (-13,-10.5) -- (-22,-14.5);\draw (-10.5,-10.5) --
(-13.5,-14.5);\draw (-5.3,-11.3) -- (-12.5,-14.5);\draw (-4,-11.3) --
(-6,-14.5);\draw (-3,-10.5) -- (5,-14.5);
\draw (1.7,-11.3) -- (-4,-14.5);\draw (3,-11.3) -- (0.5,-14.5);\draw
(4.3,-11.3) -- (11.7,-14.5);\draw (7.7,-12.3) -- (7,-14.5);\draw
(9.3,-12.3) -- (13.5,-14.5);
\draw (11,-10.5) -- (19.5,-14.5);

\node at (-26,-15.5) {$1$};\diag{-25}{-15}{6};
\node at (-17,-16) {$3$};\diagg{-16}{-15}{5}{1};

\node at (-9,-16) {$3$};\diagg{-8}{-15}{4}{2};

\node at (-2,-16) {$1$};\diagg{-1}{-15}{3}{3};

\node at (4,-16.5) {$3$};\diaggg{5}{-15}{4}{1}{1};\node at (11,-16.5)
{$2$};\diaggg{12}{-15}{3}{2}{1};\node at (17,-17)
{$1$};\diagggg{18}{-15}{3}{1}{1}{1};

\draw[thin, fill=gray,opacity=0.2] (3.5,-16)..controls +(0,12) and
+(0,12) .. (22.5,-16) .. controls +(0,-12) and +(0,-12) ..
(3.5,-16);\node at (18,-7) {$gl(2)$};
\draw[thin, fill=gray,opacity=0.2] (16,-16)..controls +(0,5) and
+(0,5) .. (21.5,-16) .. controls +(0,-5) and +(0,-5) .. (16,-16);\node
at (19,-11.5) {$gl(3)$};

\node at (-32,-0.5) {$n=1$};\node at (-32,-4.5) {$n=2$};\node at
(-32,-9.5) {$n=3$};\node at (-32,-15.5) {$n=4$};
\end{tikzpicture}
\end{center}

\setcounter{equation}{0}
\DeactivateToc
\chapter{\huge{Diagonal centralisers}}\label{chap-rac}
\ActivateToc
\addcontentsline{toc}{chapter}{\large{Chapter \thechapter. \hspace{0.2cm}Diagonal centralisers\vspace{0.3cm}}}

\setcounter{minitocdepth}{1}
\minitoc

\section*{Introduction}
\addcontentsline{toc}{section}{Introduction}

So far in our discussion of centralisers, we have been mainly discussing the centralisers of tensor products of representations:
\[V_1\otimes \dots\otimes V_n\,,\]
for a given set of representations of a Lie algebra $\mathfrak{g}$ or a quantum group $U_q(\mathfrak{g})$. In this chapter we will change gears and consider some sort of ``universal'' centraliser which will be related to the above centralisers for any choice of representations $V_1,\dots,V_n$.

To explain how this comes about, let us recall first how we perform a tensor product of Lie algebra representations. First we embed the universal enveloping algebra $U(\mathfrak{g})$ in its $n$-fold tensor product $U(\mathfrak{g})^{\otimes n}$, via the diagonal embedding:
\[\delta^{(n)}\ : \ \begin{array}{rclr}
 U(\mathfrak{g}) & \to & U(\mathfrak{g})^{\otimes n}\\[0.5em]
 g & \mapsto & g\otimes 1\otimes \dots \otimes 1+\dots+1\otimes \dots\otimes 1\otimes g & \quad (g\in\mathfrak{g})
 \end{array}\]
The morphism $\delta^{(n)}$ is defined by the above formula on the elements $g\in\mathfrak{g}$ and extended to the whole algebra $U(\mathfrak{g})$ multiplicatively. At this step, we are in $U(\mathfrak{g})^{\otimes n}$ and it remains only to naturally represent this algebra on the tensor product of $n$ representations of $U(\mathfrak{g})$.

This suggests the idea that we could study directly the centraliser at the universal level of $U(\mathfrak{g})^{\otimes n}$ before choosing the representations. This is what we will do in this chapter.

\begin{frameddefi}
The diagonal centraliser $Z_n(\mathfrak{g})$ is the centraliser in $U(\mathfrak{g})^{\otimes n}$ of the image of $\delta^{(n)}$:
\begin{align}\label{def-Zn}
 Z_n(\mathfrak{g})=\left\{X\in U(\mathfrak{g})^{\otimes n}~|~[\delta^{(n)}(g),X]=0\,, \quad
 \forall g\in U (\mathfrak{g})\right\}\ .
\end{align}
\end{frameddefi}
For a quantum group, the definition is completely similar, with the diagonal embedding $\delta^{(n)}$ replaced by the iterated coproduct $\Delta^{(n)}$. The morphism $\Delta^{(n)}$ is defined recursively by $\Delta^{(n)}=(\Delta\otimes \text{Id}\otimes\dots\otimes\text{Id})\circ\Delta^{(n-1)}$, where  $\Delta^{(2)}=\Delta$ is the coproduct.

\paragraph{Motivations.} As already indicated, a first motivation is that understanding the diagonal centraliser $Z_n(\mathfrak{g})$ would be an interesting first step towards the study of centralisers of any tensor product of representations. In fact, we have naturally a surjective morphism of algebras:
\[Z_n(\mathfrak{g})\to \text{End}_{U(\mathfrak{g})}(V_1\otimes\dots\otimes V_n)\,,\]
for any choice of  irreducible representations $V_1,\dots,V_n$. In other words, the centralisers of the representations will all be quotients of the same algebra $Z_n(\mathfrak{g})$. This is all very nice only if we can describe the algebra $Z_n(\mathfrak{g})$ in a reasonable way.

The second motivation is that for the simplest examples of such diagonal centralisers, we do find very interesting algebras. The first interesting example is when $\mathfrak{g}=sl_2$ and $n=3$. The algebra that appears here for $U(sl_2)$ is known as the Racah algebra and its deformation corresponding to $U_q(sl_2)$ is known as the Askey--Wilson algebra. The Racah and Askey--Wilson algebras arise in numerous areas of mathematics and physics. Their names come from the theory of orthogonal polynomials and they were found first as the algebras controlling the bispectral properties of the eponymous orthogonal polynomials \cite{GZ, Zhe}. In addition to the theory of orthogonal polynomials and the representation theory of $sl_2$, they have also appeared as symmetry of physical models \cite{GVZ,KMP}, in algebraic combinatorics \cite{GWH,Ter}, and are related to Kauffman bracket skein algebras (see for example \cite{CFGPRV} for a review on Askey--Wilson algebras). These algebras will be discussed in more details in the first section of this chapter. 

Finally, diagonal centralisers can be seen as natural deformations of classical algebras of invariant polynomials. Namely, as subalgebras of $U(\mathfrak{g})^{\otimes n}$ they are filtered by the degrees in the elements of $\mathfrak{g}$ and their graded quotients are algebras of invariant polynomial functions. More precisely, they are the polynomial functions on 
$\mathfrak{g}\times\dots\times\mathfrak{g}$ ($n$ times) that are invariant under simultaneous adjoint actions of $G$. For example, if $G=SL(N)$, we are talking about polynomial functions on $n$-tuples of traceless $N\times N$ matrices, invariant under simultaneous conjugation:
\begin{align}\label{eq:simultaneous_conjug}
f(M_1,\dots,M_n)
 =f(g^{-1}M_1g,\dots,g^{-1}M_ng),\ \ \ \text{for any $g\in SL(N)$.}
\end{align}
This is a classical and well studied theory, see \cite{Pr,Raz,Sib} or \cite{Dre} and
references therein. For example we know that the invariant polynomial functions are generated by the so-called polarized traces:
\begin{align}\label{eq:generatorsinvtfct}
(M_1,\dots,M_n)\mapsto Tr(M_{a_1}\dots M_{a_d})
\end{align}
where $d\geq2$ and $a_1,\dots,a_d\in\{1,\dots,n\}$. Moreover, for small $N$ (less than $5$), it is known that we can restrict the degrees $d$ to be less than $\frac{N(N+1)}{2}$ (for arbitrary $N$, it is conjectured to be true but we only know that we can take $d\leq N^2$). The ideal of relations between these generators is well understood only for $N=2$.

To summarise, the study of $Z_n(\mathfrak{g})$ is closely related to classical invariant theory, and in particular the question of finding a generating set for $Z_n(\mathfrak{g})$ corresponds to what is called the First Fundamental Theorem in invariant theory, while the problem of obtaining defining relations of $Z_n(\mathfrak{g})$ corresponds to the quest of a Second Fundamental Theorem in invariant theory. In this regard, the results that we obtain for $Z_n(sl_2)$ and $Z_2(sl_3)$ in this chapter can be seen as non-commutative analogues of the First and Second
Fundamental Theorems for the corresponding invariant polynomial functions.

\paragraph{organisation of the chapter.} The first section recalls in detail the description of $Z_3(sl_2)$ and the appearance of the Racah algebra. The $q$-version for $U_q(sl_2)$ is also mentioned. In the next section, we completely solve the problem of describing $Z_n(sl_2)$ for any $n$, defining and studying along the way an analogue of the usual Racah algebra. Amusing combinatorial considerations on the Hilbert--Poincar\'e series are given in an Appendix. The third section consists in the first attempt, as far as we know, to consider diagonal centralisers outside of the $sl_2$ case. It turns out that there is much to say about $Z_2(sl_3)$ and that we are able to give a complete description. A surprising symmetry under the Weyl group of type $E_6$ is made explicit. The final section gives an application to the problem of labelling vectors in a tensor product of $sl_3$ representations. Our algebraic description of $Z_2(sl_3)$ is used and new symmetries of this classical problem in mathematical physics are obtained.

\paragraph{References.} The publications relevant for this chapter, where more details can be found, are \cite{CPV3,CGPV,CFGPRV,CPV2,CPV1}.

\section{The Racah and the Askey--Wilson algebras}

\subsection{The Racah algebra}

Let us start with the simplest simple Lie algebra $\mathfrak{g}=sl_2$. The algebra $U(\sl_2)$ is generated by elements $e_{ij}$, $i,j\in\{1,2\}$, with the
defining relations:
$$[e_{ij},e_{kl}]=\delta_{jk}e_{il}-\delta_{li}e_{kj}\ \ \ \ \ \text{and}\ \ \ \ \ e_{11}+e_{22}=0\ .$$
We will describe the diagonal centraliser $Z_n(sl_2)$ up to the first non-trivial case $n=3$, where we will find an incarnation of the Racah algebra.

\paragraph{n=1.} The Casimir element $C$ of $U(sl_2)$ is:
\[C=\sum_{i,j=1}^2e_{ij}e_{ji}\ ,\]
and is a generator of the center of $U(sl_2)$, that is of $Z_1(sl_2)$. This concludes the (well-known) algebraic description of $Z_1(sl_2)$ as a polynomial algebra with one generator (namely, $C$).

\paragraph{n=2.} There are several options to produce elements of the diagonal centraliser $Z_2(sl_2)$, starting from the Casimir element $C$. One can use the different natural embedding of $U(sl_2)$ inside $U(sl_2)^{\otimes 2}$ (not forgetting the diagonal one):
\[C_1=C\otimes 1\,,\ \ \ \ \ \ C_2=1\otimes C\,,\ \ \ \ \ \ \ C_{12}=\delta^{(2)}(C)\ .\]
To be completely explicit before moving on, note that $C_{12}=\sum_{i,j=1}^2(e_{ij}\otimes 1+1\otimes e_{ij})(e_{ji}\otimes 1+1\otimes e_{ji})$. It turns out that $C_1,C_2,C_{12}$ commute, are algebraically independent, and generate the whole diagonal centraliser $Z_2(sl_2)$. This leads to the algebraic description of $Z_2(sl_2)$.
\begin{prop}
The diagonal centraliser $Z_2(sl_2)$ is isomorphic to a polynomial algebra in three variables (namely, $C_1,C_2,C_{12}$).
\end{prop}

\paragraph{n=3.} The first interesting situation occurs for $n=3$. Consider the 6 following ways of producing elements of $U(sl_2)^{\otimes 3}$ out of the Casimir element $C$ (we introduce notations for later use):
\begin{gather}\label{gen-Racah}
 k_1=C_1=C\otimes 1\otimes 1,\quad k_2=C_2=1\otimes C\otimes 1,\quad
 k_3=C_3=1\otimes 1\otimes C,\quad k_4=C_{123}=\delta^{(3)}(C),\nonumber\\[0.5em]
 X=C_{12}=\delta^{(2)}(C)\otimes 1,\qquad Y=C_{23}=1\otimes \delta^{(2)}(C)\ .
\end{gather}
It is easy to check that, by construction, these 6 elements belong to the diagonal centraliser $Z_3(sl_2)$ (one might notice that there ought to be also an element $C_{13}$, but it turns out that it is equal to a linear combination of the others). It is also immediate that $C_1,C_2,C_3,C_{123}$ commute with every element in $Z_3(sl_2)$. However, $C_{12}$ and $C_{23}$ do not commute and account for the non-commutativity of $Z_3(sl_2)$\footnote{The commutativity of $Z_2(sl_2)$ was related to the ``no-multiplicity'' phenomenon for 2-fold tensor products of irreducible representations of $sl_2$. On the contrary, $Z_3(sl_2)$ must be non-commutative since multiplicities occur when decomposing a 3-fold tensor product.}.

So we seek the algebraic description of $Z_3(sl_2)$ and in order to do so, we look for the commutation relation between $C_{12}$ and $C_{23}$. The rough idea is as follows. Since they do not commute, we calculate their commutator. If this new element writes in terms of the others, then we are done. Otherwise, we take this new element, and we calculate its commutator with $C_{12}$ and $C_{23}$. This could go on for a while, but it does not. We find the following relations:
\begin{equation}\label{relR1}
\begin{array}{l}
 k_1,k_2,k_3,k_4\quad \text{commute with all elements},\\[0.4em]
 \left[X,Y\right]=Z,\\[0.4em]
 \left[X,Z\right]=4\{X,Y\}+4X^2-4(k_1+k_2+k_3+k_4)X+4(k_1-k_2)(k_4-k_3),\\[0.4em]
 \left[Z,Y\right]=4\{X,Y\}+4Y^2-4(k_1+k_2+k_3+k_4)Y+4(k_3-k_2)(k_4-k_1),
\end{array}
\end{equation}
where $\{X,Y\}=XY+YX$ is the anticommutator. We find quadratic relations that we may think as generalising Lie algebra relations. Now here is the ``miracle''. These relations also appear in a seemingly very different context.

\paragraph{The Racah algebra and the Racah polynomials.} The relations (\ref{relR1}) define what is called the Racah algebra, which we will denote $R(3)$. A first approach to the Racah algebra is from the theory of orthogonal polynomials.  The
Racah polynomials form a family of bispectral classical orthogonal polynomials $\{R_n(x)\}_{n\geq 0}$ which sit on top of the Askey scheme and which are characterized
by a difference and recurrence operator (see for example \cite{KLS}). We will be very sketchy. The recurrence and difference relations look like:
\[\lambda(x)R_n(\lambda(x))=A_nR_{n+1}(\lambda(x))-(A_n+C_n)R_n(\lambda(x))+C_nR_{n-1}(\lambda(x))\,,\ \ \ \ \text{with $\lambda(x)=x(x+\gamma+\delta+1)$\,,}\]
\[n(n+\alpha+\beta+1)R_n(\lambda(x))=B(x)R_{n}(\lambda(x+1))-(B(x)+D(x))R_n(\lambda(x))+D(x)R_{n}(\lambda(x-1))\,,\]
where the coefficients $A_n,C_n,B(x),D(x)$ are explicitly given in terms of four parameters $\alpha,\beta,\gamma,\delta$. 

Now if we take the difference operators on the right hand side of the second equation and call it $X$ and if we denote $Y$ the operator of multiplication by $\lambda(x)$, we find that (as operators on polynomials) $X$ and $Y$ satisfy exactly the relations (\ref{relR1}) with $k_1,k_2,k_3,k_4$ given explicitly in terms of $\alpha,\beta,\gamma,\delta$.

\vskip .2cm
Long story short, the meaning of the Racah algebra in connection with the Racah polynomials is the following. Say we have two diagonalizable matrices $X$ and $Y$. Say also that we know how to diagonalize $X$ and we would like to diagonalize $Y$ (or vice versa). If $X$ and $Y$ satisfy the relations of the Racah algebra then the change of basis between the eigenbasis of $X$ and the eigenbasis of $Y$ will be given in terms of the Racah polynomials. In a physics language, we say that the overlaps between $X$ and $Y$ are given in terms of the Racah polynomials. We refer to \cite{DBI+,GVZ,KMP} and references therein for more details and for the other appearances of the Racah algebra in particular as symmetry algebra of superintegrable models.

What we have learnt in this short paragraph is that in representations, the ``overlaps'' between the operators $X$ and $Y$ of $Z_3(sl_2)$ will be given by Racah polynomials. But finding these overlap coefficients is a well-studied question in physics, they are called $6j$-coefficients or Racah coefficients. Thus the algebraic structure of the diagonal centraliser $Z_2(sl_3)$ explains at once why the Racah or $6j$-coefficients in the recoupling theory of $SU(2)$ are expressed in terms of Racah polynomials\footnote{Of course, this is not the historical way. Racah coefficients were calculated first (by Racah \cite{Rac}) and they were later identified with some orthogonal polynomials, hence called Racah polynomials.}.

\paragraph{Back to the diagonal centraliser $Z_3(sl_2)$.} We have not completely concluded our algebraic description of $Z_3(sl_2)$. It is not exactly isomorphic to the Racah algebra $R(3)$, even if this subtlety was somewhat overlooked in the literature. This can be seen as follows. Thanks to the commutation relations and an application of the diamond lemma (for example as in \cite{CPV2}), one can show that the Racah algebra $R(3)$ admits a nice PBW basis:
\[
 k_1^ik_2^jk_3^kk_{4}^mX^nT^pZ^q, \qquad\qquad i,j,k,m,n,p,q\in\mathbb{Z}_{\geq 0}\ .
\]
However, the images of these elements are not linearly independent in $Z_3(sl_2)$, this algebra is more complicated. In fact, using the natural notion of degree in $U(sl_2)^{\otimes 3}$, the Hilbert--Poincar\'e series of $Z_3(sl_2)$, which records the dimensions of the graded components, is:
\[
 \frac{1-t^6}{(1-t^2)^6(1-t^3)}\ .
\]
The $-t^6$ in the numerator is not compatible with the basis above. There is a hidden relation (of degree 6) satisfied in $Z_3(sl_2)$. This relation turns out to be not so bad. Indeed, there is a central element in the Racah algebra $R(3)$:
\[\Gamma=Z^2-8(XYX+YXY)+4(k_1+k_2+k_3+k_4-4)\{X,Y\}-8(k_1-k_2)(k_4-k_3)Y-8(k_3-k_2)(k_4-k_1)X\]
Now the additional relation satisfied in $Z_3(sl_2)$ (and not implied by the relations of $R(3)$) is simply fixing the values of the central element $\Gamma$:
\begin{equation}\label{relR2}
 \Gamma=8(k_1-k_2+k_3-k_4)(k_1k_3-k_2k_4)-32(k_1k_3+k_2k_4)\ .
\end{equation}
This relation can be seen as expressing $Z^2$ in terms of other elements and leads to the following basis for $Z_3(sl_2)$:
\[
 k_1^ik_2^jk_3^kk_{4}^mX^nY^pZ^q, \qquad\qquad i,j,k,m,n,p,q\in\mathbb{Z}_{\geq 0},\qquad q\in\{0,1\},
\]
This now is in accordance with the Hilbert--Poincar\'e series above. Here is the conclusion.
\begin{framedconc}
If we define the special Racah algebra $sR(3)$ with (\ref{relR1}) and (\ref{relR2}), we have that:
\[\text{The diagonal centraliser $Z_3(sl_2)$ is isomorphic to $sR(3)$.}\]
\end{framedconc}

\paragraph{The Askey--Wilson algebra.} There is a completely similar story for the quantum group $U_q(sl_2)$. The Askey--Wilson algebra $AW(3)$ is defined by generators $C_{12}$, $C_{23}$, $C_{13}$,
central elements $C_1$, $C_2$, $C_3$, $C_{123}$, together with the following defining relations:
\[\begin{array}{ll}
 \displaystyle C_{12}+\frac{[C_{23},C_{13}]_q}{q^{2}-q^{-2}}&=\displaystyle \frac{C_1C_2+C_3C_{123}}{q+q^{-1}},\\[1em]
 \displaystyle C_{23}+\frac{[C_{13},C_{12}]_q}{q^{2}-q^{-2}}&=\displaystyle \frac{C_2C_3+C_1C_{123}}{q+q^{-1}},\\[1em]
 \displaystyle C_{13}+\frac{[C_{12},C_{23}]_q}{q^{2}-q^{-2}}&=\displaystyle \frac{C_3C_1+C_2C_{123}}{q+q^{-1}},
\end{array}\]
where the $q$-commutator is defined by $[A,B]_q=qAB-q^{-1}BA$. As for the classical case, the diagonal centraliser turns out to be a quotient of $AW(3)$, and again this quotient amounts to fixing the value of a certain central element of $AW(3)$. 

The Askey--Wilson algebra has the same sort of interpretation in the context of orthogonal polynomials as the Racah algebra. We refer to \cite{CFGPRV} and references therein for more details.

\section{The diagonal centraliser $Z_n(sl_2)$}

\subsection{The Racah algebras $R(n)$}

\paragraph{Definition.} We consider the following definition of a generalisation of the Racah algebra. It is an easy exercise to check that for $n=3$, the algebra coincides (even though the presentation is slightly different) with the Racah algebra $R(3)$ of the previous section\footnote{To be precise, one has to use the relation $C_1+C_2+C_3+C_{123}=C_{12}+C_{23}+C_{13}$ connecting $C_{123}$ and $C_{13}$.}.
\begin{defi}
The Racah algebra $R(n)$ is the associative algebra with generators:
\[
 \p_{ij},\ 1\leq i\leq j\leq n\qquad \text{and}\qquad\f_{ijk},\ 1\leq i<j<k\leq n,
\]
and the defining relations are, for all possible indices $i,j,k,l,m$ in $\{1,\dots,n\}$:
\[\begin{array}{rl}
 \p_{ii}&\quad \text{is central},\\[0.4em]
 \left[\p_{ij},\p_{k\ell}\right]&=0\hspace{6em}
 \text{if both $i,j$ are distinct from $k,\ell$},\\[0.4em]
 \left[\p_{ij},\p_{jk}\right]&=2\f_{ijk},\\[0.4em]
 \left[\p_{jk},\f_{ijk}\right]&=\p_{ik}(\p_{jk}+\p_{jj})-(\p_{jk}+\p_{kk})\p_{ij},\\[0.4em]
 \left[\p_{k\ell},\f_{ijk}\right]&=\p_{ik}\p_{j\ell}-\p_{i\ell}\p_{jk},\\[0.4em]
 \left[\f_{ijk},\f_{jk\ell}\right]&=-(\f_{ij\ell}+\f_{ik\ell})\p_{jk},\\[0.4em]
 \left[\f_{ijk},\f_{k\ell m}\right]&=\f_{i\ell m}\p_{jk}-\f_{j\ell m}\p_{ik},
\end{array}\]
where in each relation all indices involved are distinct and
$\p_{ij}$ and $\f_{ijk}$ are defined by:
\[
 \p_{ij}=\p_{ji}\qquad \text{and}\qquad \f_{ijk}=-\f_{jik}=\f_{jki}\qquad
 \text{for any $i,j,k\in\{1,\dots,n\}$.}\]
\end{defi}
This is a lot of defining relations. However, the logic is clear, we give as many relations as are needed to close the algebra under commutators starting from the generators $\p_{ij}$. This is the direct generalisation of the definition of the usual Racah algebra $R(3)$. The third relation can be seen as defining the elements $\f_{ijk}$, and all the others assert that all commutators can be expressed in terms of the generators.

\paragraph{Connection with the diagonal centraliser $Z_n(sl_2)$.} Of course, the above definition would be somewhat meaningless for our purposes if these relations were not satisfied in the diagonal centraliser $Z_n(sl_2)$. So we state the first main result at once. Below, $C$ is the Casimir element of $U(sl_2)$ as in the previous section and we have used Sweedler notation 
$\delta^{(2)}(C)=\sum C_{(1)}\otimes C_{(2)}$.
\begin{framedtheo}[\cite{CGPV}]\label{thm:Racah1}
There is a surjective morphism of algebras from $R(n)$ onto $Z_n(sl_2)$ given by:
\[\p_{ii}\mapsto \frac{1}{2}C_i\,,\qquad \p_{ij}\mapsto C_{ij}-C_{i}-C_{j}\qquad\text{and}\qquad \f_{ijk}\mapsto \frac{1}{2}[C_{ij},C_{jk}]\,,\qquad\text{for $i,j,k$ all distinct,}\]
where:
\[C_{i}=1^{\otimes(i-1)}\otimes C\otimes1^{\otimes(n-i)},\qquad
 C_{ij}=\sum 1^{\otimes(i-1)}\otimes C_{(1)}\otimes1^{\otimes(j-i-1)}\otimes
                                C_{(2)}\otimes1^{\otimes(n-j)}\ .\]
\end{framedtheo}
The fact that the defining relations of $R(n)$ are satisfied by the given elements of $Z_n(sl_2)$ is a straightforward calculation (they were calculated like that in the first place). The less straightforward statement is the surjectivity. This amounts to the fact that $Z_n(sl_2)$ is generated by the intermediate Casimir elements $C_i,C_{ij}$. For this, we used results in classical invariant theory.

\paragraph{Central elements in $R(n)$.} As for $R(3)$, the defining relations of $R(n)$ are not enough for a complete description of $Z_n(sl_2)$. This is where the difficult part of the story is. Remarkably, the relations of $R(n)$ imply the existence of several families of central elements in $R(n)$, and thus in $Z_n(sl_2)$. These are the analogues of the central element $\Gamma$ of $R(3)$.

We will need the notion of determinant for matrices with non-commuting entries. If $A$
is a $n\times n$ matrix with entries $A_{i,j}$ ($1 \leq i,j \leq n$), we define the
symmetrized determinant of $A$ as follows
\begin{align*}
 \det A=\frac{1}{n!} \sum_{\rho, \sigma\in
S_{n}}sgn(\rho) sgn(\sigma)A_{\rho(1),\sigma(1)}A_{\rho(2),\sigma(2)}\dots
 A_{\rho(n),\sigma(n)},
\end{align*}
where $S_n$ is the permutation group of $n$ elements and $sgn(\sigma)$ is the signature of
$\sigma$. For commuting entries, it is the usual definition of the determinant of a matrix.
We define also the following $3\times 3$ matrix
\begin{align*}
 \p_{ijk}^{\,abc}=
 \begin{pmatrix}
  \p_{ia} & \p_{ib} & \p_{ic}\\
  \p_{ja} & \p_{jb} & \p_{jc}\\
  \p_{ka} & \p_{kb} & \p_{kc}
 \end{pmatrix}.
\end{align*}
We now introduce some elements of $R(n)$ that will play an important part.
\begin{align*}
  w_{ijk}&:={\f_{ijk}}^2 + \tfrac12\det( \p_{ijk}^{\,ijk} )\\
  &\qquad
   -\tfrac13\left( \{\p_{ij},\p_{ik}\} + \{\p_{ij},\p_{jk}\} + \{\p_{ik}\p_{jk}\}
                 + \p_{ij}\p_{kk} + \p_{ik}\p_{jj} + \p_{jk}\p_{ii} \right),
 \\[0.4em]
& x_{ijk\ell}: = \f_{ijk}\f_{jk\ell}+\tfrac12\det(\p_{ijk}^{jk\ell})
 + \tfrac12(\f_{ij\ell}+\f_{ik\ell})\p_{jk}
 - \tfrac13(\p_{ij}\p_{k\ell} + \p_{ik}\p_{j\ell} + \p_{i\ell}\p_{jk}),
\\[0.4em]
& y_{ijk\ell m}: = \f_{ijk}\f_{k\ell m} + \tfrac12\det(\p_{ijk}^{k\ell m})
               + \tfrac12(\f_{ij\ell}\p_{km} - \f_{ijm}\p_{k\ell}),
\\[0.4em]
& z_{ijk\ell mp} := \f_{ijk}\f_{\ell mp} + \tfrac12\det(\p_{ijk}^{\ell mp}),
\end{align*}
where indices $i,j,k,\ell,m,p\in\{1,\dots,n\}$ are all distinct. Here is a first summary of what we can prove using the defining relations of the algebra $R(n)$ (we will say more later about $y$ and $z$).
\begin{prop}[\cite{CGPV}]\label{pr:R6}
For all distinct $i,j,k,\ell,m,p$, the elements $w_{ijk}$, $x_{ijk\ell}$, $y_{ijk\ell m}$ and $z_{ijk\ell mp}$ are invariant under permutation of their indices, and are central in $R(n)$.
\end{prop}

\subsection{The special Racah algebra $sR(n)$ and $Z_n(sl_2)$}

Now that we have a working definition of an algebra $R(n)$, analogue of the Racah algebra $R(3)$, we can fulfil our goal to describe the diagonal centraliser $Z_n(sl_2)$. As for $n=3$, it is going to be a quotient of $R(n)$ thanks to Theorem \ref{thm:Racah1}. Remarkably it is enough to fix the values of the central elements of $R(n)$ discovered above.

What is hidden in this discussion is that we use a precise knowledge of the ideal of relations for the commutative algebra of invariants (see \cite{Dre} and references therein). In fact, we know a convenient set of generators (a Groebner basis) for this ideal, and we need to find the non-commutative version (below, we give only our results in the non-commutative case, for $Z_n(sl_2)$).

First, it turns out that some of the relations we need are already consequences of the defining relations of $R(n)$. This is summarised as follows.
\begin{prop}[\cite{CGPV}]
The following assertions are true in $R(n)$:
\begin{itemize}
\item The relations below hold for $1\leq a \leq n$ and $1\leq
i< j< k < \ell \leq n$:
\begin{align*}
 & \p_{ai} \f_{jk\ell}-\p_{aj} \f_{ik\ell}+\p_{ak} \f_{ij\ell}-\p_{a\ell} \f_{ijk}=0.
\end{align*}
\item For all distinct $i,j,k,\ell,m,p\in\{1,\dots,n\}$, we have:
\begin{align*}
 y_{ijk\ell m}=0\qquad\text{and}\qquad z_{ijk\ell mp}=0\ .
\end{align*}
\end{itemize}
\end{prop}
This is not enough, and so we really need to quotient out $R(n)$ by additional relations. Everything was designed such that these new relations amount to cancel the remaining central elements of $R(n)$.
\begin{defi}
The special Racah algebra $sR(n)$ is the quotient of $R(n)$ by:
\begin{align*}
 w_{ijk}=0\qquad\text{and}\qquad x_{ijk\ell}=0,\qquad\forall 1\leq i<j<k<\ell\leq n\ .
\end{align*}
\end{defi}

\begin{framedtheo}[\cite{CGPV}]\label{thm:Racah2}
The morphism of Theorem \ref{thm:Racah1} factors through the quotient $sR(n)$ and becomes an isomorphism:
\begin{align*}
sR(n)\cong Z_n(\sl_2)\ .
\end{align*}
\end{framedtheo}
It is quite remarkable that we only need to add to $R(n)$ the relations above in order to recover the centraliser for any value of $n$. Indeed, they use only up to 4 indices and one could have expected that we would need to quotient by elements spanning an
increasing number of indices when $n$ gets bigger. This is a simplification due to the non-commutative setting (many relations are already implied in $R(n)$) and in this regard, the non-commutative situation is simpler than the commutative one.

\section{The diagonal centraliser $Z_2(sl_3)$}

From now on, we will focus on the Lie algebra $sl_3$, for which the two-fold tensor product already produces a very interesting algebraic situation.

The algebra $U(sl_3)$ is generated by elements $e_{ij}$, $i,j\in\{1,2,3\}$, with the
defining relations:
$$[e_{ij},e_{kl}]=\delta_{jk}e_{il}-\delta_{li}e_{kj}\ \ \ \ \ \text{and}\ \ \ \ \ e_{11}+e_{22}+e_{33}=0\ .$$
The diagonal embedding $\delta^{(2)}$ from $U(sl_3)$ to $U(sl_3)\otimes U(sl_3)$ is as in the introduction and the diagonal centraliser $Z_2(sl_3)$ is defined accordingly.

\subsection{Generators of $Z_2(sl_3)$} Following the same line of thought as for $sl_2$, we first produce elements of $Z_2(sl_3)$ using the Casimir elements of $U(sl_3)$ (\emph{i.e.} generators of the center of $U(sl_3)$). There are two Casimir elements in $U(sl_3)$, one of degree 2 and one of degree 3. For symmetry reasons, we choose them as follows:
\[
C^{(2)}= \sum_{i,j=1}^3 e_{i j} e_{j i}\ \ \ \text{and}\ \ \ C^{(3)}=  \frac{1}{2}\sum_{i,j=1}^3 (e_{ij} e_{jk}e_{ki}+e_{ji} e_{kj}e_{ik})\ .
\]
They generate the center of $U(sl_3)$. Out of these two elements, we can produce 6 elements of $Z_2(sl_3)$:
\[
\begin{array}{l}
k_1=C^{(2)}\otimes 1\,,\ \ k_2=1\otimes C^{(2)}\,,\ \ k_3=\delta^{(2)}(C^{(2)})\,,\\[0.5em]
l_1=C^{(3)}\otimes 1\,,\ \ l_2=1\otimes C^{(3)}\,,\ \ l_3=\delta^{(2)}(C^{(3)})\ .
\end{array}
\] 
Now here comes the main difference with $sl_2$: these elements do not generate the diagonal centraliser $Z_2(sl_3)$. In fact, this must be since by construction these elements all commute together (and in fact must be central in $Z_2(sl_3)$). However, the diagonal centraliser $Z_2(sl_3)$ cannot be commutative\footnote{this is because multiplicities arise in the tensor product of two irreducible $sl_3$-representations.}. So we need to come up with a better way to find elements of $Z_2(sl_3)$.

\paragraph{Polarized traces.} We find inspiration from the classical construction of invariants. We define the following elements of $U(sl_3)^{\otimes 2}$:
\[T^{(a_1,\dots,a_d)}=\sum_{i_1,\dots,i_d=1}^3e_{i_2i_1}^{(a_1)}e_{i_3i_2}^{(a_2)}\dots e_{i_1i_d}^{(a_d)}\,,\ \ \ \ \ \quad\text{where $a_1,\dots,a_d\in\{1,2\}$\ .}\]
The notation used is $e_{ij}^{(1)}=e_{ij}\otimes 1$ and $e_{ij}^{(2)}=1\otimes e_{ij}$. By a straightforward exercise using the commutation relations in matrix algebras, one can show that these elements commute with $\delta^{(2)}(e_{pq})$, for any $p,q$. Therefore they belong to $Z_2(sl_3)$. They are the non-commutative analogues of the polarized traces from the introduction, and we will call them polarized traces as well.

Skipping all details, we can use some results from classical invariant theory (see \cite{Dre} and references therein) to obtain that:
\[
T^{(1,1)},\ T^{(1,2)},\ T^{(2,2)},\ T^{(1,1,1)},\ T^{(1,1,2)},\ T^{(1,2,2)},\ T^{(2,2,2)},\ T^{(1,1,2,2)}\ ,\]
will generate the whole diagonal centraliser $Z_2(sl_3)$. 
\begin{rema}
In fact, an attentive reader comparing with \cite{Dre} may notice that there seems to be one generator of degree 6 missing, such as for example $T^{(1,1,2,2,1,2)}$. Remarkably, this ``missing'' generator is not missing in the non-commutative case since it will appear in the commutator of two of the above elements. This is an instance where the non-commutative version is simpler than its commutative analogue. 
\end{rema}
The elements above $k_1,k_2,k_3,l_1,l_2,l_3$ are directly related to $T^{(1,1)},T^{(1,2)},T^{(2,2)},T^{(1,1,1)},T^{(2,2,2)}$ and $T^{(1,1,2)}+T^{(1,2,2)}$. So we introduce two additional elements:
\[
\begin{array}{l}
\displaystyle X=\frac{1}{2}(T^{(1,1,2)}-T^{(1,2,2)})+\frac{1}{3}(l_1-l_2)\ ,\\[0.8em]
\displaystyle Y=T^{(1,1,2,2)}+\frac{3}{2}(T^{(1,1,2)}+T^{(1,2,2)})-\frac{1}{12}(T^{(1,2)})^2-\frac{5}{12}T^{(1,1)}T^{(2,2)}+\frac{5}{2}T^{(1,2)}\,.
\end{array}
\]
This awkward and seemingly complicated choice of $X$ and $Y$ is dictated by symmetry considerations, that we skip, and we will be largely rewarded for such a choice in the last subsection and in the applications in the next section. Let us draw a conclusion before moving on.

\begin{framedconc}
The diagonal centraliser $Z_2(sl_3)$ is generated by:
\[k_1,\ k_2,\ k_3,\ \ l_1,\ l_2,\ l_3,\ \ \ X,\ Y\ .\]
\end{framedconc}

\subsection{A presentation of $Z_2(sl_3)$}

Now we need to find a complete set of relations for the generators above of $Z_2(sl_3)$. The result will be remarkably similar to the story of $Z_3(sl_2)$ and the Racah algebra. Let us start with the definition of a new algebra, which will be the main character of this section (and of the following).
\begin{defi}\label{def:cA2}
The algebra $\mathcal{A}$ is the algebra generated by central elements $a_2,a_5,a_6,a_8,a_9$ and generators $A,B,C$ with the following defining relations:
\begin{equation}\label{relAlg2}
\begin{array}{l}
[A,B]=C\,,\\[0.5em]
[A,C]= -6B^2 + a_2 A^2+ a_5 A+a_8\,,\\[0.5em]
[B,C]= -2A^3-a_2\{A,B\} -a_5 B+a_6A+a_9\,,
\end{array}
\end{equation}
\end{defi}
We have the following theorem:
\begin{theo}[\cite{CPV2}]\label{prop-morA}
There is a surjective morphism of algebras from $\mathcal{A}$ to $Z_2(sl(3))$ given by:
\[
A\mapsto X\,,\qquad B\mapsto Y\,,\qquad C\mapsto Z:=[X,Y]\,,
\]
for a certain choice of $a_2,a_5,a_6,a_8,a_9\in\mathbb{C}[k_1,k_2,k_3,l_1,l_2,l_3]$.
\end{theo}
There is an explicit expression for $a_2,a_5,a_6,a_8,a_9$ but we will skip it here and come back to it later. There is a nice PBW basis of the algebra $\mathcal{A}$, over the polynomials in $a_2,a_5,a_6,a_8,a_9$, which is:
\[\{A^aB^bC^c\ ,\ \ a,b,c\in\mathbb{Z}_{\geq 0}\}\ .\]
This shows that $\mathcal{A}$ cannot be isomorphic to $Z_2(sl_3)$ since this would contradict the known Hilbert--Poincar\'e series of $Z_2(sl_3)$ (see Appendix). So an additional relation must be found. Here is how this goes.

It turns out that the algebra $\mathcal{A}$ has a non-trivial central element, which is given explicitly as:
\[\Omega =x_1A+x_2B+x_3A^2+x_4\{A,B\}+x_5 B^2+x_6 ABA-A^4+4B^3+C^2\ ,\]
where $x_1= 6 a_5  + 2 a_9\,,$\ \ $x_2= -2a_6 - 2 a_8\,,$\ \ $x_3= 6 a_2 + a_6\,,$\ \ $x_4=-a_5\,,$\ \ $x_5= 8 a_2 -24\,,$\ \ $x_6= -2 a_2 +12$\,.

We can formulate our final result.
\begin{framedtheo}[\cite{CPV2}]
The diagonal centraliser $Z_2(sl_3)$ is isomorphic to the quotient of the algebra $\mathcal{A}$ by the additional relation:
\[\Omega=a_{12}\,,\ \ \ \quad \text{for some $a_{12}\in\mathbb{C}[k_1,k_2,k_3,l_1,l_2,l_3]$\,.}\]
\end{framedtheo}
Again, there is an explicit expression for $a_{12}$ which we skip here. About the notation, the subscripts in the parameters $a_2,a_5,a_6,a_8,a_9,a_{12}$ correspond to their degrees as detailed in the next subsection.

\subsection{Highest-weight specialisation of $Z_2(sl_3)$ and $E_6$ symmetry}

\paragraph{The highest-weight specialisation of $Z_2(sl_3)$.}
Consider a highest weight representation $V_{m_1,m_2}$ of $U(sl_3)$. It is parametrised by two complex numbers $m_1,m_2$, and generated as a $U(sl_3)$-module by a highest weight vector $v_{m_1,m_2}$ satisfying
\begin{equation}\label{high-weight-rep}
\begin{array}{rcl} h_{p}v_{m_1,m_2}&=& (m_p-1) v_{m_1,m_2} \quad \text{with} \quad p=1,2\,,\\[0.5em]
 e_{pq}v_{m_1,m_2}&=&0 \quad \text{with} \quad 1\leq p< q \leq 3 \,,
\end{array}
\end{equation}
where we set $h_p=e_{pp}-e_{p+1,p+1}$. For example, if $m_1,m_2\in\mathbb{Z}_{>0}$ then $V_{m_1,m_2}$ can be the finite-dimensional irreducible representation of $sl_3$ with highest-weight $(m_1,m_2)$.

In the representation $V_{m_1,m_2}$, the Casimir elements $C^{(2)}$ and $C^{(3)}$ of $U(sl_3)$ are proportional to the identity matrix. Let us denote their values by:
\[ C^{(2)}=c^{(2)}(m_1,m_2)\text{Id}_{V_{m_1,m_2}}\ \ \ \ \text{and}\ \ \ \ C^{(3)}=c^{(3)}(m_1,m_2)\text{Id}_{V_{m_1,m_2}}\ .\]
Now pick three pairs of complex numbers $(m_1,m_2)$, $(m'_1,m'_2)$ and $(m''_1,m''_2)$, and look at the tensor product of representations $V_{m_1,m_2}\otimes V_{m'_1,m'_2}$. In this space, consider all the highest-weight vectors under the diagonal action of $U(sl_3)$ corresponding to the weight $(m''_1,m''_2)$. The subspace generated by all these vectors as a $U(sl_3)$-module is denoted by:
\[M_{m_1,m_2,m_1',m_2'}^{m''_1,m''_2}\ .\] 
By construction, the diagonal centraliser $Z_2(sl_3)$ commutes with the diagonal action of $U(sl_3)$, and thus it leaves invariant the subspace $M_{m_1,m_2,m_1',m_2'}^{m''_1,m''_2}$. which thus naturally becomes a $Z_2(sl_3)$-module. 
\begin{exam}
If $V_{m_1,m_2},V_{m'_1,m'_2},V_{m''_1,m''_2}$ are finite-dimensional irreducible representations, then the subspace $M_{m_1,m_2,m_1',m_2'}^{m''_1,m''_2}$ is the isotypic component of $V_{m''_1,m''_2}$ in the tensor product $V_{m_1,m_2}\otimes V_{m'_1,m'_2}$ (the direct sum of all summands in $V_{m_1,m_2}\otimes V_{m'_1,m'_2}$ isomorphic to $V_{m''_1,m''_2}$).
\end{exam}
In the representation $M_{m_1,m_2,m_1',m_2'}^{m''_1,m''_2}$ of $Z_2(sl_3)$, the central parameters $k_1,k_2,k_3,l_1,l_2,l_3$ take definite complex values, since they are expressed in terms of the Casimir elements. Explicitly:
\[\begin{array}{c}
k_1=c^{(2)}(m_1,m_2)\,,\\[0.5em]
l_1=c^{(3)}(m_1,m_2)\,,
\end{array}\ \ \ \ \ \ \begin{array}{c}
k_2=c^{(2)}(m'_1,m'_2)\,,\\[0.5em]
l_2=c^{(3)}(m'_1,m'_2)\,,
\end{array}\ \ \ \ \ \ \begin{array}{c}
k_3=c^{(2)}(m''_1,m''_2)\,,\\[0.5em]
l_3=c^{(3)}(m''_1,m''_2)\,,
\end{array}
\]
We denote by $Z_2(sl_3)^{spec}$ the specialization of $Z_2(sl_3)$ corresponding to these values. In $Z_2(sl_3)^{spec}$, the parameters $a_i$ become polynomials in $m$'s and their subscripts correspond to their degrees:
\[a_2,a_5,a_6,a_8,a_9,a_{12}\in\mathbb{C}[m_1,m_2,m'_1,m'_2,m''_1,m''_2]\ \ \quad\text{with}\ \text{deg}(a_i)=i\ .\]
The numbers $2,5,6,8,9,12$ turn out to be the fundamental degrees of a root system of type $E_6$. Quite amazingly, this is not a mere numerical coincidence, as we will see now.

\paragraph{An action of the Weyl group of type $E_6$.} Let us consider a root system of type $E_6$ and choose the simple roots $\alpha_1,\alpha_2,\alpha_3,\alpha_4,\alpha_5,\alpha_6$ with the numeration according to the following Dynkin diagram:
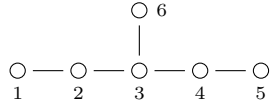
\begin{center}
\begin{tikzpicture}[scale=0.2]
\draw (1,1) circle [radius=0.5];
\node [above] at (1,-1.5) {\scriptsize{$1$}};
\draw (2,1)--(4,1);

\draw (5,1) circle [radius=0.5];
\node [above] at (5,-1.5) {\scriptsize{$2$}};
\draw (6,1)--(8,1);

\draw (9,1) circle [radius=0.5];
\node [above] at (9,-1.5) {\scriptsize{$3$}};
\draw (10,1)--(12,1);

\draw (13,1) circle [radius=0.5];
\node [above] at (13,-1.5) {\scriptsize{$4$}};
\draw (14,1)--(16,1);

\draw (17,1) circle [radius=0.5];
\node [above] at (17,-1.5) {\scriptsize{$5$}};

\draw (9,5) circle [radius=0.5];
\node [right] at (9.5,5) {\scriptsize{$6$}};
\draw (9,2)--(9,4);

\end{tikzpicture}
\end{center}
Let us associate the parameters $m_1,m_2,m'_1,m'_2,m''_1$ and $m''_2$ with the simple roots as follows:
\[
\begin{array}{l}
m_1=\alpha_1\,,\\[0.4em]
m_2=\alpha_2\,,
\end{array}\ \ \ \ \ \ \ \begin{array}{l}
m'_1=\alpha_5\,,\\[0.4em]
m'_2=\alpha_4\,,
\end{array}\ \ \ \ \ \ \ \
\begin{array}{l}
m''_1=\Theta\,,\\[0.4em]
m''_2=-\alpha_6\,,
\end{array}
\]
where $\Theta=\alpha_1+2\alpha_2+3\alpha_3+2\alpha_4+\alpha_5+2\alpha_6$ is the longest positive root of $E_6$. The choice of $(m_1,m_2)$, $(m'_1,m'_2)$ and $(m''_1,m''_2)$ corresponds to three subsystems of type $A_2$ which are pairwise orthogonal. In other words, we single out a subsystem of type $A_2\times A_2\times A_2$ in $E_6$. This is best understood by recalling that the addition of minus the longest root produces the associated affine Dynkin diagram. For $E_6$, the three subsystems of type $A_2$ in the affine Dynkin diagram appear naturally:
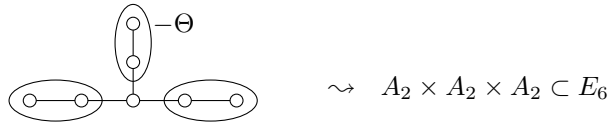
\begin{center}
\begin{tikzpicture}[scale=0.17]
\draw (1,1) circle [radius=0.5];
\draw (1.5,1)--(4.5,1);

\draw (5,1) circle [radius=0.5];
\draw (5.5,1)--(8.5,1);

\draw (9,1) circle [radius=0.5];
\draw (9.5,1)--(12.5,1);

\draw (13,1) circle [radius=0.5];
\draw (13.5,1)--(16.5,1);

\draw (17,1) circle [radius=0.5];

\draw (9,4) circle [radius=0.5];
\draw (9,1.5)--(9,3.5);

\draw (9,7) circle [radius=0.5];
\node [right] at (9.8,7) {$-\Theta$};
\draw (9,4.5)--(9,6.5);

\draw (3,1) ellipse (3.6cm and 1.6cm);
\draw (15,1) ellipse (3.6cm and 1.6cm);
\draw (9,5.5) ellipse (1.4cm and 3cm);

\node at (35,2) {$\leadsto\ \ A_2\times A_2\times A_2\subset E_6$};

\end{tikzpicture}
\end{center}
The Weyl group $W(E_6)$ of the root system of type $E_6$ acts on the set of roots by its usual reflection representation, and through the above identification, this leads to a linear action of $W(E_6)$ on $\mathbb{C}[m_1,m_2,m'_1,m'_2,m''_1,m''_2]$. It is elementary to calculate explicitly this action and we skip the resulting formulas. 

\vskip .2cm
All this may seem quite arbitrary, but here is our surprising discovery.
\begin{framedtheo}[\cite{CPV2}]
The polynomials $a_2,a_5,a_6,a_8,a_9,a_{12}$ appearing in the presentation of $Z_2(sl_3)^{spec}$ are invariant polynomials under the action of $W(E_6)$.\\
In particular, the algebra $Z_2(sl_3)^{spec}$ depends on the parameters $(m_1,m_2,m'_1,m'_2,m''_1,m''_2)$ only through their orbit under $W(E_6)$.
\end{framedtheo}
It is well-known that the subalgebra of invariant polynomials for $W(E_6)$ is generated by 6 algebraically independent polynomials of degrees, respectively, 2,5,6,8,9,12 (the fundamental degrees of $E_6$). The polynomials $a_2,a_5,a_6,a_8,a_9,a_{12}$ are in fact such a set of generators.

To be complete, the algebra $Z_2(sl_3)^{spec}$ is also invariant under the central symmetry $-\text{Id}$ (sending $X$ to $-X$ and $Y$ to $Y$) so that the full symmetry group $\pm W(E_6)$ of the root lattice is realised as a symmetry of $Z_2(sl_3)^{spec}$.

\section{Application: The missing label of $sl_3$ and its symmetry}

In this section, we fix three highest weights for $sl_3$ that we denote collectively by 
$$\tm=(m_1,m_2,m'_1,m'_2,m''_1,m''_2)\in\mathbb{Z}_{>0}\ .$$
We will use the notations:
\[
 \ell=\frac{1}{3}(m_1+2m_2+m'_1+2m'_2-m''_1-2m''_2)\ \ \ \ \text{and}\ \ \ \ n=\frac{1}{3}(2m_1+m_2+2m'_1+m'_2-2m''_1-m''_2)\ .
\]
We consider the corresponding irreducible finite-dimensional representations $V_{m_1,m_2}$, $V_{m'_1,m'_2}$ and $V_{m''_1,m''_2}$. We form the tensor product $V_{m_1,m_2}\otimes V_{m'_1,m'_2}$ and we assume that the third representation $V_{m''_1,m''_2}$ appears in the Clebsch--Gordan series:
\begin{equation}\label{tensorproduct}
V_{m_1,m_2}\otimes V_{m'_1,m'_2}=\bigoplus_{(m''_1,m''_2)} V_{m''_1,m''_2}^{\oplus d_{\tm}}\ .
\end{equation}
The multiplicity $d_{\tm}$ with which the representation $V_{m''_1,m''_2}$ appears in this decomposition is the Littlewood--Richardson coefficient. 

\paragraph{A missing label.} The multiplicity $d_{\tm}$ may be larger than one and this is where a so-called missing label problem appears. This terminology is borrowed from physics, and more precisely from quantum mechanics. The ``labels'' refer to a possibility of indexing the vector basis of the space in (\ref{tensorproduct}) by eigenvalues of some complete set of commuting observables. 

In the left-hand-side (the uncoupled basis), one can work in each factor of the tensor product individually. In this case, there is a standard way of labelling the vectors of, say, $V_{m_1,m_2}$. In physics literature, this is called the total isospin, a component of this isospin and the hypercharge. These observables come from the enveloping algebra $U(sl_3)$ (their precise meaning is not important for our considerations).

However, we would like to have labels referring to the right-hand-side (the coupled basis). That is, we consider that we lost track of the individual factors in the tensor product, but we still have access to observables from the diagonal action of $sl_3$ (total isospin, total hypercharge, etc.). Now using the diagonal embedding of the Casimir elements (in our notations, $k_3$ and $l_3$) we can separate between different irreducible representations appearing in the sum in (\ref{tensorproduct}). And using the (diagonal embedding of) the standard labelling operator we can specify a vector in the given representation $V_{m''_1,m''_2}$. But this cannot distinguish between the different but isomorphic copies of $V_{m''_1,m''_2}$ appearing in the sum. 

\paragraph{Resolution.} To resolve this ambiguity, the diagonal centraliser $Z_2(sl_3)$ may be considered. Indeed no operator coming from the diagonal embedding of $U(sl_3)$ can remove the ambiguity, so we need something from a larger algebra acting on the tensor product $V_{m_1,m_2}\otimes V_{m'_1,m'_2}$. Obviously, the tensor product $U(sl_3)\otimes U(sl_3)$ is the natural environment here. Moreover, to be compatible with the set of labels we already have, we would like to add an operator commuting with the diagonal embedding of $U(sl_3)$. This is where $Z_2(sl_3)$ makes its appearance. We can see $Z_2(sl_3)$ as the algebra containing all the possible additional labelling operators.

In the preceding section, two generators $X$ and $Y$ of $Z_2(sl_3)$ have been identified, so we shall try to use either of them as the operator providing the missing label. The action of $X$ was calculated in \cite{PST}. It follows from this calculation that $X$ can indeed distinguish between the various copies of $V_{m''_1,m''_2}$ and 
thus its eigenvalues can serve as a missing label. The second operator $Y$, which allows (through various combinations of $X$ and $Y$) to select among different missing label operators, 
was not considered so far. Moreover, the matrix obtained in \cite{PST} for $X$ did not display symmetries apart from the obvious ones. 

In this section we will provide the action of $X$ and $Y$ in such a way that symmetries will appear. By a symmetry of the missing label, we mean a transformation on the parameters $\tm$ such that the matrices for $X$ and $Y$ are transformed 
to equivalent matrices (up to a $\pm$ sign for $X$). Thus the eigenvalues will be preserved (up to a sign for $X$) under these transformations.

\subsection{Formulas for the missing label operators}\label{sec-XY}

We record the multiplicity of $V_{m''_1,m''_2}$ in the tensor product $V_{m_1,m_2}\otimes V_{m'_1,m'_2}$ with a multiplicity space $M_{\tm}$:
\[
V_{m_1,m_2}\otimes V_{m'_1,m'_2}=\dots\oplus\ \ V_{m''_1,m''_2}\otimes M_{\tm}\ \ \oplus\dots\ ,
\]
where the dimension of $M_{\tm}$ is the multiplicity $d_{\tm}$, and the space $M_{\tm}$ can be identified with the subspace of 
$V_{m_1,m_2}\otimes V_{m'_1,m'_2}$ consisting of the highest-weight vectors of weight $(m''_1,m''_2)$. In the spirit of the general ideas explained in Chapter \ref{chap-prel}, we see the space $M_{\tm}$ as a representation of the centraliser of $sl_3$, so here $Z_2(sl_3)$. In fact, as explained in the previous section, $M_{\tm}$ is a representation of the specialization $Z_2(sl_3)^{spec}$:
\[\begin{array}{rcl}Z_2(sl_3)^{spec} & \to & \text{End}(M_{\tm})\\[0.5em]
    X,Y & \mapsto & X_{\tm},Y_{\tm}\,.
    \end{array}\]
We will give the matrices of $X_{\tm}$ and $Y_{\tm}$ in a basis of $M_{\tm}$ which is carefully selected according to symmetry considerations. We will hide all details and just give the resulting matrices. 

We define
\[(\xi_1\,,\ \xi_2\,,\ \xi_3)=\left\{\begin{array}{l}
(\ell\,,\ m_2\,,\ m'_1+\ell-n)\ \ \ \text{if $\ell\leq n$\,,}\\[0.5em]
(n\,,\ m'_1\,,\ m_2+n-\ell)\ \ \ \text{if $n\leq \ell$\,,}
\end{array}\right.\ \ \ \ \text{and}\ \ \ \ (\xi_4\,,\ \xi_5\,,\ \xi_6)=(\ell-m'_2\,,\ n-m_1\,,\ 0)\,,\]
as well as the following linear combinations of these parameters
\[\Lambda=\sum_{i=1}^3\xi_i-\sum_{i=4}^6\xi_i+2|\ell-n|\ \ \ \ \text{and}\ \ \ \lambda_{\pm}=\pm\frac{1}{2}\Lambda+\frac{1}{6}\sum_{i=1}^6\xi_i\ .\]
With these notations, the size of the matrices (\emph{i.e.} the multiplicity) $d_{\tm}$ can be rewritten as follows 
\[d_{\tm}=\xi_a-\xi_b\,,\ \ \ \ \ \quad\text{where\ \ \ $\xi_a=\text{min}\{\xi_1,\xi_2,\xi_3\}\ $ and $\ \xi_b=\text{max}\{\xi_4,\xi_5,\xi_6\}$\ .}\]
The action of $X$ and $Y$ are given by the following tridiagonal matrices:
\begin{equation}\label{eq:XY}
 X_{\tm}=\left(\begin{array}{ccccc}
    a_{11} & a_{12}  &0  & \dots  & 0  \\
    a_{21} & a_{22}  & \ddots  & \ddots & \vdots  \\
     0 & \ddots  & \ddots & \ddots & 0 \\
     \vdots & \ddots & \ddots & \ddots & a_{d-1,d}\\
     0 & \dots & 0 &    a_{d,d-1}  &  a_{d,d}
 \end{array}\right)\ \ \ \ \text{and}\ \ \ \ Y_{\tm}=\left(\begin{array}{ccccc}
    b_{11} & b_{12}  &0  & \dots  & 0  \\
    b_{21} & b_{22}  & \ddots  & \ddots & \vdots  \\
     0 & \ddots  & \ddots & \ddots & 0 \\
     \vdots & \ddots & \ddots & \ddots & b_{d-1,d}\\
     0 & \dots & 0 &    b_{d,d-1}  &  b_{d,d}
 \end{array}\right)
\end{equation}
where $d=d_{\tm}$ and the off-diagonal coefficients are:
\[\begin{array}{rcl}a_{j,j+1}=(j+\xi_b-\xi_1)(j+\xi_b-\xi_2)(j+\xi_b-\xi_3) \quad & \text{and} & \quad b_{j,j+1}=a_{j,j+1}(j+\xi_b-\lambda_-)\ ,\\[0.5em]
a_{j+1,j}=(j+\xi_b-\xi_4)(j+\xi_b-\xi_5)(j+\xi_b-\xi_6) \quad & \text{and} & \quad b_{j+1,j}=a_{j+1,j}(j+\xi_b-\lambda_+)\ .
\end{array}\]
To write the diagonal coefficients, we set:
\[E_{j,k}=\sum_{i=1}^6x_i(j)^k\,,\ \ \ \ \text{where $x_i(j)=(\xi_i-j-\xi_b+\frac{1}{2})$\ for $i\in\{1,\dots,6\}$\ .}\]
Then we have:
\[a_{jj}=-\frac{1}{108}\Bigl(\frac{7}{2}E_{j,1}^3-18 E_{j,1}E_{j,2}+18 E_{j,3}\Bigr)-\frac{1}{24}E_{j,1}\Bigl(\Lambda^2+2\Bigr)\ ,\]
\[b_{jj}=\frac{1}{288}\Bigl(\frac{5}{2}E_{j,1}^4+32 E_{j,1}E_{j,3}+6\bigl(E_{j,2}^2-3E_{j,1}^2E_{j,2}-4E_{j,4})+6E_{j,2}(\Lambda^2+2)-3E_{j,1}^2(\Lambda^2-2)-\frac{3}{2}\Lambda^4+6\Lambda^2-36\Bigr)\ .\]
The off-diagonal coefficients of $X_{\tm}$ and $Y_{\tm}$ take a very simple form whereas their diagonal coefficients are more involved (though their chosen form allows to extract easily some symmetries).

\vskip .2cm
This explicit form allows to show that the eigenvalues of $X_{\tm}$ are all distinct, and as well for $Y_{\tm}$ (an easy general argument for tridiagonal matrices with non-zero off-diagonal coefficients).
Therefore both sets of eigenvalues can serve as labels for the different copies of $V_{m''_1,m''_2}$ appearing in the tensor product $V_{m_1,m_2}\otimes V_{m'_1,m'_2}$.
In this sense each of them provides a solution of the ``missing label'' problem. 

We note that there does not seem to exist exact formulas, in the general case, for the eigenvalues of $X_{\tm}$ or $Y_{\tm}$. 
In fact, diagonalizing $X_{\tm}$ in some examples immediately puts $Y_{\tm}$ in a complicated form (and vice versa). So it seems that the above simultaneous tridiagonal form for $X_{\tm}$ or $Y_{\tm}$ is the simplest presentation that can be achieved, and one may find it quite remarkable that such a simple form exists.

\begin{rema}
The matrices $X_{\tm}$ and $Y_{\tm}$ provide in particular a representation of the diagonal centraliser $Z_{2}(sl_3)$. Remarkably, there exists a geometric way to obtain a 
set of infinite-dimensional representations of this algebra by using the polytope associated to the root system $E_6$. Finite-dimensional representations can be extracted from the infinite-dimensional ones and the representations above can be identified among this set of representations \cite{CPV3}. This looks intriguingly similar to the usual Verma module construction for Lie algebra. 
\end{rema}

\subsection{Symmetry}

Given two sets of parameters $\tm$ and $\tm'$, it is of interest to determine when the two pairs of matrices $(X_{\tm},Y_{\tm})$ and $(X_{\tm'},Y_{\tm'})$ are simultaneously conjugated by the same change of basis. If this is so the eigenvalues of $X$ and of $Y$ (the possible missing labels) will be invariant when $\tm$ is replaced by $\tm'$.
Since $X$ and $Y$ generate, up to the central elements, the whole centraliser, any missing label operator taken from the centraliser $Z_2(sl_3)$ will have the same eigenvalues in the representation associated to $\tm$ and in the representation associated to $\tm'$. 

If $(X_{\tm},Y_{\tm})$ and $(X_{\tm'},Y_{\tm'})$ are conjugated, we call the corresponding transformation of the parameters $\tm$ into $\tm'$ a ``symmetry'' of the missing label operators (in fact we will consider the slightly more general situation where $X_{\tm}$ and $X_{\tm'}$ are conjugated up to a sign).

We have natural candidates for (linear) symmetry operations on the parameters $\tm$. These are the transformations from the extended Weyl group $\pm W(E_6)$ described in the preceding section. From this rather large group (of order $103680$), we extract the subgroup of symmetries of the missing label operators (we skip the explicit description as a subgroup of $\pm W(E_6)$. More details can be found in \cite{CPV3}).

\begin{framedtheo}[\cite{CPV3}]\label{theo-sym}
On the parameters in $\tm$, the symmetries of the missing label are described as transformations of the following arrangement:
\begin{equation}
 \label{eq:thm1}
\left[\begin{array}{ccc|ccc}
m_1 & m'_1 & m''_2                             \ \ \ & \ \ \ m_1+\ell-n & m'_1+\ell-n & m''_2+\ell-n  \\[0.5em]
m'_1+m'_2-\ell & m_1+m_2-\ell & n      \ \ \  & \ \ \ m'_1+m'_2-n & m_1+m_2-n & \ell \\[0.5em]
m_2+n-\ell & m'_2+n-\ell & m''_1+n-\ell \ \ \ & \ \ \ m_2 & m'_2 & m''_1
\end{array}\right]
\end{equation}
which are combinations of the following operations:\\
$\bullet$ simultaneously on the left and right $3\times 3$ squares: permutations of the lines, permutations of the columns and transposition of matrices;\\
$\bullet$ exchange of the left and right $3\times 3$ squares.
\end{framedtheo}
\noindent To be precise, the transformations in the theorem transform $X_{\tm}$ to a matrix equivalent to $\pm X_{\tm}$, the sign being determined by the number of transpositions of lines and of columns, and they transform $Y_{\tm}$  to a matrix equivalent to $Y_{\tm}$. Note that the two $3\times 3$ squares above are magic squares, in the sense that all the sums over lines, columns, diagonals or antidiagonals give the same result. 

\paragraph{Some remarks.} $\bullet$
The theorem has a very practical application. Given a set of parameters, 
the symmetry can be used to obtain a new set of parameters such that $\ell$ is the minimum of the 18 numbers in the arrangement. 
In this case, the explicit formulas for $X_{\tm}$ and $Y_{\tm}$ are easier to apply, since the parameters $\xi$'s are always the same, and the parameter $\xi_b$ is $0$.

\vskip .1cm
$\bullet$ There are some ``obvious'' symmetries of the missing label operator, obtained by permuting the representations and taking the duals. They form altogether a group of order 12. The symmetry group above is larger than that (of order 144) so we actually have interesting additional ``hidden'' symmetries. The structure of the group is in fact a semidirect product: 
\[\Bigl((S_3\times S_3)\rtimes \mathbb{Z}/2\mathbb{Z}\Bigr)\times \mathbb{Z}/2\mathbb{Z}\ ,\]
where we can see, successively, the permutations of the lines, the permutations of the columns, the transposition and the exchange of the two squares. This explains the order $144=6\times 6\times 2\times 2$ of the group.

\vskip .1cm
$\bullet$ The exchange of the left and right $3\times 3$ squares commutes with all other symmetries, so the whole symmetry group is the direct product of a group of order 72 
with $\bZ/2\bZ$. Rather surprisingly, the subgroup of order 72 acting simultaneously on both $3\times 3$ squares has the same structure as the group describing the symmetries of the Clebsch--Gordan coefficients, 
or $3j$-symbols, for $sl_2$.

\vskip .1cm
$\bullet$ The multiplicity $d_{\tm}$ is the minimum (if non-zero) of the $18$ numbers appearing in the arrangement above. Therefore, the condition for $\tm$ to correspond to a non-zero multiplicity in the Clebsch--Gordan series is that every number in this array be a positive integer. Of course a symmetry of the missing label preserves 
in particular the size of the matrices $X_{\tm}$ and $Y_{\tm}$, and thus is in particular a symmetry of the number $d_{\tm}$. Indeed it is clear that the transformations of the theorem are symmetries of 
the Littlewood--Richardson coefficient (they permute the 18 elements).

\begin{figure}[htb]
\begin{center}
\begin{tikzpicture}[scale=0.55]
\draw[fill] (1,1) circle [radius=0.08];
\node [above] at (1,1) {$m_2+n-\ell$};
\draw[fill] (1,5) circle [radius=0.08];
\node [above] at (1,5) {$m'_1$};
\draw[fill] (1,9) circle [radius=0.08];
\node [above] at (1,9) {$n$};
\draw[fill] (3,11) circle [radius=0.08];
\node [above] at (3,11) {$m'_1+m'_2-n$};
\draw[fill] (3,7) circle [radius=0.08];
\node [above] at (3,7) {$m_2+m'_2-\ell$};
\draw[fill] (3,3) circle [radius=0.08];
\node [above] at (3,3) {$m'_2$};

\draw[fill] (5,1) circle [radius=0.08];
\node [above] at (5,1) {$m_1+m_2-\ell$};
\draw (5,5) circle [radius=0.08];
\node [above] at (5,5) {$m_1+m'_1-n$};
\draw[fill] (5,9) circle [radius=0.08];
\node [above] at (5,9) {$m_1$};
\draw[fill] (7,11) circle [radius=0.08];
\node [above] at (7,11) {$m'_1-n+\ell$};
\draw[fill] (7,7) circle [radius=0.08];
\node [above] at (7,7) {$m_2$};
\draw[fill] (7,3) circle [radius=0.08];
\node [above] at (7,3) {$\ell$};

\draw[fill] (9,1) circle [radius=0.08];
\node [above] at (9,1) {$m''_2$};
\draw[fill] (9,5) circle [radius=0.08];
\node [above] at (9,5) {$m'_1+m'_2-\ell$};
\draw[fill] (9,9) circle [radius=0.08];
\node [above] at (9,9) {$m'_2+n-\ell$};
\draw[fill] (11,11) circle [radius=0.08];
\node [above] at (11,11) {$m''_1$};
\draw[fill] (11,7) circle [radius=0.08];
\node [above] at (11,7) {$m_1+m_2-n$};
\draw[fill] (11,3) circle [radius=0.08];
\node [above] at (11,3) {$m_1-n+\ell$};
\draw[thin] (0,0)--(0,12)--(12,12)--(12,0)--(0,0);

\foreach \y in {0,4,8}
\foreach \x in {0,4,8} 
{\draw [fill, opacity=0.2] (\x,\y)--(\x+1,\y+1)--(\x,\y+2)--(\x,\y);
\draw [fill, opacity=0.2] (\x+1,\y+1)--(\x+2,\y+2)--(\x+2,\y)--(\x+1,\y+1);
\draw [fill, opacity=0.2] (\x+2,\y+2)--(\x+3,\y+3)--(\x+2,\y+4)--(\x+2,\y+2);
\draw [fill, opacity=0.2] (\x+3,\y+3)--(\x+4,\y+4)--(\x+4,\y+2)--(\x+3,\y+3);
}

\end{tikzpicture}

\caption{Symmetries of the missing labels seen as affine transformations on a torus. \label{fig}}
\end{center}
\end{figure}
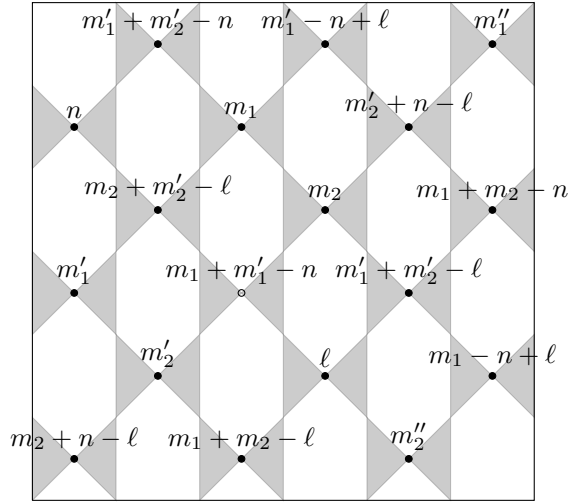

\vskip .1cm
$\bullet$ Another way to describe this symmetry group consists in drawing the $18$ elements of the array above on a torus as displayed in Figure \ref{fig}. The opposite sides of the square are identified.
The $144$ transformations of Theorem \ref{theo-sym} correspond to the affine transformation leaving stable the set of dots of this picture. 
Given one of the 18 dots, the transformations fixing this dot are the symmetries of a square of which this dot is the center (a dihedral group of order 8). 
The full symmetry is the semi-direct product of one of these dihedral groups with the subgroup of translations leaving the set of dots invariants. These translations are generated, for example, by two translations of order 6: say, $t_1$ sending $m_1$ to $m_2$, and $t_2$ sending $m_1$ to $m_2+m'_2-\ell$. These two translations satisfy $t_1^3=t_2^3$, so that the translation subgroup is of order $18$. This checks with $144=8\times 18$.

Moreover, if a symmetry maps gray zones to gray zones then it corresponds to a transformation $\tm\mapsto \tm'$ such that $X_{\tm}$ and $X_{\tm'}$ are conjugated. Otherwise, it corresponds to a case where $X_{\tm}$ and $-X_{\tm'}$ are conjugated.

\section{Conclusion and outlook}\label{sec:concl}

According to the general philosophy of this thesis which, roughly, claims that we should look at centralisers in order to discover interesting algebras, we have been looking at some diagonal centralisers for the simplest Lie algebras. The case of $Z_n(sl_2)$ is now well-understood and a first situation outside of $sl_2$, namely $Z_2(sl_3)$, was also studied in details. In both cases, the algebras that appear are expected to be useful in various contexts in mathematics and mathematical physics. One of these expected applications to the missing label of $sl_3$ was described in detail in the last section.

We note that the relations of the generalisation $R(n)$ of the Racah algebra appeared in \cite{DBG+} in the study of a superintegrable system in $n$ dimensions. They also have been recently observed in physical models \cite{LMZ}
and the study of their relation to multivariate Racah polynomials has been
initiated \cite{DBVV} but is still under active investigation.

An important generalisation concerns the quantum group $U_q(sl_2)$ for which the $3$-fold tensor product leads to the Askey--Wilson algebra.  Many attempts \cite{DBDC+} to generalise this result to $n$-fold tensor products have yielded relations of this centraliser but did not give all the defining relations. To make connections with knot theory, we emphasise that the diagonal centraliser $U_q(sl_2)$ for the $3$-fold tensor product was shown to be isomorphic to the Kauffman skein bracket algebra of the $4$-punctured sphere (see \cite{CFGPRV} for example). The case of the $n$-fold tensor product was recently shown to correspond exactly to the $n+1$-punctured sphere \cite{CL}. A $q$-generalisation of $R(n)$ and of its special quotient leading to a presentation of the centraliser is work in progress. 

The algebra appearing in this chapter, for $Z_2(sl_3)$, seems also to be relevant in the context of superintegrable systems \cite{CDO+}. More on the mathematical side, it turns out that this algebra belongs to the class of Calabi--Yau algebra \cite{EG,Gin} and making connections with the results of \cite{ELOR}, it is also a particular case of spherical symplectic reflection algebra (and the type $E_6$ appears in this setting also). These connections remain quite mysterious. Furthermore, a $q$-deformation of our description of $Z_2(sl_3)$ corresponding to $U_q(sl_3)$ is not known.
 
Classical results about invariant polynomials have been shown to be closely related to the description of diagonal centralisers. We have now more or less exhausted the situations where the commutative algebra of invariant polynomials is explicitly known. During our studies, we have encountered the interesting phenomenon that in the non-commutative setting, less generators are needed and more relations are implied from the commutation relations. This may well be an example where the non-commutative world can shed some light on its commutative analogue. In any case, a study of the non-commutative algebra generated by the polarized traces will be necessary for the diagonal centralisers and it could be that it will turn up to be useful also for the study of the commutative algebra of invariant polynomials.

Last but not least, both this chapter and the preceding one should be combined in some sense. When representing $Z_n(sl_2)$ on a tensor product of representations, we recover the centralisers studied in the preceding chapter. This means that there are connections, unknown at this time, between the fused Hecke algebras and the Racah/Askey--Wilson algebra and their generalisations. For a 3-fold tensor product, we have a conjecture stating that the centraliser, for $sl_2$, is a quotient of the Racah algebra $R(3)$ \cite{CPV1} (see also \cite{CVZ} for the $q$-deformed case). The generalisation of these ideas to the case of the $n$-fold tensor product is desirable. 

We can also ask the following interesting questions from the point of view of the braid group:  If one adds the Askey--Wilson relations to the braid group, is it sufficient to obtain finite-dimensional algebras? If yes, do we obtain this way a presentation of the centralisers of tensor products? In case of a 3-fold tensor product, it seems that we have some positive answers to these questions. This is work in progress.

\section*{Appendix \thechapter.A $\ $The Hilbert--Poincaré series of $Z_n(sl_2)$ and $Z_2(sl_3)$}
\addcontentsline{toc}{section}{Appendix \thechapter.A $\ $The Hilbert--Poincaré series of $Z_n(sl_2)$ and $Z_2(sl_3)$}

For more information on Hilbert--Poincaré series of graded algebras, we refer to
\cite{Sta1}. Denote by $Z_n(sl_2)_{\leq k}$ the subspace of elements of degree less or equal to $k$ (the degree is in the generators of $U(sl_2)^{\otimes n}$). The Hilbert--Poincar\'e series of $Z_n(sl_2)$
records the dimensions of these increasing sequence of subspaces:
\begin{align}\label{}
F_n(t)=\sum_{k\geq0}\dim\bigl(Z_n(sl_2)_{\leq k}/Z_n(sl_2)_{<k}\bigr)t^k,
\end{align}
From the point of view of classical invariant theory, the coefficient in front of $t^k$ is the dimension of homogeneous invariant polynomials of degree $k$.

\paragraph{An explicit formula.} Several formulas, using various approaches, have been obtained for the Hilbert--Poincar\'e
series $F_n(t)$ (see references in \cite{Dre,For,Tera}). The
following formula seems to be new.
\begin{prop}[\cite{CGPV}] Let $n\geq2$ and set $r=n-2$. The Hilbert--Poincar\'e
series of $Z_n(sl_2)$ is:
\begin{align*}
 F_n(t)=\frac{P_r(t)}{(1-t^2)^{3(r+1)}} ,
\end{align*}
where the numerator is given by:
\begin{align*}\label{}
 P_r(t)=(1+t)^r\sum_{k=0}^{2r}(-1)^ka_kt^k,\qquad\text{where}\
 \left\{\begin{array}{l}
  a_{2k}=\binom{r}{k}^2,\\[0.5em]
  a_{2k+1}=\binom{r}{k}\binom{r}{k+1}.
 \end{array}\right.
\end{align*}
\end{prop}
We make a few comments about this formula with a flavour of classical algebraic geometry. First, the exponent $3(r+1)$ appearing in the denominator of $F_n(t)$ is the Krull, or
Gelfand--Kirillov, dimension of the algebra of invariant polynomial functions (see
\cite{For,Kir}). Using Hochster--Roberts theorem from general invariant theory \cite{HoRo}, there is a set of $3(r+1)$ algebraically independent elements (a system of parameters)
$\theta_1,\dots,\theta_{3(r+1)}$ such that the algebra is a free module of finite
dimension over the polynomial subalgebra $\mathbb{C}[\theta_1,\dots,\theta_{3(r+1)}]$. The
form $F_n(t)$ above with the positivity of the numerator $P_r(t)$ (see below) strongly suggests that the system of parameters could consist of $3(r+1)$ elements of degrees $2$. 

If it were to be the case, it would be very interesting since then $P_r(1)$ would be the dimension of the algebra over
$\mathbb{C}[\theta_1,\dots,\theta_{3(r+1)}]$. Moreover, the different monomials in $P_r(t)$ would indicate in which degrees the elements of a basis over $\mathbb{C}[\theta_1,\dots,\theta_{3(r+1)}]$ would have to be found.

Finally, the palindromic property of the numerator $P_r(t)$ in the formula above shows
directly that the Hilbert--Poincar\'e series satisfies the functional equation:
\begin{align*}\label{}
 F_n(t^{-1})=(-1)^{(n-1)}t^{3n}F_n(t).
\end{align*}
This recovers a known property which is related to a property, called being Gorenstein, for the algebra of
invariant polynomial functions (see \cite{Dre} and references therein).

\paragraph{Some combinatorics.} We have obtained an expression for the Hilbert--Poincar\'e series of $Z_n(sl_2)$ of the
form:
\begin{align*}\label{}
 F_n(t)=\frac{(1+t)^rQ_r(t)}{(1-t^2)^{3(r+1)}} ,\qquad
 \text{where $Q_r(t)=\sum_{k=0}^{2r}(-1)^ka_kt^k$},
\end{align*}
and the coefficients $a_k$ are given in the proposition above. It is perhaps not so
surprising that the coefficients of the various polynomials involved show some connections
with well-studied combinatorial objects of ``Catalan'' flavour.

\paragraph{The polynomial $Q_r(t)$.} The coefficient $a_k$ in the polynomial $Q_r(t)$
counts the number of symmetric Dyck paths of semi-length $2r+1$ with $k+1$ peaks (see
A088855 in \cite{OEIS}). Their expression with binomial coefficients corresponds to
choosing a certain number of peaks and troughs in the first $r$ steps of the paths.

In fact, the polynomial $Q_r(t)$ is a $t$-deformation of the well-known Catalan number,
that is, the value of $Q_r(t)$ at $t=1$ is the $r$-th Catalan number:
\begin{align*}\label{}
 Q_r(1)=c_r=\binom{2r}{r}-\binom{2r}{r+1}.
\end{align*}
It is not difficult to give a combinatorial proof that the alternating sum of the
$a_k$'s is equal to the Catalan number $c_r$ (the number of Dyck paths of length $2r$). Maybe it is also enjoyable to see it as follows. From its explicit expression, the polynomial $Q_r(t)$ is equal to the constant term of the following Laurent
polynomial in $x$:
\begin{align*}
 Q_r(t)=\bigl[(1-x^2)(1+tx^{2})^{r}(1+tx^{-2})^{r}\bigr]_0.
\end{align*}
For $t=1$, this combines nicely and gives:
\begin{align*}\label{}
\bigl[(1-x^2)(1+x^{2})^{r}(1+x^{-2})^{r}\bigr]_0=\bigl[(1-x^2)(x+x^{-1})^{2r}\bigr]_0
    =c_r\ .
\end{align*}
In this sense, the polynomial $Q_r(t)$ is a natural $t$-deformation of the $r$-th Catalan
number.

\paragraph{The numerator $P_r(t)$.} The numerator of the Hilbert--Poincar\'e series of
$Z_n(sl_2)$ is $P_r(t)=(1+t)^rQ_r(t)$. Its coefficients are
all positive and we will show this explicitly.

First, from what we have said above about $Q_r(t)$, it follows that $P_r(t)$ is a
$t$-deformation of the number $2^rc_r$, that is, its value at $t=1$ is $P_r(1)=2^rc_r$.
This number counts several classes of combinatorial objects (see A151374 in \cite{OEIS}). The $t$-deformation giving $P_r(t)$
can be expressed similarly as before as:
\begin{align*}
 P_r(t)=\bigl[(1-x^2)(1+t)^r(1+tx^{2})^{r}(1+tx^{-2})^{r}\bigr]_0.
\end{align*}
Now regrouping the terms with an $r$-th power gives the following expression:
\begin{align*}
 P_r(t) =  \bigl[(1-x^2)\bigl(1+t^3+(t+t^2)(1+x^2+x^{-2})\bigr)^{r}\bigr]_0 = \displaystyle\sum_{k=0}^rR_k\binom{r}{k}(1+t^3)^{r-k}(t+t^2)^k,
\end{align*}
where the positive integer $R_k$ is the Riordan number, one of the closest relative of the
Catalan number (see A005043 in
\cite{OEIS}). They are given by either one of the following equalities:
\begin{align*}
 R_n=\bigl[(1-x^2)(1+x^2+x^{-2})^n\bigr]_0=\sum_{i=0}^n(-1)^{n-i}\binom{n}{i}c_i.
\end{align*}
The formula above for $P_r(t)$ has the advantage to show explicitly that it has
positive coefficients. So $P_r(t)$ is a $t$-deformation with positive coefficients of
$2^rc_r$ and therefore should be given by an interesting statistics on a certain set of
$2^rc_r$ objects.

\paragraph{The diagonal centraliser $Z_2(sl_3)$.} Here is the Hilbert--Poincar\'e series of $Z_2(sl_3)$:
\[\frac{1+t^{6}}{(1-t^2)^3(1-t^3)^4(1-t^4)}\ .\]
We see in the denominators the degrees of the generators $Z_2(sl_3)$. The numerators shows that there is an additional element of degree 6 (we denoted it $Z$) whose square is expressed in terms of the generators. We do not know a nice uniform formula for $Z_n(sl_3)$.

\setcounter{equation}{0}
\DeactivateToc
\chapter{\huge{Yokonuma--Hecke algebras and link invariants}}\label{chap-yok}
\ActivateToc
\addcontentsline{toc}{chapter}{\large{Chapter \thechapter. \hspace{0.2cm}Yokonuma--Hecke algebras and link invariants\vspace{0.3cm}}}

\setcounter{minitocdepth}{1}
\minitoc

\section*{Introduction}
\addcontentsline{toc}{section}{Introduction}

In this chapter, we will consider Yokonuma--Hecke algebras. They belong to the class of algebras studied in this thesis, since their are centraliser algebras, although not of tensor products of representations of quantum groups. As their names indicate, they have been introduced by Yokonuma \cite{Yok}, and they belong to the world of Hecke algebras, in the sense of Chapter \ref{chap-prel}, Section \ref{sec-Hecke}. Indeed they are centraliser algebras of some permutation representation of a group.

In fact, from the point of view of centralisers of permutation representations, they are one of the most natural generalisation of the usual Iwahori--Hecke algebras. Recall that the usual Iwahori--Hecke algebra is obtained as centraliser algebra of  the permutation representation of a finite group $G$ of Lie type with respect to a Borel subgroup. To define Yokonuma--Hecke algebras, one simply needs to replace the Borel subgroup by a maximal unipotent subgroup of $G$. For $G=\text{GL}_n (\mathbb{F}_{q})$, the maximal unipotent subgroup consists of upper-triangular matrices with only $1$'s on the diagonal.

The Yokonuma--Hecke algebras have been revived more recently, when their natural description has been transformed into a  simple presentation with generators and relations \cite{Juy1}. From this presentation, one can observe that the Yokonuma--Hecke algebra of $G=\text{GL}_n (\mathbb{F}_{q})$ (sometimes called the Yokonuma--Hecke algebra of type $A$) is a deformation of the group algebra of the complex reflection group of type $G(d,1,n)$, where $d=q-1$. In particular, it contains a homomorphic image of the braid group on $n$ strands, and thus the whole business of constructing Markov traces and producing invariants for links could start \cite{Juy2}.

In this chapter, as usual, we will forget about the origin as a centraliser algebra, and in particular that $q$ is the power of a prime number and $d=q-1$. What we will call the Yokonuma--Hecke algebra, and denote $Y_{d,n}$, is the generic algebra defined for any positive integer $d$ and depending on an indeterminate $q$ (in fact we will even add another indeterminate $u$ since it is more convenient with respect to Markov traces and link invariants).

\paragraph{The story of link invariants associated to Yokonuma--Hecke algebras.} In \cite{Juy2}, Juyumaya introduced on $Y_{d,n}$ an analogue of the Ocneanu trace of the Iwahori--Hecke algebra of type $A$. This trace was subsequently used by Juyumaya and Lambropoulou to produce invariants for knots and links \cite{JL3}. At this point, the question of whether these invariants were related (that is, equivalent, stronger, weaker or simply different) to the HOMFLYPT polynomial arises \cite{ChLa}. At first sight, they have no reason to be the same. However, the construction is so similar to the usual construction of the HOMFLYPT polynomial that connections were suspected. Somehow, it would be a nice surprise if such an extension of the usual construction were to produce a genuinely stronger invariant.

The first examples \cite{CJKL} seemed to indicate that the Juyumaya--Lambropoulou invariant was topologically equivalent to HOMFLYPT, but these first examples were only for knots (links with one connected component). For links, the situation was different and an example was soon found \cite{CJKL} where two different links were distinguished by the new invariant while they had the same HOMFLYPT polynomial. Thus the new invariant does contain more information than the HOMFLYPT polynomial. It remained to understand what was this new information because it was still suspected that this new invariant was somehow connected to the HOMFLYPT polynomial. The conclusion of this story was reached when a complete description of the Juyumaya--Lambropoulou invariant was obtained: it can be expressed in terms of the HOMFLYPT polynomials of sublinks and the linking numbers between different connected components (\cite[Appendix]{CJKL} and \cite{PW1}). This explains immediately that no difference was seen for knots, and this also explains why some links were distinguished: they were distinguished by a subtle combination of the linking numbers and the HOMFLYPT polynomials for their sublinks.

With the conclusion of the story, one may feel a little bit disappointed since in fact the new invariant was no more than a combination of other known invariants. However, the whole story was interesting and it is certainly a non-trivial matter to realise that this combination of known invariants can be obtained through Markov traces on the Yokonuma--Hecke algebras. The Yokonuma--Hecke algebras also opened several perspectives related to knots and links that will be briefly mentioned in the conclusion of this chapter.

\paragraph{The algebraic approach.} In this chapter, we will explain our algebraic approach to the story summarised above. In a nutshell, one would like to construct Markov traces on the Yokonuma--Hecke algebras $Y_{d,n}$ and somehow relate them to the usual Markov trace on the Hecke algebra. 

The program goes as follows. First of all, precise algebraic connections between $Y_{d,n}$ and the usual Hecke algebra should be understood. It turns out that one can prove an explicit isomorphism theorem between $Y_{d,n}$ and a direct sum of matrix algebras over tensor products of Hecke algebras. Complicated as it seems, this purely algebraic result is the key to the following developments.

Using the isomorphism theorem, the classification of Markov traces on $\{Y_{d,n}\}_{n\geq 1}$ is obtained allowing to produce a family of link invariants. Moreover, the Markov traces are explicitly expressed in terms of the usual Markov trace on the Hecke algebra. We note that the classification allows to be sure that we do not miss any invariant, and indeed, the Juyumaya--Lambropoulou invariant can be identified among the obtained set of invariants.

Using at full power the algebraic tools coming from the isomorphism, the family of invariants is organised and a basic set is identified, from which all the others can be obtained by linear combinations. Finally, the main formula giving these basic invariants in terms of HOMFLYPT polynomials and linking numbers is obtained thereby concluding the story.

\paragraph{References.} The publications relevant for this chapter, where more details can be found, are \cite{PW1,JP2,P3,JP1,mCP2,mCP1}.

\section{Isomorphism theorem for Yokonuma--Hecke algebras}

\subsection{Yokonuma--Hecke algebras}

\paragraph{Definition.} We define the Yokonuma--Hecke algebra $Y_{d,n}$, over the ring $\C[u^{\pm 1},q^{\pm 1}]$, as the algebra with generators:
$$g_1,g_2,\ldots,g_{n-1}, t_1,\ldots, t_n,$$
and relations:
\begin{equation}\label{rel-def-Y}\begin{array}{rclcl}
g_ig_j & = & g_jg_i && \mbox{for all $i,j=1,\ldots,n-1$ such that $\vert i-j\vert > 1$,}\\[0.1em]
g_ig_{i+1}g_i & = & g_{i+1}g_ig_{i+1} && \mbox{for  all $i=1,\ldots,n-2$,}\\[0.1em]
t_it_j & =  & t_jt_i &&  \mbox{for all $i,j=1,\ldots,n$,}\\[0.1em]
g_it_j & = & t_{s_i(j)}g_i && \mbox{for all $i=1,\ldots,n-1$ and $j=1,\ldots,n$,}\\[0.1em]
t_j^d   & =  &  1 && \mbox{for all $j=1,\ldots,n$,}\\[0.2em]
g_i^2  & = & u^2 + u(q-q^{-1})e_{i}  g_i && \mbox{for  all $i=1,\ldots,n-1$,}
\end{array}
\end{equation}
where, for all $i=1,\ldots,n-1$, $s_i$ denotes the transposition $(i,i+1)$ and
$$e_i :=\frac{1}{d}\sum\limits_{0\leq s\leq d-1}t_i^s t_{i+1}^{-s}\ .$$
The elements $e_i$ are idempotents and $g_ie_i=e_ig_i$. The generators $g_i$ are invertible,
$g_i^{-1} = u^{-2} g_i - u^{-1} (q-q^{-1}) e_i\,,$ so the subalgebra generated by $g_i$ is a quotient of the braid group algebra.

If $u$ and $q$ are specialized to $1$, we recover the defining relations of a group, the complex reflection group denoted $G(d,1,n)$. This group is the natural semidirect product $\bigl(\mathbb{Z}/d\mathbb{Z}\bigr)^n\ltimes S_n$. Thus the algebra $Y_{d,n}$ is a deformation of the group algebra of $G(d,1,n)$, and this deformation is flat in the following sense. 

If we set $g_w=g_{i_1}\ldots g_{i_r}$, for $w\in S_n$ and $s_{i_1} \ldots s_{i_r}$ a reduced expression for $w$, then we have a basis of $Y_{d,n}$ given by
 \[
 \{ t_1^{k_1} \ldots t_n^{k_n} \tg_w \ |\ w\in S_n,\ k_1,\ldots, k_n \in \mathbb{Z}/d\mathbb{Z}\}\ .
\]
The rank of  $Y_{d,n}$ is therefore equal to $d^n n!$, which is the order of $G(d,1,n)$.

For $d=1$, we have $t_i=e_i=1$ and $Y_{1,n}$ is simply the usual Hecke algebra.

\paragraph{Representations.} In the generic case ($u$ and $q$ are indeterminate), the irreducible representations of $Y_{d,n}$ are parametrised by $d$-tuples of partitions of total size $n$. They can be constructed explicitly in a seminormal form, analogue to the construction in Chapter \ref{chap-prel}, Appendix. We refer to \cite{mCP1}. The representation theory of $Y_{d,n}$ can also be understood from the isomorphism theorem presented below.

\paragraph{Affine Yokonuma--Hecke algebras.} We add generators $X_1,\dots,X_n$ to $Y_{d,n}$, which are required to be invertible, and with additional defining relations:
\begin{equation}\label{def-aff2}
\begin{array}{rclcl}
X_iX_j & =  & X_jX_i  && \mbox{for all $i,j=1,\ldots,n$,}\\[0.1em]
X_it_j & =  & t_jX_i  && \mbox{for all $i,j=1,\ldots,n$,}\\[0.1em]
g_iX_ig_i & = & X_{i+1} && \mbox{for all $i=1,\ldots,n-1$,}\\[0.1em]
g_iX_j & = & X_{j}g_i && \mbox{if $j\neq i,i+1$.}
\end{array}
\end{equation}
For $d=1$, this is the well-known affine Hecke algebra of type GL. The affine Yokonuma--Hecke algebra was originally introduced from the point of view of the representation theory of $Y_{d,n}$ \cite{mCP1} and subsequently studied in more details in \cite{mCP2} where it was used to construct invariants for links in the solid torus. The isomorphism theorem and the consequences for link invariants below generalise to the affine Yokonuma--Hecke algebras \cite{P4}.

From the point of view of centralisers of induced representations of finite or $p$-adic groups, a question naturally arises when looking at the following table:
\begin{center}
\begin{tabular}{|c|c|}
\hline & \\[-0.4em]
$\text{Ind}_H^G\mathbf{1}_H$ & centraliser algebra \\[0.5em]
\hline & \\[-0.2em]
$G=GL_n(\mathbb{F}_q)$ and $H=$Borel subgroup & Hecke algebra \\[0.5em]
$G=GL_n(\mathbb{Q}_p)$ and $H=$Iwahori subgroup & affine Hecke algebra \\[0.5em]
$G=GL_n(\mathbb{F}_q)$ and $H=$unipotent subgroup & Yokonuma--Hecke algebra \\[0.5em]
??? & affine Yokonuma--Hecke algebra\\
\hline
\end{tabular}
\end{center}
Of course, the answer should be that the affine Yokonuma--Hecke algebra corresponds to  $G=GL_n(\mathbb{Q}_p)$ and $H$ is the analogue of the unipotent subgroup. Such a subgroup is called pro-p-Iwahori subgroup and the resulting centraliser algebra, called the pro-p-Iwahori--Hecke algebra, does coincide with the affine Yokonuma--Hecke algebra \cite{CS}.

\subsection{Isomorphism theorem}

In this chapter, the quadratic equation for the generators of the Hecke algebra is taken to be $\si_i^2=u^2+u(q-q^{-1})\si_i$. For the statement below, the convention is that $H_0=H_1=\mathbb{C}[u^{\pm 1},q^{\pm 1}]$.

A composition $\mu$ of $n$ with $d$ parts is a $d$-tuple of non-negative integers $\mu=(\mu_1,\dots,\mu_d)$ summing to $n$. We denote $\mu\models_d n$. Note that some of the parts may be $0$.
 \begin{framedtheo}[\cite{JP1}]\label{theo-iso}  
Over $\C[u^{\pm1},q^{\pm 1}]$, we have an explicit isomorphism of algebras: 
\[Y_{d,n}\cong \bigoplus_{\mu \models_d n } \operatorname{Mat}_{m_{\mu}} (H_{\mu_1}\otimes \dots\otimes H_{\mu_d} )\ , \ \ \ \ \ \ \text{where $m_{\mu}=\displaystyle \frac{n!}{\mu_1!\mu_2!\dots\mu_d!}\ .$}\]
 \end{framedtheo}
This result is in fact a special case of a result for unipotent Hecke algebras \cite{Lus}, but one important feature of the statement above is that we have an explicit description of the isomorphism. It is crucial for what follows, we will use it to answer very concrete questions. Note the nice equality of dimensions:
\[d^n.n!=\sum_{\mu \models_d n }\bigl(\frac{n!}{\mu_1!\mu_2!\dots\mu_d!}\bigr)^2\mu_1!\mu_2!\dots\mu_d!=n!\sum_{\mu \models_d n }\frac{n!}{\mu_1!\mu_2!\dots\mu_d!}\ .\]

\section{Classification of Markov traces on Yokonuma--Hecke algebras}\label{sec-mark}

In this section, we give the complete classification of Markov traces on $Y_{d,n}$, which is obtained using the explicit isomorphism from the previous theorem. As in the Hecke algebra situation, from now on, we extend the ground ring to $\mathbb{C}[u^{\pm1},q^{\pm1},(q-q^{-1})^{-1}]$. We denote $\{\tau_n\}_{n\geq 1}$ the unique Markov trace on the family of Hecke algebras $\{H_n\}_{n\geq 1}$ normalized by $\tau_1(1)=1$ (see Chapter \ref{chap-prel}).

The rough idea behind the following theorem is as follows. We work with the direct sum of matrix algebras obtained in Theorem \ref{theo-iso}. On each of this matrix algebra, there is a natural trace to consider: first apply the usual trace of a matrix; then we are in a tensor product $H_{\mu_1}\otimes \dots\otimes H_{\mu_d}$ and we apply the Markov trace of the Hecke algebra on each factor: $\tau_{\mu_1}\otimes\dots\otimes \tau_{\mu_d}$.

The above procedure gives obviously a trace function on $Y_{d,n}$, but it does not give a Markov trace. However, the Markov traces on $Y_{d,n}$ are obtained by summing some of these trace functions.

First, for a composition $\mu\models_d n$, we define its support $\text{supp}(\mu)$ as the subset of $\{1,\dots,d\}$ for which the parts $\mu_i$ are non-zero. Formally, we set:
\[\text{supp}(\mu)=\{i\in\{1,\dots,d\}\ \text{such that}\ \mu_i\neq 0\}\ .\]
\begin{framedtheo}\label{theo-classif}
Let $S$ be a non-empty subset of $\{1,\dots,d\}$. The family of maps $\{\rho^S_n\}_{n\geq1}$ is a Markov trace on $\{Y_{d,n}\}_{n\geq1}$ where
\begin{equation}\label{rho-mu0}
\rho^S_{n}=\sum_{\substack{\mu\models_d n\\ \text{supp}(\mu)=S}}(\tau_{\mu_1}\otimes\dots\otimes \tau_{\mu_d})\circ\operatorname{Tr}_{\Mat_{m_{\mu}}}\,.
\end{equation}
Moreover, as $S$ varies, the Markov traces $\{\rho^S_n\}_{n\geq1}$ form a basis of the space of Markov traces on $\{Y_{d,n}\}_{n\geq1}$.
\end{framedtheo}
In the formula above, it is understood that in each summand indexed by $\mu\models_d n$, we act on the corresponding algebra $\Mat_{m_{\mu}}(H_{\mu_1}\otimes \dots\otimes H_{\mu_d})$.

\begin{exam}\label{exam-markov}
Here we show graphically the supports of the three Markov traces for $d=2$. Each pair below represents the matrix algebra associated to the corresponding composition.
$$\begin{array}{ccccccccccc}
(n=1)&&&&& \underline{\text{$(1,0)$}} && \mathbf{(0,1)} && \\[0.5em]
(n=2)&&&& \underline{\text{$(2,0)$}} && (1,1) && \mathbf{(0,2)} &\\[0.5em]
(n=3)&&& \underline{\text{$(3,0)$}} && (2,1) && (1,2) && \mathbf{(0,3)} \\[0.5em]
\vdots && \reflectbox{$\ddots$} && \reflectbox{$\ddots$} &&\vdots && \ddots&&\ddots
\end{array}$$
A first Markov trace is living on the underlined summands. Another one is living on the boldfaced summands. These two are simply the usual Markov traces on the two chains of Hecke algebras embedded in $\{Y_{d,n}\}_{n\geq 1}$. They are already known from the usual Hecke algebra (the case $d=1$). The third one (on the remaining summand) is the only new Markov trace appearing for $d=2$. We will see below that there is in fact only one new Markov trace appearing for each $d$.
\end{exam}

\section{Link invariants from Yokonuma--Hecke algebras}\label{sec-inv}

Now that we have obtained a complete description of the Markov traces for $Y_{d,n}$, we will use them to deduce invariants for knots and links.

\paragraph{Definition of the invariants.} Usually, to define invariants, we would simply combine the Markov trace with the natural map from $B_n$ to $Y_{d,n}$ ($\si_i\mapsto g_i$). Nevertheless, here we will use a slightly more general procedure.

Let $\gamma$ be another indeterminate, and extend again our ring of coefficients to $\C[u^{\pm1},q^{\pm1},(q-q^{-1})^{-1},\gamma^{\pm1}]$ to include this new indeterminate. We define the following map, for which it is not hard to show that it gives a morphism from $B_n$ to $Y_{d,n}$:
\[
\pi_{\gamma,n}\ :\ \ \sigma_i\mapsto \bigl(\gamma+(1-\gamma)e_i\bigr)g_i\ \ (i=1,\dots,n-1)\,.
\]
For $\gamma=1$, this is simply the natural map $\si_i\mapsto g_i$. It turns out that we can use the more general map $\pi_{\gamma,n}$ to construct invariants.
\begin{defi}
Let $d>0$ and $S\subset\{1,\dots,d\}$. For a link $L$, let $\beta_L$ be a braid in $B_n$ closing to it, and define:
\[T_d^S(L)=\rho^S_n\bigl(\pi_{\gamma,n}(\beta_L)\bigr)\ ,\]
where $\{\rho^S_n\}_{n\geq1}$ is the Markov trace on $Y_{d,n}$
\end{defi}
\begin{framedtheo}[\cite{JP1}]
$T_d^S(L)$ is an invariant of the link $L$ for any $d>0$ and any $S\subset\{1,\dots,d\}$, 
\end{framedtheo}
Of course, most of the work for this result was done when proving that $\{\rho^S_n\}_{n\geq1}$ is a Markov trace on $Y_{d,n}$. The only remaining thing was to handle the fact that we used the more general map $\pi_{\gamma,n}$, and this is not difficult. We emphasize that it is important to consider the parameter $\gamma$ since it interpolates between various conventions used in the literature, and thus the invariants above contain all invariants previously considered, and in particular the Juyumaya--Lambropoulou invariant (see \cite{P3}). The parameter $\gamma$ also clarifies the final formula for the invariants below.

\paragraph{Description of the invariants.} Now that we have a rather large set of invariants, it is time to organise them and to see if they can be compared with the invariant coming from the Hecke algebra, namely, the HOMFLYPT polynomial. The key tool that we have is the explicit isomorphism and all results below rely on this algebraic theorem. Also, we have been selecting in fact a basis of the space of Markov traces on $\{Y_{d,n}\}_{n\geq 1}$, and the choice we made is justified by the relatively clean formulas that we present below.

\vskip .2cm
First of all, there is a big reduction on the number of invariants to actually consider. This relies on a comparison of the various invariants for different $d$. As already hinted at in Example \ref{exam-markov}, we basically obtain that there is only one new invariant for each new value of $d$.

Indeed, let $S\subset\{1,\dots,d\}$ and denote by $d'$ the size of $S$: $d'=|S|$. What we can prove is the following \cite{P3}:
\[T_d^S(L)=T_{d'}^{\{1,\dots,d'\}}(L)\ \ \ \ \ \text{for any link $L$.}\]
Thus it is enough to consider, for each $d>0$, the full subset $S=\{1,\dots,d\}$.

\vskip .2cm
Now we can state our final result on the invariants, describing them entirely in terms of linking numbers and HOMFLYPT polynomials of sublinks. Below $P(L)$ is the HOMFLYPT polynomial of a link. A similar result has been proved independently by Lickorish using diagrammatic techniques, for some specialisations of the parameter $\gamma$ (see \cite[Appendix]{CJKL}).
\begin{framedtheo}[\cite{PW1}]\label{theo-final-YH}
Let $L$ be a link and $d>0$. We have
\begin{equation*}
T_d(L)=d!\sum_{\{L_1,\dots,L_d\}}(u\gamma)^{\ell(L_1,\dots,L_d)}P(L_1)\dots P(L_d),
\end{equation*}
where the sum is over the set of all complete families $\{L_1,\dots,L_d\}$ of $d$ distinct non-empty sublinks of $L$.
\end{framedtheo}
The total linking number $\ell(L_1,\dots,L_d)$ for a family of links is:
\[\ell(L_1,\dots,L_d):=\sum_{1\leq i<j\leq d}\text{lk}(L_i,L_j) \ ,\]
where $\text{lk}(L_i,L_j)$ is the linking number between the links $L_i$ and $L_j$. For a link $L$, a sublink is a union of some of its disconnected components, and a complete family $\{L_1,\dots,L_d\}$ of sublinks for $L$ means that $L$ is the union of 
$L_1,\dots,L_d$.

Note that the subfamily in the sum must consist of non-empty and distinct sublinks of $L$. It means that if $d$ is strictly larger than the number of connected components of $L$, then the invariant $T_d(L)$ is $0$. For a link with $N$ connected components, it is thus enough to consider the finite set of invariants:
\[T_1(L),\dots,T_N(L)\ .\]
In particular for a knot $K$ (that is, $N=1$), all invariants from the Yokonuma--Hecke algebra reduce to $T_1(K)$, which is simply equal to the HOMFLYPT polynomial.

\section{Conclusion and outlook}

Looking at the final formula in Theorem \ref{theo-final-YH}, one soon realises a posteriori that this could serve as a definition of the polynomial $T_d(L)$ for a link $L$. With this definition, it is obviously an invariant. What is less clear while looking at the formula above is that it can be provided by a Markov trace on the Yokonuma--Hecke algebra. In fact using the Markov trace on $Y_{d,n}$ can be seen as an alternative way to calculate this invariant. 

One can imagine defining other invariants for links by similar formulas using other invariants than the HOMFLYPT polynomial, and then ask for the algebras producing them via a Markov trace. As explained in this chapter, the decomposition of the invariants over subfamilies of links is the topological shadow of the decomposition of the Yokonuma--Hecke algebra appearing in the isomorphism theorem. So for these other algebras corresponding to other invariants, one would expect also an analogue of the isomorphism theorem. This has been studied in full details for the Jones polynomial \cite{Chl} where the algebra appearing is called a framization of the Temperley--Lieb algebra. One could consider also for instance the $2$-variable Kauffman polynomial, which should lead to a framization of the BMW algebra. One can play this game with any invariant, such as the Reshetikhin--Turaev invariants coming from quantum groups.

Other algebras connected to the Yokonuma--Hecke algebras have been studied. This includes the braids and ties algebra \cite{AJ} which has been generalised for any complex reflection group \cite{Ma3,Ma4}. The cyclotomic and affine versions of the Yokonuma--Hecke algebra give rise to invariants of links in the solid torus and the whole story of this chapter repeats for these invariants \cite{P4,mCP2,mCP1}.

In fact, to complete the picture of the Yokonuma--Hecke algebra in the light of this thesis, the connection with centralisers of quantum groups is missing. One can expect that there is one, being roughly speaking the ``Schur--Weyl translation'' of the isomorphism theorem involving Hecke algebras. An action of the Yokonuma--Hecke algebra on tensor products of vector spaces is indeed available \cite{ERH}. And by the way, since we are here, we note that there is a solution of the Yang--Baxter equation inside the Yokonuma--Hecke algebra, given by
\[
g_i(u)=g_i+(q-q^{-1})\frac{e_i}{u-1}\ .
\]
The Yang--Baxter equation was one of the main characters of this thesis so this is a good way to conclude.

\chapter*{\huge{List of publications}}
\addcontentsline{toc}{chapter}{\large{List of publications\vspace{0.3cm}}}
\mtcaddchapter

\newcommand{\etalchar}[1]{$^{#1}$}


\begin{enumerate}

\itemsep-0.1em

\bibitem{P5} L. Poulain d'Andecy, \emph{Fusion for the braid group and the Yang--Baxter equation}, Winter Braids Lecture Notes, Vol. 7 (2020), Course no III, p. 1--49.

\bibitem{CPV3} N. Cramp\'e, L. Poulain d'Andecy, L. Vinet, \emph{The missing label of $su_3$ and its symmetry}, Commun. Math. Phys. (2023). https://doi.org/10.1007/s00220-022-04596-3

\bibitem{CGPV} N. Cramp\'e, J. Gaboriaud, L. Poulain d'Andecy, L. Vinet, \emph{Racah algebras, the centralizer $Z_n(sl_2)$ and its Hilbert--Poincaré series}, Annales de l'Institut Henri Poincaré 23 (2022), no. 7, 2657--2682.

\bibitem{CFGPRV}  N. Cramp\'e, L. Frappat, J. Gaboriaud, L. Poulain d'Andecy, E. Ragoucy, L. Vinet, \emph{The Askey-Wilson algebra and its avatars}, J. Phys. A 54 (2021), no. 6, Paper No. 063001, 32 pp. (2020).

\bibitem{CPV2} N. Cramp\'e, L. Poulain d'Andecy, L. Vinet, \emph{A Calabi--Yau algebra with $E_6$ symmetry and the Clebsch-Gordan series of $sl(3)$}, J. Lie Theory 31 (2021), no. 4, 1085--1112.

\bibitem{publi-CP2} N. Cramp\'e, L. Poulain d'Andecy, \emph{Baxterisation of the fused Hecke algebra and R-matrices with $gl(N)$-symmetry},  Lett. Math. Phys. 111 (2021), no. 4, Paper No. 92, 21 pp.

\bibitem{publi-CP1} N. Cramp\'e, L. Poulain d'Andecy, \emph{Fused braids and centralizers of tensor representations of $U_q(gl_N)$}, Algebr. Represent. Theor. (2022). https://doi.org/10.1007/s10468-022-10116-7

\bibitem{PRo} L. Poulain d'Andecy, S. Rostam, \emph{Morita equivalences for cyclotomic Hecke algebras of type B and D}, Bull. Soc. Math. France 149 (2021), no. 1, 179--233.

\bibitem{CPV1} N. Cramp\'e, L. Poulain d'Andecy, L. Vinet, \emph{Temperley--Lieb, Brauer and Racah algebras and other centralizers of $\mathfrak{su}(2)$},  Trans. Amer. Math. Soc. 373 (2020), no. 7, 4907--4932.

\bibitem{PRu2} L. Poulain d'Andecy, R. Walker, \emph{Affine Hecke algebras of type D and generalisations of quiver Hecke algebras},  J. Algebra 552 (2020), 1--37.

\bibitem{PRu1} L. Poulain d'Andecy, R. Walker, \emph{Affine Hecke algebras and generalisations of quiver Hecke algebras for type B}, Proc. Edinb. Math. Soc. (2) 63 (2020), no. 2, 531--578.

\bibitem{PW1} L. Poulain d'Andecy, E. Wagner, \emph{The HOMFLYPT polynomials of sublinks and the Yokonuma-Hecke algebras}, Proc. Roy. Soc. Edinburgh Sect. A 148 (2018), no. 6, 1269--1278.

\bibitem{JP2} N. Jacon, L. Poulain d'Andecy, \emph{Clifford theory for Yokonuma--Hecke algebras and deformation of complex reflection groups}, J. London Math. Soc. (2) 96 (2017) 501--523.

\bibitem{P4} L. Poulain d'Andecy, \emph{Young tableaux and representations of Hecke algebras of type ADE}, J. Comb. Algebra 1 (2017), 371--423.

\bibitem{P3} L. Poulain d'Andecy, {\it Invariants for links from classical and affine Yokonuma--Hecke algebras}, Algebraic modeling of topological and computational structures and applications, 77--95, Springer Proc. Math. Stat., 219 (2017).

\bibitem{P2} L. Poulain d'Andecy, \emph{Fusion formulas and fusion procedure for the Yang-Baxter equation}, Algebr. Represent. Theory 20 (2017), no. 6, 1379--1414.

\bibitem{JP1} N. Jacon, L. Poulain d'Andecy, {\it An isomorphism theorem for Yokonuma--Hecke algebras and applications to link invariants}, Mathematische Zeitschrift 283 (2016), no. 1, 301--338.

\bibitem{mCP2} M. Chlouveraki, L. Poulain d'Andecy, {\it Markov traces on cyclotomic and affine Yokonuma--Hecke algebras}, Int. Math. Res. Not. IMRN (2016), no. 14, 4167--4228.

\bibitem{OP7} O. Ogievetsky, L. Poulain d'Andecy, {\it Induced representations and traces for chains of affine and cyclotomic Hecke algebras}, 
Journal of Geometry and Physics, 87 (2015), 354--372.

\bibitem{mCP1} M. Chlouveraki, L. Poulain d'Andecy, {\it Representation theory of the Yokonuma--Hecke algebra}, 
Advances in Mathematics, 259 (2014), no. 10, 134--172.

\bibitem{P1} L. Poulain d'Andecy, {\it Fusion procedure for wreath products of finite groups by the symmetric group},  
 Algebras and Representation Theory, 17 (2014), no. 3, 809--830.

\bibitem{OP6} O. Ogievetsky, L. Poulain d'Andecy, {\it Fusion procedure for cyclotomic Hecke algebras}, SIGMA, 10 (2014), 039, 13 pages.

\bibitem{OP5} O. Ogievetsky, L. Poulain d'Andecy, {\it Alternating subalgebras of Hecke algebras and alternating subgroups of braid groups}, Communications in Algebra, 42 (2014), no. 5, 1921--1936.

\bibitem{OP4} O. Ogievetsky, L. Poulain d'Andecy, {\it Fusion procedure for Coxeter groups of type B and complex reflection groups G(m,1,n)}, Proceedings of the AMS, 142 (2014), no. 9, 2929--2941.

\bibitem{OP3} O. Ogievetsky, L. Poulain d'Andecy, {\it On representations of complex reflection groups G(m,1,n)}, Theoretical and Mathematical Physics, 174 (2013), no. 1, 95--108.

\bibitem{OP2} O. Ogievetsky, L. Poulain d'Andecy, {\it Alternating subgroups of Coxeter groups and their spinor extensions}, Journal of Pure and Applied Algebra, 217 (2013), no. 11, 2198--2211.

\bibitem{OP1} O. Ogievetsky, L. Poulain d'Andecy, {\it On representations of cyclotomic Hecke algebras}, Modern Physics Letters A 26 (2011), no. 11, 795--803.

\end{enumerate}

\begin{thebibliography}{MMM{\etalchar{+}}21}

\addcontentsline{toc}{chapter}{\large{Bibliography\vspace{0.3cm}}}

\itemsep-0.1em

\bibitem[APV]{APV} P. Abramenko, J. Parkinson, H. Van Maldeghem, \emph{A classification of commutative parabolic Hecke algebras}, J. Algebra 385 (2013), 115--133.

\bibitem[AJ]{AJ} F. Aicardi, J. Juyumaya, \emph{Markov trace on the algebra of braids and ties}, Mosc. Math. J. 16 (2016), no. 3, 397--431.

\bibitem[Ale]{Ale} J.W. Alexander, \emph{A lemma on systems of knotted curves}, Proc. Nat. Acad. Sci. U.S.A. 9(3) (1923) 93.


\bibitem[AP]{AP} A. Appel, T. Przezdziecki, \emph{Generalized Schur--Weyl dualities for quantum affine symmetric pairs and orientifold KLR algebras}, (2022) arXiv:2204.04123.

\bibitem[AV]{AV} A. Appel, B. Vlaar, \emph{Universal K-matrices for quantum Kac--Moody algebras}, (2020)
arXiv:2007.09218.


\bibitem[Ar1]{Ar1} E. Artin, \emph{Theorie der Z\"opfe}, Abh. Math. Sem. Univ. Hamburg 4
(1925), 47--72.

\bibitem[Ar2]{Ar2} E. Artin, \emph{Theory of braids}, Ann. of Math. (2) 48 (1947), 101--126.

\bibitem[AYP]{AYP} H. Au-Yang and J.H.H. Perk, \emph{Onsager's star-triangle equation: Master
key to integrability}, Advanced Studies in Pure Mathematics 19 (1989) 57--94.


\bibitem[Bat]{Bat} M.T. Batchelor, \emph{The importance of being integrable: Out of the paper, into the lab}, Int. J.
Mod. Phys. B 28.18 (2014) 1430010.

\bibitem[BF]{BF} M.T. Batchelor and A. Foerster, \emph{Yang--Baxter integrable models in experiments: from condensed matter to ultracold atoms},  J. Phys. A: Math. Theor 49.17 (2016): 173001.

\bibitem[Bax]{Bax} R. J. Baxter, \emph{Exactly solved models in statistical mechanics}, 
Academic Press (1982).

\bibitem[Bir]{Bir} J. Birman, \emph{Braids, links, and mapping class groups}, Annals of Mathematics Studies, No. 82. Princeton University Press, Princeton, N.J.; University of Tokyo Press, Tokyo, 1974. ix+228 pp.

\bibitem[BW]{BW} J. S. Birman and H. Wenzl, \emph{Braids, link polynomials and a new algebra}, Trans. Amer. Math. Soc. 313 (1989), 249--273.

\bibitem[Boh]{Bo} F. Bohnenblust, \emph{The algebraical braid group},  Ann. of Math. (2) 48 (1947), 127--136.


\bibitem[BDVO]{BDVO} C. Bowman, M. De Visscher, R. Orellana, \emph{The partition algebra and the Kronecker coefficients}, Trans. Amer. Math. Soc. 367 (2015), no. 5, 3647--3667.

\bibitem[Br]{Br} R. Brauer, \emph{On algebras which are connected with the semisimple continuous groups,} Ann. of Math. (2) 38 (1937), no. 4, 857--872.


\bibitem[Bro]{Bro} W. Brown, \emph{An algebra related to the orthogonal group,}
Michigan Math. J. 3 (1955), 1--22.

\bibitem[BDK]{BDK} J. Brundan, R. Dipper, A. Kleshchev, \emph{Quantum linear groups and representations of $GL_n(F_q)$},
Mem. Amer. Math. Soc. 149 (2001), no. 706, viii+112 pp.

\bibitem[BK]{BK} J. Brundan, A. Kleshchev, \textit{Blocks of cyclotomic Hecke algebras and Khovanov--Lauda algebras}, Invent. Math. 178(3) (2009) 451--484.

\bibitem[Bum]{Bum} D. Bump, \emph{Automorphic forms and representations},
Cambridge Studies in Advanced Mathematics, 55. Cambridge University Press, Cambridge, (1997) xiv+574 pp.

\bibitem[BZ]{BZ} G. Burde, H. Zieschang, \emph{Knots,} Second edition. De Gruyter Studies in Mathematics, 5. Walter de Gruyter \& Co., Berlin, 2003. xii+559 pp.

\bibitem[CP1]{CP1} V. Chari and A. Pressley, \emph{Quantum affine algebras}, Commun. Math. Phys. 142 (1991) 261--283.

\bibitem[CP2]{CP} V. Chari and A. Pressley, \emph{A guide to quantum groups}, Cambridge University Press (1995).

\bibitem[CP3]{CP2} V. Chari and A. Pressley, \emph{Minimal affinizations of representations of quantum groups: the simply laced case}, J. Algebra 184.1 (1996) 1--30.

\bibitem[Cha]{Ch} E. Chavli, \emph{Universal deformations of the finite quotients of the braid group on 3 strands}, J. Algebra 459 (2016), 238--271.



\bibitem[Chl]{Chl} M. Chlouveraki, \emph{From the Framisation of the Temperley--Lieb algebra to the Jones polynomial: an algebraic approach}, Knots, Low-Dimensional Topology and Applications, Springer PROMS 284 (2019), 247--276.

\bibitem[CJKL]{CJKL} M. Chlouveraki, J. Juyumaya, K. Karvounis, S. Lambropoulou (with an Appendix by W.B.R. Lickorish), {\em Identifying the invariants for classical knots and links from the Yokonuma--Hecke algebra}, Int. Math. Res. Not. IMRN (2020), no. 1, 214--286.

\bibitem[ChLa]{ChLa} M. Chlouveraki and S. Lambropoulou, \emph{The Yokonuma--Hecke algebras and the HOMFLYPT polynomial}, J. Knot Theory Ramifications 22,  No. 14 (2013), 1350080.

\bibitem[CS]{CS} M. Chlouveraki, V. Sécherre, \emph{The affine Yokonuma--Hecke algebra and the pro-p-Iwahori--Hecke algebra}
Math. Res. Lett. 23 (2016), no. 3, 707--718.

\bibitem[CL]{CL} J. Cooke, A. Lacabanne, \emph{Higher rank Askey--Wilson algebras as skein algebras}, arXiv:2205.04414.

\bibitem[CDO+]{CDO+} F. Correa, M. A. del Olmo, I. Marquette, J. Negro, \emph{Polynomial algebras from $su(3)$ and a quadratically superintegrable model on the two sphere}, J. Phys. A 54 (2021), no. 1

\bibitem[CDVM1]{CDVM1} A. Cox, M. De Visscher, P. Martin, \emph{The blocks of the Brauer algebra in characteristic zero}, Represent. Theory 13 (2009), 272--308.

\bibitem[CDVM2]{CDVM2} A. Cox, M. De Visscher, P. Martin, \emph{A geometric characterisation of the blocks of the Brauer algebra}, J. Lond. Math. Soc. (2) 80 (2009), no. 2, 471--494.


\bibitem[CVZ]{CVZ} N. Crampé, L. Vinet, M. Zaimi, \emph{Temperley--Lieb, Birman--Murakami--Wenzl and Askey--Wilson algebras and other centralizers of $U_q(sl_2)$}, Ann. Henri Poincaré 22 (2021), no. 10, 3499--3528.

\bibitem[Cur]{Cur} C. W. Curtis, \emph{On Lusztig's isomorphism theorem for Hecke algebras}, J. Algebra 92 (1985), no. 2, 348--365.

\bibitem[Cur2]{Cur2} C. W. Curtis, \emph{Pioneers of representation theory: Frobenius, Burnside, Schur, and Brauer},
History of Mathematics, 15. American Mathematical Society, Providence, RI; London Mathematical Society, London (1999).


\bibitem[DBDC+]{DBDC+}H. De Bie, H. De Clercq, W. van de Vijver, \emph{The higher rank q-deformed Bannai--Ito and Askey--Wilson algebra}, Comm. Math. Phys. 374 (2020), no. 1, 277--316.

\bibitem[DBG+]{DBG+} H. De Bie, V. Genest, W. van de Vijver, L. Vinet, \emph{A higher rank Racah algebra and the $\mathbb{Z}_2^n$ Laplace--Dunkl operator}, J. Phys. A 51 (2018), no. 2, 025203, 20 pp.

\bibitem[DBI+]{DBI+}H. De Bie, P. Iliev, W. van de Vijver, L. Vinet, \emph{The Racah algebra: an overview and recent results}, Lie groups, number theory, and vertex algebras, 3--20, Contemp. Math., 768, Amer. Math. Soc., [Providence], RI, (2021).

\bibitem[DBVV]{DBVV} H. De Bie, W. van de Vijver, \emph{A discrete realization of the higher rank Racah algebra},
Constr. Approx. 52 (2020), no. 1, 1--29.


\bibitem[DG]{DG} S. Doty, A. Giaquinto, \emph{Presenting Schur algebras}, Int. Math. Res. Not. 2002, no. 36, 1907--1944.

\bibitem[Dre]{Dre} V. Drensky, \emph{Computing with matrix invariants}, Mathematica Balkanika 21, 141(172) (2007).

\bibitem[Dri1]{Dr-QG} V.G. Drinfeld, \emph{Hopf algebras and the Yang--Baxter quantum equation,} Dokl. Akad. Nauk SSSR. Vol. 283. No. 5. (1985).

\bibitem[Dri2]{Dri} V.G. Drinfeld, \emph{Quantum groups}, in Proceedings of the International Congress of Mathematicians (A.M. Gleason, ed.), Amer. Math. Soc. (1986) 798--820.

\bibitem[EMTW]{EMTW} B. Elias, S. Makisumi, U. Thiel, G. Williamson, \emph{Introduction to Soergel bimodules},
RSME Springer Series, 5. Springer, Cham, (2020) 588 pp.

\bibitem[EK]{EK} N. Enomoto, M. Kashiwara, \emph{Symmetric crystals and affine Hecke algebras of type B}, Proc. Japan Acad. 82 Ser. A no. 8 (2006)  131--136.

\bibitem[ERH]{ERH} J. Espinoza, S. Ryom-Hansen, \emph{Cell structures for the Yokonuma--Hecke algebra and the algebra of braids and ties}, J. Pure Appl. Algebra 222 (2018), no. 11, 3675--3720.

\bibitem[EFK]{EFK} P. Etingof, I. Frenkel, A.A. Kirillov, \emph{Lectures on representation theory and Knizhnik--Zamolodchikov equations}, No. 58, American Mathematical Soc. (1998).

\bibitem[EG]{EG} P. Etingof, and V. Ginzburg, 
\emph{Noncommutative del Pezzo surfaces and Calabi--Yau algebras}, J. Eur. Math. Soc 12 (2010) 1371--1416.

\bibitem[ELOR]{ELOR} P. Etingof, S. Loktev, A. Oblomkov, L. Rybnikov, \emph{A Lie-theoretic construction of spherical symplectic reflection algebras}, Transformation Groups 13 (2008) 541--556.

\bibitem[Fad]{Fad} L.D. Faddeev, \emph{How Algebraic Bethe Ansatz works for integrable model}, arXiv:hep-th/9605187

\bibitem[FST]{FST}  L.D. Faddeev, E.K. Sklyanin, L.A. Takhtadzhyan, 
\emph{Quantum inverse problem method I,} Theor. and Math. Phys. 40 (1979) 86.

\bibitem[For]{For} E. Formanek, \emph{The invariants of $n\times n$ matrices}, Invariant theory, 18--43,
Lecture Notes in Math., 1278, Springer, Berlin, 1987.

\bibitem[FR]{FR} I. Frenkel,  N. Reshetikhin, \emph{Quantum affine algebras and holonomic difference equations,} Commun. Math. Phys. 146 (1992), 1--60.

\bibitem[FYH+]{HOMFLY} P. Freyd, D. Yetter, J. Hoste, W.B.R. Lickorish, K. Millett, A. Ocneanu, \emph{A new polynomial invariant of knots and links}, Bull. Amer. Math. Soc. (N.S.) 12 (1985), no. 2, 239--246.

\bibitem[Ful]{Ful} W. Fulton, \emph{Young tableaux. With applications to representation theory and geometry.} London Mathematical Society Student Texts, 35. Cambridge University Press, Cambridge, 1997. x+260 pp.

\bibitem[FH]{FH} W. Fulton, J. Harris, \emph{Representation theory: a first course}, Springer (1991).

\bibitem[GWH]{GWH} S. Gao, Y. Wang, B. Hou, \emph{The classification of Leonard triples of Racah type},
Linear Algebra Appl. 439 (2013), no. 7, 1834--1861.

\bibitem[GJ]{GJ} M. Geck, N. Jacon, \emph{Representations of Hecke algebras at roots of unity},
Algebra and Applications, 15. Springer-Verlag London, Ltd., London, (2011) xii+401 pp.

\bibitem[GM]{GM} M. Geck, G. Malle, \emph{The character theory of finite groups of Lie type.
A guided tour.} Cambridge Studies in Advanced Mathematics, 187. Cambridge University Press, Cambridge, (2020).

\bibitem[GP]{GP} M. Geck, G. Pfeiffer, \emph{Characters of finite Coxeter groups and Iwahori-Hecke algebras}. London Mathematical Society Monographs, 21, Oxford University Press, 2000.

\bibitem[GVZ]{GVZ} V. Genest, L. Vinet, and A. Zhedanov, \emph{The Racah algebra and superintegrable models},
Journal of Physics: Conference Series 512, 012011 (2014).

\bibitem[Gin]{Gin} V. Ginzburg, \emph{Calabi-Yau algebras}, arXiv:math/0612139.

\bibitem[GRS]{GRS} C. Gomez, M. Ruiz-Altaba and G. Sierra, \emph{Quantum groups in two-dimensional physics},  Cambridge University Press (1996).

\bibitem[GW]{GW} R. Goodman and R. Wallach, \emph{Symmetry, representations, and invariants}, Springer (2009).

\bibitem[GZ]{GZ} Ya. A. Granovskii, A.S. Zhedanov, \emph{Nature of the symmetry group of the $6j$-symbol,}
JETP 67:1982--1985 (1988).

\bibitem[Gre]{Gr} J. A. Green, \emph{Polynomial representations of $GL_n$}, Lecture Notes in Mathematics, Vol. 830, Springer, 1980.

\bibitem[HR]{HR} T. Halverson, A. Ram, \emph{Partition algebras}, European J. Combin. 26 (2005), no. 6, 869--921

\bibitem[HW]{HW} P. Hanlon, D. Wales,
\emph{On the decomposition of Brauer's centralizer algebras,}
J. Algebra 121 (1989), no. 2, 40--445.

\bibitem[Haw1]{Haw1} T. Hawkins, \emph{Emergence of the theory of Lie groups. An essay in the history of mathematics 1869--1926}, Sources and Studies in the History of Mathematics and Physical Sciences. Springer-Verlag, New York (2000).

\bibitem[Haw2]{Haw2} T. Hawkins, \emph{The mathematics of Frobenius in context. A journey through 18th to 20th century mathematics}, Sources and Studies in the History of Mathematics and Physical Sciences. Springer, New York (2013).

\bibitem[Hec]{Hec} E. Hecke, \emph{\"Uber Modulfunktionen und die Dirichletschen Reihen mit Eulerscher Produktentwicklung. I},
Math. Ann. 114 (1937), no. 1, 1--28.

\bibitem[Her]{He} D. Hernandez, \emph{Advances in R-matrices and their applications (after Maulik-Okounkov, Kang-Kashiwara-Kim-Oh,...),} (2017) arXiv:1704.06039.

\bibitem[HoRo]{HoRo} M. Hochster, J. Roberts, \emph{Rings of invariants of reductive groups acting on regular rings are Cohen--Macaulay}, Advances in Math. 13 (1974), 115--175.

\bibitem[Hoe]{Ho} P. Hoefsmit, \emph{Representations of Hecke algebras of finite groups with BN-pairs of classical type}, Ph. D. thesis, University of British Columbia (1974). 

\bibitem[HL]{HL} R. B. Howlett, G. I. Lehrer, \emph{Induced cuspidal representations and generalised Hecke rings},
Invent. Math. 58 (1980), no. 1, 37--64.

\bibitem[Hum]{Hum} J. E. Humphreys, \emph{Reflection groups and Coxeter groups},
Cambridge Studies in Advanced Mathematics, 29. Cambridge University Press, Cambridge, (1990) xii+204 pp.


\bibitem[IMO]{IMO} A.P. Isaev, A.I. Molev, O.V. Ogievetsky, \emph{A new fusion procedure for the Brauer algebra and evaluation homomorphisms}, Int. Math. Res. Not. IMRN 2012, no. 11, 2571--2606.



\bibitem[IO1]{IO1} A. Isaev and O. Ogievetsky, \emph{On representations of Hecke algebras}, Czech. Journ. Phys. 55 No. 11 (2005) 1433--1441.

\bibitem[IO2]{IO2} A. Isaev, O. Ogievetsky, \emph{On Baxterized solutions of reflection equation and integrable chain models}, 
Nucl. Phys. B 760[PM] (2007) 167--183.

\bibitem[Iwah]{Iwah} N. Iwahori, \emph{On the structure of a Hecke ring of a Chevalley group over a finite field},
J. Fac. Sci. Univ. Tokyo Sect. I 10 (1964), 215--236.

\bibitem[IM]{IM} N. Iwahori, H. Matsumoto, \emph{On some Bruhat decomposition and the structure of the Hecke rings of p-adic Chevalley groups}, Inst. Hautes Études Sci. Publ. Math. No. 25 (1965), 5--48.

\bibitem[Iwa]{Iwa} H. Iwaniec, \emph{Topics in classical automorphic forms},
Graduate Studies in Mathematics, 17. American Mathematical Society, Providence, RI (1997). xii+259 pp.

\bibitem[JS]{JS} J. C. Jantzen, G. M. Seitz, \emph{On the representation theory of the symmetric groups},
Proc. London Math. Soc. (3) 65 (1992), no. 3, 475--504.

\bibitem[Jim1]{Ji-QG} M. Jimbo, \emph{A $q$-difference analogue of $U(g)$ and the Yang--Baxter equation}, Lett. Math. Phys. 10 (1985) 63--69.

\bibitem[Jim2]{Ji86} M. Jimbo, \emph{A q-Analogue of $U(\mathfrak{gl}(N+1))$, Hecke algebra, and the Yang--Baxter equation}, Lett. Math. Phys. 11 (1986), 247--252.

\bibitem[Jim3]{Ji-int}M. Jimbo, \emph{Introduction to the Yang--Baxter equation,} in ``Braid Group, Knot Theory And Statistical Mechanics'' 9 (1991) 111.

\bibitem[Jim4]{Ji-ed} M. Jimbo (Editor), \emph{Yang-Baxter equation in integrable systems}, Vol. 10. World Scientific, (1990).

\bibitem[JM]{JM} M. Jimbo, T. Miwa, \emph{Algebraic analysis of solvable lattice models}, American Mathematical Soc. Vol. 85 (1994).

\bibitem[Jo1]{Jo1} V.F.R. Jones, \emph{Hecke algebra representations of braid groups and link polynomials}, Annals of Math. 126 (1987), no. 2, 335--388.

\bibitem[Jo2]{Jo2} V.F.R. Jones, \emph{On a Certain Value of the Kauffman Polynomial}, Commun. Math. Phys. 125(1989) 459.

\bibitem[Jo3]{Jo3} V.F.R. Jones, \emph{Baxterization}, Int. J. Mod. Phys.B 4 (1990) 701.

\bibitem[Jo4]{Jo4} V.F.R. Jones, \emph{The Potts model and the symmetric group}, Subfactors (Kyuzeso, 1993), 259--267, World Sci. Publ., River Edge, NJ, 1994.

\bibitem[Juy1]{Juy1} J.~Juyumaya, \emph{Sur les nouveaux g\'en\'erateurs de l'alg\`ebre de Hecke H(G,U,1)}, J.~Algebra 204 (1998) 49--68.

\bibitem[Juy2]{Juy2} J.~Juyumaya, \emph{Markov trace on the Yokonuma--Hecke algebra}, J.~Knot Theory Ramifications 13 (2004) 25--39.

\bibitem[JL3]{JL3} J. Juyumaya, S. Lambropoulou, \emph{An adelic extension of the Jones polynomial}, M. Banagl, D. Vogel (eds.)
The mathematics of knots, Contributions in the Mathematical and Computational Sciences, Vol. 1, (2011) Springer.

\bibitem[KMP]{KMP} G. Kalnins, W. Miller Jr., and S. Post, \emph{Contractions of 2D 2nd Order Quantum Superintegrable
Systems and the Askey Scheme for Hypergeometric Orthogonal Polynomials},
SIGMA Symmetry Integrability Geom. Methods Appl. 9 (2013), Paper 057, 28 pp.

\bibitem[Kam]{Kam} S. Kamada, \emph{Braid and knot theory in dimension four},
Mathematical Surveys and Monographs, 95. American Mathematical Society, Providence, RI (2002).

\bibitem[KK]{KK} S.-J. Kang, M. Kashiwara, \textit{Categorification of highest weight modules via Khovanov--Lauda--Rouquier algebras}. Invent. Math. 190 (2012) 699-742.

\bibitem[KKK]{KKK}S.-J. Kang, M. Kashiwara, M. Kim, \emph{Symmetric quiver Hecke algebras and R-matrices of quantum affine algebras}, Invent. Math. 211 (2018), no. 2, 591--685.

\bibitem[KM]{KM} M. Kashiwara, V. Miemietz, \emph{Crystals and affine {H}ecke algebras of type D}, Proc. Japan Acad. 83 Ser. A Math. Sci. no. 7  (2007)  135--139.

\bibitem[Kas]{Kas} C. Kassel, \emph{Quantum groups}, Springer (1995).

\bibitem[KRT]{KRT} C. Kassel, M. Rosso, V. Turaev, \emph{Quantum groups and knot invariants},
Panoramas et Synth\`eses 5. Soci\'et\'e Math\'ematique de France, Paris (1997).

\bibitem[KT]{KT} C. Kassel, V. Turaev, \emph{Braid groups}, Graduate Texts in Mathematics, Vol. 247, Springer, 2008.

\bibitem[Kau]{Kau} L.H. Kauffman, \emph{An invariant of regular isotopy}, Trans. Amer. Math.
Soc. 318 (1990), no. 2, 417--471.

\bibitem[KL]{KL} D. Kazhdan, G. Lusztig, \emph{Representations of Coxeter groups and Hecke algebras},
Invent. Math. 53 (1979), no. 2, 165--184.

\bibitem[KhLa1]{KhLa1} M. Khovanov, A. D. Lauda, \textit{A diagrammatic approach to categorification of quantum groups I}. Represent. Theory 13 (2009) 309--347.

\bibitem[KhLa2]{KhLa2} M. Khovanov, A. D. Lauda, \textit{A diagrammatic approach to categorification of quantum groups II}. Trans. Amer. Math. Soc. 363 (2011) 2685--2700.

\bibitem[Kir]{Kir} A.A. Kirillov, \emph{Certain division algebras over a field of rational functions}, Funct Anal Its Appl 1   (1967) 87--88.

\bibitem[KS]{KS} A. Klimyk, K. Schm\"udgen, \emph{Quantum groups and their representations}, Springer, 2012.


\bibitem[Knu]{Knu} D. Knuth, \emph{Permutations, matrices, and generalized Young tableaux}, Pacific J. Math. 34(3) (1970) 709--727.

\bibitem[KLS]{KLS} R. Koekoek, P. A. Lesky, R. F. Swarttouw, \emph{Hypergeometric orthogonal polynomials and their q-analogues},
With a foreword by Tom H. Koornwinder. Springer Monographs in Mathematics. Springer-Verlag, Berlin, (2010).


\bibitem[KBI]{KBI} V.E. Korepin, N.M. Bogoliubov, A.G. Izergin, \emph{Quantum inverse scattering method and
correlation functions}, Cambridge University Press (1993).


\bibitem[KR]{KR} P.P. Kulish, N.Y. Reshetikhin, \emph{Quantum linear problem for the sine-Gordon equation and higher representations}, J. Soviet Math. 23(4) (1983) 2435--2441.


\bibitem[Lam]{Lam} S. Lambropoulou, \emph{Solid torus links and Hecke algebras of B-type}, Proceedings of the Conference on Quantum Topology (Manhattan, KS, 1993), 225--245, World Sci. Publ., River Edge, NJ, (1994).


\bibitem[LMZ]{LMZ} D. Latini, I. Marquette, Y.-Z. Zhang, \emph{Embedding of the Racah algebra $R(n)$ and superintegrability},
Ann. Physics 426 (2021), Paper No. 168397, 18 pp.

\bibitem[LR]{LR} R. Leduc, A. Ram, \emph{A ribbon Hopf algebra approach to the irreducible representations of centralizer algebras: the Brauer, Birman--Wenzl, and type A Iwahori--Hecke algebras,} Adv. Math. 125 (1997), no. 1, 1--94.

\bibitem[Leh]{Leh} G. Lehrer, \emph{A survey of Hecke algebras and the Artin braid groups}, Braids (Santa Cruz, CA, 1986), 365--385,
Contemp. Math., 78, Amer. Math. Soc., Providence, RI (1988).

\bibitem[LZ1]{LZ1} G. Lehrer, R. Zhang, \emph{Strongly multiplicity free modules for Lie algebras and quantum groups}, Journal of Algebra 306(1) (2006) 138--174.

\bibitem[LZ2]{LZ2} G. Lehrer, R. Zhang, \emph{A Temperley--Lieb analogue for the BMW algebra}, in Representation theory of algebraic groups and quantum groups. Birkh\"auser Boston, (2010) 155--190.

\bibitem[Lus1]{Lus2} G. Lusztig, \emph{Characters of reductive groups over a finite field},
Annals of Mathematics Studies, 107. Princeton University Press, Princeton, NJ, (1984) xxi+384 pp.

\bibitem[Lus2]{Lus1} G. Lusztig, \emph{Hecke algebras with unequal parameters},
CRM Monograph Series, 18. American Mathematical Society, Providence, RI, (2003) vi+136 pp.

\bibitem[Lus3]{Lus} G.~Lusztig, \emph{Character sheaves on disconnected groups. VII.} Represent. Theory 9 (2005), 209--266.

\bibitem[Mac1]{Mac1} G.W. Mackey, \emph{The scope and history of commutative and noncommutative harmonic analysis},
History of Mathematics, 5. American Mathematical Society, Providence, RI; London Mathematical Society, London (1992).

\bibitem[Ma1]{Ma1} I. Marin, \emph{The cubic Hecke algebra on at most 5 strands}, J. Pure Applied Algebra 216 (2012), 2754--2782

\bibitem[Ma2]{Ma3} I. Marin, \emph{Artin groups and Yokonuma--Hecke algebras},
Int. Math. Res. Not. IMRN 2018, no. 13, 4022--4062.

\bibitem[Ma3]{Ma4} I. Marin, \emph{Lattice extensions of Hecke algebras}, J. Algebra 503 (2018), 104--120.

\bibitem[Ma4]{Ma2} I. Marin, \emph{A maximal cubic quotient of the braid algebra I}, J. Algebra (2020).

\bibitem[Mark]{Mark} A.A. Markov, \emph{\"Uber die freie Aquivalenz der geschlossner Zopfe,} Rec. Soc. Math. Moscou 1 (1935): 73--78.

\bibitem[Mar1]{Mar1} P. Martin, \emph{Potts models and related problems in statistical mechanics}, Series on Advances in Statistical Mechanics, 5. World Scientific Publishing Co., Inc., Teaneck, NJ, 1991.

\bibitem[Mar2]{Mar2} P. Martin, \emph{Temperley--Lieb algebras for nonplanar statistical mechanics --- the partition algebra construction}, J. Knot Theory Ramifications 3(1994), 51--82.

\bibitem[Mar3]{Mar3} P. Martin, \emph{The structure of the partition algebras}, J. Algebra 183(1996), 319--358.

\bibitem[Mar4]{Mar4} P. Martin, \emph{The partition algebra and the Potts model transfer matrix spectrum in high dimensions}, J. Phys. A:Math. Gen. 33 (2000), 3669--3695.

\bibitem[Mar5]{Mar5} P. Martin, \emph{The decomposition matrices of the Brauer algebra over the complex field,} Trans. Amer. Math. Soc. 367 (2015), no. 3, 179--1825.

\bibitem[MMcA]{MMcA} P.P. Martin, D.S. McAnally, \emph{On commutants, dual pairs and nonsemisimple algebras from statistical mechanics}, Infinite analysis, Part A, B (Kyoto, 1991), 675--705, Adv. Ser. Math. Phys., 16, World Sci. Publ., River Edge, NJ, (1992).

\bibitem[MR]{MR} P. Martin, G. Rollet \emph{ The Potts model representation and a Robinson--Schensted correspondence for the partition algebra}, Compositio Math. 112 (1998), no. 2, 237--254.


\bibitem[McG]{McG} J.B. McGuire, \emph{Study of exactly solvable one-dimensional $N$-body problems}, J. Math. Physics, 5 (1964) 622--636.


\bibitem[MRR]{MRR} A. Morin-Duchesne, J. Rasmussen and D. Ridout, \emph{Boundary algebras and Kac modules for logarithmic minimal models}, Nucl. Phys. B899, 677 (2015)

\bibitem[Mor]{Mor} H. Morton, \emph{A basis for the Birman--Wenzl algebra}, math.QA/1012.3116, (2010).

\bibitem[MT]{MT} H. Morton, P. Traczyk, \emph{Knots and algebras}, Zaragoza 201 (1990): 220.

\bibitem[Mu1]{Mu1} J. Murakami, \emph{The Kauffman polynomial of links and representation theory}, Osaka J.
Math. 24 (1987) 745--758.

\bibitem[Mu2]{Mu2} J. Murakami, \emph{Solvable lattice models and algebras of face operators}, Adv. Studies in
Pure Math. 19 (1989) 399--415

\bibitem[Na1]{Na1} M. Nazarov, \emph{Young's orthogonal form for Brauer's centralizer algebra},
J. Algebra 182 (1996), no. 3, 664--693.


\bibitem[OEIS]{OEIS} OEIS Foundation Inc., \emph{The On-Line Encyclopedia of Integer Sequences}, 
http://oeis.org (2021)

\bibitem[Ogg]{Ogg} A. Ogg, \emph{Modular forms and Dirichlet series}, W. A. Benjamin, Inc., New York-Amsterdam (1969) xvi+173 pp.

\bibitem[OV]{OV} A. Okounkov and A. Vershik, \emph{A new approach to representation theory of symmetric groups II},
Selecta Math (New series) 2 No. 4 (1996) 581--605.

\bibitem[Ons]{Ons} L. Onsager, \emph{Crystal Statistics I. A Two-Dimensional Model with an Order-Disorder Transition}, Phys. Rev. 65 (1944) 117--149

\bibitem[PAY]{PAY} J.H.H. Perk and H. Au-Yang, \emph{Yang--Baxter Equation},
in Encyclopedia of Mathematical Physics, eds. J.-P. Fran\c{c}oise, G.L. Naber and Tsou S.T., 
2006, Vol. 5, pp. 465--473.

\bibitem[PST]{PST} Z. Pluha\v{r}, Yu. F. Smirnov and V. N. Tolstoy,
\textsl{Clebsch--Gordan coefficients of SU(3) with simple symmetry properties,}
J. Phys. A 19 (1986) 21.

\bibitem[Pro]{Pr} C. Procesi, \emph{The invariant theory of $n\times n$ matrices}, Adv. Math. 19 (1976) 306--381.

\bibitem[PT]{PT} J.H. Przytycki, P. Traczyk, \emph{Invariants of links of Conway type},
Kobe J. Math. 4 (1988), no. 2, 115--139.

\bibitem[Rac]{Rac} G. Racah, \emph{Theory of Complex Spectra II}, Physical Review. 62 (9-10): 438--46 (1942).

\bibitem[Ram]{Ra} A. Ram, \emph{Characters of Brauer's centralizer algebras}, Pacific J. Math. 169 (1995), no. 1, 173--200.

\bibitem[Ram2]{Ra2} A. Ram, \emph{Skew shape representations are irreducible}, Combinatorial and Geometric Representation Theory (Seoul, 2001), Contemporary Mathematics 325 (AMS, Providence, RI, 2003) 161--189.

\bibitem[Raz]{Raz} Yu.P. Razmyslov, \emph{Trace identities of full matrix algebras over a field of characteristic zero}
(Russian), Izv. Akad. Nauk SSSR, Ser. Mat. 38 (1974) 723--756. Translation: Math. USSR,
Izv. 8 (1974) 727--760.

\bibitem[Res]{Res} N. Yu. Reshetikhin, \emph{Quantized universal enveloping algebras, the Yang--Baxter equations
and invariants of links, I and II}, LOMI preprints E-4-87 and E-17-87, Leningrad (1987).


\bibitem[Rou]{Rou} R. Rouquier, \textit{2-Kac--Moody algebras} (2008) arXiv:0812.5023.

\bibitem[Rut]{Rut} D. Rutherford , \emph{Substitutional analysis}. Edinburgh University Publications, Science and Mathematics No. 1., University Press,  Edinburgh (1948).

\bibitem[Sch]{Sch} I. Schur, \emph{\"Uber eine Klasse von Matrizen, die sich einer gegebenen Matrix zuordnen lassen}, in Gesammelte Abhandlungen. Band I. Springer-Verlag, Berlin-New York, 1973. xv+491 pp.

\bibitem[SVV]{SVV} P. Shan, M. Varagnolo, E. Vasserot, \emph{Canonical bases and affine {H}ecke algebras of type {$D$}}. Adv. Math. 227 no. 1 (2011) 267--291.

\bibitem[Shi]{Shi} G. Shimura, \emph{Introduction to the arithmetic theory of automorphic functions}, Reprint of the 1971 original. Publications of the Mathematical Society of Japan, 11. Kanô Memorial Lectures, 1. Princeton University Press, Princeton, NJ, (1994) xiv+271 pp.

\bibitem[Sib]{Sib} K. S. Sibirskii, \emph{Algebraic invariants of a system of matrices}, Sibirsk. Mat. Z. 9 (1968) 152--164.

\bibitem[Skl]{Skl} E.K. Sklyanin, \emph{Some algebraic structures connected with the Yang--Baxter equation,} , Funct. Anal. Appl, 16(4) (1982) 263-270.

\bibitem[Sol1]{Sol1} M. Solleveld, \emph{Endomorphism algebras and Hecke algebras for reductive p-adic groups}, arXiv:2005.07899 (2020).

\bibitem[Sol2]{Sol2} M. Solleveld, \emph{Affine Hecke algebras and their representations},
Indag. Math. (N.S.) 32 (2021), no. 5, 1005--1082.


\bibitem[Sta1]{Sta1} R. Stanley, \emph{Hilbert functions of graded algebras},
Advances in Math. 28 (1978), no. 1, 57--83.

\bibitem[Sta2]{St} R. Stanley, \emph{Enumerative Combinatorics, vol. 2}, Cambridge Stud. Adv. Math (1999).

\bibitem[Tera]{Tera} Y. Teranishi, \emph{The ring of invariants of matrices}, Nagoya Math. J. 104 (1986), 149--161.

\bibitem[Ter]{Ter} P. Terwilliger, \emph{Two linear transformations each tridiagonal with respect to an eigenbasis of the other}, Linear Algebra Appl. 330 (2001), no. 1--3, 149--203.

\bibitem[Tur]{Tur} T.G. Turaev, \emph{The Yang-Baxter equation and invariants of links}, Invent. Math., 92 (1988)
527--553.

\bibitem[VV]{VV} M. Varagnolo, E. Vasserot,  \textit{Canonical bases and affine Hecke algebras of type B}. Invent. Math. 183(3) (2011) 593--693.

\bibitem[Wen]{Wen} H. Wenzl, \emph{On the structure of Brauer's centralizer algebras},
Ann. of Math. (2) 128 (1988), no. 1, 173--193.

\bibitem[Wey1]{We1} H. Weyl, \emph{Generalized Riemann matrices and factor sets,}
Ann. of Math. (2) 37 (1936), no. 3, 709--745.

\bibitem[Wey2]{We2} H. Weyl, \emph{The classical groups, their invariants and representations}, Princeton University Press (1946).

\bibitem[Wey3]{Wey3} H. Weyl, \emph{The theory of groups and quantum mechanics}, Courier Corporation (1950).

\bibitem[Yan1]{Yang1} C. N. Yang, \emph{Some exact results for the many-body problem in one dimension with repulsive $\delta$-function interaction,}
Phys. Rev. Lett. 19 (1967) 1312.

\bibitem[Yan2]{Yang2} C. N. Yang, \emph{$S$-matrix for the one-dimensional $n$-body problem with repulsive or attractive $\delta$-function interaction,}
Phys. Rev. 168 (1968) 1920. 

\bibitem[Yok]{Yok} T. Yokonuma, \emph{Sur la structure des anneaux de Hecke d'un groupe de Chevalley fini}, C. R. Acad. Sci. Paris Ser. I Math. 264 (1967) 344--347.

\bibitem[You]{Yo} A. Young, \emph{On quantitative substitutional analysis} (sixth and eigth), Proc. London Math. Soc. 34(2) (1931) 196--230 and 37 (1934), 441--495. 

\bibitem[ZZ]{ZZ} A. B. Zamolodchikov, A. B. Zamolodchikov, 
\emph{Factorized S-matrices in two dimensions as the exact solutions of certain relativistic quantum field theory models,}
Annals of Physics 120 (1979) 253--291.

\bibitem[ZGB]{ZGB} R.B. Zhang, M.D. Gould, A.J. Bracken, 
\emph{Quantum group invariants and link polynomials,} Comm. Math. Phys. 137 (1991), no. 1, 13--27.

\bibitem[Zhe]{Zhe} A. Zhedanov,
\emph{``Hidden symmetry'' of Askey--Wilson polynomials,}
Theoret. and Math. Phys. 89 (1991), no. 2, 1146--1157.


\bibitem[ZJ2]{ZJ2} P. Zinn-Justin, \emph{The trigonometric E8 R-matrix}, Lett. Math. Phys. 110 (2020), no. 12, 3279--3305.

\end{thebibliography}
\end{document}